# CATEGORICAL ASPECTS OF TOPOLOGICAL QUANTUM FIELD THEORIES

BRUCE BARTLETT

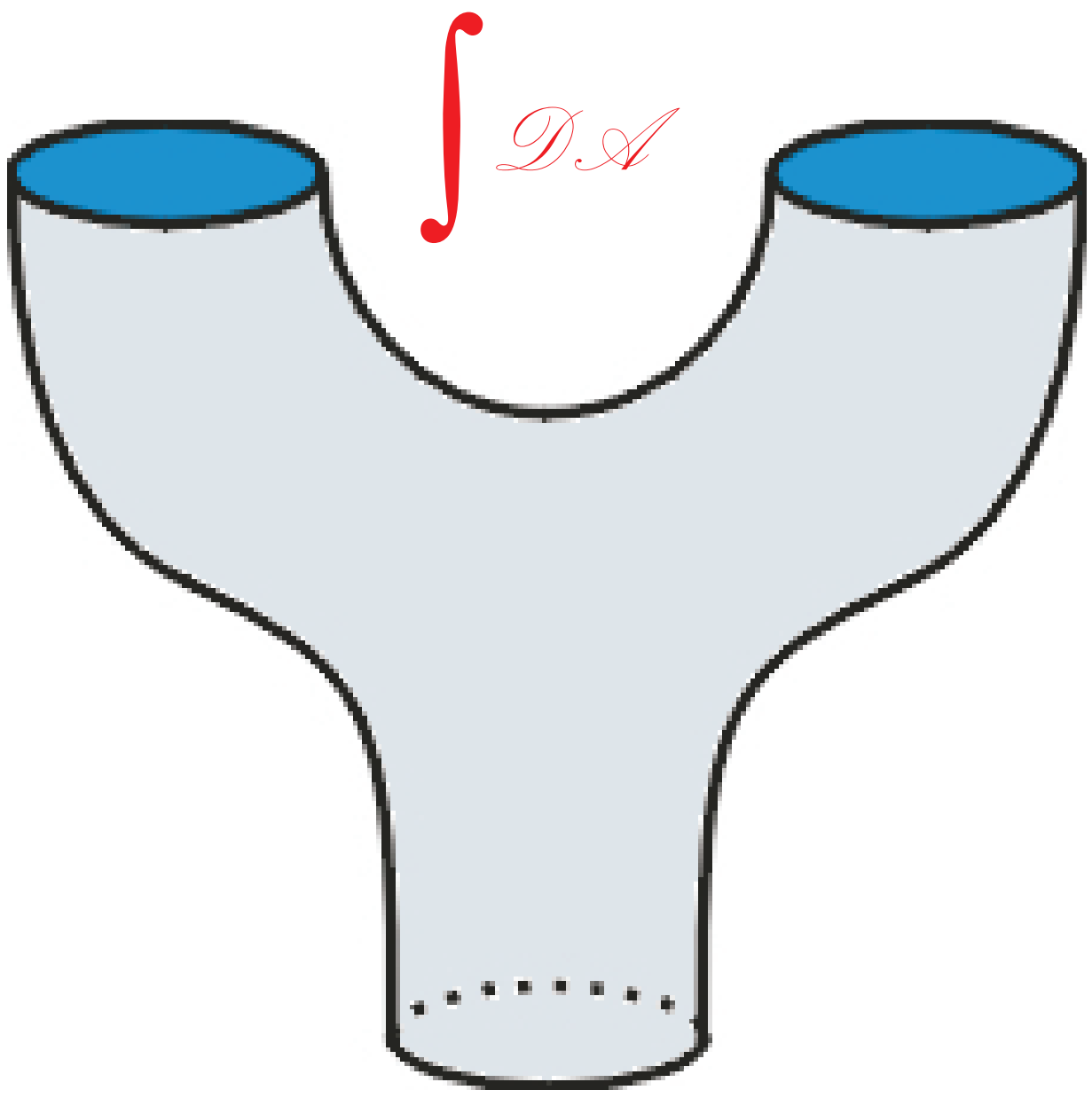

MASTERS THESIS, UTRECHT UNIVERSITY

SUPERVISORS : DR. S. VANDOREN, PROF. I. MOERDIJK

SEPTEMBER 2005

# Abstract

This thesis provides an introduction to the various category theory ideas employed in topological quantum field theory (TQFT). These theories can be viewed as symmetric monoidal functors from topological cobordism categories into the category of vector spaces. In two dimensions, they are classified by Frobenius algebras. In three dimensions, and under certain conditions, they are classified by modular categories. These are special kinds of categories in which topological notions such as braidings and twists play a prominent role. There is a powerful graphical calculus available for working in such categories, which may be regarded as a generalization of the Feynman diagrams method familiar in physics. This method is introduced and the necessary algebraic structure is graphically motivated step by step.

A large subclass of two-dimensional TQFTs can be obtained from a lattice gauge theory construction using triangulations. In these theories, the gauge group is finite. This construction is reviewed, from both the original algebraic perspective as well as using the graphical calculus developed in the earlier chapters. This finite gauge group toy model can be defined in all dimensions, and has a claim to being the simplest non-trivial quantum field theory. We take the opportunity to show explicitly the calculation of the modular category arising from this model in three dimensions, and compare this algebraic data with the corresponding data in two dimensions, computed both geometrically and from triangulations. We use this as an example to introduce the idea of a quantum field theory as producing a tower of algebraic structures, each dimension related to the previous by the process of categorification.

# Uittreksel

Deze scriptie vormt een inleiding in het gebruik van de categorie theorie in topologische kwantumvelden theorie (TQFT). TQFT's kunnen worden beschouwd als symmetrische monoidale functors van topologische cobordisme-categorieën naar de categorie van vectorruimten. Twee dimensionale TQFT's worden geclassificeerd door Frobenius algebras. Drie dimensionale TQFT's worden onder bepaalde voorwaarden geclassificeerd door modulaire categorieën. Dit is een specifieke klasse categorieën waarin begrippen als braidings en twists een prominente rol spelen. Er is een krachtige grafische methode voor handen om te rekenen in zulke categorieën, die opgevat kan worden als een generalizatie van the Feynman-diagrammethode, bekend uit de natuurkunde. Deze methode wordt ingevoerd en de vereiste algebraische structuren worden stap voor stap grafisch gemotiveerd.

Een grote subklasse van tweedimensionale TQFT's kan ontstaan uit triangulaties van rooster ijktheorieën. In deze theorieën is de ijkgroep eindig. We bekijken constructies enerzijds vanuit het algebraische perspectief, en anderzijds vanuit de eerder in de scriptie ontwikkelde grafische rekenmethode. Het eenvoudige model met de eindige ijkgroep kan in alle dimensies worden gedefinieerd, en is mogelijk de meest eenvoudige niet-triviale kwantumveldentheorie. We maken van de gelegenheid gebruik om expliciet de berekening van de modulaire categorie afkomstig van dit drie dimensionale model te laten zien. We vergelijken de algebraische data met de uitkomst van het twee dimensionale model. Dit laatste model wordt zowel berekend met de geometrische als de triangulatiemethode. Aan de hand van dit voorbeeld kan men zien dat een kwantumveldentheorie een toren van algebraische structuren produceert, waar elke dimensie gerelateerd is aan de vorige via zogenaamde categorificatie.

## Acknowledgements

First and foremost, I would like to thank the Julius School of Physics and Astronomy for giving me the opportunity to take part in the International Masters Programme at Utrecht. In this regard special thanks go to Leonie Silkens and Dr. Frank Witte for their endless assistance, hospitality and friendly demeanour. Secondly, I would like to sincerely thank my two supervisors, Dr. Stefan Vandoren and Prof. Ieke Moerdijk. Dr. Vandoren was a big help in my attempt to understand how topological string theory and conformal theory ideas fit into the 'topological quantum field theory' framework, always encouraging me to emphasize the connection with physics, stressing the need to write more clearly and concisely (a task in which I failed miserably) and patiently bearing with me as the write-up took longer than expected! A warm and special mention must go to Prof. Moerdijk, who generously gave many hours of his time in explaining (and asking *me* to explain!) the various notions in this thesis, and for his patience and graciousness in taking the time to carefully listen and read my often longwinded and laborious arguments. Thirdly I would like to thank John Baez, Aaron Lauda, and Michael Müger : John for his inspirational 'This weeks finds' series and his 'quantum gravity seminar' lecture notes; both John and Aaron for their numerous patient answers to my questions; and Michael for his help in exorcising my demons with respect to orientations and cobordisms. Fourthly I would like to thank my friends : Sheer El-Showk for his kind encouragement and discussions, Izak Snyman for being an attentive and friendly 'sounding-board' (!), Ygor Geurts for pointing out mistakes and always being encouraging and helpful, my room 301 friends Wessel, Gerben, Tristan, Gijs, Sasha and Jaap for their camaraderie and humour, Lotte Hollands for answering my questions about topological string theory, and finally my 'Cambridgelaan colleagues' Banban, Richard, Dima, Dimos, Miles, Ari and Wouter for their endless friendship and encouragement.

# Contents

# Chapter 1

# Introduction

One of the great themes of contemporary theoretical physics is the quest to combine general relativity and quantum mechanics. This endeavour is beset with a multitude of well-known profound and vexing difficulties. Many of these difficulties may be traced back to the venerable old puzzles and paradoxes which have plagued quantum mechanics since its inception, and have stubbornly persisted with us for almost a century. Although none of the current approaches to quantum gravity have even come close to attaining a consistent and physically sensible final theory, one feels that there is, surely, some element of truth in all of them, and that they have brought us closer to our goal, if not in foundational understanding, then at least in developing the appropriate tools. Reading papers about string theory, loop quantum gravity, supergravity, and so on inevitably reminds one of Alice's experience of reading 'Jabberwocky' for the first time:

*'Twas brillig, and the slithy toves*
*Did gyre and gamble in the wade.*
*All mimsy were the borogroves,*
*And the mome raths outgrabe.*

> *"It seems very pretty," she said when she had finished it, "but its* rather *hard to understand!" (You see she didn't like to confess, even to herself, that she couldn't make it out at all.) "Somehow it seems to fill my head with ideas - only I don't exactly know what they are!"*

- taken from ''Through the Looking-Glass' by Lewis Carroll, 1896.

This thesis adopts a new and subtle philosophy towards quantum gravity, and quantum field theory in general [9]. In this approach it is believed that the mathematical subject of *category theory* could shed much light on modern theoretical physics, and especially on the quantum world. This may sound more plausible if one remembers that category theory achieved precisely this in the mathematical world when it illuminated large tracts of modern mathematics. One is struck by the similarity between Alice's experience, and the following quotation about category theory [59]:

> *There are some ideas you simply could not think without a vocabulary to think them. And the language that they introduced made huge swaths of modern mathematics possible.*

The fundamental tenet of this new philosophy is that *quantum theory will make more sense when regarded as a theory of spacetime, and that one can only see this from a category-theory perspective [10]* - in particular, one that de-emphasizes the primary role of the category of sets and functions. We shall shortly attempt to explain what this all means.

This approach has been inspired by the work of mathematicians such as Segal, Atiyah, Baez and Freed, and of theoretical physicists such as Witten, Dijkgraaf, Fuchs and Moore (to name only a rough selection of those whose ideas feature in this thesis), who have exposed deep relationships between topology, quantum field theory, general relativity and category theory. These are, of course, grandiose and ambitious claims, so let us therefore explain the particular type of relationship between the first pair of topics, topology and quantum field theory, with which we shall be concerned in this work. Soon it will become clear how general relativity and category theory naturally enter the picture!

## 1.1 Chern-Simons theory

Consider then the following path integral, which serves as an excellent motivating example and in addition might serve to set at ease those minds which prefer to encounter a few concrete and familiar symbols before embarking into more abstract terrain:

$$J_k(K) = \langle K \rangle_k \equiv \int_{\mathcal{A}} \mathcal{D}A \, \mathrm{Tr}\left(P \exp \oint_K A\, ds\right) \exp \frac{ik}{4\pi} \int_{S^3} \mathrm{Tr}(A \wedge dA + \frac{2}{3} A \wedge A \wedge A) \tag{1.1}$$

The sequence of symbols appearing in (1.1) was first written down in 1989 by Witten [95] in a paper entitled 'Quantum Field Theory and the Jones Polynomial', and is one of the ideas which led to him receiving the Fields Medal for mathematics in 1990. It marked more or less the beginning of the subject, 'Topological Quantum Field Theory'. What is going on in Eqn. (1.1)? Simply put, it is a quantum field theory interpretation of the Jones polynomial invariant $J(K)$ for knots $K$ embedded in $S^3$, which presents the value of the Jones polynomial as the expectation value of the holonomy of the connection around the knot (the Wilson loop operator), where the action is the Chern-Simons action,

$$S[A] = \exp \frac{1}{4\pi} \int_{S^3} \mathrm{Tr}(A \wedge dA + \frac{2}{3} A \wedge A \wedge A). \tag{1.2}$$

Here $A$ is the pullback of a SU(2)-valued connection on a trivial principal SU(2)-bundle over $S^3$ - a $\mathfrak{su}(2)$ valued 1-form on $S^3$ - or in physics terms, a SU(2) gauge field. $k$ is an integer, which we translate into a root of unity $q$ in the complex plane by setting $q^{\frac{1}{2}} = \exp(\frac{\pi i}{k+2})$. The trace is taken in the fundamental representation of $\mathfrak{su}(2)$. Once the trace has been taken, we are left with an ordinary 3-form on $S^3$, that is, a volume form. In other words, the Chern-Simons action (unlike the Yang-Mills action, for example) supplies its own volume form and is defined without reference to a metric. In co-ordinates the action would take the following form,

$$S[A] = \exp \frac{1}{4\pi} \int_{S^3} d^3x \epsilon^{\mu\nu\alpha} \mathrm{Tr}\left[A_\mu(x)\partial_\nu A_\alpha(x) + \frac{2}{3} A_\mu(x) A_\nu(x) A_\alpha(x)\right]. \tag{1.3}$$

where the $A_\mu$ are matrix valued vector fields, $A_\mu(x) = A^\alpha_\mu(x)\sigma^\alpha/2i$.

Eqn. (1.1) may be placed into context as follows [9]. In 1984 the mathematician Vaughan Jones announced the discovery of a new link invariant, which soon led to a bewildering profusion of generalizations. Given a knot $K$ embedded in $S^3$, that is, a smooth embedding of $S^1$ into $S^3$, the Jones polynomial outputs a polynomial in $q^{\frac{1}{2}}$. Two knots cannot be smoothly deformed into each other (that is, they are not equivalent) if their Jones polynomials differ, although the converse is not true. For this breakthrough, Jones was also awarded the Fields Medal, together with Witten, in 1990.

These polynomials have nothing *a priori* to do with geometry or quantum field theory. Jones discovered them in the course of some investigations into subfactor theory, the study of how various Von-Neumann algebras can fit together. Soon it was clear, however, that these new invariants were intimately related to conformal field theory. Atiyah, [7] however, conjectured that there should be an intrinsically 3-dimensional definition of these invariants using gauge theory. Witten's Eqn. (1.1) confirmed this conjecture. It is remarkable since it amounts to the statement that one may compute the Jones polynomial by considering all connections $A$ on $S^3$: For each connection one computes a certain number (the holonomy around the knot or Wilson loop), and then one sums over all connections. Since the action is dependent on the level $k$, the expectation values output polynomials in $q^{\frac{1}{2}} = \exp(\frac{\pi i}{k+2})$ - these are precisely the Jones polynomials. This was Witten's achievement : a beautiful geometric, field theory interpretation of an object (the Jones polynomial) which manifestly concerned itself with geometric structures (knots in $S^3$) yet had a purely algebraic definition. One might say that Jones *discovered* his polynomial, while Witten *explained* it.

We shall have far more to say about Eqn. (1.1) in later chapters. For now, we should content ourselves with observing a few of its salient features. Firstly note that *it is the quantum nature of Witten's path integral representation of the Jones polynomial - the weighted sum over* all *holonomies around the knot - that has produced the interesting mathematics.* That is, once again it is Feynman's ubiquitous path integral which forms the spiritual background of the entire field of topological quantum field theory. Since the path integral has stubbornly resisted all attempts to be made completely rigorous, we can unfortunately, at this point in mathematical history, view beautiful formulas such as Eqn. (1.1) only heuristically. Mathematicians *have* however, found settings in which Witten's argument can be made rigorous - not via path integrals, but via quantum groups and modular categories. We shall encounter these in the chapters to come. Sadly, none are as charming as Witten's simple formula, but they are profound and elegant in their own right, and also reveal important concepts that are not manifest in Eqn (1.1).

## 1.2 Functorial view of quantum field theory

The fundamental importance of the path integral suggests that it might be enlightening to simplify things somewhat by stripping away the knot observable $K$ and studying only the bare partition functions of the theory, considered over arbitrary spacetimes. That is, consider the path integral

$$Z(M) = \int \mathcal{D}A \exp\left(i \int_M S[A]\right) \tag{1.4}$$

where $M$ is an arbitrary closed 3d manifold, that is, compact and without boundary, and $S[A]$ is the Chern-Simons action (1.2). Immediately one is struck by the fact that, since the action is topological, the number $Z(M)$ associated to $M$ should be a topological invariant of $M$. This is a remarkably efficient way to produce topological invariants! Recall that, unlike closed 2d manifolds which are classified by their genus, the classification of closed 3d manifolds is highly non-trivial, and is still an unsolved problem. We are reminded of one of the Clay Institute's 'Millennium Problems' [20], that is, the famous conjecture of Poincaré that there are no 'fake spheres':

**1.2.1 Poincaré Conjecture.** *If $M$ is a closed 3-manifold, whose fundamental group $\pi_1(M)$, and all of whose homology groups $H_i(M)$ are equal to those of $S^3$, then $M$ is homeomorphic to $S^3$.*

One therefore appreciates the simplicity of the quantum field theory approach to topological invariants, which runs as follows.

1. Endow the space with extra geometric structure in the form of a connection (alternatively a field, a section of a line bundle, an embedding map into spacetime, ...).
2. Compute a number from this manifold-with-connection (the action).
3. Sum over all connections.

This may be viewed as an extension of the general principle in mathematics that one should classify structures by the various kinds of extra structure that can live on them. Indeed, the Chern-Simons lagrangian was originally introduced in mathematics in precisely this way. Chern-Weil theory provides access to the cohomology groups (that is, topological invariants) of a manifold $M$ by introducing an arbitrary connection $A$ on $M$, and then associating to $A$ a closed form $f(A)$ (for instance, via the Chern-Simons lagrangian), whose cohomology class is, remarkably, independent of the original arbitrary choice of connection $A$. Quantum field theory takes this approach to the extreme by being far more ambitious; it associates to a connection $A$ the actual numerical value of the action (usually obtained by integration over $M$) - this number certainly depends on the connection, but field theory atones for this by summing over all connections.

Quantum field theory is however, in its path integral manifestation, far more than a mere machine for computing numbers associated with manifolds. There is *dynamics* involved, for the natural purpose of path integrals is not to calculate bare partition functions such as (1.4), but rather to express the *probability amplitude for a given field configuration to evolve into another.* Thus one considers a 3d manifold $M$ (spacetime) with boundary components $\Sigma_1$ and $\Sigma_2$ (space), and considers $M$ as the evolution of space from its initial configuration $\Sigma_1$ to its final configuration $\Sigma_2$:

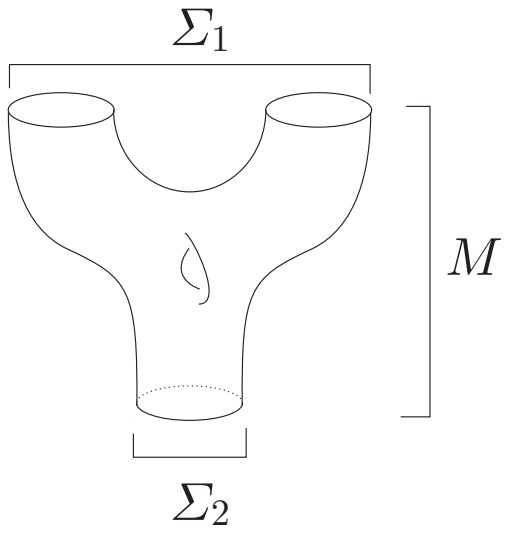

This is known mathematically as a *cobordism* from $\Sigma_1$ to $\Sigma_2$. To a 2d closed manifold $\Sigma$ we associate the space of fields $\mathcal{A}(\Sigma)$ (in this case, connections) living on $\Sigma$. A physical state $\Psi$ corresponds to a functional on this space of fields. This is the Schrodinger picture of quantum field theory: if $A \in \mathcal{A}(\Sigma)$, then $\Psi(A)$ represents the probability that the state known as $\Psi$ will be found in the field $A$. Such a state evolves with time due to the dynamics of the theory; $\Psi(A) \to \Psi(A,t)$. The space of states has a natural basis, which consists of the delta functionals $\hat{A}$ - these are the states satisfying $\langle \hat{A} | \hat{A}' \rangle = \delta(A - A')$. Any arbitrary state $\Psi$ may be expressed as a superposition of these basis states. The path integral instructs us how to compute the time evolution of states, by first expanding them in the $\hat{A}$ basis, and then specifying that the amplitude for a system in the state $\hat{A}_1$ on the space $\Sigma_1$ to be found in the state $\hat{A}_2$ on the space $\Sigma_2$ is given by:

$$\boxed{\langle \hat{A}_2 | U | \hat{A}_1 \rangle = \int_{A|_{\Sigma_1} = A_1}^{A|_{\Sigma_2} = A_2} \mathcal{D}\mathcal{A} \exp iS[A]} \tag{1.5}$$

This equation is the fundamental formula of quantum field theory, and is the basis for most of what is written in this thesis. It might not be out of place therefore to repeat its instructions in words : 'Perform a weighted sum over all possible fields (connections) living on spacetime that restrict to $A_1$ and $A_2$ on $\Sigma_1$ and $\Sigma_2$ respectively'. This formula constructs the *time evolution operator* $U$ associated to the cobordism $M$.

In this way we see that, at the very heart of quantum mechanics and quantum field theory, is a formula which associates to every space-like manifold $\Sigma$ a Hilbert space of fields $\mathcal{A}(\Sigma)$, and to every cobordism $M$ from $\Sigma_1$ to $\Sigma_2$ a time evolution operator $U(M) : \Sigma_1 \longrightarrow \Sigma_2$. To specify a quantum field theory is nothing more than to give rules for constructing the Hilbert spaces $\mathcal{A}(\Sigma)$ and the rules (correlation functions) for calculating the time evolution operators $U(M)$. *This is precisely the statement that a quantum field theory is a functor from the cobordism category* ***nCob*** *to the category of Hilbert spaces* ***Hilb*** *!* Before we recall precisely what a 'category' and a 'functor' are, it should be pointed out that the latter observation was the fundamental insight of Segal [80] in his paper 'The definition of Conformal Field Theory'.

Some readers may need more convincing at this point and demand to know the procedure of how to express, say, ordinary four-dimensional $\phi^4$ scalar field theory in this language. The answer is simple : the 'space-like' manifolds $\Sigma_t$ are the time-slices $\mathbb{R}^3 \times t$ , and the spacetimes $M_{[t_1,t_2]}$ are simply segments of Minkowski space, $M_{[t_1,t_2]} = \mathbb{R}^3 \times [t_1, t_2]$. $\phi^4$ theory is nothing more than the process of associating, to each $\Sigma_t$, the space of functionals on the fields living on $\Sigma_t$, and to each spacetime $M_{[t_1,t_2]}$ the time evolution operator $U(t_1, t_2)$. Where has the topology gone, one asks? This is buried in the assumption that $\phi^4$ theory takes place in a flat, topologically trivial background space. A general quantum field theory, such as a theory of quantum gravity, obviously should not make this assumption. Of course, in what follows it is normally assumed that both space and spacetime are compact. This is a well-known procedure, in all approaches to quantum field theory, since it provides a few technical advantages.

Appendix A comprises a quick-fire introduction to the language of category theory, but let us recall for the sake of the reader, the definition. A *category* $\mathcal{C}$ consists of a collection of objects, a collection of arrows $f : a \longrightarrow b$ from any object $a$ to any object $b$, a rule for composing arrows $f : a \longrightarrow b$ and $g : b \longrightarrow c$ to obtain an arrow $gf : a \longrightarrow c$, and for each object A an identity arrow $1_a : a \longrightarrow a$. These must satisfy the associative law $f(gh) = (fg)h$ and the left and right unit laws $1_a f = f$ and $f 1_a = f$ whenever these composites are defined. In many cases, the objects of a category are best thought of as *sets equipped with extra structure*, while the morphisms are *functions preserving the structure*. However, this is neither true for the category of Hilbert spaces nor for the category of cobordisms!

The fundamental idea of category theory is to consider the 'external' structure of the arrows between objects instead of the 'internal' structure of the objects themselves - that is, the actual elements inside an object - if indeed, an object is a set at all : it need not be, since category theory waives its right to ask questions about what is *inside* an object, but reserves its right to ask how one object is *related to* another.

A *functor* $F : \mathcal{C} \longrightarrow \mathcal{D}$ from a category $\mathcal{C}$ to another category $\mathcal{D}$ is a rule which associates to each object $a$ of $\mathcal{C}$ an object $b$ of $\mathcal{D}$, and to each arrow $f : a \longrightarrow b$ in $\mathcal{C}$ a corresponding arrow $F(f) : F(a) \longrightarrow F(b)$ in $\mathcal{D}$. This association must preserve composition and the units, that is, $F(fg) = F(f)F(g)$ and $F(1_a) = 1_{F(a)}$.

Armed with this definition, let us quickly consider the three examples of categories which are relevant to this introduction.

1. **Set** is the category whose objects are sets, and whose arrows are the functions from one set to another.
2. **nCob** is the category whose objects are closed $(n-1)$-dimensional manifolds $\Sigma$, and whose arrows $M : \Sigma_1 \longrightarrow \Sigma_2$ are cobordisms, that is, $n$-dimensional manifolds having an input boundary $\Sigma_1$ and an output boundary $\Sigma_2$. [1]
3. **Hilb** is the category whose objects are Hilbert spaces and whose arrows are the bounded linear operators from one Hilbert space to another.

The 'new philosophy' alluded to earlier amounts to the following observation due to John Baez [10]: The last two categories, **nCob** and **Hilb**, resemble each other far more than they do the first category, **Set**! If we loosely regard general relativity or geometry to be represented by **nCob**, and quantum mechanics to be represented by **Hilb**, then perhaps many of the difficulties in a theory of quantum gravity, and indeed in quantum mechanics itself, arise due to our silly insistence of *thinking of these categories as similar to* **Set**, *when in fact the one should be viewed in terms of the other.* That is, the notion of points and sets, while mathematically acceptable, might be highly unnatural to the subject at hand!

Indeed, in a remarkable recent expository article, Baez [10] shows that quantum chestnuts such as entanglement and the no-cloning theorem are easily understood if one looks at **Hilb** through the eyes of **nCob**, while they appear paradoxical in they language of **Set**. Moreover, a large amount of quantum information theory has also very recently been put into this framework [75], achieving many simple explanations of facts which required a long stream of linear algebra calculations in the past.

## 1.3 Historical background

We have arrived at a description of a field theory as a functor from a cobordism category to the category of Hilbert spaces. Originally it was Segal who first promoted this idea in 1989, in the case of conformal field theory. Segal created an unfinished manuscript about it which circulated through the mathematical community for years [2]. Thankfully it has finally been published [80], although still in 'unfinished form'. Segal begins the published version with the comment,

> "The manuscript that follows was written fifteen years ago. On balance, though, conformal field theory has evolved less quickly than I expected, and to my mind the difficulties that kept me from publishing the paper are still not altogether elucidated."

The cobordism category for conformal field theory is indeed quite complicated since it has far more resolution than in the topological case. Cobordisms are only considered equivalent if they are conformally equivalent, that is, if they have the same complex structure. This leads one to consider moduli spaces, etc. It was Atiyah who decided to reformulate Segal's viewpoint in terms of a precise set of axioms for a 'topological field theory' [7]. This forms the background for the material presented in the following chapters.

When one says the words 'topological quantum field theory' to an academic audience, various things are understood by various people. Most mathematicians will probably understand the term to refer to the Atiyah definitions. Physicists, on the other hand, are probably more familiar with an alternative 'working definition' of a topological field theory.

---

[1] The cautious reader is assured that this category will be more rigourously defined in Sec. 2.1.1.

[2] This despite the strong handwritten words on the front cover, "Do Not Copy"!

In physics, one distinguishes between a topological field theory of *Schwarz* or *Chern-Simons* type and of *Witten* or *cohomological* type. In Schwarz type theories the action is *explicitly* independent on the metric,

$$\frac{\delta S}{\delta g_{\mu\nu}} = 0. \tag{1.6}$$

Examples of such theories are Chern-Simons theories (as we have described it here) and $BF$ theories. The name 'Schwarz' is associated with Chern-Simons since it was he who first suggested, at least in the physics community, that the Jones polynomial may be related to Chern-Simons theory [78]. In Witten type topological field theories, the action and the stress energy tensor are BRST exact forms so that their functional averages are zero. The topological observables in these theories therefore form cohomological classes. In four dimensions, such theories involving Yang-Mills theories provide a (heuristic) field-theoretic representation for Donaldson invariants.

Unfortunately, many of these theories are not *truly* topological since, at some stage or another, either in their formulation or in the renormalization process, a hidden dependence on a metric or a complex structure creeps in. One can already see a simple case of this phenomenon with Chern-Simons theory. When one renormalizes the path integral (1.1) it is found that one must include a *framing* of the knot, something which was not immediately obvious from the action[3]. One quick-fire way to determine if a physics theory is *truly* topological is to ask whether the Hilbert space of states is finite dimensional. If this is the case, then the theory is probably (but not always) an Atiyah style TQFT. If it is not, then the theory cannot be a true topological field theory, as we shall see in Chapter 3.

## 1.4 Layout of this thesis

This thesis runs as follows. In chapter 2, we introduce the language of monoidal categories, and a graphical notation for doing calculations with them. This graphical calculus is used as a guiding principle in order to determine the necessary additional structure on the monoidal categories we shall require. In this way duality, braiding and twist operations on monoidal categories are introduced. We end the chapter by applying this language to a familiar physical context : feynman diagrams.

In chapter 3, we introduce the formal definition of topological quantum field theory, and discuss some general properties of these theories in all dimensions. Then we specialize to two dimensions, and offer three different proofs of the result that such theories are classified by commutative Frobenius algebras. Next we consider adding 'open strings' to our framework, and discuss the work of Segal and Moore which classifies D-branes in this heavily simplified topological context.

Chapter 4 is about actual physical models which are topological field theories in two dimensions. We begin by looking at the language of gauge theory from a functorial perspective, and specialize this to the case when the 2d surfaces are equipped with a triangulation. We then introduce the Dijkgraaf-Witten finite group model, which is valid in any number of dimensions, and which forms a central theme of the thesis. We solve the model in two dimensions from first principles, as well as from triangulations. Finally we consider Yang-Mills theory in two dimensions, which can be regarded as a direct generalization of the model to the case when the group $G$ is continuous.

In chapter 5, we study three dimensional TQFT's. We begin by extending the language of ribbon categories introduced in chapter 2 in order to define modular categories. We give a flavour of how calculations work in these categories, and prove the Verlinde formula which is valid in all modular categories. Next we describe how to construct a 3d TQFT armed with the data of a modular category, using surgery on links. We also show an algorithm for producing these categories from ordinary monoidal categories, called the quantum double construction. Next, we outline some fundamental ideas of Daniel Freed which unifies the approach to TQFT via modular categories (eg. quantum groups) and the approach via path integrals (eg. the Chern-Simons action). In this framework, the modular category is actually *produced* from a path integral. We then illustrate these ideas with the Dijkgraaf-Witten model in three dimensions, showing how the modular category is derived. Remarkably, the passage from the 2d theory to the 3d theory corresponds to taking the quantum double of the category of representations of the group. Finally we show explicitly the modular data arising from the representation theory of the associated quantum group. This allows us to make contact with all the work from the previous chapters, by comparing the dimensions of certain Hilbert spaces associated with this modular category, computed via the Verlinde formula, with similar formulas obtained in the 2d case in chapter 4.

[3]The need for a framing arises because divergent integrals are produced from the perturbation expansion of equations such as 1.1. It is remarkable that in the mathematical framework of modular categories, one also requires a framing, but for apparently unrelated reasons.

## 1.5 Important topics left out

Due to lack of space and insufficient expertise, many important topics specifically related to the ideas in this thesis have been left out, or only touched on briefly. These include Tillman's redefinition of a modular category a 'categorified Frobenius algebra' - the image of the circle under a 2-functor $Z_2 : \mathbf{2Cob_2} \rightarrow \mathbf{Hilb_2}$ [88]; Lauda's development of this idea for open strings in the plane [1, 2]; Baez's method of obtaining 3d TQFT's from the categorification of the Fukama-Hosono-Kawai method for 2d TQFT's [15]; Kauffman's development of Chern-Simons field theory calculations and a geometric explanation for why the Chern-Simons action gives knot invariants [47]; Turaev's definition of a homotopy quantum field theory as an embedded version of a topological quantum field theory in an auxiliary space $X$ and his explicit examples of such theories in dimensions 2 and 3 [91, 92]; Bunke, Turner and Willerton's ideas on the relationship between 2 dimensional homotopy quantum field theories and gerbes[4] [100]; Baez and his collaborators' work on higher gauge theory expressed in terms of higher category theory [16]; Moore and Seiberg's method of viewing the data in a modular category via its Hom spaces [64, 18] and Fuchs, Runkel and Schweigert's framework for expressing boundary conditions in conformal field theory as Frobenius algebras inside the modular category associated to the chiral data[5] [35, 36, 37, 38, 39].

[4]Gerbes describe parallel transport of *strings* around a manifold, so they are a higher dimensional version of a connection. The $B$ fields in the $RR$ sector of string theory appear to be candidates for such a description [82].

[5]See Sec. 5.1.1.

# Chapter 2

# Categories, Cobordisms and Feynman Diagrams

## 2.1 Categories

A *category*[1] $\mathcal{C}$ consists of a collection of objects, $\mathrm{Ob}(\mathcal{C})$, and a collection of arrows (also called morphisms), $\mathrm{Ar}(\mathcal{C})$. Each arrow $f$ has a source object and a target object which is denoted as $f : x \to y$. Besides the collection of objects and arrows, there is given a rule for composing arrows; that is, if $f : x \to y$ and $g : y \to z$ then there must be an arrow $g \circ f : x \to z$, and this composition is required to be associative, $f \circ (g \circ h) = (f \circ g) \circ h$. Finally, each object $x$ comes equipped with at least an identity arrow $\mathrm{id} : x \to x$ such that $\mathrm{id} \circ f = f \circ \mathrm{id} = f$ whenever it makes sense. All the information in a category is contained in the data (objects, arrows, composition rules of arrows).

A category can be viewed as a special kind of directed graph, where one has the ability to 'compose' the arrows:

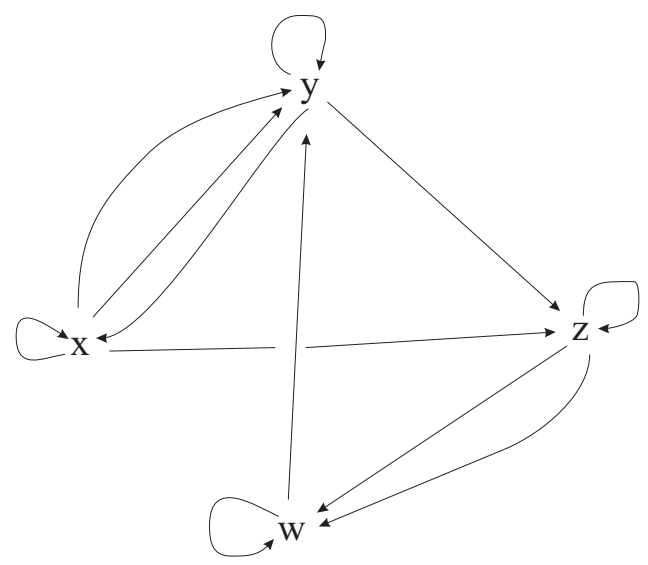

This theme of representing categories graphically is very useful and is the subject of this chapter. Physicists should find themselves at home with this notation - they are nothing but Feynman diagrams!

The archetypal category is **Set** which has sets $X$ for objects, functions between sets $f : X \to Y$ for arrows, and ordinary function composition for composition of arrows. This example can be misleading because the only information owned by an object $x$ in a category is its label '$x$'. Objects have no substructure and it does not make sense to ask whether a certain element is 'contained' in $x$. Thus one should view the objects (sets) in **Set** merely as labels, and forget any elements they might have contained. The principle to remember is that *all the information resides in the arrows* - in the case of **Set**, one can recover the elements from the arrows by noticing that each arrow $1 \xrightarrow{f} X$ from some fixed one element set 1 to an arbitrary set $X$ defines an 'element' of $X$.

[1]The mysterious word *collection* in the definition above should be understood as a *class*, that is, a weaker notion of set capable of dealing sensibly with notions like 'The class of all sets'. This leads one to distinguish *small* categories, where the class of objects actually do form a set, from *large* categories, where they do not. These distinctions are important from a viewpoint of logic but are not relevant(yet!) for physics.

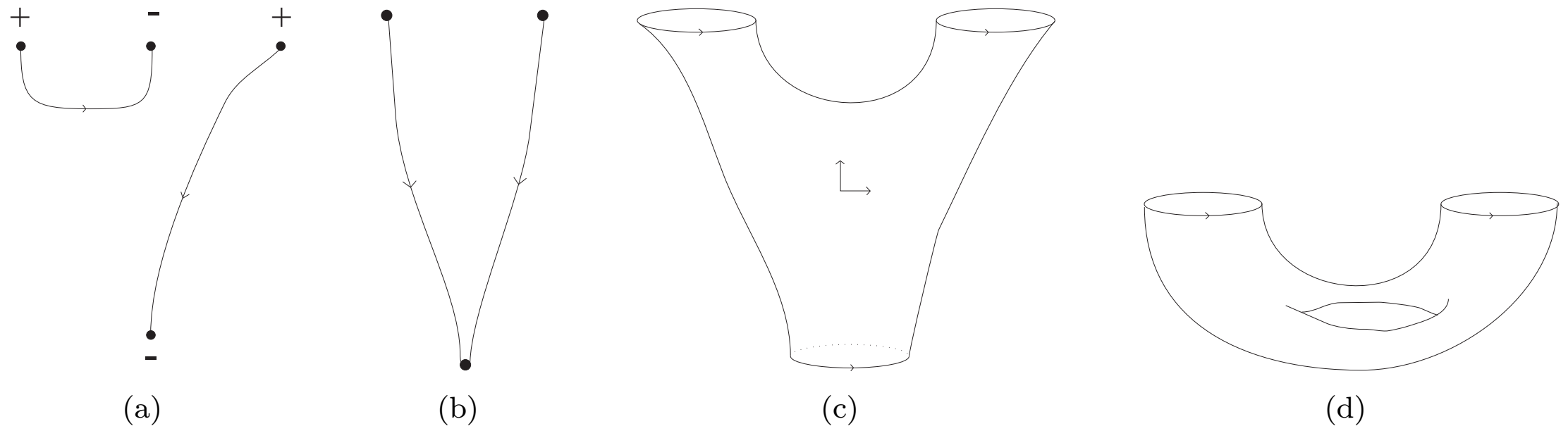

Figure 2.1: (a) A cobordism from $(+,-,+)$ to $(-)$. (b) This is not a cobordism since it is not a valid 1-manifold with boundary. (c) The 'pair of pants' - a cobordism from two circles to one circle. (d) A cobordism from two circles to the empty set.

### 2.1.1 Cobordisms

However, objects and arrows in general categories can be very different from sets and functions. The category which quantum field theory concerns itself with is called **nCob**, the '$n$ dimensional cobordism category', and the rough definition is as follows. Objects are oriented closed (that is, compact and without boundary) $(n-1)$-manifolds $\Sigma$, and arrows $M : \Sigma \to \Sigma'$ are compact oriented $n$-manifolds $M$ which are *cobordisms* from $\Sigma$ to $\Sigma'$. Composition of cobordisms $M : \Sigma \to \Sigma'$ and $N : \Sigma' \to \Sigma''$ is defined by gluing $M$ to $N$ along $\Sigma'$.

Let us fill in the details of this definition. Let $M$ be an oriented $n$-manifold with boundary $\partial M$. Then one assigns an induced orientation to the connected components $\Sigma$ of $\partial M$ by the following procedure. For $x \in \Sigma$, let $(v_1, \ldots, v_{n-1}, v_n)$ be a positive basis for $T_x M$ chosen in such a way that $(v_1, \ldots, v_{n-1}) \in T_x\Sigma$. It makes sense to ask whether $v_n$ points inward or outward from $M$. If it points inward, then an orientation for $\Sigma$ is defined by specifying that $(v_1, \ldots, v_{n-1})$ is a positive basis for $T_x\Sigma$. If $M$ is one dimensional, then $x \in \partial M$ is defined to have positive orientation if a positive vector in $T_x M$ points into $M$, otherwise it is defined to have negative orientation.

Let $\Sigma$ and $\Sigma'$ be closed oriented $(n-1)$-manifolds. An *cobordism* from $\Sigma$ to $\Sigma'$ is a compact oriented $n$-manifold $M$ together with smooth maps

$$\Sigma \xrightarrow{i} M \xleftarrow{i'} \Sigma' \tag{2.1}$$

where $i$ is a orientation preserving diffeomorphism of $\Sigma$ onto $i(\Sigma) \subset \partial M$, $i'$ is an orientation reversing diffeomorphism of $\Sigma'$ onto $i'(\Sigma') \subset \partial M$, such that $i(\Sigma)$ and $i'(\Sigma')$ (called the in- and out-boundaries respectively) are disjoint and exhaust $\partial M$. Observe that the empty set $\phi$ can be considered as an $(n-1)$-manifold.

Examples and non-examples of cobordisms are shown in Fig. 2.1. One can view $M$ as *interpolating* from $\Sigma$ to $\Sigma'$. An important property of cobordisms is that they can be glued together. Let $M : \Sigma_0 \to \Sigma_1$ and $M' : \Sigma_1 \to \Sigma_2$ be cobordisms,

$$\Sigma_0 \xrightarrow{i_0} M \xleftarrow{i_1} \Sigma_1, \qquad \Sigma_1 \xrightarrow{i'_1} M' \xleftarrow{i_2} \Sigma_2. \tag{2.2}$$

Then we can form a composite cobordism $M' \circ M : \Sigma_0 \to \Sigma_2$ by gluing $M$ to $M'$ using $i'_1 \circ i_1^{-1} : \partial M \to \partial M'$:

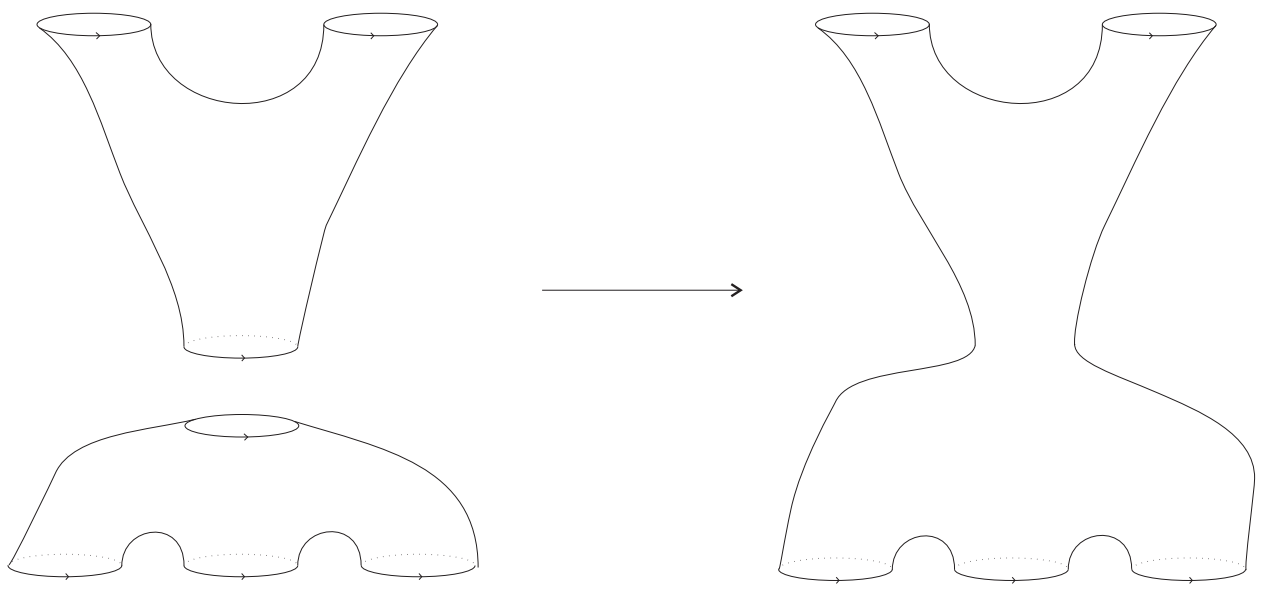

### 2.1.2 Further examples of categories

Having defined cobordisms, we can now list some important examples of categories which we shall encounter again and again in this thesis. Each of them in fact possess additional structure, with which they shall be equipped in stages over the course of this chapter.

a) Given a mathematical structure on a set, one can always consider the category where objects are sets possessing this structure, and arrows are functions preserving the structure. For instance, we have the categories **Grp** (objects are groups, arrows are group homomorphisms), *Ab* (objects are abelian groups, arrows are group homomorphisms), **Vect** (objects are vector spaces, arrows are linear maps), **FinVect** (objects are finite dimensional vector spaces, arrows are linear maps), ....

b) The category **Hilb**. Objects are Hilbert spaces, and arrows are bounded linear maps. This does not quite fit into the previous series since arrows need not be unitary. Nevertheless we shall see how choosing general bounded linear maps is a more natural choice.

c) Given a fixed set-with-structure $X$, one can consider the category of *representations* of $X$ as linear operators on vector spaces. For instance, if $G$ is a group, then we can consider $\text{Rep}(G)$ (Objects are representations $(\rho, V)$ of $G$ where $\rho : G \to \text{End}(V)$ is a representation on some vector space $V$. Arrows $f : (\rho, V) \to (\rho', W)$ are intertwining maps, i.e. linear maps $V \to W$ such that $\rho' f = f\rho$). Similarly, one can consider $\text{Rep}(A)$ where $A$ is an algebra, $\text{Rep}(\mathfrak{g})$ where $g$ is a Lie algebra, or $\text{Rep}(U_q(\mathfrak{g}))$ where $U_q(\mathfrak{g})$ is the '$q$-deformed universal enveloping algebra of $\mathfrak{g}$'. These representation categories are central in this thesis; in fact, we shall see how the whole machinery of Feynman diagrams (such as those in QED) can be viewed as taking place inside $\text{Rep}(G)$, where $G = SO(3,1) \times U(1)$!

d) For a given $n$, one can construct the (smooth) cobordism category **nCob**. Objects are closed, oriented $(n-1)$-manifolds $\Sigma$, and arrows $M : \Sigma \to \Sigma'$ are cobordisms. In order to make this a well-defined category with identity arrows, we must quotient out diffeomorphic cobordisms. Specifically, let $M$ and $M'$ be cobordisms from $\Sigma$ to $\Sigma'$:

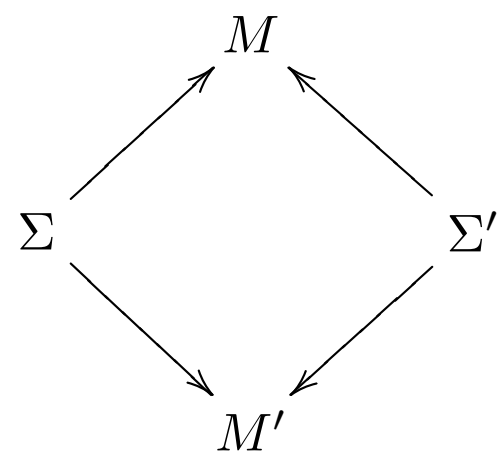

Then they are considered equivalent if there is an orientation preserving diffeomorphism $\psi : M \xrightarrow{\sim} M'$ making the following diagram commute:

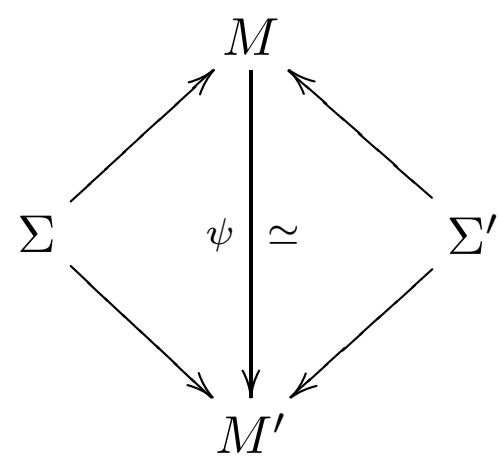

After identifying equivalent cobordisms, the 'cylinder' $\Sigma \times [0,1]$ functions as the identity arrow for $\Sigma$.

The cobordism category is a geometric category which captures the way $n$-manifolds glue together. Note that all the manifolds involved are abstract and are not embedded in some ambient space. So far we have only defined the smooth cobordism category, but there are many variations of this idea. For instance, one could use (topological, triangulated) manifolds, and identify cobordisms when they are (homeomorphic, homeomorphic preserving triangulation) - these would define the (topological, triangulated) cobordism categories. This thesis focuses on two and three dimensions, where luckily every topological manifold can be equipped with a smooth structure in essentially a unique way, and similarly triangulated and smooth manifolds are essentially equivalent, so that (smooth) **nCob**, (topological) **nCob** and (triangulated) **nCob** coincide for $n = 2, 3$. This is not the case in four dimensions, where it is well-known that some topological 4d manifolds cannot carry a smooth structure at all, while there are uncountably many different smooth structures on $R^n$, for example, and also that triangulated(also called piecewise-linear) manifolds are *not* always equivalent to smooth or even topological ones.

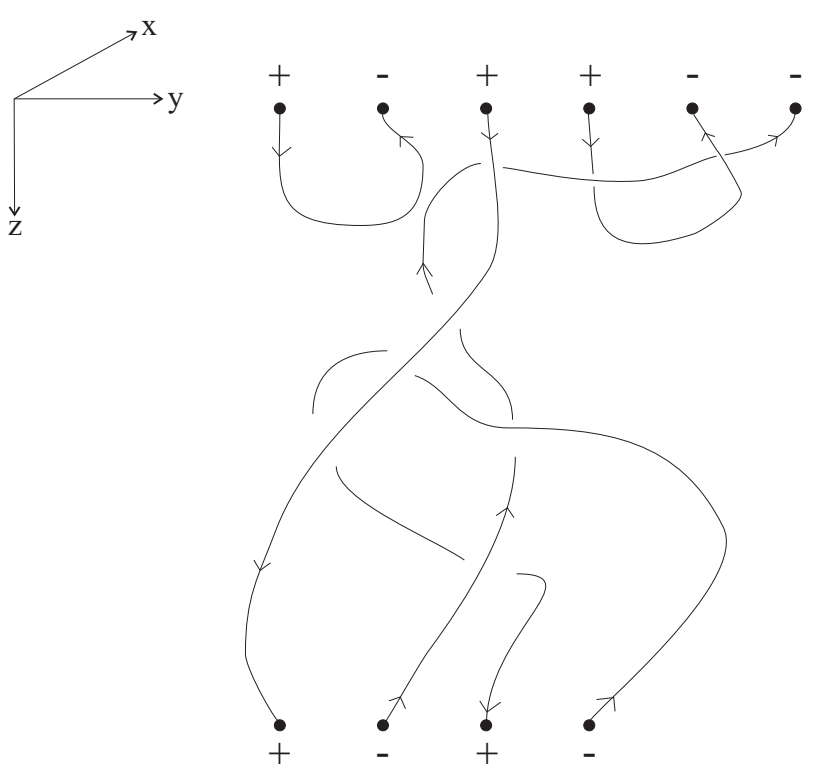

Figure 2.2: A typical morphism in the **Tangle** category. Note the orientation of the axes.

There is another cobordism category $\mathcal{M}$ which Graeme Segal introduced in the mathematical foundations of conformal field theory. $\mathcal{M}$ has the same objects as **2Cob**, except that cobordisms are now equipped with a complex structure, and are considered equivalent only up to a diffeomorphism preserving this structure. For two 1-manifolds $\Sigma, \Sigma'$, the space of all morphisms from $\Sigma$ to $\Sigma'$ has a natural smooth structure and is the moduli space familiar from conformal field theory.

e) There is a category called **Tangle** which can be viewed as **1Cob** embedded in $\mathbb{R}^3$. The objects $p$ are finite sequences of $+$ and $-$, eg. $p = (+,+,+,-,+,-,-)$. An arrow $T : p_1 \to p_2$ is an embedded oriented 1d submanifold $T \subset \mathbb{R}^2 \times [0,1]$ with $\partial T = p_1 \cup p_2$ where $p_1$ (the in boundary) and $p_2$ (the out boundary) are interpreted as sequences of $+$ or $-$ marked points living at $(0, n, 0)$ and $(0, m, 1)$ respectively. Moreover, the curves that leave or enter $+$ points are oriented downwards near the point, while curves that leave or enter $-$ signs are oriented upwards near the point. For instance, Fig. 2.2 shows a morphism from $(+,-,+,+,-,-)$ to $(+,-,+,-)$.

Morphisms are composed by gluing the out-boundaries to the in-boundaries, and rescaling the result in the $z$ direction by a factor of $\frac{1}{2}$. Morphisms $M : p_1 \to p_2$ and $M' : p_1 \to p_2$ are considered equivalent if there is an isotopy taking the one into the other, i.e. a 1-parameter smooth family of diffeomorphisms $\psi_t : \mathbb{R}^2 \times [0,1] \to \mathbb{R}^2 \times [0,1]$ such that $\psi_0 = \text{id}$ and $\psi_1(M) = M'$. The identity arrow on each object $p$ is simply a collection of straight lines:

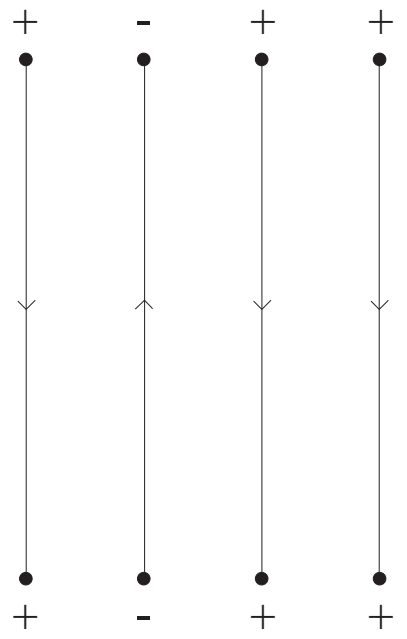

## 2.2 Monoidal categories

We have seen how categories can be interpreted in a graphical way as directed graphs with a composition operation on the arrows. We now want to take this up more seriously, using the **Tangle** category as our paradigm. Suppose we have a category $\mathcal{C}$, and we wish to draw pictures of the morphisms in $\mathcal{C}$ as a kind of enriched tangle diagram where we label the curves with objects of $\mathcal{C}$ and have allowed for the action of various morphisms $f \in \mathrm{Ar}(\mathcal{C})$ in $\mathcal{C}$ (drawn as labeled coupons) along the way:

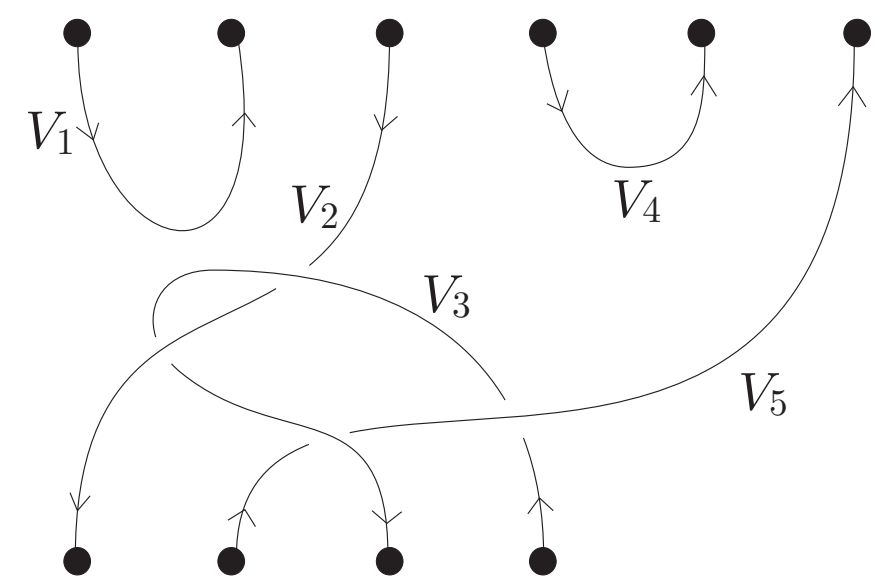

We now ask : What additional structure is required on $\mathcal{C}$ in order to make these graphical pictures precise?

Algebraic structures in mathematics are based on the idea of a binary operation on a set, that is, the sum or product of two elements living inside the set. The most general structure of this form is a monoid, that is, a set $X$ armed with an associative binary operation and a unit element. For example, groups, rings, modules and algebras are all monoids. Similarly we seek to define a binary operation on a category. This is most conveniently achieved by a process known as *categorification*, a fundamental algorithm which converts any structure on sets into a structure on categories. The idea is to express the axioms of structure at hand as a bunch of commuting diagrams of arrows inside **Set**, and then to interpret these diagrams as diagrams inside **Cat**.

In our case, this means noticing that the axioms for a monoid $X$, which would traditionally be phrased as

> For all $(a, b) \in X \times X$ there is defined a product $a \Box b$ such that $(a \Box b) \Box c = a \Box (b \Box c)$ for all $a, b, c \in X$, and there is also a unit element $1 \in X$ such that $1 \Box a = a \Box 1 = a$ for all $a \in X$.

can be equivalently expressed 'diagrammatically' as the existence of two functions,

$$\mu : X \times X \longrightarrow X, \qquad \eta : 1 \longrightarrow X \tag{2.3}$$

where 1 is a singleton set, such that the following diagrams commute:

$$\begin{array}{ccccc} & & X \times X \times X & & \\ & \swarrow^{\mu \times \mathrm{id}_X} & & \searrow^{\mathrm{id}_X \times \mu} & \\ X \times X & & & & X \times X \\ & \searrow_{\mu} & & \swarrow_{\mu} & \\ & & X & & \end{array} \tag{2.4}$$

$$\begin{array}{ccc} 1 \times X & \xrightarrow{\eta \times \mathrm{id}_X} & X \times X \\ & \searrow & \downarrow \mu \\ & & X \end{array} \qquad \begin{array}{ccc} X \times X & \xleftarrow{\mathrm{id}_X \times \eta} & X \times 1 \\ \mu \downarrow & \swarrow & \\ X & & \end{array} \tag{2.5}$$

The advantage of phrasing it in this latter way is that such diagrams make sense in any category. Hence we define a *monoidal* (or *tensor*) category as a category $\mathcal{C}$ together with a functor[2] $\otimes : \mathcal{C} \times \mathcal{C} \longrightarrow \mathcal{C}$,[3] and

[2] According to common practice, we depict the functor $\mu$ as $\otimes$

[3] The product of two categories $\mathcal{C} \times \mathcal{D}$ is defined in the cartesian way. Objects are pairs $(a_1, a_2)$ where $a_1 \in \mathcal{C}$ and $b \in \mathcal{D}$, and similarly maps are pairs $(f_1, f_2)$ where $f_1 \in \mathcal{C}$ and $f_2 \in \mathcal{D}$. Composition is defined pair-wise.

$\eta : \hat{1} \to \mathcal{C}$,[4] satisfying the associativity axiom and the unit object axiom:

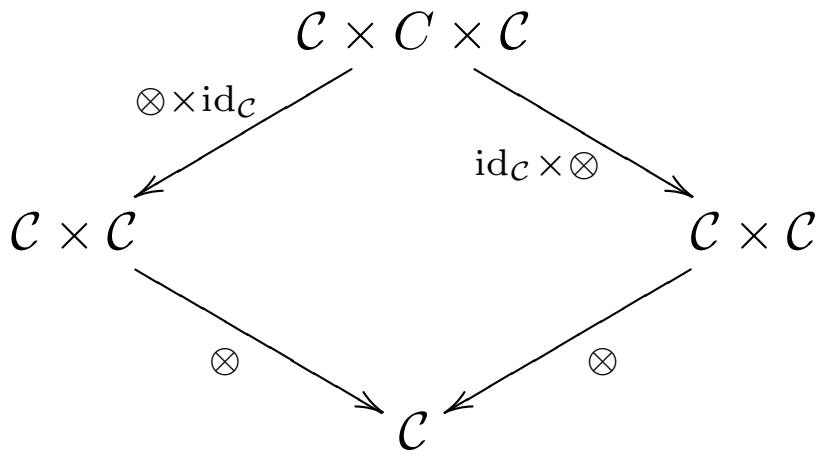

(2.6)

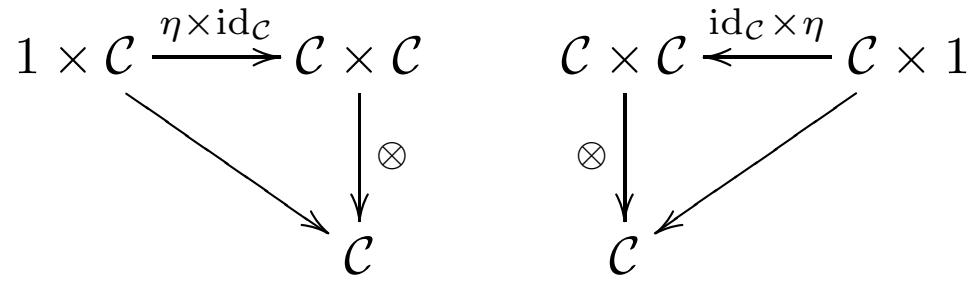

(2.7)

Let us spell it out. To each pair of objects $a, b$ of $\mathcal{C}$ is associated an object $a \otimes b$ of $\mathcal{C}$ and to each pair of morphisms $a \xrightarrow{f} a'$, $b \xrightarrow{g} b'$ a morphism $a \otimes a' \xrightarrow{f \otimes g} b \otimes b'$ [5], satisfying the following identities:

| | | |
|---|---|---|
| $(a \otimes b) \otimes c = a \otimes (b \otimes c)$ | Associativity diagram | (2.8) |
| $1 \otimes a = a = a \otimes 1$ | Unit diagram | (2.9) |
| $(f \circ f') \otimes (g \circ g') = (f \otimes g) \circ (f' \otimes g')$ | Functoriality | (2.10) |
| $\mathrm{id}_a \otimes \mathrm{id}_b = \mathrm{id}_{a \otimes b}$ | Functoriality | (2.11) |

What we are defining here is properly called a *strict* monoidal category, since $(a \otimes b) \otimes c$ is actually equal to $a \otimes (b \otimes c)$ and similarly $1 \otimes a$ and $a \otimes 1$ are actually equal to $a$. The weaker alternative (which is mathematically, and philosophically, more correct) is to only require associativity and the unit to hold up to a coherent isomorphism. 'Coherent' means that these isomorphisms should satisfy commuting diagrams of their own, in order for the whole package to be globally consistent. Luckily MacLane has proved [58] that all monoidal categories (in the weaker sense) are equivalent to strict ones[6], so that we will glibly ignore these issues here[7].

A good example of a monoidal category is **Vect**, with $\otimes$ the ordinary tensor product of vector spaces (and of linear maps). Another example is $\mathrm{Rep}(G)$, where $\otimes$ is the tensor product of representations of groups. In **nCob**, $\otimes$ is the disjoint union of manifolds.

### 2.2.1 Tannaka duality

This immediately begs the question, 'Can the group $G$ be recovered from the monoidal category $\mathrm{Rep}(G)$?'. For compact topological groups the beautiful answer is *yes*. This is known as Tannaka duality (for a good introduction, see [45]). This duality exchanges the picture of a compact Lie group as a smooth manifold with a group multiplication, for a discrete (since $G$ is compact) picture in terms of objects (reps) and arrows (intertwining maps); see Fig. 2.3.

These considerations lead one to the theory of *quantum groups* - algebras whose category of representations looks very similar to $\mathrm{Rep}_G$ [84]. As one would suspect, quantum groups are highly relevant in topological quantum field theory.

[4] The category $\hat{1}$ is the category consisting of a single object and single identity arrow

[5] We agree to denote the singleton category 1 and its image $\eta(1)$ in $\mathcal{C}$ by the same symbol 1.

[6] Indeed, the commuting diagrams which the associativity and unit isomorphisms must satisfy were designed precisely in order to ensure this.

[7] Nevertheless, the associativity and unit isomorphisms are an important and fascinating structure which reflect interesting topological properties [15].

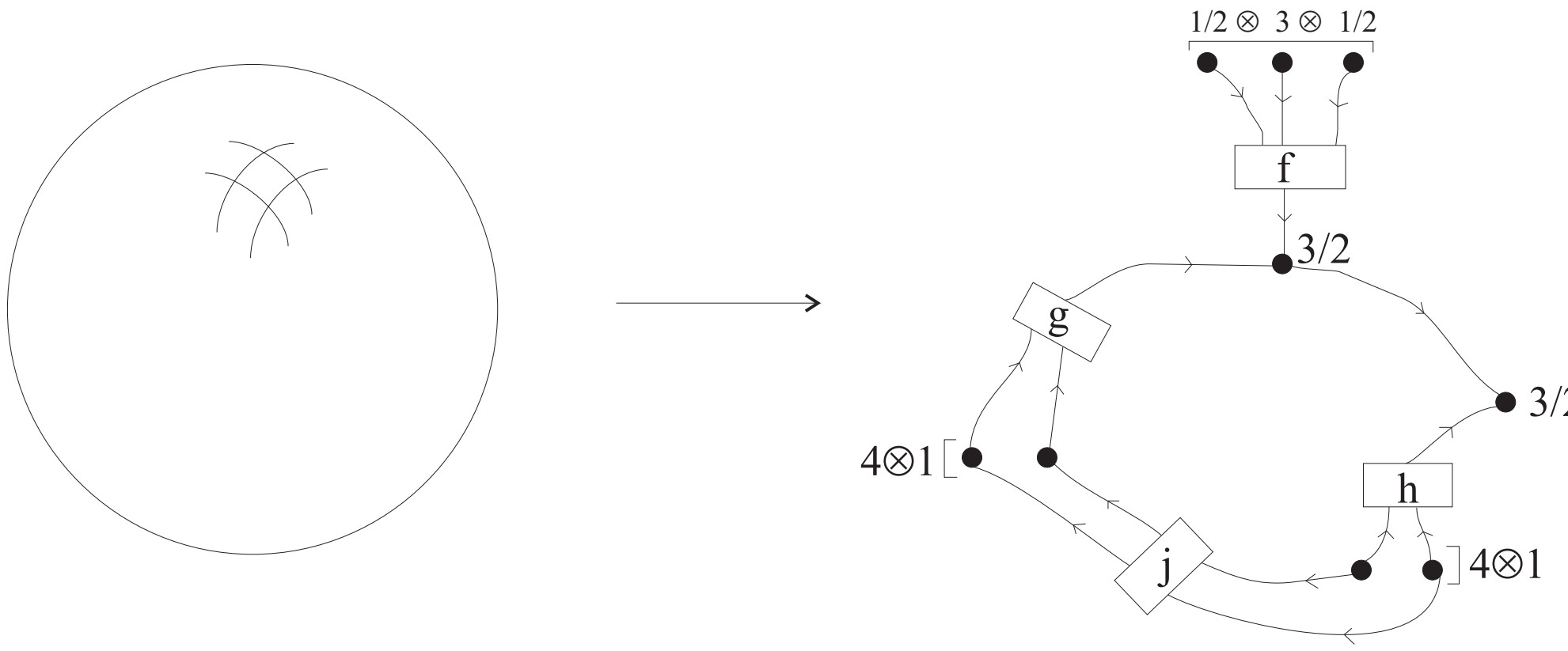

Figure 2.3: Tannaka duality turns the group $SU(2) \simeq S^3$ into its discrete category of reps, Rep($G$). Each rep is completely reducible to a direct sum of spin $\frac{j}{2}$ irreps. The morphisms $f, g, h, j$ represent intertwining maps. It is remarkable that these two pictures contain the same information.

## 2.3 Graphical calculus for monoidal categories

There is a beautifully intuitive way of drawing morphisms in monoidal categories. It can be traced back to Penrose's graphical notation for tensor calculus[8]. To draw, say, $U \otimes V \otimes W \xrightarrow{f} P \otimes Q$ simply line the source objects along the top, the target objects along the bottom, and use a rectangular 'coupon' to represent $f$:

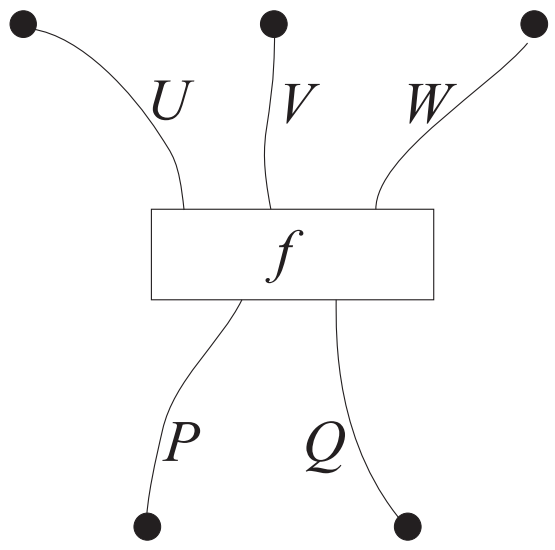

The object labels $U$, $V$, etc. are viewed as labeling the *curve*, not the end points[9]. In this notation, composition of arrows is the vertical dimension (the flow of time runs from top to bottom). For $f : U \to V$ and $g : V \to W$ we write $g \circ f$ in the following two ways:

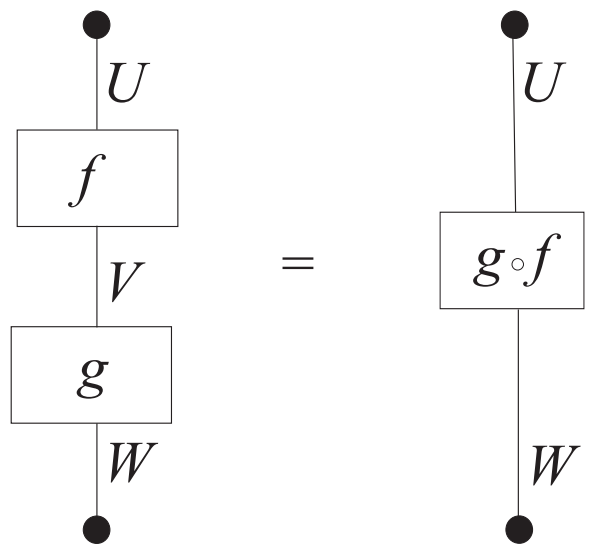

The identity arrow $\mathrm{id}_V$ on an object $V$ is drawn as:

[8]'In an appendix we outline an alternative and equivalent diagrammatic notation which is very valuable for use in private calculations' - taken from *Spinors and space-time*, R. Penrose and W. Rindler (1984).

[9]This is a clever way to accommodate duality, which we will define shortly.

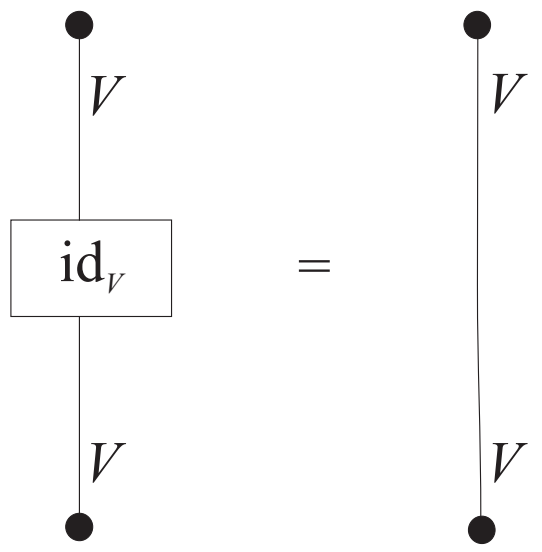

The tensor product supplies the horizontal dimension. For $f : U \to V$ and $f' : U' \to V'$ we write the tensor product $f \otimes g : U \otimes U' \to V \otimes V'$ as:

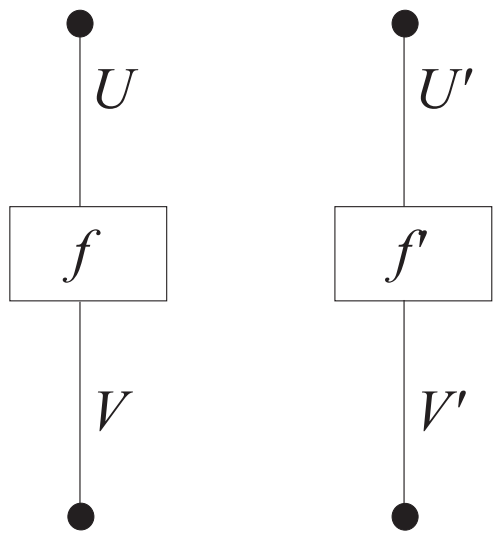

The unit object 1 is left out in drawings (in order to ensure that $V \otimes 1 = 1 \otimes V$), and morphisms such as $f : V \otimes W \to 1$ and $g : 1 \to V$ are drawn as:

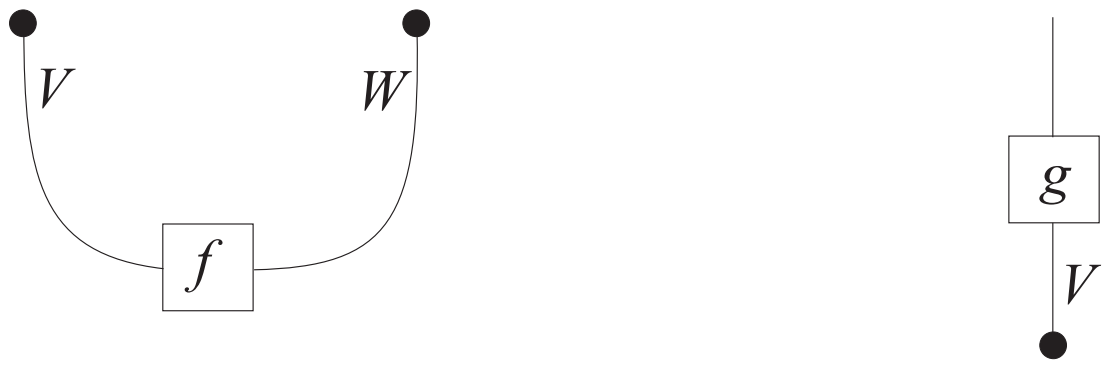

We can decompose a large complicated diagram into a composition of sections where only one morphism acts at a time, by drawing a sequence of horizontal lines ('time slices')[10] through the diagram:

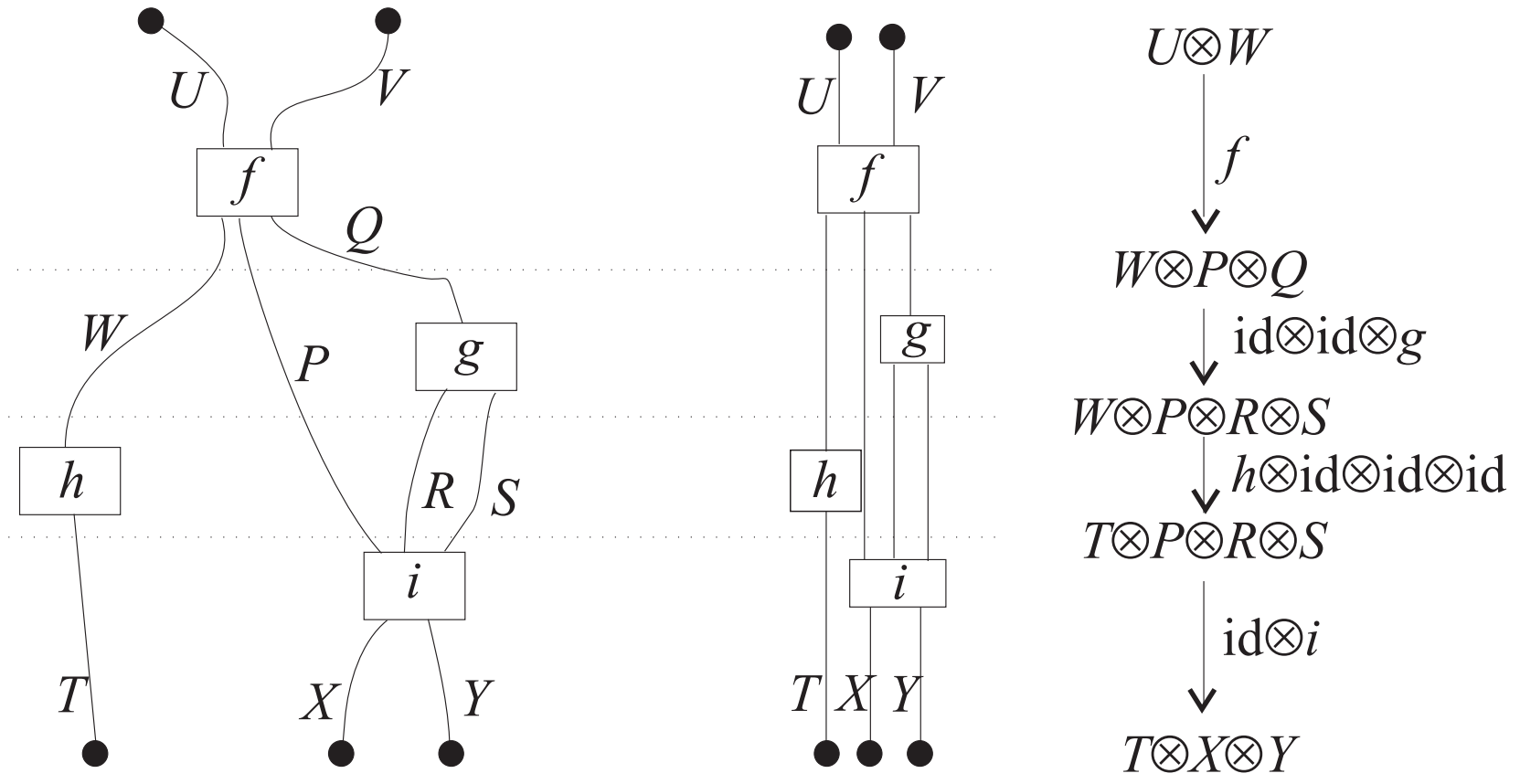

The resulting morphism obtained is in this case

$$(\mathrm{id} \otimes i) \circ (h \otimes \mathrm{id} \otimes \mathrm{id}) \circ (\mathrm{id} \otimes \mathrm{id} \otimes g) \circ f, \qquad (2.12)$$

and is called the *evaluation* of the diagram.

[10] The resemblance with Feynman's time slicing procedure for deriving the path integral is not accidental, and there are some uncanny similarities between the two approaches.

### 2.3.1 Deformation of diagrams and duality

In order to ensure that the evaluation of the diagram is invariant under an arbitrary planar deformation keeping the endpoints fixed, we need to introduce *duality* into our monoidal category. Since a monoidal category is defined so that it closely resembles **Vect** with the usual tensor product of vector spaces, one may ask what other properties of the tensor product in **Vect** can be carried over to a general monoidal category. The first step in this programme is to categorize the notion of duals of vector spaces - how can one generalize this to an arbitrary monoidal category?

Suppose that to each object $V$ of $\mathcal{C}$ there is associated an object $V^*$ of $\mathcal{C}$ and two morphisms

$$i_V : 1 \longrightarrow V \otimes V^*, \quad e_V : V^* \otimes V \longrightarrow 1, \tag{2.13}$$

such that the following two composites

$$V = 1 \otimes V \xrightarrow{i_V \otimes \mathrm{id}_V} V \otimes V^* \otimes V \xrightarrow{\mathrm{id}_V \otimes e_V} V \otimes 1 = V \tag{2.14}$$

$$V^* = V^* \otimes 1 \xrightarrow{\mathrm{id}_{V^*} \otimes i_V} V^* \otimes V \otimes V^* \xrightarrow{e_V \otimes \mathrm{id}_{V^*}} 1 \otimes V^* = V^* \tag{2.15}$$

are equal to $\mathrm{id}_V$ and $\mathrm{id}_{V^*}$ respectively. Then this structure is called a *duality* on $\mathcal{C}$. Note that, at least for finite dimensional vector spaces, this is equivalent to the usual notion of the dual space as the space of linear functionals on $V$. Indeed, if $(e_i)$ is a basis for $V$, with dual basis $(e^i)$ for $V^*$, then the maps defined by

$$i(1) = \sum_i e_i \otimes e^i, \quad e(e^i \otimes e_j) = \delta^i_j \tag{2.16}$$

satisfy (2.14) and (2.15). This is a basis-independent construction!

We incorporate duality into the graphical notation by adding arrows to the string diagrams. The general rule is, that at any horizontal slice through the diagram, downward pointing strings refer to the labeled objects as usual but upward pointing strings refer to their duals. For instance, here are four equivalent ways to depict a morphism $f : V \to W^*$:

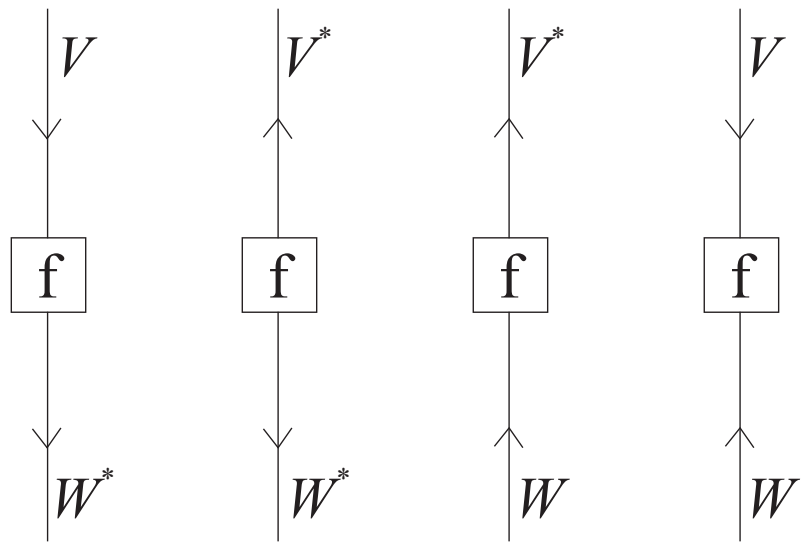

The duality morphisms $i_V$ and $e_V$ are drawn as follows:

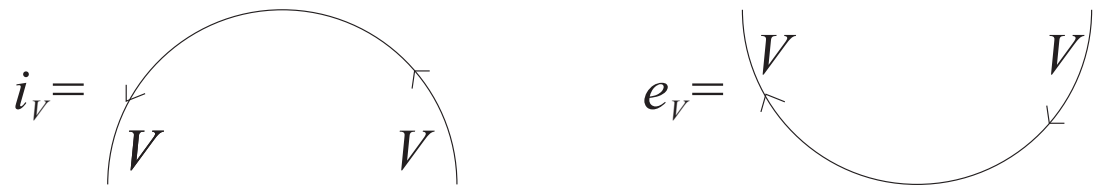

(2.17)

We now see why a duality structure is necessary in order to allow for deformations of string diagrams. Eqs. (2.14) and (2.15) tell us that 'snakes can be straightened out':

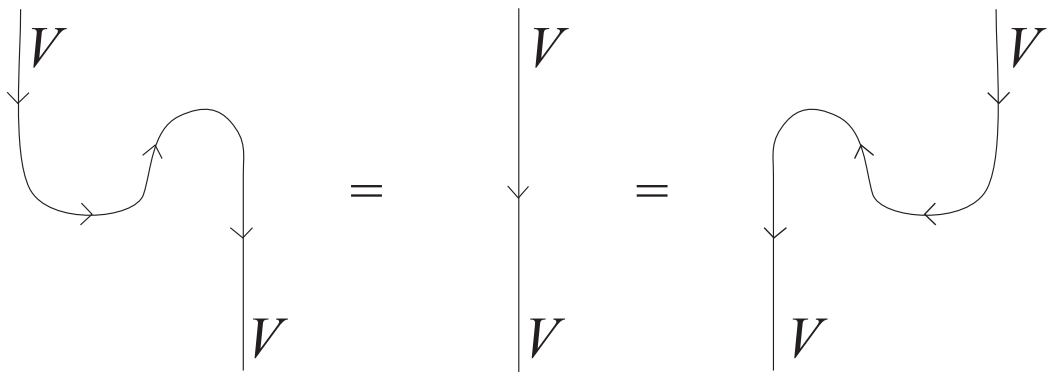

In this way, the graphical calculus is invariant under deformation. Let us record this important fact:

**2.3.1 Theorem.** **[87]** *Let $(\mathcal{C}, \otimes, 1)$ be a monoidal category with duals. Then the evaluation of a string diagram into a morphism in $\mathcal{C}$ is invariant under planar isotopy which keeps the endpoints fixed.*

To illustrate the use of the graphical calculus, consider the following eery lemma:

**2.3.2 Eckmann-Hilton Argument.** **[27]** *Let $(\mathcal{C}, \otimes, 1)$ be a monoidal category. Then the monoid $Hom_{\mathcal{C}}(1,1)$ is abelian.*

*Graphical proof.* Let us represent $f, g \in \mathrm{Hom}_{\mathcal{C}}(1,1)$ as follows:

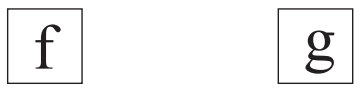

There are no endpoints, since they are morphisms from the unit object to itself. Thus planar isotopy interchanges $g \circ f$ for $f \circ g$:

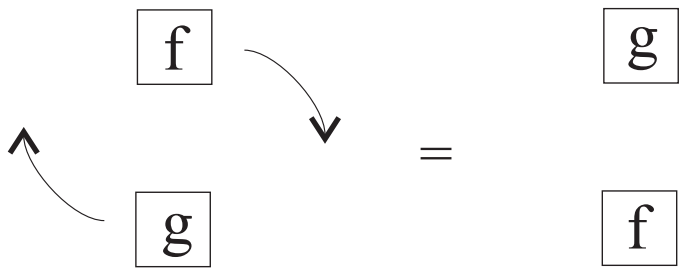

□

*Algebraic proof.*

$$\begin{aligned}
g \circ f =& (g \otimes \mathrm{id}_1) \circ (\mathrm{id}_1 \otimes f) && \text{(1 is the unit object)} && (2.18)\\
=& (g \circ \mathrm{id}_1) \otimes (\mathrm{id}_1 \circ f) && \text{(Functoriality of } \otimes) && (2.19)\\
=& g \otimes f && (\mathrm{id}_1 \text{ is the identity arrow}) && (2.20)\\
=& (\mathrm{id}_1 \circ g) \otimes (f \circ \mathrm{id}_1) && (\mathrm{id}_1 \text{ is the identity arrow}) && (2.21)\\
=& (f \otimes \mathrm{id}_1) \circ (\mathrm{id}_1 \otimes g) && \text{(Functoriality of } \otimes) && (2.22)\\
=& f \circ g && \text{(1 is the unit object)} && (2.23)
\end{aligned}$$

□

### More on Duals

If a category has duals, then they are unique up to a unique isomorphism compatible with the $i_V$ and $e_V$ maps, i.e. for any two duals $(V^*_{(1)}, e_{(1)}, i_{(1)})$ and $(V^*_{(2)}, e_{(1)}, i_{(1)})$ of an object $V$, there is a unique isomorphism $\phi : V^*_{(1)} \xrightarrow{\sim} V^*_{(2)}$ which makes the appropriate diagrams commute. Simply choose $\phi$ to be the composition:

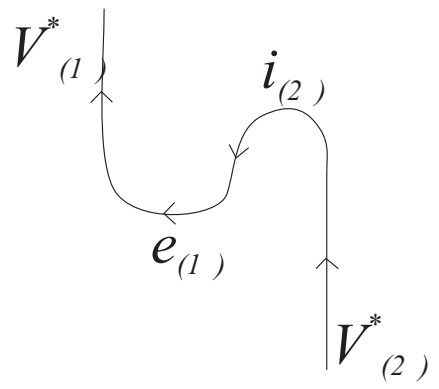

Amongst other things, this means that $(V \otimes W)^* = W^* \otimes V^*$, which makes $e_{V \otimes W}$ and $i_{V \otimes W}$ fit into our graphical notation nicely:

$i_{V\otimes W}=$

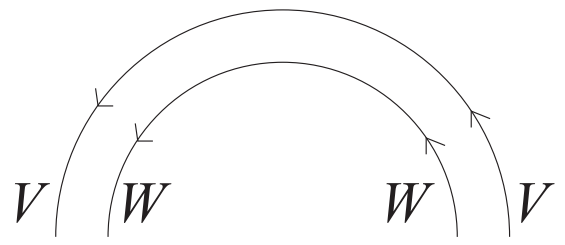

$e_{V\otimes W}=$

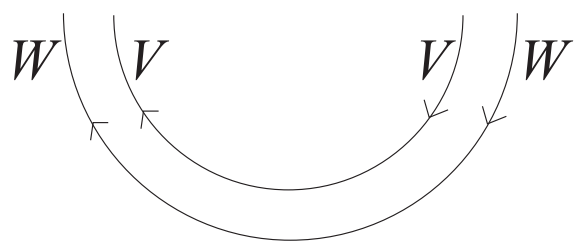

The duality structure also allows us to define duals of *arrows* and not just of objects. Given a morphism $f : U \to V$ one can define $f^* : V^* \to U^*$ in the following way:

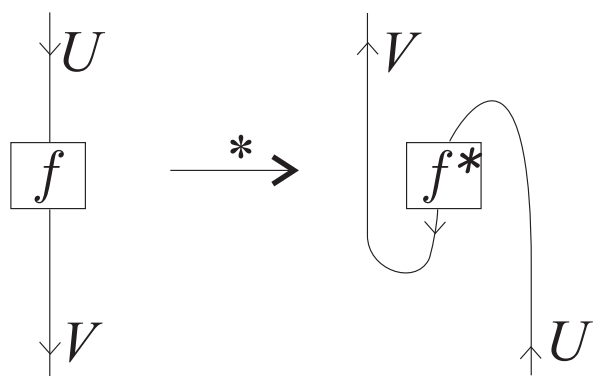

(2.24)

### 2.3.2 Braidings

So Theorem (2.3.1) tells us that a monoidal category with duals is sufficient to do graphical calculus invariant under planar deformations, but how about three dimensional deformations? We will clearly have to contend with braidings ('crossings' of the strings):

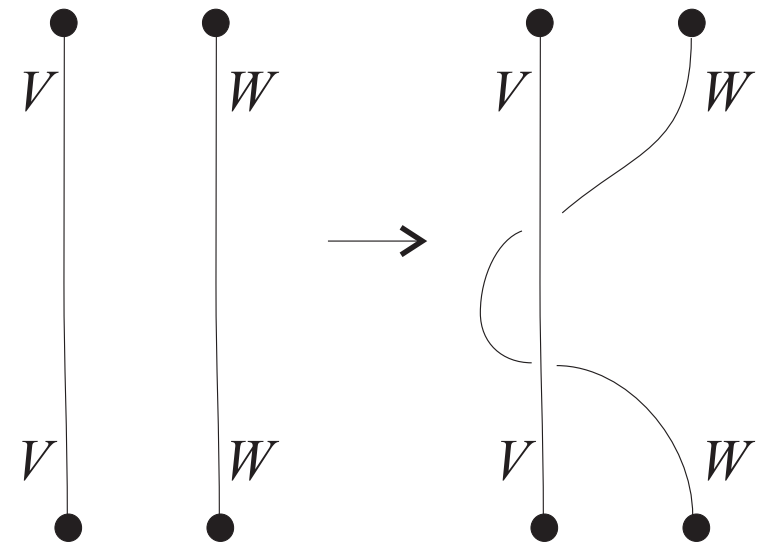

There are two types of braidings, the 'over' and the 'under', which cannot be deformed into each other (by a deformation we mean one limited to the region between the top and bottom planes, and equal to the identity on these planes). It is clear what additional structure we need to place on our category. For every pair $V, W$ of objects there should be assigned a braiding isomorphism $\sigma_{V,W} : V \otimes W \xrightarrow{\sim} W \otimes V$:

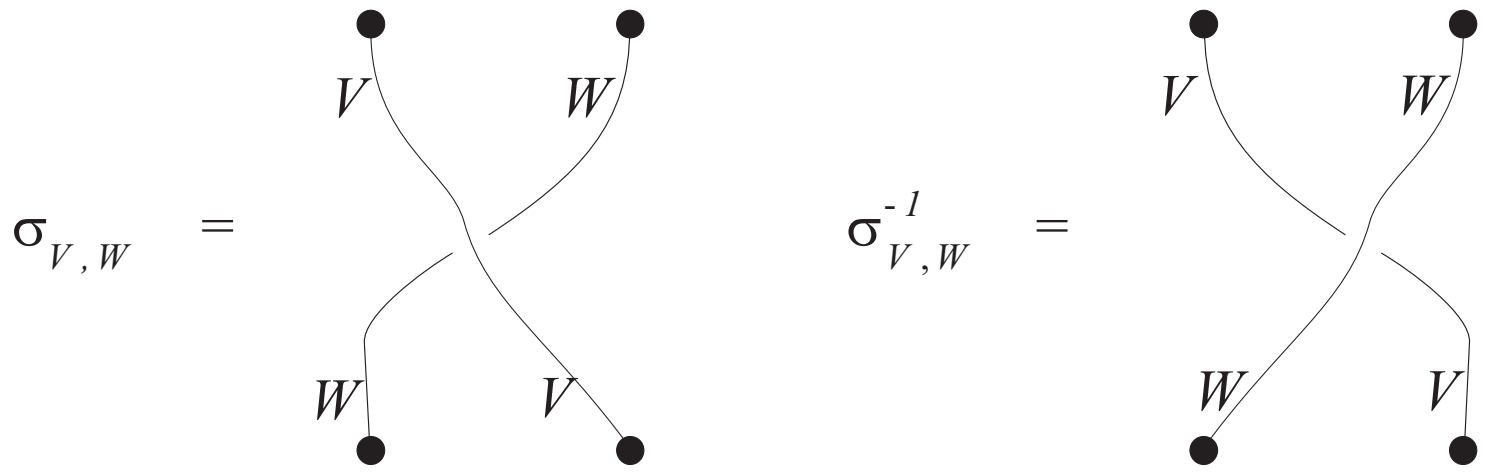

It must be an isomorphism, since clearly $\sigma^{-1}_{V,W}$ is indeed the inverse of $\sigma_{V,W}$:

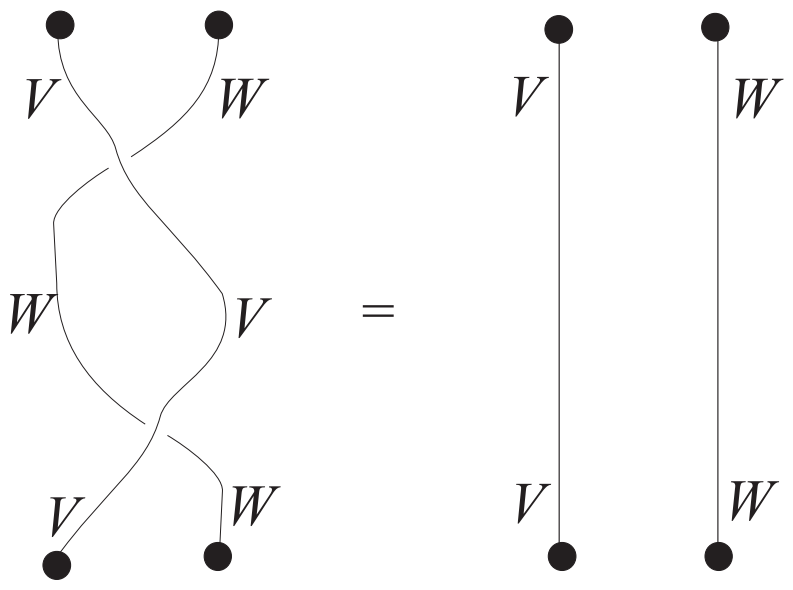

This is known as the 1st Reidemeister move. The $\sigma_{V,W}$'s are not just arbitrary isomorphisms; they will have to satisfy (at least) the following relation:

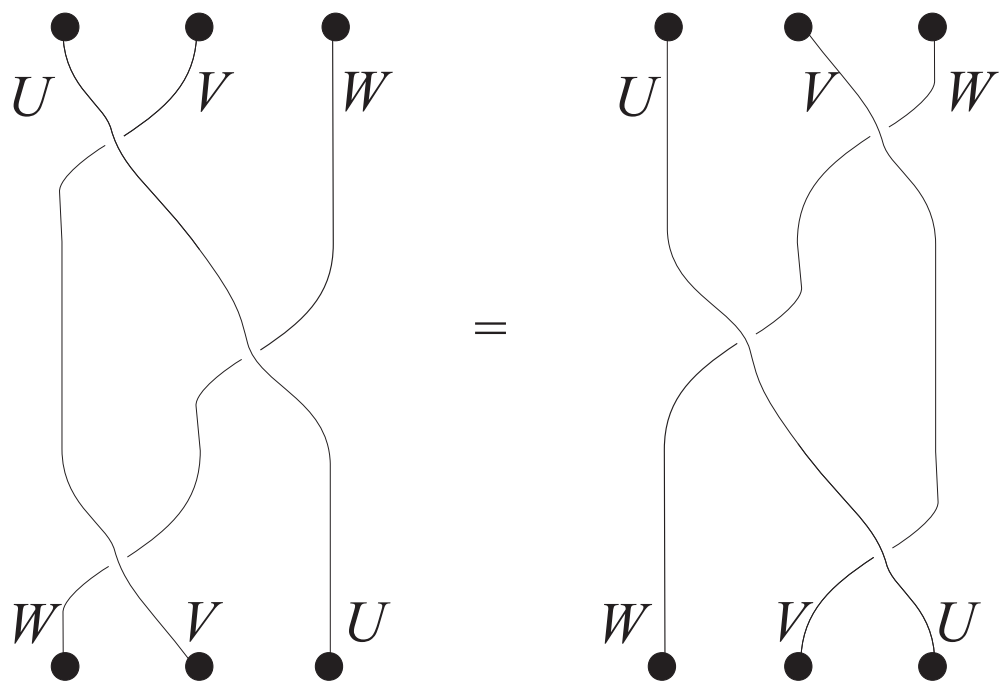

That is,

$$(\sigma_{W,V} \otimes \mathrm{id}_U) \circ (\mathrm{id}_V \otimes \sigma_{U,W}) \circ (\sigma_{U,V} \otimes \mathrm{id}_W) = (\mathrm{id}_W \otimes \sigma_{U,V}) \circ (\sigma_{U,W} \otimes \mathrm{id}_V) \circ (\mathrm{id}_U \otimes \sigma_{V,W}) \tag{2.25}$$

This is known as the 2nd Reidemeister move. It is nothing but the celebrated Yang-Baxter equation in statistical physics!

Let $(\mathcal{C}, \otimes, 1)$ be a monoidal category. A *braiding* for $\mathcal{C}$ consists of a natural family of isomorphisms

$$\sigma_{V,W} : V \otimes W \xrightarrow{\sim} W \otimes V \tag{2.26}$$

where $V, W$ run over all objects of $\mathcal{C}$, such that for any three objects $U, V, W$ we have

$$\begin{aligned} \sigma_{U,V\otimes W} &= (\mathrm{id}_V \otimes \sigma_{U,W})(\sigma_{U,V} \otimes \mathrm{id}_W), && (2.27)\\ \sigma_{U\otimes V,W} &= (\sigma_{U,W} \otimes \mathrm{id}_V)(\mathrm{id}_U \otimes \sigma_{V,W}). && (2.28) \end{aligned}$$

The word 'natural' means here that for all $f : V \longrightarrow V'$ and $g : W \longrightarrow W'$ it does not matter whether the braiding is applied before or after the maps, i.e. the following diagram commutes:

$$\begin{array}{ccc} V \otimes W & \xrightarrow{f\otimes g} & V' \otimes W' \\ \downarrow{\scriptstyle \sigma_{V,W}} & & \downarrow{\scriptstyle \sigma_{V',W'}} \\ W \otimes V & \xrightarrow{g\otimes f} & W' \otimes V' \end{array} \tag{2.29}$$

Or graphically:

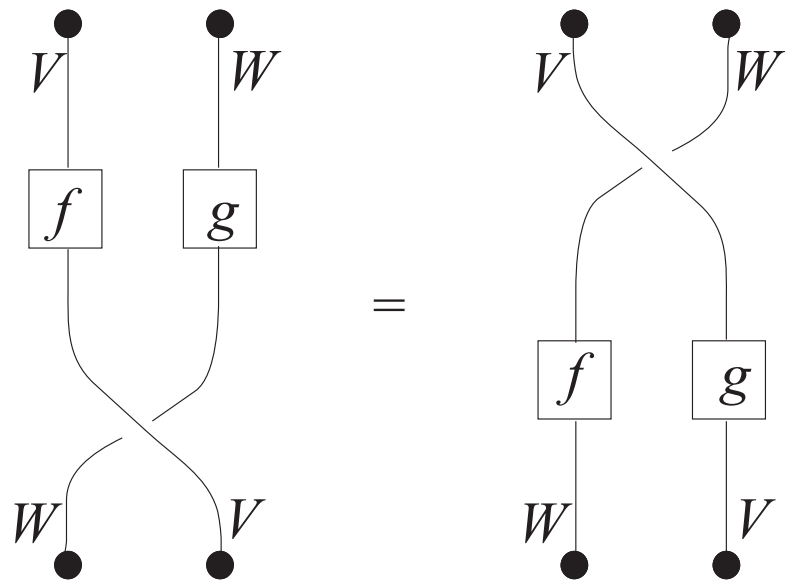

We will represent braidings between involving tensor products, such as $\sigma_{U,V\otimes W}$, by:

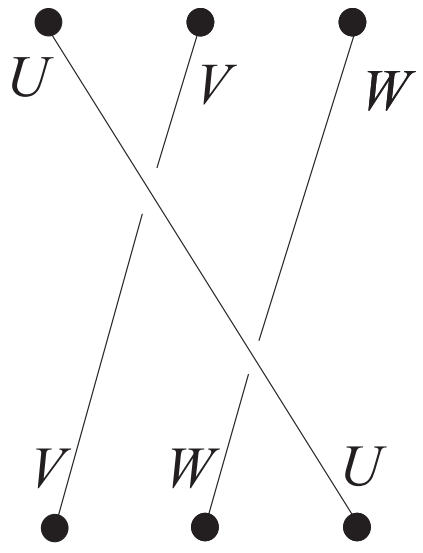

Eqs. (2.27) and (2.28) relate the braiding to the tensor product in a bilinear way:

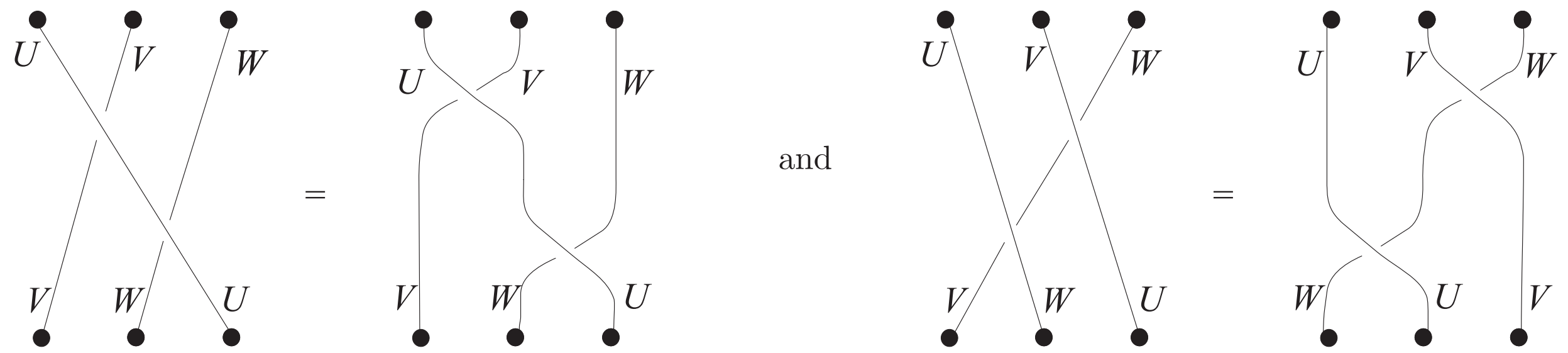

The conditions of naturality and bilinearity are strong enough to prove the 2nd Reidemeister move:

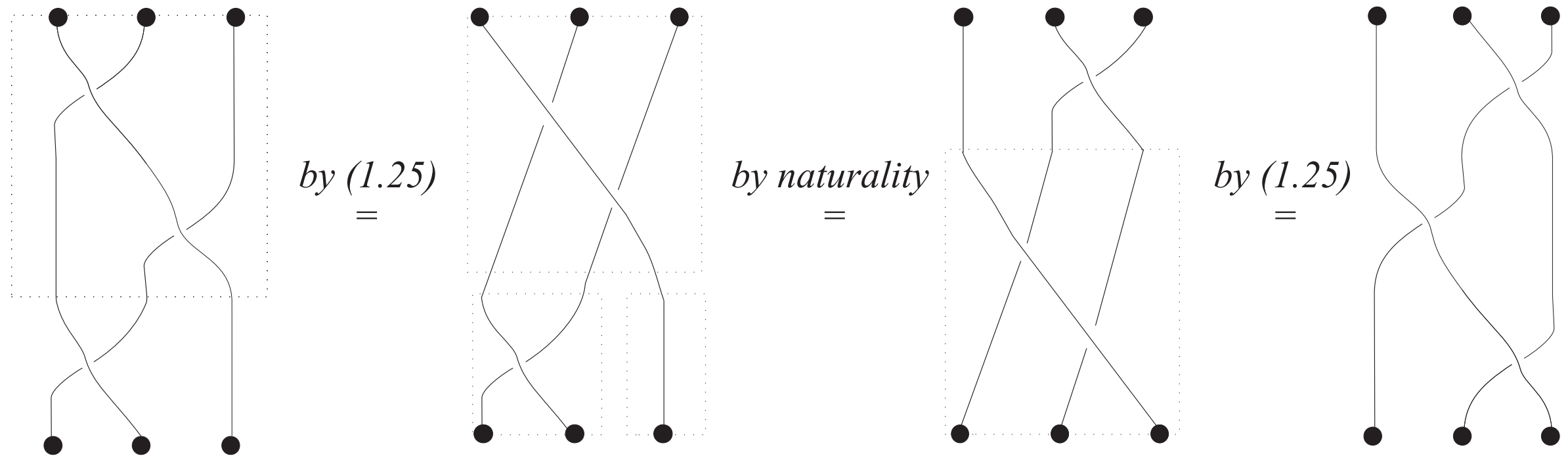

### 2.3.3 Twists, traces and dimensions

A braiding is not yet enough to ensure that the graphical calculus be deformation invariant. There is a 3rd Reidemeister move[11]:

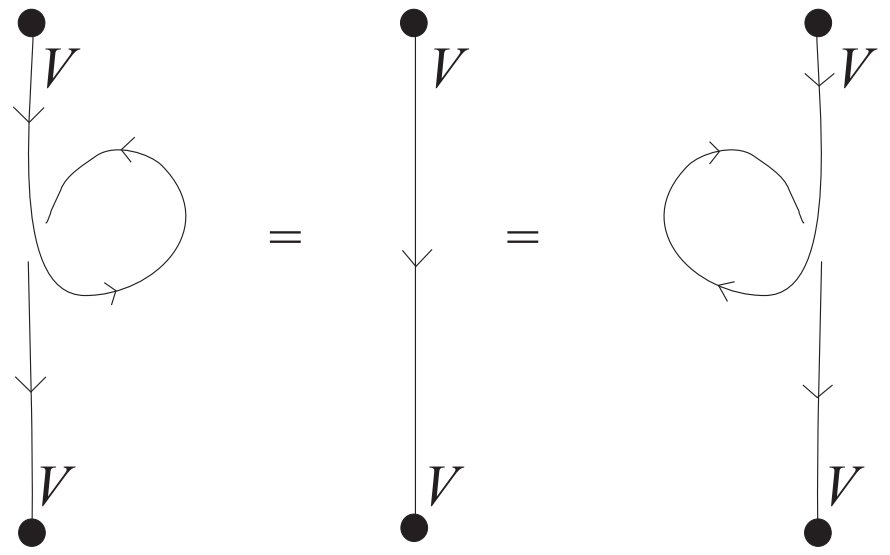

In our language, this will clearly be a relation between the duality structure and the braiding structure. The trouble is that we cannot construct such a figure since we are missing two diagram components : to the standard duality structure (2.17) we need to somehow add the arrows $i'_V : 1 \to V^* \otimes V$ and $e'_V : V \otimes V^* \to 1$:

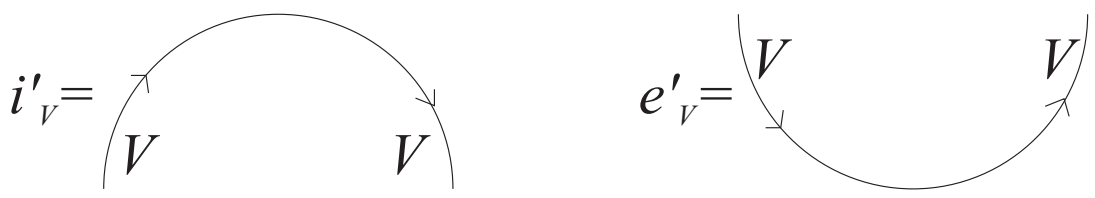

(2.30)

[11]Strictly speaking, the Reidemeister moves deal with unoriented tangles, but we have oriented them in order to demonstrate the connection with the duality.

There is one obvious way to add these arrows - use the braiding to swap the inputs of $i_V$ and $e_v$:

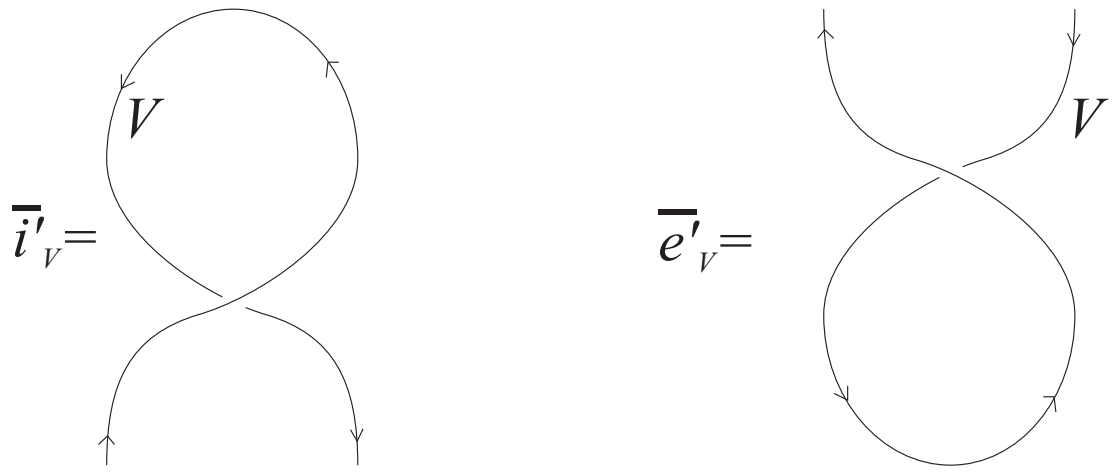

(2.31)

These maps certainly satisfy the expected duality snake relation. For instance, here is a proof that

$$V \xrightarrow{\mathrm{id}\otimes\overline{i'}_V} V \otimes V^* \otimes V \xrightarrow{\overline{e'}_V\otimes\mathrm{id}} V \tag{2.32}$$

is the identity:

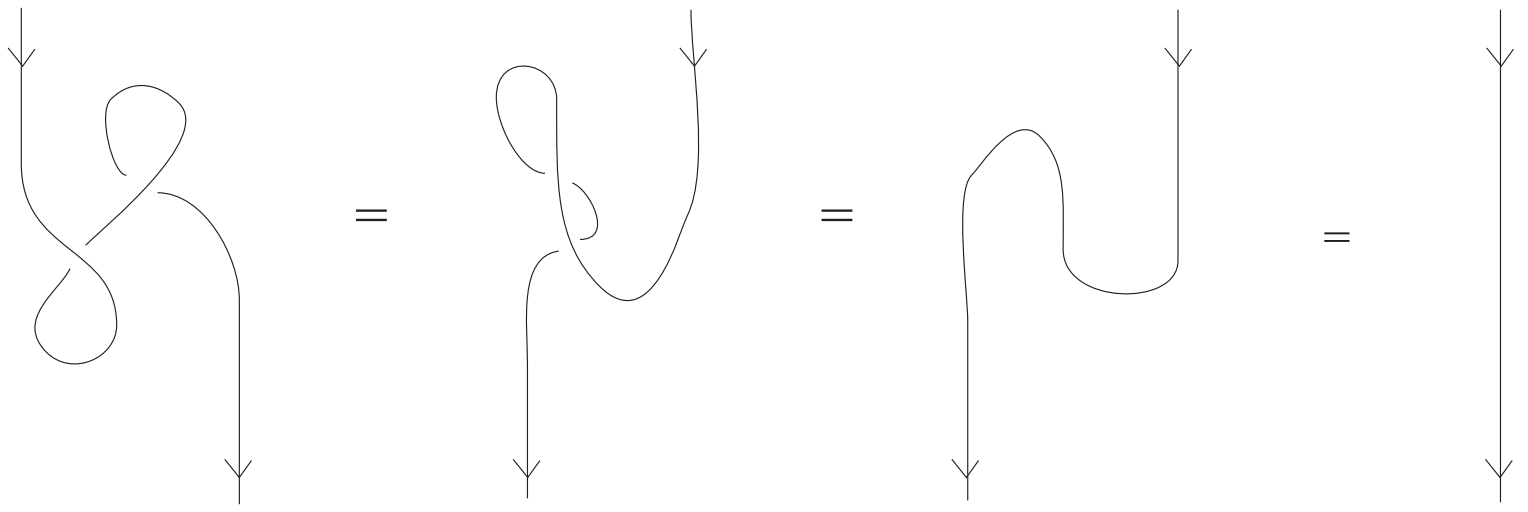

In the first and second steps we used naturality and properties (2.27) and (2.28) of the braiding, while in the third step we used the fact that $i_V$ and $e_V$ are dual maps.

The trouble with $\overline{i'}_V$ and $\overline{e'}_V$ is that they will not allow for a nice theory of traces and dimensions in our category. Suppose $V$ is an object and $f : V \to V$. We define the *trace* of f, denoted $\mathrm{tr} f \in \mathbb{C}$, as:

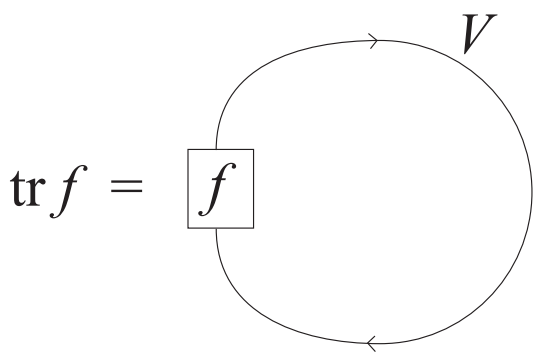

These are the '1-loop' functions of our category. Note that we are using $\overline{i'}_V$ and $\overline{e'}_V$ in these diagrams via the shorthand notation (2.30). We define the *dimension* of $V$ to be $\dim V = \mathrm{tr}\,\mathrm{id}_V$:

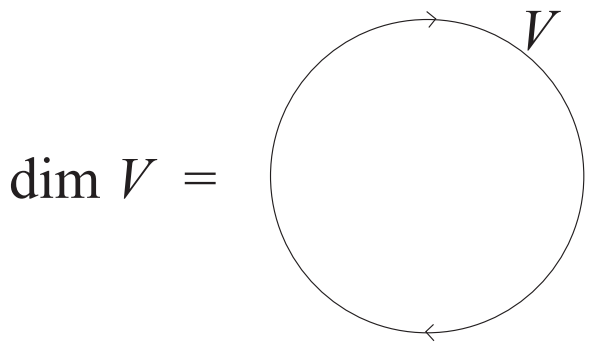

We would like these traces and dimensions to behave just like those in **Vect**. Specifically, we want:

$$\text{a) } \mathrm{tr}(f \otimes g) = \mathrm{tr} f\, \mathrm{tr} g, \quad \text{b) } \mathrm{tr}(f^*) = \mathrm{tr} f, \quad \text{c) } \mathrm{tr}(fg) = \mathrm{tr}(gf). \tag{2.33}$$

In particular,

$$\text{(A) } \dim(V \otimes W) = \dim V \dim W, \quad \text{(B) } \dim V^* = \dim V. \tag{2.34}$$

Let us try to obtain, say, (A) from our current definitions. (A) says that one should be able to pull the circles apart:

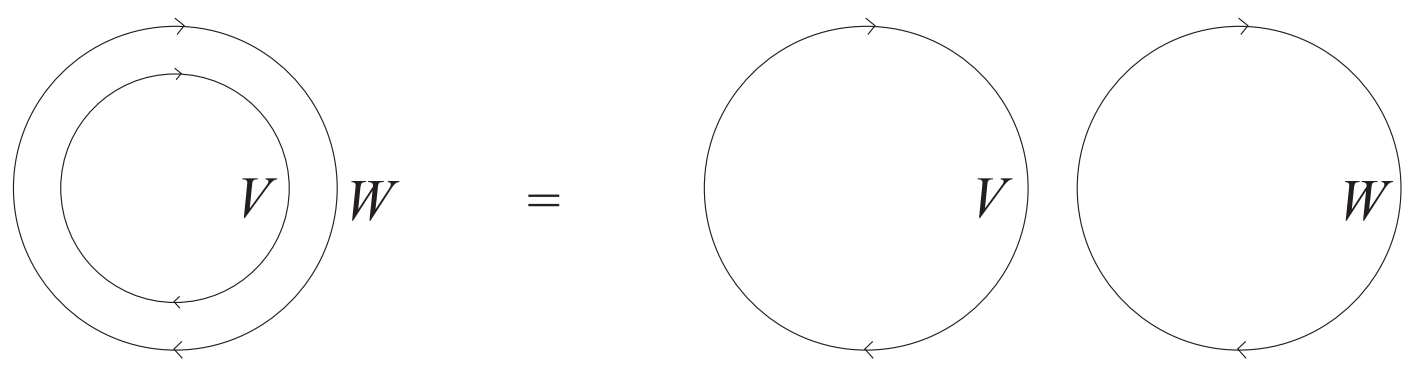

Yet our definition (2.31) has somehow linked them together:[12]

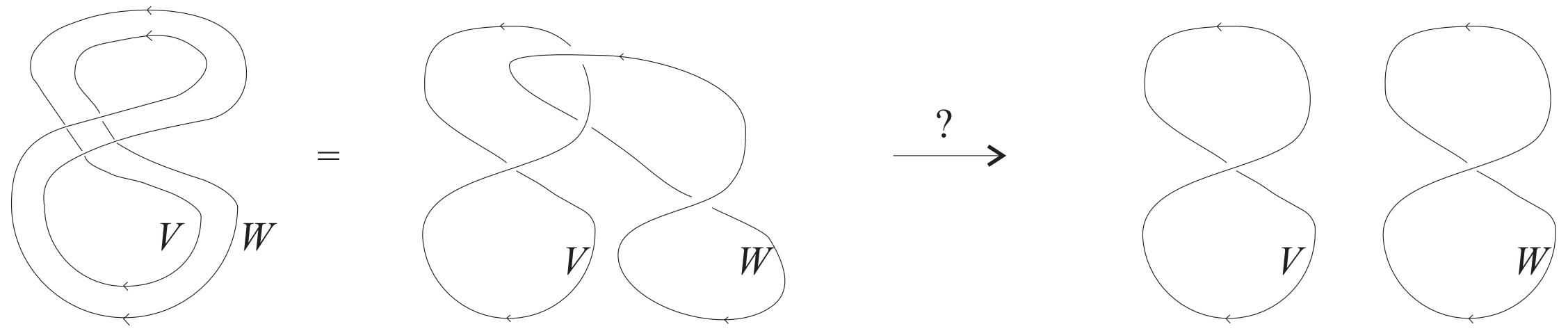

It is clear what has happened; (2.31) has introduced an unwanted twist. To get rid of it, we postulate the existence of a natural family of isomorphisms

$$\theta = \{\theta_V : V \xrightarrow{\sim} V\} \tag{2.35}$$

for all objects $V$, which interacts with the monoidal structure, braiding and duality via:

$$\theta_{V\otimes W} = \sigma_{W,V}\sigma_{V,W}(\theta_V \otimes \theta_W), \tag{2.36}$$
$$(\theta_V \otimes \mathrm{id}_{V^*})i_V = (\mathrm{id}_V \otimes \theta_{V^*})i_V, \tag{2.37}$$
$$\theta_1 = \mathrm{id}. \tag{2.38}$$

A braided rigid monoidal category with a compatible twist is called a *ribbon* category. We interpret $\theta$ as an actual $2\pi$ twist, so that instead of strings we should really be drawing ribbons. For ease of use, we will continue to employ string notation[13]. Here is the graphical expression for (2.36):

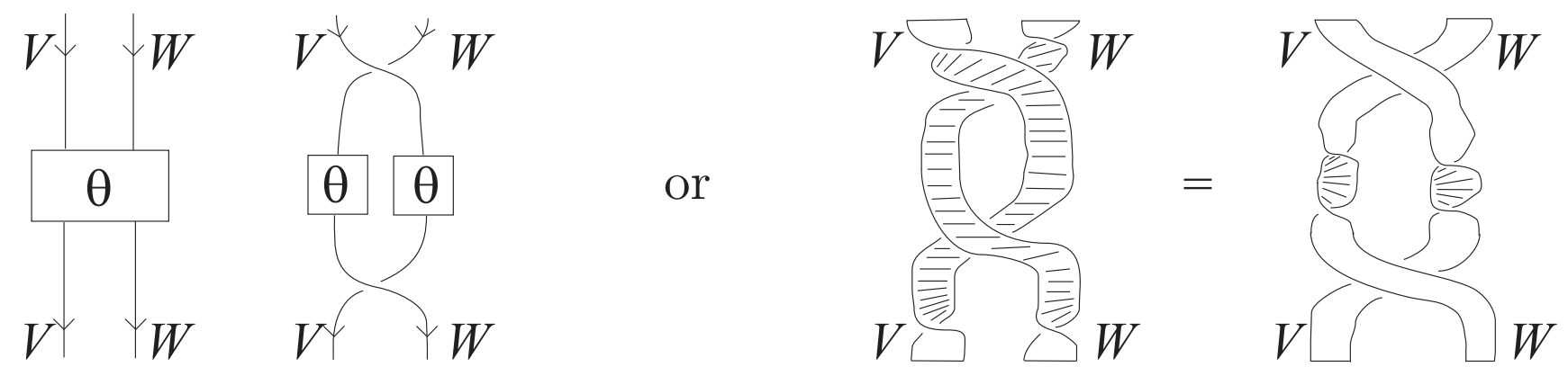

(2.37) tells us we can move the twist around a dual map:

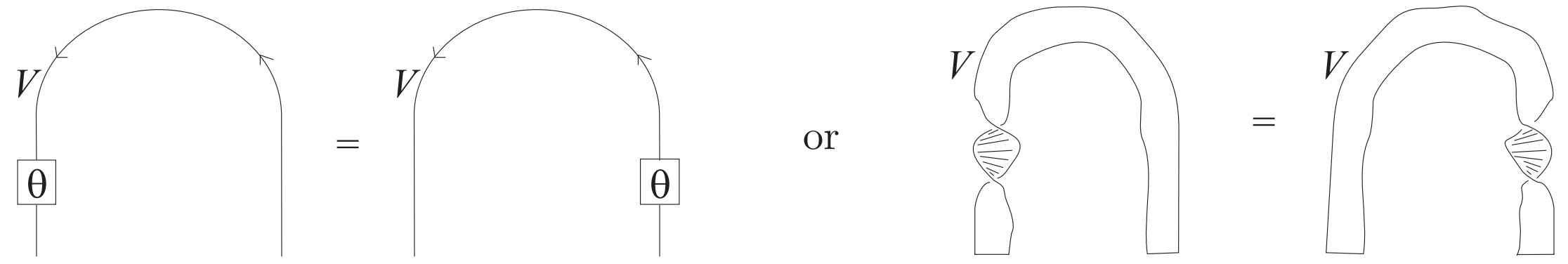

We now define $i'_V$ and $e'_V$ as:

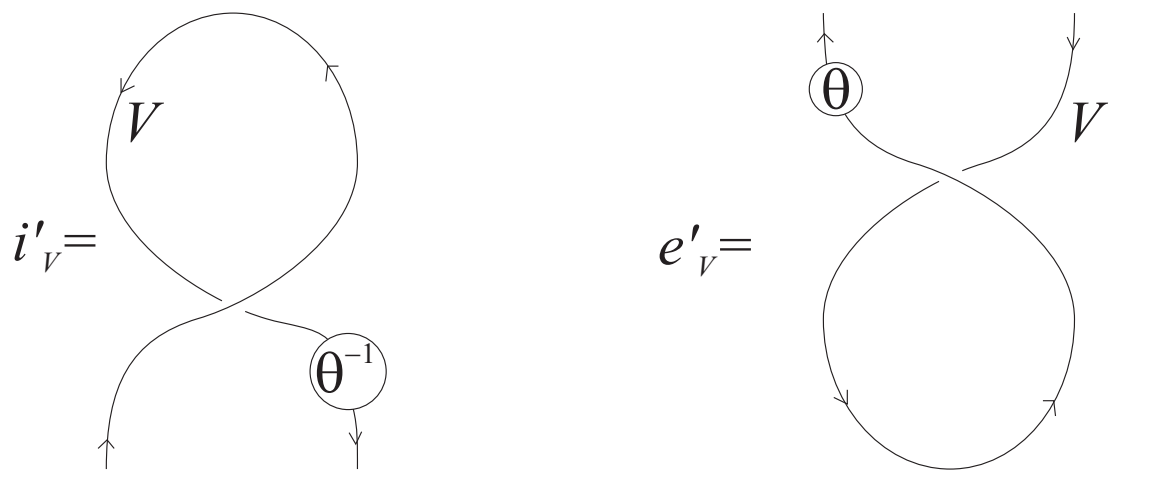

[12]The reader is urged to play with these diagrams, using the axioms!

[13]This is the notation one uses when drawing these diagrams quickly by hand.

Properties (2.36)-(2.37) of the twist enable us to pull apart the loops in dimensions and traces:

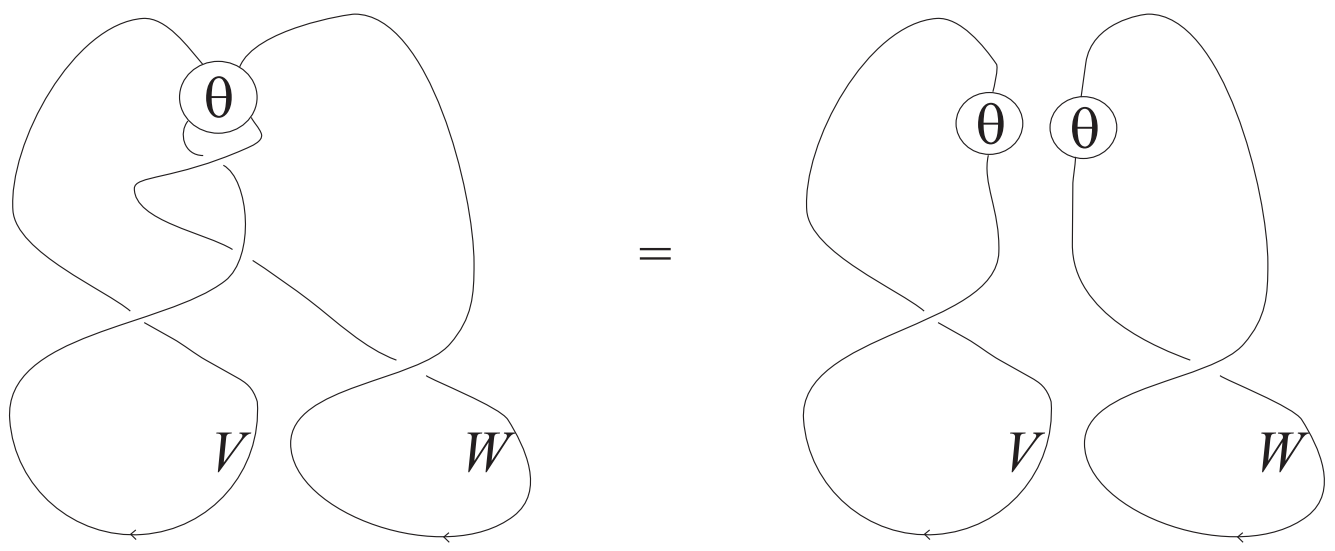

We began Section 2.3.2 by asking what additional structure is needed on a rigid monoidal category in order to turn the 2d graphical calculus into a 3d graphical calculus. We have arrived at the answer : one needs a ribbon category. This is the content of the following beautiful and powerful theorem, which lies at the very foundation of this thesis and all of 3d topological quantum field theory:

**2.3.3 Theorem (Reshetikhin and Turaev [89]).** *Let $(\mathcal{C}, \otimes, 1, \sigma, \theta)$ be a ribbon category. Then the evaluation of a ribbon diagram into a morphism in $\mathcal{C}$ is invariant under 3d isotopy*[14].

In other words, the graphical calculus in terms of ribbons is entirely consistent, and is a powerful way to do computations inside ribbon categories : simply deform the diagrams appropriately.

## 2.4 Semisimple ribbon categories

The categories which we will work with have, in addition to the tensor product, a direct sum of objects. The interaction between the two are known as the *fusion rules* familiar from conformal field theory. We are primarily interested in $\mathbb{C}$-linear abelian categories. We shall call a category $\mathcal{C}$ $\mathbb{C}$-linear if:

- $\mathrm{Hom}_{\mathcal{C}}(U, V)$ is a vector space for all objects $U, V$ in $\mathcal{C}$ and such that composition is a linear operation.
- There exists a zero object 0 in $\mathcal{C}$ such that $\mathrm{Hom}_{\mathcal{C}}(0, V) = \mathrm{Hom}_{\mathcal{C}}(V, 0) = \{0\}$ for all $V$ in $\mathcal{C}$.
- Finite direct sums[15] exist in $\mathcal{C}$.

An *abelian* category is a category where one can use the notion of kernel and cokernel of a morphism in the same way as in **Vect** [16]. An object $U$ in an abelian category is called *simple* if any injection $V \hookrightarrow U$ is either 0 or an isomorphism. In other words, $U$ does not contain nontrivial 'subobjects'. An abelian category $\mathcal{C}$ is called *semisimple* if any object $V$ is isomorphic to a direct sum of simple ones:

$$V = \bigoplus_{i \in I} N_i V_i \tag{2.39}$$

where $V_i$ are simple objects and $I$ indexes the isomorphism classes of nonzero simple objects in $\mathcal{C}$, $N_i$ are positive integers and only a finite number of $N_i$ are nonzero. The decomposition above is known as the *fusion* rules. It is a version of Schur's lemma that the Hom-sets between simple objects have a simple structure:

$$\mathrm{Hom}(V_i, V_j) = \begin{cases} \mathbb{C} & \text{if } i = j, \\ 0 & \text{otherwise.} \end{cases} \tag{2.40}$$

In a semisimple abelian category $\mathcal{C}$, one can construct the *Grothendieck group* $K(\mathcal{C})$ as the quotient of the free abelian group on the set of isomorphism classes of objects in $C$ modulo the relations $\langle V \oplus W \rangle = \langle V \rangle + \langle W \rangle$.

[14] As always, these isotopies take place in the box-like region between the end intervals and are the identity on the end intervals

[15] Direct sums can be defined abstractly in categories as coproducts. Namely, let $\mathcal{C}$ be a category and let $\{X_j : j \in J\}$ be a indexed family of objects in $\mathcal{C}$. The coproduct of the set $\{X_j\}$ is an object $X$ together with a collection of morphisms $i_j : X_j \to X$ (called injections) which satisfy a universal property: for any object $Y$ and any collection of morphisms $f_j : X_j \to Y$, there exists a unique morphism $f$ from $X$ to $Y$ such that $f_j = f \circ i_j$.

[16] For the precise definition, see eg. [94].

Here $\langle V \rangle$ denotes the isomorphism class of $V$. When $C$ is also rigid monoidal, then $K(\mathcal{C})$ becomes a ring - the *Grothendieck ring* of $\mathcal{C}$ - by defining $\langle V \rangle \langle W \rangle = \langle V \otimes W \rangle$. $K(\mathcal{C})$ is an associative ring with unit $\langle 1 \rangle$. If $\mathcal{C}$ is braided, then this ring is obviously commutative, since $\langle V \otimes W \rangle = \langle W \otimes V \rangle$ via the braiding $\sigma_{V,W}$. The Grothendieck ring does not 'see' the twist $\theta$, if there is any.

## 2.5 Hermitian and unitary ribbon categories

A Hermitian monoidal category is a category equipped with a contravariant functor[17] † which is the identity on objects but sends each morphism $f : V \to W$ to a morphism $f^\dagger : W \to V$, satisfying

$$f^{\dagger\dagger} = f, \tag{2.41}$$

$$(f \otimes g)^\dagger = f^\dagger \otimes g^\dagger, \tag{2.42}$$

$$(fg)^\dagger = g^\dagger f^\dagger. \tag{2.43}$$

The axioms imply that $\mathrm{id}_V^\dagger = \mathrm{id}_V$. A Hermitian ribbon category $\mathcal{C}$ is a ribbon category where † satisfies

$$\sigma_{V,W}^\dagger = \sigma_{V,W}^{-1}, \quad \theta_V^\dagger = \theta_V^{-1}, \quad i_V^\dagger = e'_V, \quad e_V^\dagger = i'_V. \tag{2.44}$$

If $\mathcal{C}$ is $\mathbb{C}$-linear, then we require that $f \to f^\dagger$ be an antilinear operation.

One should not confuse the duals $f^* : W^* \to V^*$ (2.24) with the daggers $f^\dagger : W \to V$:

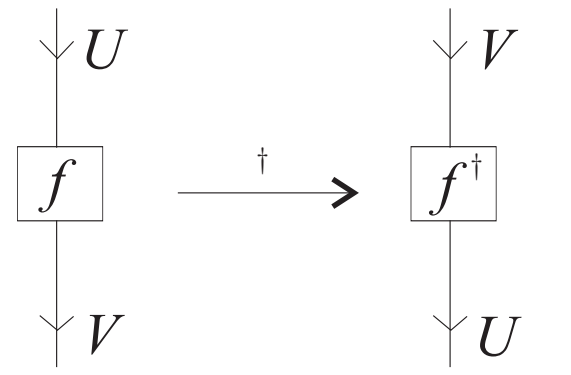

(2.45)

**Hilb** is a Hermitian category, since for any bounded linear operator $T : H \to H'$ we can define the Hermitian adjoint $T^\dagger : H' \to H$ to be given by

$$\langle T^\dagger \psi, \phi \rangle_H = \langle \psi, T\phi \rangle_{H'}. \tag{2.46}$$

Since time evolution is (normally) represented in quantum mechanics by a unitary operator (which satisfies $U^\dagger = U^{-1}$, we see that Hermitian categories are categories equipped with a time-reversal operation $t \to -t$. In contrast, the $*$ operation is a kind of space inversion $x \to -x$.

In a $\mathbb{C}$-linear Hermitian ribbon category, one can define a non-degenerate hermitian inner product on $\mathrm{Hom}(V, W)$ by

$$(f, g) = \mathrm{tr}(\overline{f} g) \tag{2.47}$$

which satisfies the usual relations

$$(f, \lambda g) = \lambda (f, g) \tag{2.48}$$

$$(f, g) = \overline{(g, f)} \tag{2.49}$$

$$(f, f) = 0 \text{ if and only if } f = 0. \tag{2.50}$$

We shall call $\mathcal{C}$ *unitary* if the inner product is positive-definite, i.e. $(f, f) \geq 0$ for all $f$ in $\mathcal{C}$. For instance, **Hilb** is unitary. Physically speaking, in unitary categories the corresponding quantum field theories are positive definite - the ground state energy is bounded from below.

## 2.6 Feynman diagrams and Rep($G$)

This section is devoted to the following profound and beautiful idea connecting the preceding category-flavoured mathematics with physics: *Feynman diagrams are nothing but the graphical calculus applied to the ribbon category* $Rep_f(G)$*, where* $G$ *is the symmetry group of the physical theory.*

[17] John Baez [13] calls them $*$-categories, which clashes with the present author's (possibly peculiar) convention of using † for adjoints and reserving $*$ for duals of linear maps.

### 2.6.1 $\text{Rep}_f(G)$ as a ribbon category

Recall from Section 2.1.2 that the category $\text{Rep}_f(G)$ has finite dimensional representations of $G$ as objects, and intertwining maps as arrows. If we think of the arrows as physical processes, then the requirement that they be intertwining is physically the requirement of *covariance* with respect to $G$. $\text{Rep}_f(G)$ is monoidal, with the usual tensor product of representations,

$$(\rho_1 \otimes \rho_2)(g) = \rho_1(g) \otimes \rho_2(g). \tag{2.51}$$

The unit 1 is provided by the trivial representation $\rho_{\text{trivial}}$. If $(\rho, V)$ is a representation of $G$ on a vector space $V$, then its dual $(\rho^*, V^*)$ is a representation of $G$ on $V^*$ defined by

$$(\rho^*(g)(f))\,(v) = f(\rho(g^{-1})v) \tag{2.52}$$

The duality structure is provided by the familiar vector space morphisms $i_V : 1 \to V \otimes V^*$ and $e_V : V^* \otimes V \to 1$ defined by setting $i_V(1) = \text{id}_V$ (this uses the identification of $V \otimes V^*$ with $\text{End}(V)$) and $e_V$ the evaluation map $f \otimes v \to f(v)$. We should check that these are intertwiners. $i_v$ is clearly an intertwiner. $e_V$ is also an intertwiner; one computes that $e_V \circ [\rho^* \otimes \rho(g)] = e_V = \rho_{\text{trivial}}(g) \circ e_V$ for all $g \in G$.

The braiding depends on the *statistics* of the particles one is describing. If $(\rho, V)$ and $(\rho', W)$ are representations, and $v \in V$ and $w \in W$, then the braiding is defined as:

$$\sigma_{V,W}(v \otimes w) = \begin{cases} w \otimes v & \text{for bosons,} \\ -w \otimes v & \text{for fermions.} \end{cases} \tag{2.53}$$

These rules are familiar - swapping two fermions picks up a minus sign in Feynman diagrams, while swapping two bosons should leave the diagram unchanged. One may ask if there are quantum field theories with a more exotic interchange law than (2.53). The answer is *yes* - for example, there are the 2d theories involving anyons (see [101] and references therein), particles with fractional statistics. In general this involves replacing the symmetry group $G$ by a *quantum group* $U_q(G)$ - a kind of deformation of the original group involving a deformation parameter $q = \exp(i\theta)$ (where $\theta$ may be complex) such that as $q \to 1$ one recovers the original group $G$. We shall treat these theories later.

The twists $\theta_V$ are trivial, $\theta_V = \text{id}$. The category $\text{Rep}_f(G)$ is clearly $\mathbb{C}$-linear and abelian. If $G$ is compact, then the simple objects are precisely the irreducible representations, since all representations are completely reducible. Thus $\text{Rep}_f(G)$ is a semisimple ribbon category, whose Grothendieck ring $K(\text{Rep}_f(G))$ is precisely the Verlinde algebra characterizing how tensor products of representations of $G$ decompose into sums of irreducible ones.

### 2.6.2 CPT Theorem

We interpret the vector space $V$ in a representation $(\rho, V)$ as the state space for a particle. Its dual $V^*$ is interpreted as the state space for the corresponding antiparticle. The graphical calculus now takes on a physical meaning : downwards pointing arrows correspond to particles, and upwards pointing arrows correspond to antiparticles. Since time travels downwards in our diagrams, this is precisely the prescription Feynman taught us[18]. The duality maps $i_V$ and $e_V$ for a representation on $V$ are literally interpreted as the creation and annihilation of a particle-antiparticle pair:

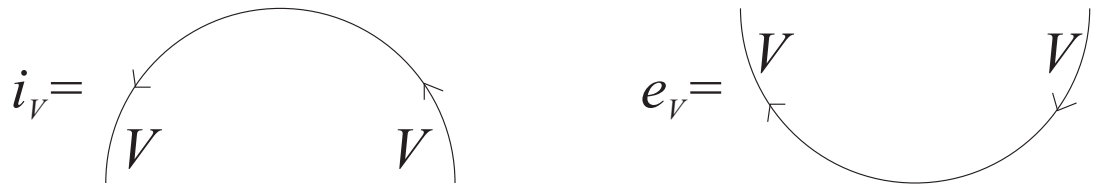

Besides the dual $(\rho^*, V^*)$, can also form the *conjugate* representation $(\bar{\rho}, \bar{V})$. The vector space $\bar{V}$ is defined to be conjugate to $V$ in the sense that addition of vectors is unchanged, but scalar multiplication is sent to its complex conjugate:

$$v \to \overline{v} \tag{2.54}$$

$$\overline{v} + \overline{w} = \overline{v + w} \tag{2.55}$$

$$\lambda\overline{v} = \overline{\overline{\lambda} v} \tag{2.56}$$

[18] According to Feynman, antiparticles are particles running backwards in time.

The operation $V \to V^*$ is a contravariant functor, while $V \to \overline{V}$ is a covariant one.For a map $T : V \to W$, we already know the definition of $T^* : W^* \to V^*$. Similarly $\overline{T} : \overline{V} \to \overline{W}$ is defined by $\overline{T}\overline{v} = \overline{Tv}$. Given $V$, one can form

$$V, V^*, \overline{V}, \overline{V}^*, \overline{V^*}. \tag{2.57}$$

These fit into a schematic diagram:

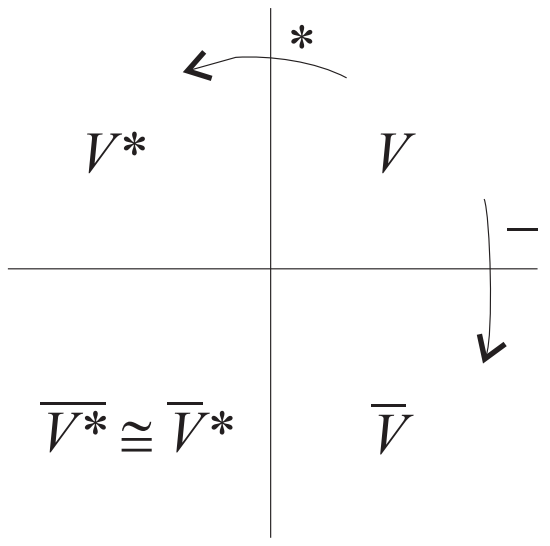

There are three discrete symmetries a quantum field theory might have :

- C - Charge conjugation : replacing particles by their antiparticles.
- P - Parity : invert space, $(x, y, z) \to (-x, -y, -z)$.
- T - Time reversal : reflect time, $t \to -t$.

The celebrated CPT theorem says that only the combination CPT is assured to be a symmetry of a relativistic quantum field theory. If we choose

$$C \leftrightarrow \ ^*, \quad PT \leftrightarrow \ \bar{}. \tag{2.58}$$

then we can prove a 'baby CPT' theorem for the category $\mathrm{Rep}_f(G)$!

**2.6.1 'CPT' Lemma.** *If $(\rho, V)$ is unitary, then $(\rho^*, V^*)$ is unitarily equivalent to $(\overline{\rho}, \overline{V})$. Thus $(\rho, V)$ is equivalent to 'CPT' $(\rho, V) = (\overline{\rho}^*, \overline{V}^*)$.*

*Proof.* The correspondence

$$v \to \langle v, \cdot \rangle, \tag{2.59}$$

thought of as a linear map from $\bar{V} \to V^*$, does the job. □

### 2.6.3 Conservation of energy-momentum from intertwiners

The symmetry group of flat space-time is the Poincaré group - the semidirect product of the Lorentz group with the spacetime translation group:

$$P = SO(3,1) \ltimes \mathbb{R}^4. \tag{2.60}$$

Let us focus on the translations. The irreducible representations of $\mathbb{R}$ are simply the functions

$$\rho_k(t) = e^{kt} \tag{2.61}$$

for any $k \in \mathbb{C}$, and $\rho_k$ is equivalent to $\rho_{k'}$ iff $k = k'$. The dual and conjugate reps turn out to be

$$\rho_k^* \simeq \rho_{-k}, \quad \overline{\rho_k} \simeq \rho_{\overline{k}}. \tag{2.62}$$

The unitary reps are those where $k$ is purely imaginary, so we can write

$$\rho_k(t) = e^{iEt} \in U(1) \tag{2.63}$$

We interpret $E$ as the energy. Since all the $\rho_k$ are irreducible, there are no intertwiners from $\rho_k$ to $\rho_{k'}$ unless $k = k'$, when multiplication by an arbitrary complex number is allowed. This is precisely energy conservation! To underline this point, we first calculate the tensor product of two reps:

$$\rho_k \otimes \rho_l \simeq \rho_{k+l} \tag{2.64}$$

Thus we conclude that in an interaction process:

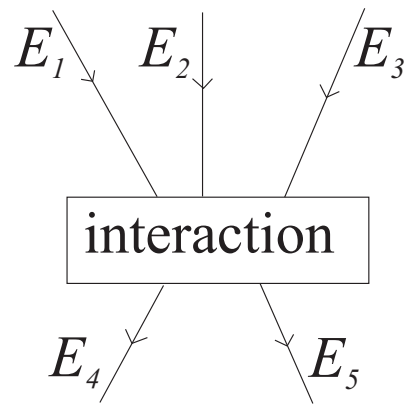

We must have total energy conservation in order for the interaction to be nonzero:

$$E_1 + E_2 + E_3 = E_4 + E_5 \tag{2.65}$$

The unitary irreps of $\mathbb{R}^4$ are of the form:

$$\rho^k = \rho_{k_0} \otimes \rho_{k_1} \otimes \rho_{k_2} \otimes \rho_{k_3} \tag{2.66}$$

where $(k_0, k_1, k_2, k_3) \in \mathbb{R}^4$ is interpreted as the energy-momentum. In precisely the same way as before, we see that an interaction

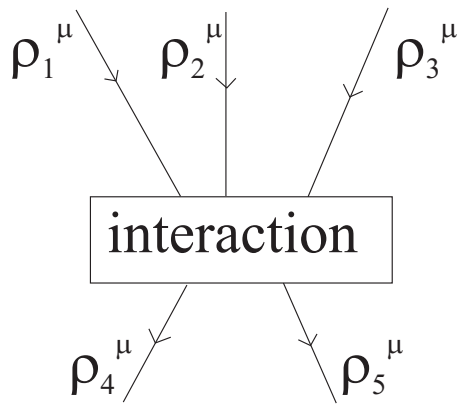

is zero unless we have conservation of energy momentum,

$$\rho_1^\mu + \rho_2^\mu + \rho_3^\mu = \rho_4^\mu + \rho_5^\mu. \tag{2.67}$$

### 2.6.4 Simple example

Consider the baby field theory where fields are simply column vectors living in an $N$ dimensional vector space $V$, with the action

$$S = z^\dagger A z + \lambda (z^\dagger z)^2 = \sum_{i,j}^{N} z_i^* A_{ij} z_j + \lambda \sum_{i,j}^{N} z_i^* z_i z_j^* z_j. \tag{2.68}$$

We interpret $V$ as the fundamental representation of $U(N)$ which acts via $z \to Uz$. The action is $U(N)$ invariant providing $U^\dagger A U = A$, or put differently, if $A : V \to V$ is an intertwiner. Since $U(N)$ is compact, and $V$ is finite dimensional, this means that $A$ must be proportional to the identity (Schur's lemma). Ignoring this for now, we note that the machinery of QFT instructs us to calculate processes such as:

$$\langle z_i^* z_j^* z_k z_l \rangle = \int dz^* dz \, z_i^* z_j^* z_k z_l e^{-z^\dagger A z - \lambda (z^\dagger z)^2}. \tag{2.69}$$

These processes can be calculated order for order in $\lambda$. That is, we first recall Wick's theorem which gives a formula for the $n$-point functions in the non-interacting theory as a sum over all pairs of elementary contractions:

$$\langle z_1^* z_2^* \cdots z_n^* z_1 z_2 \cdots z_n \rangle = \sum_{\text{contractions}} \langle z_{i_1}^* z_{j_1} \rangle_0 \langle z_{i_2}^* z_{j_2} \rangle_0 \cdots \langle z_{i_n}^* z_{j_n} \rangle_0 = \sum_{\text{contractions}} A^{-1}_{i_1 j_1} A^{-1}_{i_2 j_2} \cdots A^{-1}_{i_n j_n}. \tag{2.70}$$

Thus, the following first order process,

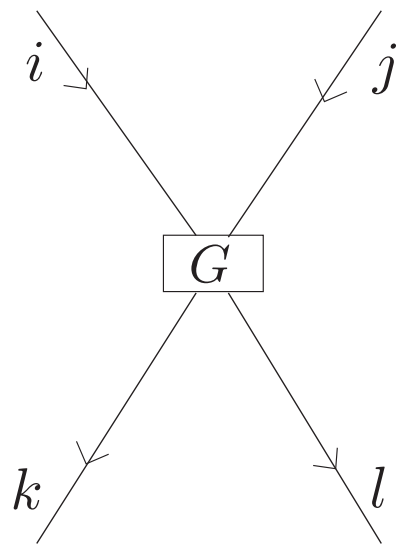

defines an intertwiner

$$G : V \otimes V \to V \otimes V, \tag{2.71}$$

given in terms of components by $e_i \otimes e_j \to G_{ij}^{kl} e_k \otimes e_l$ where $(e_1, \ldots, e_N)$ is a basis for $V$, and

$$G_{ij}^{kl} = \lambda \sum_{p,q} A_{ip}^{-1} A_{jq}^{-1} A_{pk}^{-1} A_{ql}^{-1}. \tag{2.72}$$

Higher order processes are made up of these first-order building blocks, and can thus be computed in the standard graphical calculus way. For instance, the following second order-process,

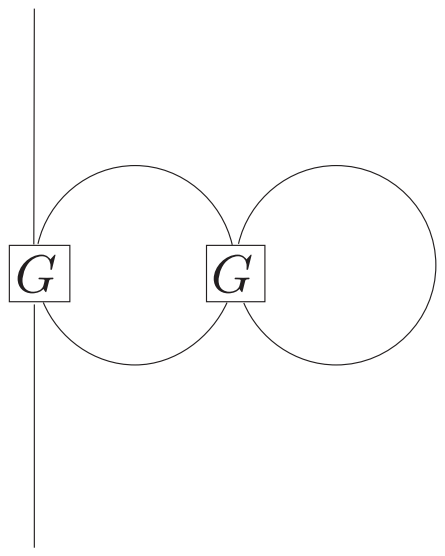

translates into the sequence of maps:[19]

$$V \xrightarrow{\mathrm{id}\otimes i_V \otimes i_V} V \otimes V \otimes V \otimes V \otimes V \xrightarrow{G\otimes G} V \otimes V \otimes V \otimes V \otimes V \xrightarrow{\mathrm{id}\otimes e_V \otimes e_V} V \tag{2.73}$$

Notice from (2.72) that in fact, $G = A^{-1} \otimes A^{-1}$, so that in this simple finite-dimensional case we have not generated any new 'interesting' intertwiners. The point of this example is simply to demonstrate how the language of perturbative quantum field theory in terms of actions and feynman rules is translated into our graphical calculus language inside representation categories.

### 2.6.5 Feynman diagrams in QED

In this language, to specify a quantum field theory is to specify a symmetry group $G$, as well as a few basic intertwiners (interactions) inside Rep($G$). These are the Feynman rules. Feynman diagrams are then built up from these basic interactions. One could therefore say that the physics takes place inside the subcategory of Rep($G$) of all diagrams built up from these interactions. Strictly speaking, one needs the representations to be finite dimensional, since the loops contribute a term dim($V$), as we have seen in Section 2.3.3. Let us ignore this for the moment and consider quantum electrodynamics. The external symmetry group is the Poincaré group $P$. There are two kinds of physically relevant irreps of $P$ : the massive irreps are classified by pairs of numbers $(m, j)$ where $m > 0$ is the mass and $j = 0, \frac{1}{2}, 1, \frac{3}{2}, \ldots$ is the spin. The massless irreps have $m = 0$ but have helicity specified by numbers $(j, 0)$ or $(0, j)$ since they have a handedness. In QED, there are two kinds of particles (irreps):

$$\text{photons} \to \quad m = 0, \quad \text{helicity} = (1,0) \oplus (0,1), \tag{2.74}$$

$$\text{electrons/positrons} \to \quad m = m_e \sim 0.511 MeV/c^2, \quad \text{spin} = \frac{1}{2}. \tag{2.75}$$

[19]Here we are using the basis $\{e_i\}$ to identify $V$ with $V^*$, so that $i_V$ is a map $i_V : 1 \to V \otimes V$.

The photon representation consists of solutions of the Maxwell equations, while the electron/photon representations consists of solutions to the Dirac equation. Both these spaces are infinite dimensional. We draw these reps as:

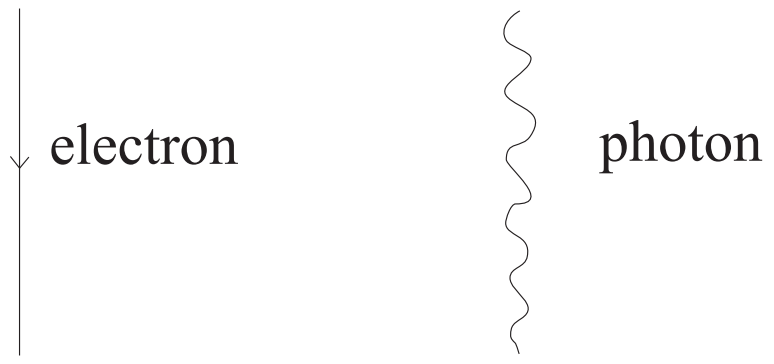

Note that the photon rep is self-dual, so that it equals its own antiparticle. The basic intertwiner (Feynman rule) is:

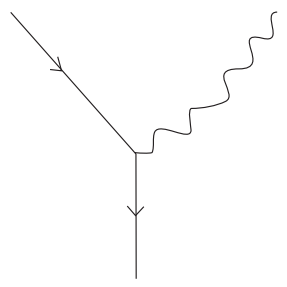

(2.76)

Using (2.76), one can build up arbitrarily complicated Feynman diagrams:

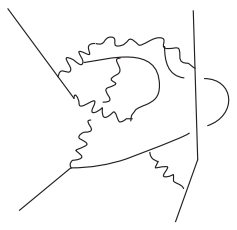

(2.77)

Unfortunately the loops contribute infinities to these expressions, which requires renormalization. Nevertheless the advantage of infinite dimensional representations is that the propagators need not be trivial, since intertwiners of irreducible representations need not be scalars.

### 2.6.6 A new look at quantum field theory

We have motivated that, *perturbatively speaking* (that is, in terms of Feynman diagrams), a quantum field theory is nothing but a finitely generated subcategory $QFT$ of the ribbon category $\mathrm{Rep}(G)$. $QFT$ is generated by the fundamental particles (irreducible representations) and all possible combinations of the fundamental interactions coming from the Feynman rules (intertwiners). Thus one is led to consider generalized quantum field theories $QFT'$ living inside arbitrary Hermitian[20] ribbon categories $\mathcal{R}$, where the braiding and twist need not be trivial. In fact, Fuchs *et. al.* [35, 36, 37, 38, 39] have shown that a fully fledged conformal field theory (with boundary conditions) is precisely such a gadget! (See Sec. 5.1.1). This is the reason why we have spent so much time developing ribbon categories in this chapter. Now Tannaka-Krein duality tells us that one can recover the group $G$ from the category $\mathrm{Rep}(G)$. Majid has shown [57] that this is just a special case of a more general reconstruction theorem. In a certain sense, *every* ribbon category is a category of representations - in the general case not of a group, but of a *quantum group.* When we do quantum field theory in ribbon categories, we are replacing the symmetry group by a quantum group.

We shall encounter this phenomenon in Chern-Simons theory, where the Lie group $G$ is replaced by its quantum deformation, $U_q(\mathfrak{g})$. The reason for this deformation of the underlying symmetry group, as one passes from the classical to the quantum theory, has not been altogether elucidated, and remains an interesting problem. In Witten's approach [98], three dimensional Chern-Simons theory defines a two dimensional conformal field theory on the boundary, the *Wess-Zumino-Witten* (WZW) model. The corresponding affine lie algebra $\hat{\mathfrak{g}}$ of the WZW model defines, for each $k \in \mathbb{Z}_+$, a category $\mathcal{C}_k(\hat{\mathfrak{g}})$ of integrable modules of level $k$, and these categories are modular [18].

On the other hand, in Turaev's approach [89], one deforms the lie algebra $\mathfrak{g}$ into a quantum group $U_q(\mathfrak{g})$, where $q = e^{\frac{\pi i}{k}}$, which for $k \in \mathbb{Z}_+$ is a root of unity. The representation categories $\mathrm{Rep}(U_q(\mathfrak{g}))$ of these quantum groups are also modular, and are the starting point in Turaev's approach.

[20] We need them to be Hermitian in order to order to accommodate time reversal, a very important physical concept. We should also require them to be unitary, in order to have bounded ground state energies. But even some 'ordinary' QFT's break this rule!

It is an important fact, proved by Finkelberg [29], that these two modular categories are equivalent (we adopt here the formulation taken from [18]):

**2.6.2 Finkelberg, 1996.** *The category $\mathcal{C}_k(\hat{\mathfrak{g}})$ is equivalent to the category Rep$(U_q(\mathfrak{g}))$ as a modular tensor category for $q = k + h^\vee$, where $h^\vee$ is the dual Coxeter number for $\mathfrak{g}$.*

Despite this theorem, the relationship between the Witten and Turaev approaches is still not completely understood.

Ordinary Lie groups are the symmetry groups of manifolds. Quantum groups are the symmetry groups of noncommutative spaces - deformed, noncommutative versions of the commutative algebra of functions on a manifold. Thus the process of passing from $QFT$ to $QFT'$ is associated with the philosophy of noncommutative geometry, a recent trend in physics. Oeckl has shown [70] how to define a path integral formalism based on Gaussian integration in general braided categories, including a generalization of Wick's theorem. Bosonic and Fermionic path integrals and Feynman rules can be recovered as special cases.

# Chapter 3

# 2d Topological Quantum Field Theories

## 3.1 Definition of topological quantum field theory

Consider again the fundamental formula of quantum field theory (Eqn. 1.5):

$$\boxed{\langle \hat{A}_2|U|\hat{A}_1\rangle = \int_{A|_{\Sigma_1}=A_1}^{A|_{\Sigma_2}=A_2} \mathcal{D}\mathcal{A} \exp iS[A]} \tag{1.5}$$

We have now developed the appropriate mathematical language to understand this equation. It tells us that a quantum field theory is a symmetric monoidal functor $Z : \mathbf{nCob_{metric}} \to \mathbf{Vect}$. It is possible that the theory is actually independent of the metric, and we thus define a *topological* quantum field theory (TQFT) as a symmetric monoidal functor

$$\tilde{Z} : \mathbf{nCob} \to \mathbf{Vect}. \tag{3.1}$$

A TQFT is thus a *representation* of **nCob** in terms of vector spaces and linear maps. In particular, we can consider a closed $n$-manifolds $M$ as a cobordism $M : \phi \to \phi$, which under the TQFT goes into a map $Z(M) : \mathbb{C} \to \mathbb{C}$, which is simply a complex number. That is, TQFT's give us topological invariants - these are simply the partition functions of the quantum field theory.

Eqn. (3.1) is the most elegant definition of a TQFT, but Atiyah's original definition [7], which is practically equivalent to (3.1) and is given below, is still often used in practice because it avoids the use of cobordisms and works directly with the manifolds. Cobordisms are indeed elegant but they come at a price of degeneracy - the same $n$-manifold $M$ with boundary $\partial M$ can be viewed as a cobordism in many different ways, depending on how one chooses the input and output boundaries in $\partial M$. Atiyah's definition is also closer to the kind of language used in conformal field theory [64]. Nevertheless, once all the details have been checked, Eqn. (3.1) remains the best overall viewpoint.

### 3.1.1 Atiyah's definition

Atiyah [7] defined a $d$-dimensional topological quantum field theory (TQFT) $Z$, as consisting of the following data:

a) A vector space $Z(\Sigma)$ associated to each $(d-1)$ dimensional closed manifold $\Sigma$.

b) A vector $Z(M) \in Z(\partial M)$ associated to each oriented $d$-dimensional manifold $M$ with boundary $\partial M$.

This data is subject to the following axioms, which we state briefly and expand upon below:

(a) $Z$ is *functorial* with respect to orientation preserving diffeomorphisms of $\Sigma$ and $M$.

(b) $Z$ is *involutory*, i.e. $Z(\Sigma^*) = Z(\Sigma)^*$ where $\Sigma^*$ is $\Sigma$ with opposite orientation and $Z(\Sigma)^*$ is the dual vector space of $Z(\Sigma)$.

(c) $Z$ is *multiplicative*.

(d) $Z(\phi) = \mathbf{C}$, where $\phi$ is interpreted as an empty $(d-1)$-dimensional closed manifold.

(e) $Z(\phi) = 1$, where $\phi$ is interpreted as an empty $d$-dimensional manifold.

These axioms meant to be understood as follows. The functoriality axiom (a) means first that an orientation preserving diffeomorphism $f : \Sigma \to \Sigma'$ induces an isomorphism $Z(f) : Z(\Sigma) \to Z(\Sigma')$ and that $Z(gf) = Z(g)Z(f)$ for $g : \Sigma' \to \Sigma''$. Also if $f$ extends to an orientation preserving diffeomorphism $M \to M'$, with $\partial M = \Sigma$, $\partial M' = \Sigma'$, then $Z(f)$ takes the element $Z(M)$ to $Z(M')$. The involutory axiom (b) is clear, but it also demonstrates that one may consider TQFT's with target an arbitrary symmetric rigid monoidal category, since this is the only structure used here. The multiplicative axiom (c) asserts first that, for disjoint unions,

$$Z(\Sigma_1 \cup \Sigma_2) = Z(\Sigma_1) \otimes Z(\Sigma_2). \tag{3.2}$$

Moreover if $\partial M_1 = \Sigma_1 \cup \Sigma_3^*$, $\partial M_2 = \Sigma_2 \cup \Sigma_3^*$ and $M = M_1 \cup_{\Sigma_3} M_2$ is the manifold obtained by gluing together the common $\Sigma_3$-component:

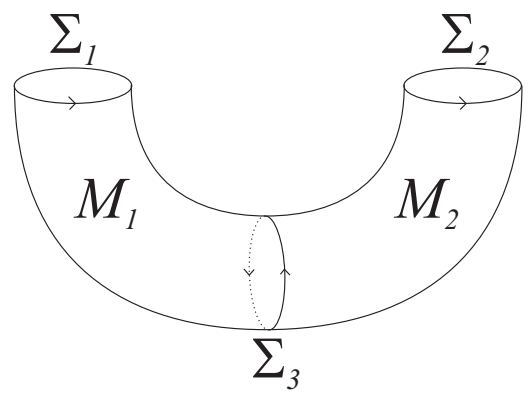

Then we require:

$$Z(M) = \langle Z(M_1), Z(M_2) \rangle \tag{3.3}$$

where $\langle\, ,\, \rangle$ denotes the natural pairing from the duality map,

$$Z(\Sigma_1) \otimes Z(\Sigma_3)^* \otimes Z(\Sigma_3) \otimes Z(\Sigma_2) \xrightarrow{\text{id} \otimes e_V \otimes \text{id}} Z(\Sigma_1) \otimes Z(\Sigma_2). \tag{3.4}$$

This is a very powerful axiom which implies that $Z(M)$ can be computed (in many different ways) by 'cutting $M$ in half' along any $\Sigma_3$. These are the gluing (or sewing) rules and 'lego game' familiar from conformal field theory. Axioms (d) and (e) are to prevent trivial theories.

### 3.1.2 Remarks on Atiyah's definition

Note that the word 'cobordism' does not enter into Atiyah's definition - it works directly with the manifolds and not on any 'man-made' interpretation of them. It is also physics oriented - the vector $Z(M) \in Z(\partial M)$ is the *vacuum state* defined by $M$. We can make contact between Atiyah's definition and (3.1) by noting that a cobordism $M : \Sigma_1 \to \Sigma_2$ can be interpreted as a decomposition of $\partial M$ into two components so that

$$\partial M = \Sigma_2 \cup \Sigma_1^*. \tag{3.5}$$

Then Atiyah's map $Z$ assigns an element $Z(M) \in Z(\Sigma_2)^* \otimes \Sigma_1 = \text{Hom}(Z(\Sigma_1), Z(\Sigma_2))$. In this way we get a functor[1] $\tilde{Z} : \mathbf{nCob} \to \mathbf{Vect}$, and the axioms (a)-(e) effectively state that $\tilde{Z}$ is a symmetric monoidal functor. Conversely, given a symmetric monoidal functor $\tilde{Z} : \mathbf{nCob} \to \mathbf{Vect}$, we get an Atiyah style TQFT Z by interpreting all our cobordisms $M : \Sigma_1 \to \Sigma_2$ as cobordisms $M : \phi \to \Sigma_1^* \amalg \Sigma_2$ and hence obtaining a map $\tilde{Z} : \mathbf{C} \to \tilde{Z}(\Sigma_1 \amalg \Sigma_2)$ which is the same as giving a vector $\tilde{Z}(1) \in \tilde{Z}(\Sigma_1 \amalg \Sigma_2)$.

When expressed as a cobordism, Atiyah's multiplicative axiom (c) shows that, for a cylinder $\Sigma \times I$, the linear map

$$Z(\Sigma \times I) \in \text{End}(Z(\Sigma)) \tag{3.6}$$

[1]It is not quite a functor until we have settled the subtlety regarding the 'identity' $Z(\Sigma \times I)$; see below.

is an idempotent $\sigma$ (it squares to 1), and more generally acts as the identity on the subspace of $Z(\Sigma)$ spanned by all elements $Z(M)$ with $\partial M = \Sigma$. If we replace $Z(\Sigma)$ by its image under $\sigma$, it is easy to see that the axioms are still satisfied. This is usually assumed as a further non-triviality axiom:

$$Z(\Sigma \times I) \text{ is the identity.} \tag{3.7}$$

With this proviso, Atiyah's definition is equivalent to (3.1), and we shall denote both by $Z$ from now on.

### 3.1.3 Formal properties of TQFT's

**nCob as a rigid, symmetric monoidal category**

**nCob** is a very interesting geometric category. It is at the very least a symmetric, rigid monoidal category, for all dimensions $n$. The tensor product of two $(n-1)$-dimensional manifolds $\Sigma_1$ and $\Sigma_2$ is just their disjoint union[2] $\Sigma_1 \sqcup \Sigma_2$, represented schematically as:

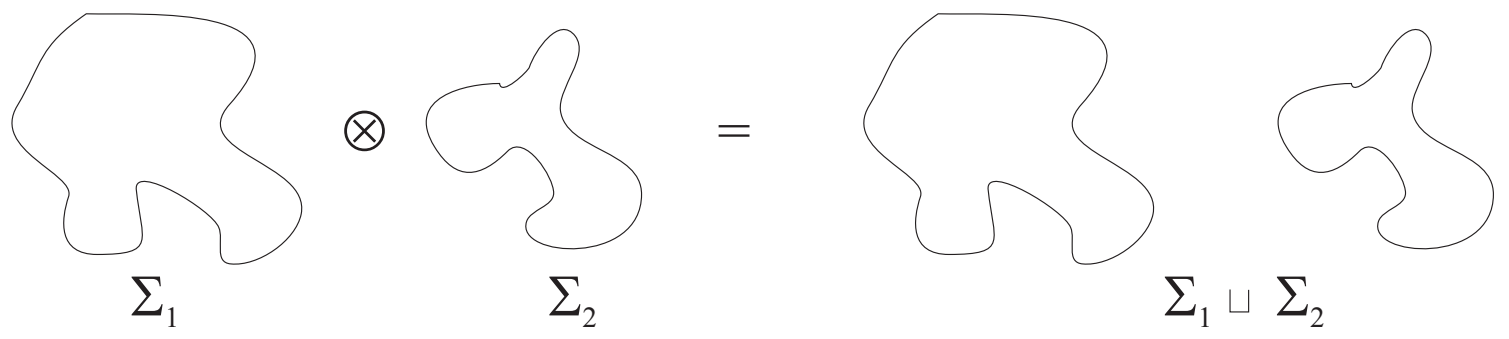

The unit 1 is the empty (n-1)-manifold $\phi$. The symmetry $\sigma_{\Sigma_1,\Sigma_2} : \Sigma_1 \sqcup \Sigma_2 \to \Sigma_2 \sqcup \Sigma_1$ is the cobordism which, as a manifold, is equal to $(\Sigma_1 \sqcup \Sigma_2) \times I$, but with injection maps which swap $\Sigma_1$ and $\Sigma_2$ in the output. It is drawn equivalently in the following two ways (we will stick to the former):

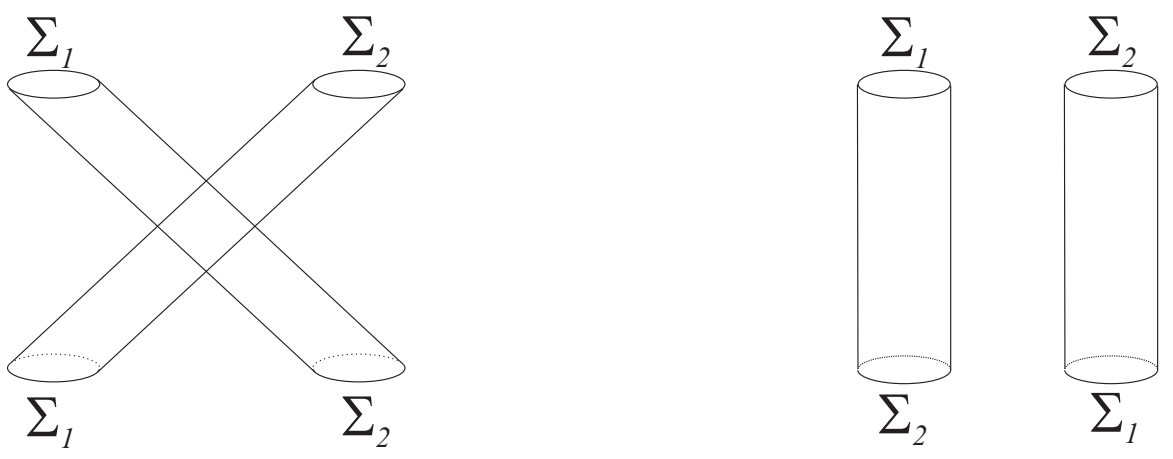

Recall that these manifolds are abstract and not embedded, it is misleading to distinguish over and under braidings. $\Sigma^*$ is defined as $\Sigma$ but with the reverse orientation. By viewing the cylinder $\Sigma \times I$ (with its canonical product orientation) as a cobordism from $\phi \to \Sigma \sqcup \Sigma^*$ or as a cobordism from $\Sigma^* \sqcup \Sigma \to \phi$ we obtain the duality maps $i_\Sigma$ and $e_\Sigma$ respectively:

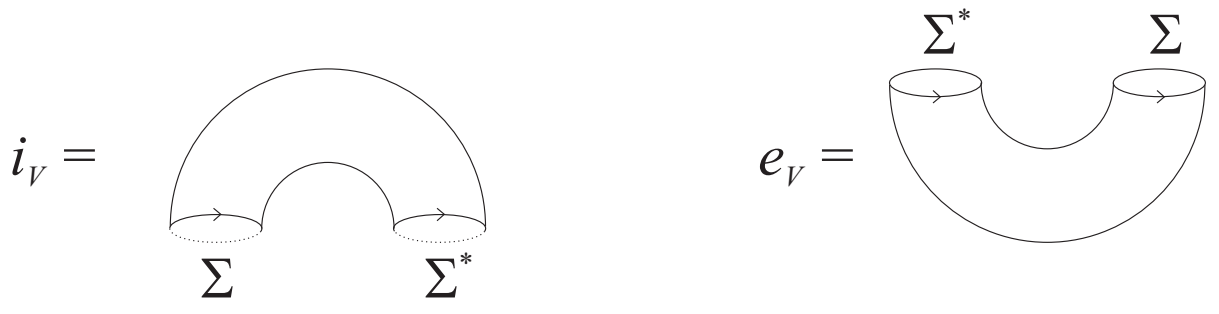

The snake relations follow since gluing three cylinders together is still a cylinder, eg.:

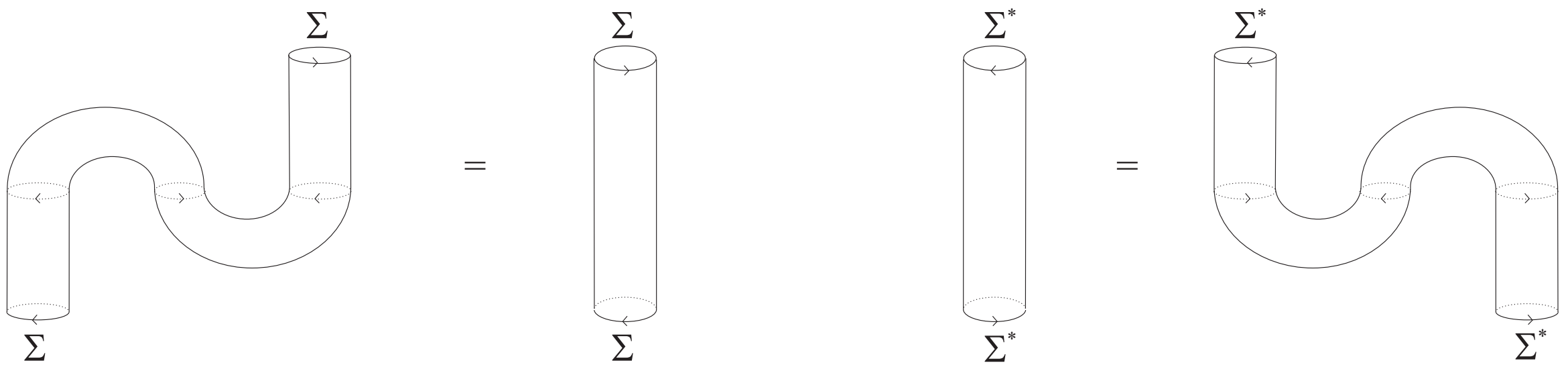

[2]The disjoint union of two sets, say $A = \{v, w, x, y\}$ and $B = \{x, y, z\}$ is defined by attaching labels and then forming the ordinary union, i.e. $A \sqcup B = \{v_A, w_A, x_A, y_A, x_B, y_B, z_B\}$.

Note that these duality maps are defined in all dimensions; drawing them as actual cylinders and circles is just schematic.

**nCob** is also a Hermitian category. The hermitian transpose of a cobordism $M : \Sigma_1 \to \Sigma_2$ is simply the cobordism $M^\dagger : \Sigma_2^* \to \Sigma_1^*$ obtained by reversing the orientation on $M$ (and hence also on its boundary $\partial M$):

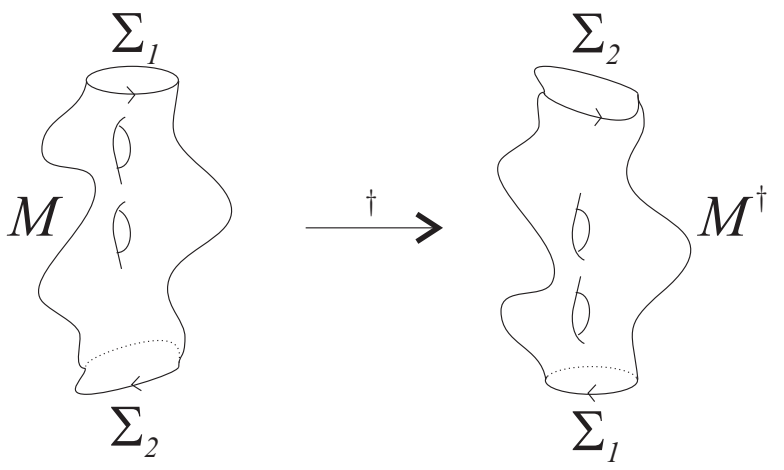

This is the time-reversal operation in **nCob**.

### Action of the mapping class group

Given a diffeomorphism $f : \Sigma_1 \xrightarrow{\sim} \Sigma_2$, we can construct a canonical cobordism $M_f : \Sigma_1 \to \Sigma_2$ by the cylinder construction:

$$\Sigma_1 \xrightarrow{\text{id}} M \times I \xleftarrow{f^{-1}} \Sigma_2 \tag{3.8}$$

We have the following important result.

**3.1.1 Lemma.** *Two diffeomorphisms* $\Sigma_1 \underset{f'}{\overset{f}{\longrightarrow}} \Sigma_2$ *induce the same cobordism class* $M : \Sigma_1 \to \Sigma_2$ *if and only if they are homotopic.*

*Proof.* Suppose we have two homotopic diffeomorphisms $\Sigma_1 \underset{f'}{\overset{f}{\longrightarrow}} \Sigma_2$. Let $F$ be the homotopy, i.e. $F$ is a map from the cylinder $\Sigma_1 \times I$ to $\Sigma_1$ which agrees with $f$ on one end of the cylinder and with $f'$ on the other. It is not hard to convince oneself that this is simply another notation for the statement that $M_f$ and $M_f'$ are equivalent cobordisms. □

For a closed manifold $\Sigma$, its *mapping class group* $\Gamma(\Sigma)$ is defined as the group of isotopy classes of diffeomorphisms $\Sigma \xrightarrow{\sim} \Sigma$. This lemma shows that the assignment $[f] \to M_f$, considered as a map $\Gamma(\Sigma) \to \mathrm{Hom}_{\mathbf{nCob}}(\Sigma, \Sigma)$, is well-defined and injective. Thus we see that a TQFT $Z$ : **nCob** $\to$ **Vect** gives a representation of the mapping class group for each $(n-1)$-manifold $\Sigma$. This is a very important property of TQFT's.

### Hermitian TQFT's

A Hermitian TQFT is one that satisfies $Z(M^\dagger) = Z(M)^\dagger$. In particular this means that the topological invariants $Z(M)$ change to their conjugates under orientation reversal. Unless all values are real then these theories can 'detect' conjugation. Such a theory is also called a *unitary* TQFT[3], since the amplitude for the time reversed process is simply the conjugate of the original process. All the TQFT's of interest in physics have this property.

An equivalent description of a Hermitian TQFT, is that it is a TQFT where *the inner product on the (kinematical) Hilbert spaces is determined by the dynamics of the fields themselves.* That is, the inner product on the space $Z(\Sigma)$ (which need have no relation to the dynamics of the theory) is precisely equal to the map $Z($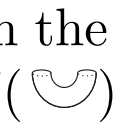$)$.

[3]This concept is compatible with the notion of 'unitary' defined in Section 2.5 if we convert the Hom spaces in **nCob** into free vector spaces over the cobordisms.

**Physical picture of TQFT's**

Let us remind the sceptical reader of the way in which the abstract, categorical notion of a TQFT as a unitary symmetric monoidal functor relates to fundamental physical principles[4]. The fact that $Z$ is a *functor* (i.e. $Z(M'M) = Z(M')Z(M)$ and $Z(\Sigma \times I) = \mathrm{id}$) means that the passage of time corresponding to the cobordism $M$, followed by the passage of time corresponding to the cobordism $M'$, has the same effect as the passage of time corresponding to the cobordism $M'M$. Moreover a passage of time in which no topology change occurs has no effect at all on the state of the universe - TQFT's have no local degrees of freedom.

The fact that $Z$ is *monoidal* means that the space of states corresponding to two non-interacting systems is given by the tensor product of the space of states associated to each individual systems, a familiar rule in quantum mechanics.

The fact that $Z$ is a *symmetric* functor is related to the *statistics* of the particles we are dealing with. It means that if we interchange two states $\phi$ and $\psi$, then on the vector space level this corresponds to sending $\phi \otimes \psi \to \psi \otimes \phi$. In other words, we are dealing with *bosons*. We can accommodate fermions by allowing for grassmann numbers in our vector spaces. Mathematically, all this means is that we pass from the category of ordinary vector spaces $V$ to the category **grVect** consisting of graded vector spaces $V = \bigoplus_n V_n$ where the grading is the fermionic degree. The symmetry map, instead of $\phi \otimes \psi \to \psi \otimes \phi$, is now given by

$$\phi \otimes \psi \to (-1)^{pq} \psi \otimes \phi, \tag{3.9}$$

where $\deg(\phi) = p$ and $\deg(\psi) = q$. That is, a fermionic TQFT is a symmetric monoidal functor $Z : \mathbf{nCob} \to \mathbf{grVect}$. What about more exotic statistics, like anyons? In this case, we need to consider an embedded version of **nCob** like **TubeCob** (which we shall define later), which is no longer symmetric but genuinely braided. We would then consider a TQFT as a *braided* monoidal functor $Z : \mathbf{TubeCob} \to \mathbf{Vect}_q$ where the number $q$ is the anyonic parameter which labels the way in which the braiding is implemented in $\mathbf{Vect}_q$:

$$\phi \otimes \psi \to e^{2\pi i q} \psi \otimes \phi. \tag{3.10}$$

The fact that $Z$ is *unitary* means that the time evolution operator corresponding a time-reversed process should be the conjugate $T^\dagger$ of the time evolution operator corresponding to the original process. This does not mean that time evolution operators are unitary operators! In fact, it has been realized by many authors that unitary time evolution is not a built-in feature of quantum theory but rather a consequence of specific assumptions about the nature of spacetime [9, 6]. In a Hermitian category, we can define a morphism $f : x \to y$ to be unitary if $f^\dagger f = 1_x$ and $f f^\dagger = 1_y$. In **Hilb**, this is the usual definition of unitary maps. In **nCob**, it is clear that cobordisms which are cylinders over some space manifold $\Sigma$ will be unitary. For dimensions $n \leq 3$ these are in fact *all* the unitary cobordisms, but interestingly enough, there are non-cylinder-like unitary cobordisms in dimensions $n \geq 4$. A unitary TQFT will thus (at least) map unitary morphisms in **nCob** to unitary operators in **Hilb** - in other words, if there is no topology change, then time evolution is unitary. In general though, time evolution is given by a non-unitary operator.

## 3.2 Closed and open strings

In two dimensions, a TQFT is similar to a string theory. More precisely, it behaves like a *closed* string theory, because the objects of **2Cob** are circles. To deal with open strings, we should define a variant of **2Cob** which has disjoint unions of intervals as objects, and 'cobordisms between intervals' as morphisms. We shall call this category **OCob**, and a working definition is as follows. Objects are integers $n$ representing the number of open strings serving as input or output. A morphism $M : n_1 \to n_2$ is a planar submanifold $M \subset \mathbb{R} \times [0, 1]$ such that $\partial M = I_n \times \{0\} \cup I_m \times \{1\}$, where

$$I_k = [1 - \frac{1}{3}, 1 + \frac{1}{3}] \cup \cdots \cup [k - \frac{1}{3}, k + \frac{1}{3}]. \tag{3.11}$$

Cobordisms are considered equivalent up to isotopy, i.e. smooth deformations of $\mathbb{R} \times [0, 1]$ (which are identity on the boundary) taking the one into the other. Composition is defined by gluing the output of the first

[4]I am indebted to Baez's interesting treatment of this subject in [9].

cobordism to the input of the second, and rescaling appropriately. We shall often use artistic license to place the inputs and outputs more symmetrically:

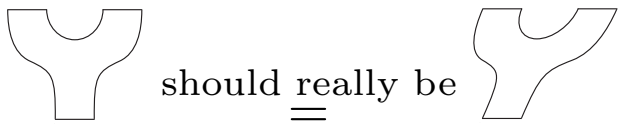

(3.12)

The tensor product $M \otimes M'$ of two cobordisms is defined by placing $M'$ to the right of $M$, adjusting its input and output slots accordingly:

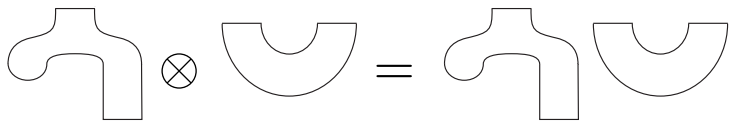

(3.13)

The unit object is 0. These constructions define **OCob** as a monoidal category. In physics terms, **OCob** is the category of worldsheets of an open string embedded (no self-intersections) in the plane. To get a feeling for **OCob**, here are examples of cobordisms from $2 \to 3$ and $3 \to 4$ respectively:

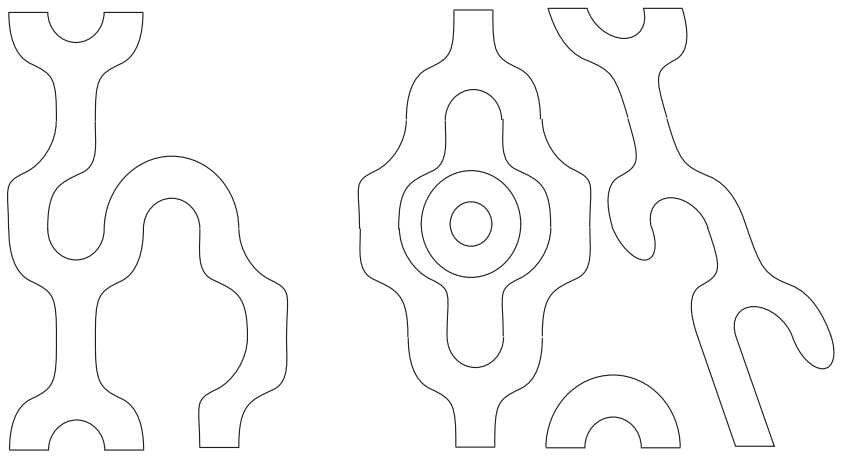

(3.14)

**OCob** is very similar to **2Cob**, the difference being that the former is embedded in the plane, while the latter is abstract with labeled boundaries, and hence can accommodate the twist map:

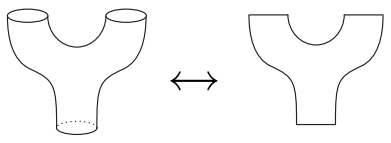 etc., but 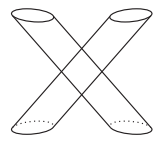 has no counterpart in **OCob**. (3.15)

Thus **2Cob** is a symmetric monoidal category, but **OCob** is not. We shall occasionally use the terms 'closed string theory' to mean an ordinary TQFT (symmetric monoidal functor) $Z : \mathbf{2Cob} \to \mathbf{Vect}$, while an 'open string theory' shall mean a monoidal functor $Z : \mathbf{OCob} \to \mathbf{Vect}$.

## 3.3 Frobenius Algebras

Since connected 2d manifolds are completely classified by their genus and the number of boundary circles, there is a hope to completely classify 2d TQFT's. Indeed this is the case, the main result being:

**3.3.1 Theorem.** *To give an open string theory is equivalent to giving a Frobenius algebra A inside* ***Vect****. To give a closed string theory is equivalent to giving a commutative Frobenius algebra B inside* ***Vect****.*

The algebra $A$ ($B$) is defined on the vector space which is the image under $Z$ of the interval $I$ (circle $S^1$). To prove that a open/closed string theory defines a Frobenius algebra on these vector spaces is easy, especially after one reformulates the definition of a Frobenius algebra in a categorical or 'topological' way. To prove the converse, that every Frobenius algebra arises as $Z(I)$ or $Z(S^1)$ for some open/closed TQFT $Z$ is the more interesting result[5]. To the author's knowledge, there are three different ways of proving this fact, and it is instructive to compare the different approaches.

The first and perhaps most modern way (elegantly set forth in Kock's recent book [50]) is to express **2Cob** and **OCob** using generators and relations, and to use a recent result of Quinn and Abrams [76, 3, 4] which formulates the axioms for a Frobenius algebra in exactly the same way. The second way (which can be found in [18]) is to use the Atiyah-style definition of a TQFT, where the burden of proof is to show that, given a Frobenius algebra $A$, one can define the vectors $Z(M) \in Z(\partial M)$ in a consistent way, i.e. the definition is independent of the cutting of $M$ into smaller pieces (this is called *consistency of the sewing* in conformal field theory). The third way (which to the author's knowledge has never been explicitly written down before, but has been implicitly suggested eg. by Moore [63]) is to take advantage of the fact that it is relatively harmless to consider 2d cobordisms as embedded inside $\mathbb{R}^3$. Then one can extend the graphical calculus ideas from Chapter 1 to show that, given a Frobenius algebra $A$, there is a 'deformation invariance' theorem exactly analogous to Theorem 2.3.3.

[5]Indeed, even the renowned mathematician Graeme Segal admitted not very long ago that he didn't know of an illuminating proof (see [79]).

### 3.3.1 Definition

Frobenius algebras are classical algebras that were once, shamefully, called 'Frobeniusean algebras' in honour of the Prussian mathematician Georg Frobenius[67, 68, 69]. They have many equivalent definitions; but before we list them it is worthwhile to record the following fact.

**3.3.2 Lemma.** *Suppose $A$ is an arbitrary vector space equipped with a bilinear pairing $( \, , \, ) : A \otimes A \to \mathbb{C}$. Then the following are equivalent:*

*(a) $A$ is finite dimensional and the pairing is nondegenerate; i.e. $A$ is finite dimensional and the map $A \to A^*$ which sends $v \to (v, \cdot)$ is an isomorphism.*

*(b) $A$ is self dual in the rigid monoidal sense; i.e. there exists a copairing $i : \mathbb{C} \to A \otimes A$ which is dual to the pairing $e : A \otimes A \to \mathbb{C}$ given by $e(a, b) = \epsilon(ab)$.*

*Proof.* (a) $\Rightarrow$ (b). Choose a basis $(e_1, \ldots, e_n)$ of $A$. Then by assumption the functionals $(e_i, \cdot)$ are a basis for $A^*$. Then there exist vectors $e^1, \ldots e^n$ in $A$ such that $(e_i, e^j) = \delta^i_j$. Define the copairing $i$ by setting $1 \to \sum_i e_i \otimes e^i$. Then a general vector $v = \lambda^i e_i$ goes through the composite $V \xrightarrow{i \otimes \mathrm{id}} V \otimes V \otimes V \xrightarrow{e \otimes id} V$ as:

$$v = \lambda^i e_i \to \lambda^i e_j \otimes e^j \otimes e_i \to \lambda^i e_j (e^j, e_i) = \lambda_i e_i = v. \tag{3.16}$$

Similarly, $w = \lambda_i e^i$ goes through the composite $V \xrightarrow{\mathrm{id} \otimes i} V \otimes V \otimes V \xrightarrow{e \otimes \mathrm{id}} V$ as:

$$w = \lambda_i e^i \to \lambda_i e^i \otimes e_j \otimes e^j \to \lambda_i (e^i, e_j) e^j = \lambda_i e^i = w. \tag{3.17}$$

(b) $\Rightarrow$ (a). The copairing $i$ singles out a vector in $A \otimes A$ by $1 \to \sum_i^n e_i \otimes e^i$ for some vectors $e_i, e^i \in A$ and some number $n$ (note that we have not used finite dimensionality here). Now take an arbitrary $v \in A$ and send it through the composite $V \xrightarrow{i \otimes \mathrm{id}} V \otimes V \otimes V \xrightarrow{e \otimes id} V$:

$$v \to e_i \otimes e^i \otimes v \to e_i (e^i, v) \tag{3.18}$$

By assumption this must be equal to $v$. This shows that $(e_1, \ldots, e_n)$ spans $A$, so $A$ is finite dimensional. Now we show that $v \to (v, \cdot)$ is injective, and hence an isomorphism. Suppose $(v, \cdot)$ is the zero functional. Then in particular $(v, e^i) = 0$ for all $i$. But these scalars are exactly the coordinates in the 'basis' $(e_1, \ldots, e_n)$, so that $v = 0$. $\square$

This lemma translates the algebraic notion of nondegeneracy into category language, and from now on we shall use the two meanings interchangeably. It also makes explicit that a nondegenerate pairing allows one to construct, from a basis $(e_1, \ldots, e_n)$ for $A$, a corresponding *dual basis* $(e^1, \ldots, e^n)$, which satisfies $e(e_i, e^j) = \delta^j_i$, and which can be recovered from the decomposition $i(1) = \sum_i e_i \otimes e^i$.

**3.3.1 Defn.** *A Frobenius algebra is:*

*(a) A finite dimensional algebra $A$ equipped with a nondegenerate form (also called trace) $\epsilon : A \to \mathbb{C}$.*

*(b) A finite dimensional algebra $(A, \beta)$ equipped with a pairing $\beta : A \otimes A \to \mathbb{C}$ which is nondegenerate and associative[6].*

*(c) A finite dimensional algebra $(A, \gamma)$ equipped with a left algebra isomorphism to its dual $\gamma : A \to A^*$.*

[6] Recall that a pairing $(\cdot, \cdot)$ on an algebra is called associative when $(ab, c) = (a, bc)$ for all $a, b, c \in A$.

Observe that if $A$ is an algebra, then there is a one-to-one correspondence between forms $\epsilon : A \to \mathbb{C}$ and associative bilinear pairings $(\cdot,\cdot) : A \otimes A \to \mathbb{C}$. Given a form, define the pairing by $(a, b) = \epsilon(ab)$, this is obviously associative. Given the pairing, define a form by $\epsilon(a) = (1, a) = (a, 1)$; these are equal since the pairing is associative. This establishes the equivalence of (a) and (b). It is not difficult to show the equivalence of (a) and (c). However, (a) and (b) are more convenient for our purposes, since due to Lemma 3.3.2 they can be expressed entirely in terms of commuting diagrams inside **Vect**.

### 3.3.2 Graphical notation

We can take advantage of this last comment by developing a graphical notation, in the spirit of Chapter 1, for the various maps involved. We draw $\mathrm{id} : A \to A$, $\mu : A \otimes A \to A$, $\eta : \mathbb{C} \to A$, $\epsilon : A \to \mathbb{C}$, $i : \mathbb{C} \to A \otimes A$ and $e : A \otimes A \to 1$ as follows:

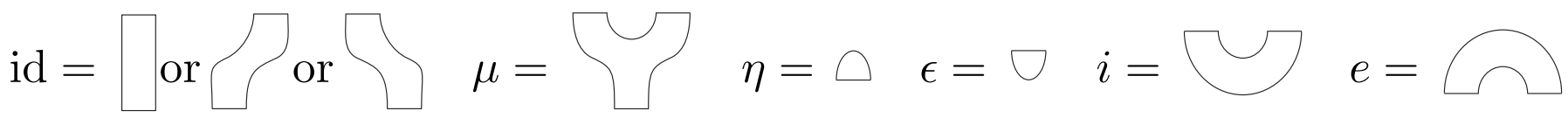

(3.19)

Then Defn. 3.3.1 (a) says that a Frobenius algebra $A$ is a vector space equipped with maps

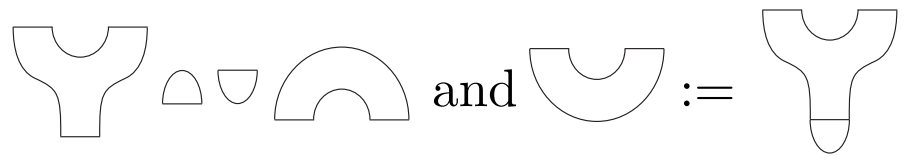

(3.20)

such that

(Associativity)

(Unit)

(Nondegeneracy)

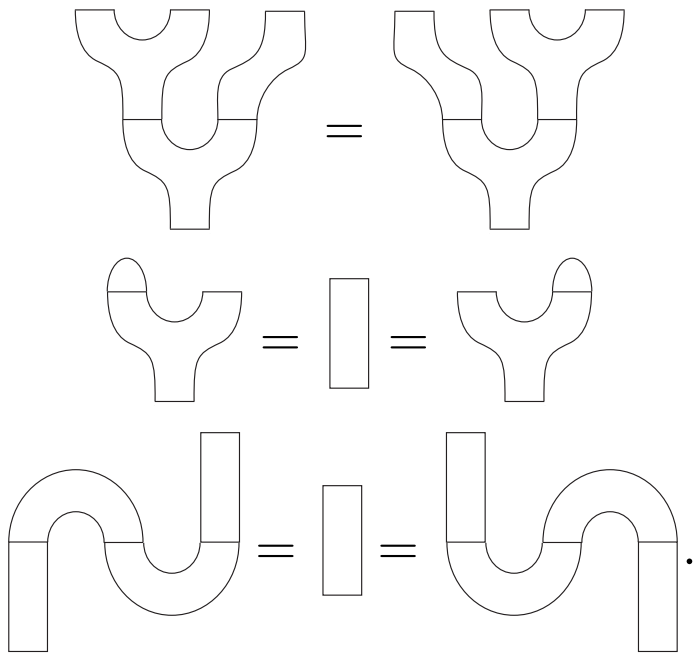

(3.21)

When we speak about commutative Frobenius algebras, then we shall used closed string notation. Thus, if we denote the symmetry which sends $v \otimes w \to w \otimes v$ in **Vect** by $\sigma$,

$\sigma =$ 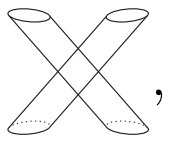 ,

(3.22)

then a commutative Frobenius algebra $B$ is a vector space equipped with maps

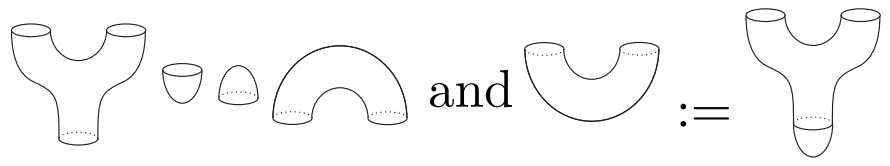

(3.23)

such that

(Associativity) 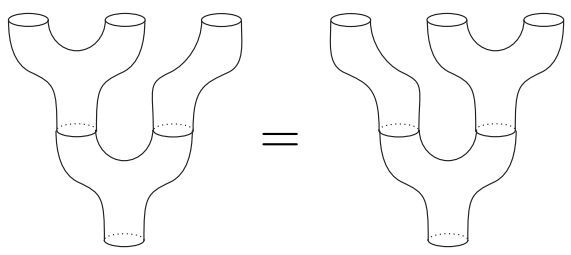

(Unit) 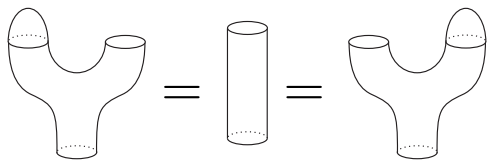

(Nondegeneracy) 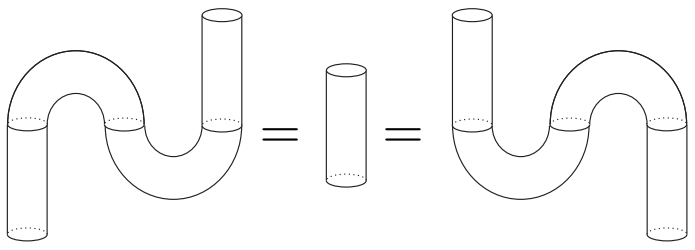

(Commutativity) 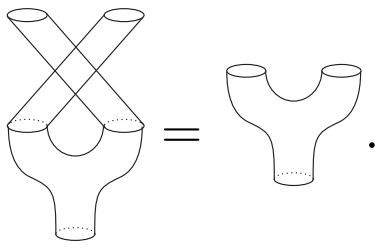

(3.24)

Expressing things this way, we immediately obtain the first part of Theorem 3.3.1:
**Theorem 3.3.1 (⇒)** : *An open string theory $Z$ gives a Frobenius algebra $Z(I)$ on the image of the interval $I$. A closed string theory $Z$ gives a commutative Frobenius algebra $Z(S^1)$ on the image of the circle $S^1$.*

*Proof.* We simply observe that all the relations above do indeed hold in **OCob** and **2Cob**, so that they must be preserved by the functor $Z$. □

### 3.3.3 Examples of Frobenius Algebras

- **Matrix algebras.** Here $A$ is the space $Mat_N(\mathbb{C})$ of $N \times N$ matrices. The trace $\epsilon$ is given by the actual trace Tr. This is nondegenerate because if $\mathrm{Tr}(BC) = 0$ for all $C$, then $B = 0$. The dual of $B_{ij}$ is $B_{ji}$, since $\mathrm{Tr}(B_{ij}Bji) = 1$.

- **Finite group algebras.** If $G$ is a finite group, then we form the group algebra $\mathbb{C}[G]$ as the vector space with basis $\{g\}, g \in G$, and with algebra multiplication given by $g \cdot h = gh$. A vector $f = \sum_g f_g g \in \mathbb{C}[G]$ can be regarded as a function $f : G \to \mathbb{C}$, i.e. $f \in C(G)$. The product in $\mathbb{C}[G]$, when expressed in terms of $C(G)$, becomes the convolution product,

$$(f * f')(g) = \sum_h f_g f_{h^{-1}g} \tag{3.25}$$

The trace, in group element notation, is $\epsilon(g) = \delta_g^e$ where $e$ is the identity of the group, while in function notation it is $\epsilon(f) = f(e)$, although one usually sets $\epsilon(f) = \frac{1}{|G|} f(e)$ to ensure that the total volume of the group is unity. The dual of the basis $(g)$ is $(g^{-1})$. The beautiful Peter-Weyl theorem tells us that this example is simply a special case of the former, since $C(G)$ splits up as a direct sum of matrix algebras,

$$C(G) = \bigoplus_{\rho \in \hat{G}} \mathrm{End}(V_\rho), \tag{3.26}$$

where $\hat{G}$ is the (finite) set of irreducible representations of $G$. The map $\mathrm{End}(V_\rho) \to C(G)$ is given by

$$A \to \mathrm{Tr}(A\rho(\cdot)), \tag{3.27}$$

under which $\epsilon \to |G|\mathrm{Tr}$.

- **Class functions on a group.** An important commutative Frobenius algebra is the center of the algebra of functions on the group, $Z(C(G))$. The Peter-Weyl theorem tells us that this space is spanned by the identity matrices $1_{V_\rho}$ for each representation $\rho$, which translate back via (3.27) into the characters[7] of the representations, $\chi_\rho \in C(G)$. The same Peter-Weyl theorem also tells us that the characters span the space of class functions on the group $G$, so we conclude that

$$Z(C(G)) = C_{\text{class}}(G) := \{f : G \to \mathbb{C} : f(h^{-1}gh) = f(g)\}. \tag{3.28}$$

On the other hand, $C_{\text{class}}(G) \cong Z(\mathbb{C}[G])$ since the latter has a basis $\{e_\alpha\}$ labeled by conjugacy classes $\alpha$ of $\mathbf{G}$, where

$$e_\alpha = \sum_{g \in \alpha} g. \tag{3.29}$$

The isomorphism sends $e_\alpha \to f_\alpha$ where $f_\alpha$ is equal to 1 on $\alpha$ and zero elsewhere. The multiplication (fusion rules) computes as $e_\alpha e_\beta = N^\gamma_{\alpha\beta} e_\gamma$ where

$$N^\gamma_{\alpha\beta} = |\{h \in \beta : gh^{-1} \in \alpha \ \text{ for } g \in \gamma\}|. \tag{3.30}$$

It is natural to adopt the same normalization for the Frobenius form as before, i.e.

$$\epsilon(f) = \frac{1}{|G|} f(e). \tag{3.31}$$

This Frobenius algebra forms a central part of this thesis.

- **2d Landau Ginzburg Models.** These are field theories described in terms of a set of chiral $N = 2$ superfields $x_i$ and $\overline{x}_i$ for $i = 1 \dots n$. The action takes the form

$$S = \int d^z d^4\theta K(x_i, \overline{x}_i) + \int d^2 z d^2\theta W(x_i) + c.c. \tag{3.32}$$

where $K(x_i, \overline{x}_i)$ is called the kinetic term and $W(x_i)$ is a holomorphic function called the superpotential. This theory is a conformal field theory whose chiral primary ring turns out to be

$$\mathcal{R} = \mathbb{C}[x_1, \dots, x_n]/(dW) \tag{3.33}$$

where $\mathbb{C}[x_1, \dots, x_n]$ is the ring (algebra) of polynomials in $x_1, \dots, x_n$ and $(dW)$ is the ideal generated by the derivatives $\{W_i = \frac{\partial W}{\partial x_1}\}$ of $W$. The trace $\epsilon$ of a polynomial $g$ is defined as the residue obtained by integrating around a ball $B$ around the origin:

$$g \to \int_B \frac{g(x_1, \dots, x_n) dx_1 \wedge \dots \wedge dx_n}{W_1(x) \cdots W_n(x)} \tag{3.34}$$

When $W(x_i)$ is a generic nth order polynomial then $\mathcal{R}$ is a finite dimensional algebra called the *nth perturbed minimal model*. In this case it turns out that[8][9]

$\mathcal{R}$ is semisimple $\Leftrightarrow$ the zeros of $W_i$ are simple $\Leftrightarrow$ the infrared fixed point of the system is a massive theory.

- **Cohomology rings.** For $M$ a compact, oriented manifold, the De Rham cohomology $H^*(M) = \oplus_{i=0}^n H^i(M)$ is an algebra under the wedge product. It is graded in the sense that if $\alpha \in H^p(M)$ and $\beta \in H^q(M)$ then $\alpha \wedge \beta \in H^{p+q}(M)$. The trace $\epsilon$ is integration over $M$ (with respect to a chosen volume form). The corresponding bilinear form $(\cdot,\cdot) : H^*(M) \otimes H^*(M) \to \mathbb{R}$ is known as the intersection form. The fact that it is nondegenerate is precisely Poincaré duality. If $(e_{p,i})$ is a basis for $H^*(M)$ with $e_{p,i} \in H^p(M)$, then the dual basis are the Poincaré duals $e^i_{N-p}$. This algebra is not commutative but it is graded commutative, $\alpha \wedge \beta = (-1)^{pq} \beta \wedge \alpha$.

  If $M$ is in addition a symplectic manifold, then one can define a deformation $Q$ of $H^*(M)$ known as *quantum cohomology*, which is also a Frobenius algebra. Quantum cohomology is inspired from string theory where it has to do with Gromov-Witten invariants and $\sigma$ models [41, 96].

[7] Recall that the characters are defined as $\chi_\rho(g) = \text{Tr}(\rho(g))$.

[8] See Section 3.3.4 for a definition of semisimple Frobenius algebras.

[9] Recall that a zero $a$ of a complex function $f$ is called *simple* if $f$ can be written as $f(z) = (z-a)g(z)$ where $g$ is a holomorphic function such that $g(a)$ is not zero.

### 3.3.4 Semisimple Frobenius algebras

A semisimple algebra is usually defined as an algebra containing no nontrivial nilpotent ideals. Wedderburn's theorem informs us that a semisimple algebra is isomorphic to a direct sum of matrix algebras. In particular, this means that semisimple commutative Frobenius algebras are isomorphic to a direct sum of the complex numbers, i.e.

$$A = \oplus_i \mathbb{C} e_i \tag{3.35}$$

where $e_i e_j = \delta_{ij}$, so the only degrees of freedom are the dimension of the algebra and the traces $\epsilon_i := \epsilon(e_i) \neq 0$. An alternative interpretation of Eqn. (3.35) is that a semisimple commutative algebra is one where the 'fusion rules' $e_i * e_j = \sum_k \mu_{ij}^k e_k$ are diagonalizable. The only way they can fail to be diagonalizable is if some of the multiplication matrices $(\mu_i) : e_j \to e_i e_j$ are nilpotent.

As an aside, Abrams has proved the following elegant theorem [3] relating semisimplicity to the graphical notation:

**3.3.3 Theorem [Abrams .]** *A Frobenius algebra A is semisimple iff the handle element*

$\omega :=$ 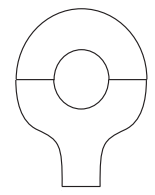

*is invertible.*

If a 2d TQFT is unitary, then the handle operator $H : A \to A$ which sends $v \to \omega v$ is Hermitian, and hence diagonalizable. Durhuus and Jonsson have used this idea [26] to classify 2d unitary TQFT's as those of the form (3.35) where $e_i$ are forced to be positive real numbers. Thus a unitary 2d TQFT results in an exceedingly simple algebra!

Of the examples listed above, matrix algebras, group algebras and the space of class functions on a group are all semisimple. Landau-Ginsburg models are generically semisimple when they are massive theories, as discussed above. Cohomology rings are *not* semisimple, because almost the entire algebra is nilpotent - simply wedge a form with itself enough times, and it will be zero. *Quantum* cohomology is, however, generically semisimple.

### 3.3.5 Proof of Theorem 3.3.1 (I)

**Coalgebras.**

An algebra can be expressed in terms of commuting diagrams as a triple $(A, \mu, \eta)$ consisting of a vector space $A$, a multiplication map $\mu : A \otimes A \to A$ and a unit map $\eta : \mathbb{C} \to A$ satisfying:

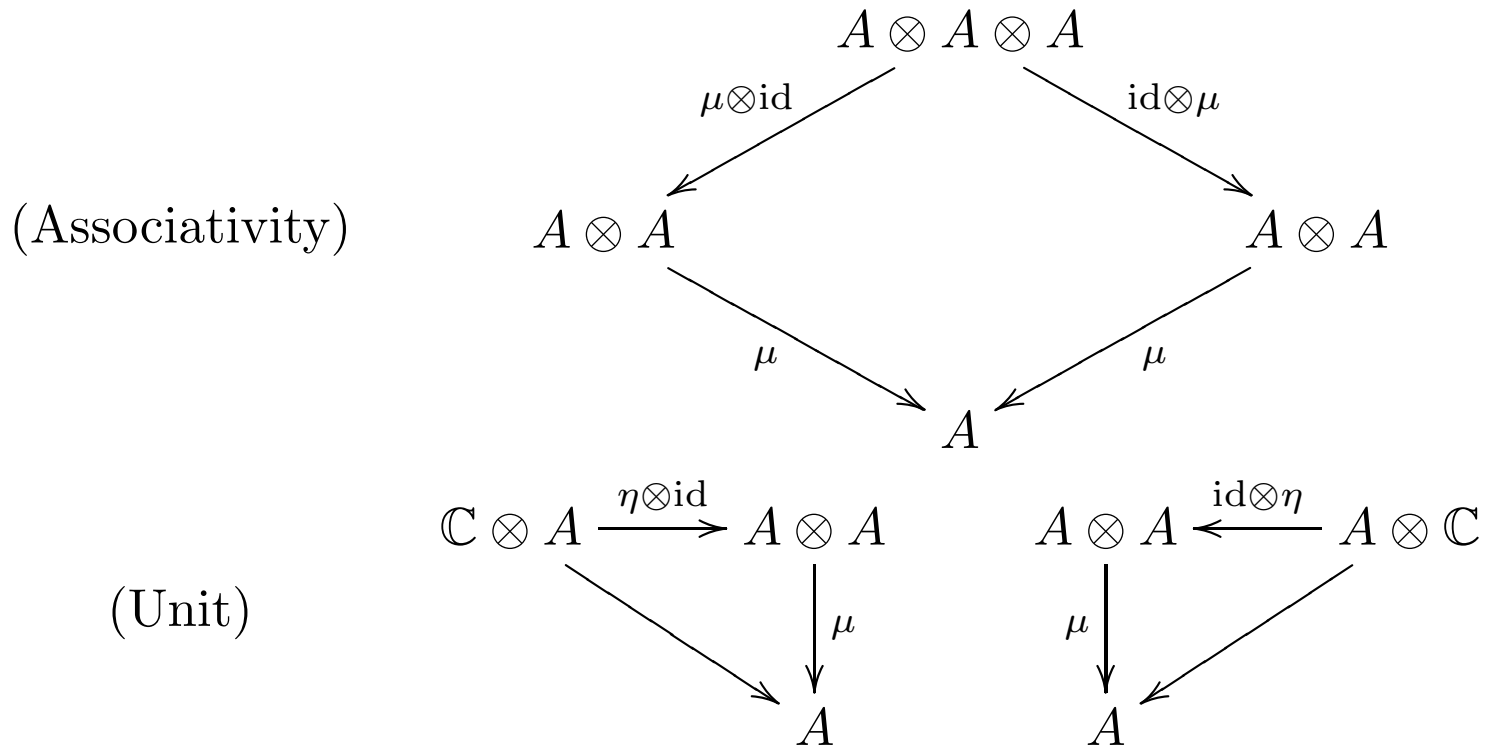

The multiplication $\mu$ sends $a \otimes b \to ab$ and the unit $\eta$ sends $1_{\mathbb{C}} \to 1_{\mathbb{A}}$. Similarly, a *coalgebra* is defined as the structure obtained when all the arrows above are reversed, i.e. it is a triple $(A, \Delta, \epsilon)$ consisting of a vector space $A$, a *co*multiplication map $\Delta : A \to A \otimes A$ and a *co*unit $\epsilon : A \to \mathbb{C}$ satisfying:

(Coassociativity)

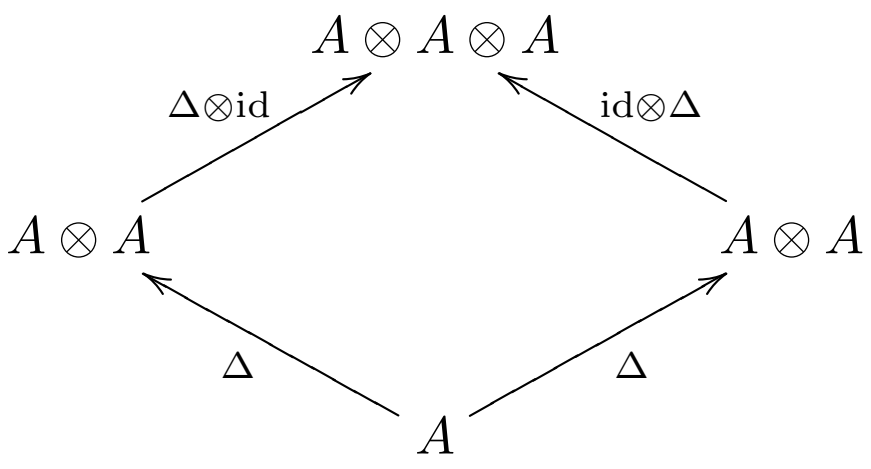

(Counit)

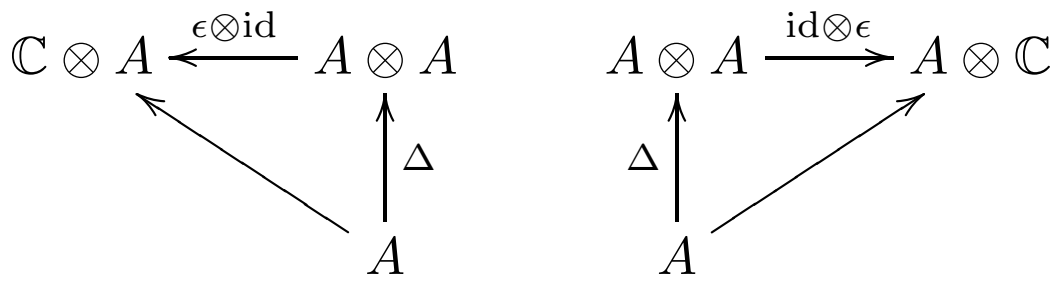

The great thing about phrasing algebras and coalgebras in this way is that they make sense in any monoidal category, allowing us to define *monoids* and *comonoids* as objects in monoidal categories equipped with morphisms having the properties listed above.

**Another definition of a Frobenius algebra**

We can now formulate the following more modern definition of a Frobenius algebra [76, 3, 4].

**3.3.4 Theorem.** *The following are equivalent:*

*(a) A Frobenius algebra (i.e. an algebra A equipped with a nondegenerate trace* $\epsilon : A \to \mathbb{C}$*).*

*(b) A finite-dimensional vector space A equipped with maps* $\mu : A \otimes A \to \mathbb{C}$, $\Delta : \mathbb{C} \to A \otimes A$, $\eta : \mathbb{C} \to A$ *and* $\epsilon : A \to \mathbb{C}$ *such that*

$$\begin{aligned} \mu(\eta \otimes \mathrm{id}) = \mathrm{id} = \mu(\mathrm{id} \otimes \eta) \qquad & \text{(Unit)} \\ (\epsilon \otimes \mathrm{id})\Delta = \mathrm{id} = (\mathrm{id} \otimes \epsilon)\Delta \qquad & \text{(Counit)} \\ (\mathrm{id} \otimes \mu)(\Delta \otimes \mathrm{id}) = \Delta\mu = (\mu \otimes \mathrm{id})(\mathrm{id} \otimes \Delta) \qquad & \text{(Frobenius condition)} \end{aligned} \tag{3.36}$$

This theorem is most transparent in graphical notation. We draw the comultiplication $\Delta : A \to A \otimes A$ as:

$\Delta =$ 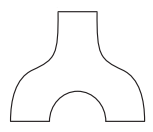 (3.37)

Then the theorem becomes:

**3.3.5 Theorem.** *The following are equivalent:*

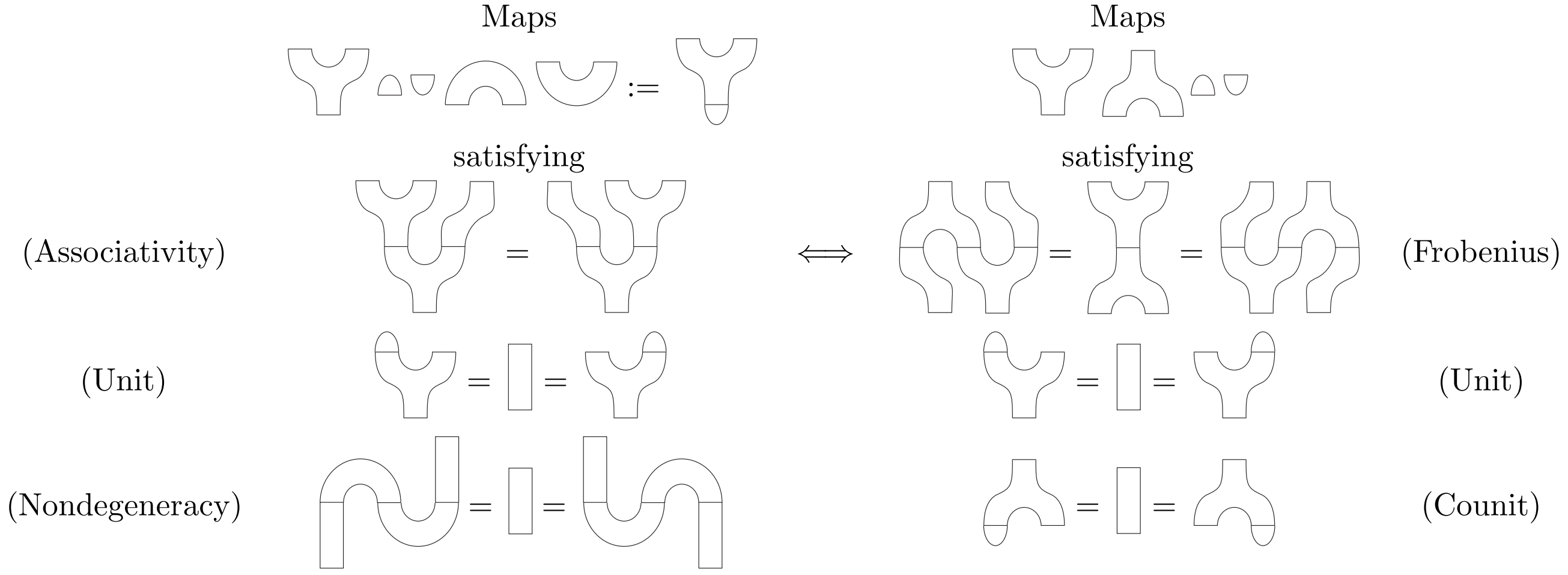

*Proof.* The proof is matter of fun and games, using the lego building blocks and the rules given above. For example, given the left-hand definition, we define a comultiplication as:

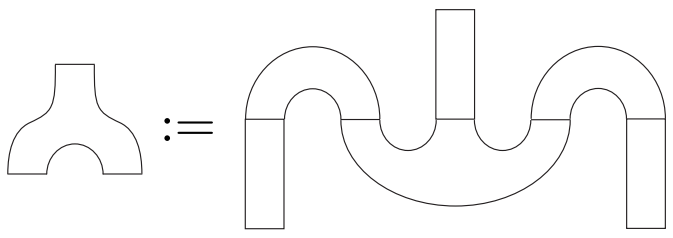

(3.38)

where we have introduced the 3-point function as:

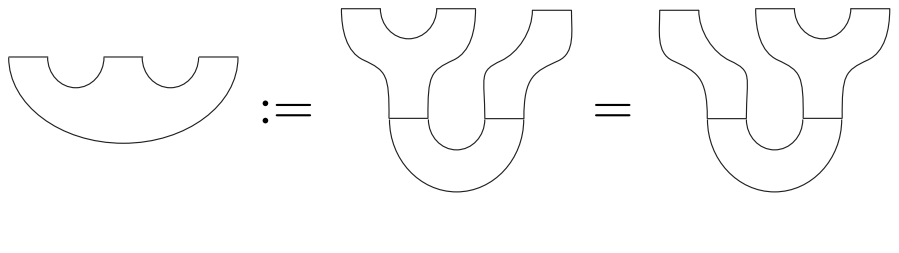

(3.39)

The fact that these two definitions are the same is a matter of writing them out and using associativity:

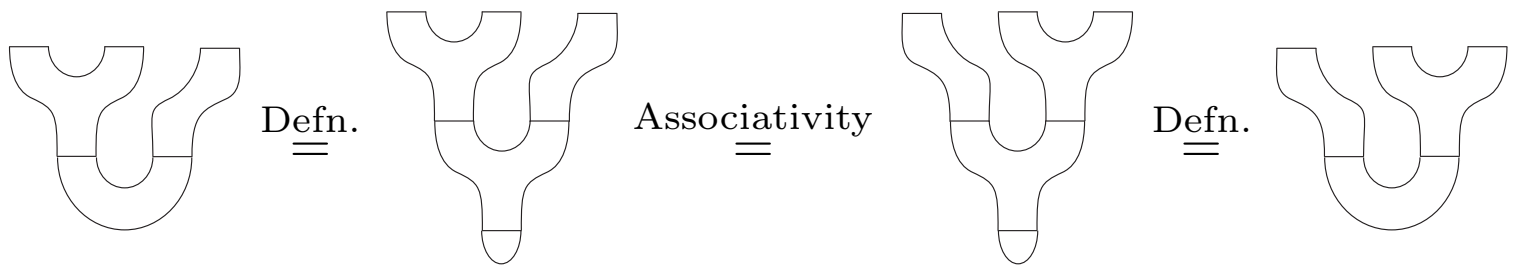

(3.40)

One then plays around a bit more and proves that (1) 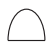 is a unit for 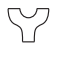, (2) 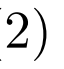 is a counit for 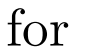, and (3) that 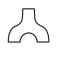 and 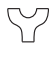 satisfy the Frobenius relation.

Given the right-hand definition, one shows, using the Frobenius relation, that 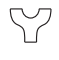 is associative, in other words, we have an algebra. Next one defines the pairing  and the copairing 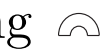 as:

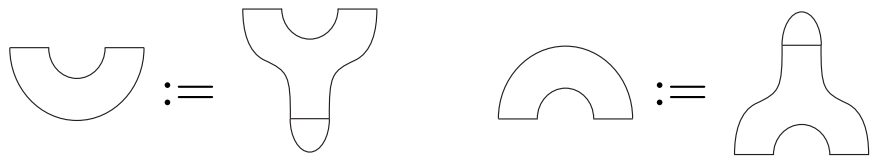

(3.41)

Some more lego games establishes the nondegeneracy of  and . This completes the proof. □

The fact that we can prove results using exclusively lego games means that they make sense in any monoidal category. Thus one defines a *Frobenius object* $(A, \mu, \Delta, \eta, \epsilon)$ in a general monoidal category as an object equipped with maps $\mu : A \otimes A \to A$, $\Delta : A \to A \otimes A$, $\eta : 1 \to A$ and $\epsilon : A \to 1$ satisfying the unit, counit and Frobenius conditions. For instance, Aaron Lauda has shown how a Frobenius object in **Cat** turns out to be related to an ambidextrous adjunction [2].

### OCob and 2Cob in terms of generators and relations

To say that a monoidal (resp. symmetric monoidal) category $\mathcal{C}$ is generated by a bunch of generators $G$ and relations $R$ means that every arrow in $\mathcal{C}$ can be obtained by composing and tensoring the arrows in $G$, and every equality in $\mathcal{C}$ can be obtained as a consequence of the relations in $R$ (resp. together with naturality of the twist map). Now we have,

**3.3.6 Lemma.** *The monoidal category **OCob** is generated by the following cobordisms:*

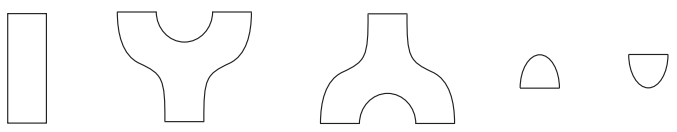

(3.42)

*subject to the relations:*

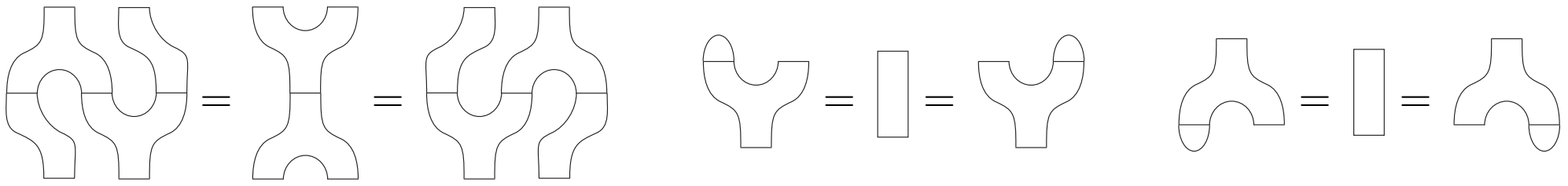

(3.43)

*Similarly, the symmetric monoidal category* ***2Cob*** *is generated by the following cobordisms:*

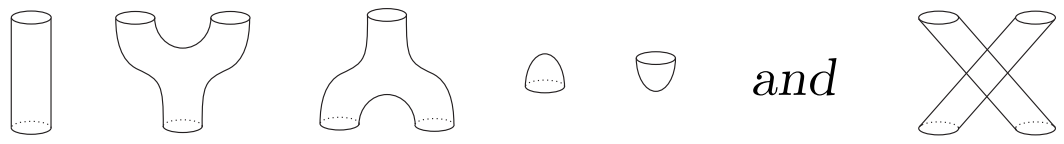

(3.44)

*subject to the relations:*

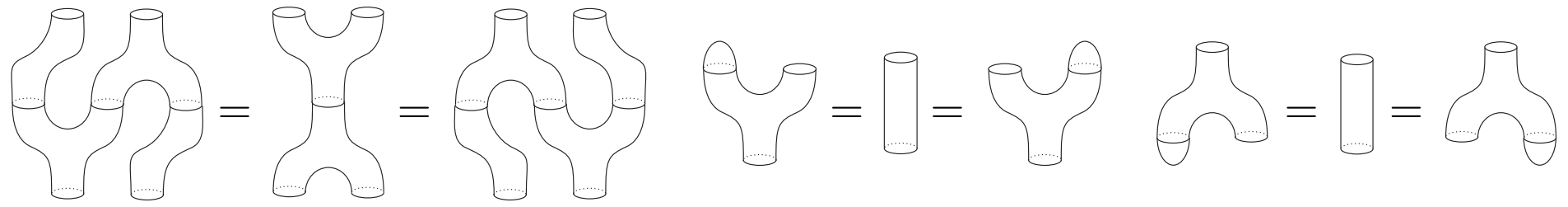

(3.45)

*together with the twist map relations:*[10]

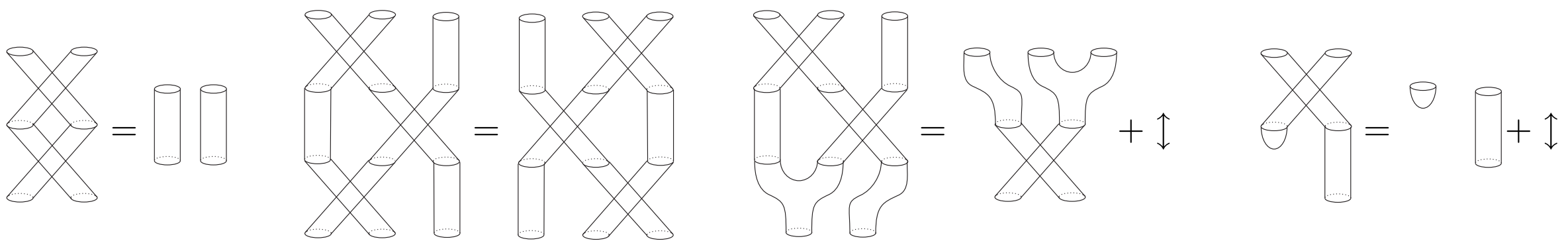

*Proof.* A two dimensional connected surface is classified by its number of boundaries and its genus (For planar surfaces like **OCob**, the role of genus is played by planar holes living inside the surface). Thus we can build up a normal form for surfaces as follows:

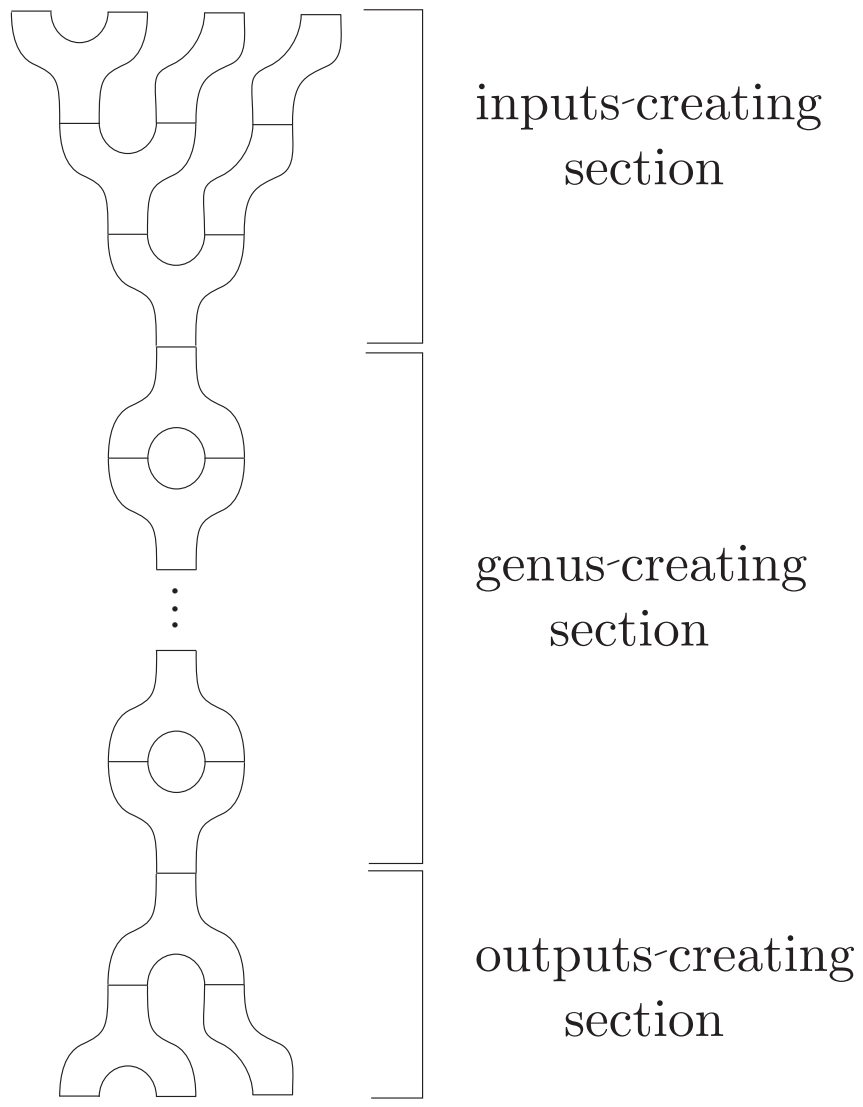

(3.46)

If the input or output is empty, just use the caps 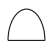 and 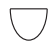 appropriately. A disconnected cobordism in **OCob** is simply a tensor product of cobordisms like (3.46). In **2Cob** we have to be a bit more careful because there may be twist maps involved. Luckily in such a case we can insert a double copy of the twist,

[10]It is technically not necessary to list these relations, since they are part of the axioms for a symmetric monoidal category, and hence implied. We list them for clarity.

eg:

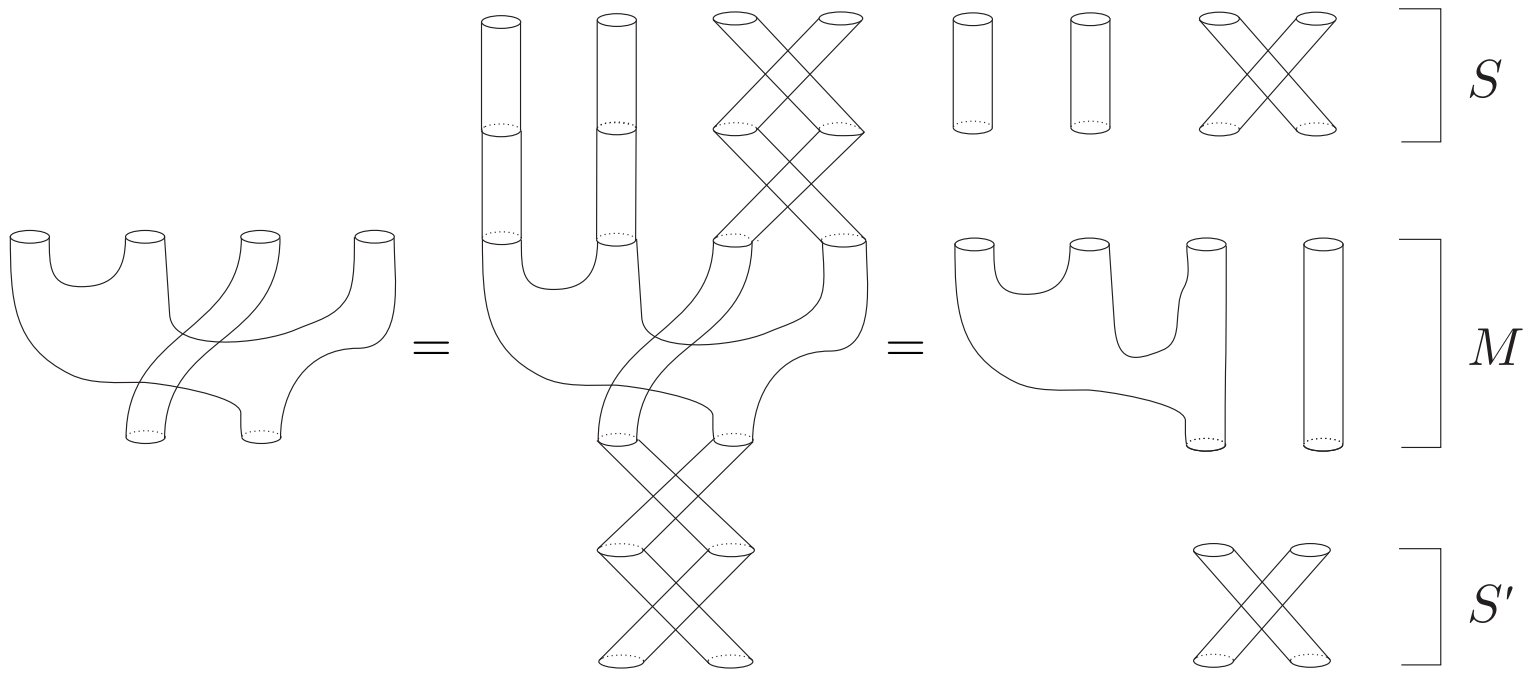

(3.47)

Applying this trick, we see that every disconnected cobordism $M'$ factors as $S'MS$ where $M$ is connected and $S$ is made up of twist maps. Thus these morphisms do indeed generate **OCob** and **2Cob**. To show that the relations are sufficient, firstly note that the relations (3.43) and (3.45) automatically imply the relations in Theorem 3.3.5. Now suppose we have a cobordism $M$, and we want to use the relations to bring it into normal form. Firstly put a Morse function[11] $f : M \to I$ on $M$ (For **OCob** this is literally the projection onto the $y$ axis (deforming the cobordism if necessary to ensure the height function is nondegenerate), while for **2Cob** one may appeal to Morse theory which ensures that one always exists). In the neighborhood of a critical point, the surface looks like (remember, time runs down the page):

Local minimum : 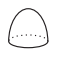 Local maximum : 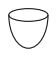 Saddle point I : 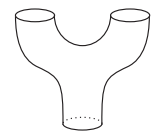 Saddle point II :

The Morse function is a timeline which chops up time into intervals $[t_i, t_{i+1}]$ which contain only one critical point, and can be represented by an elementary cobordism:

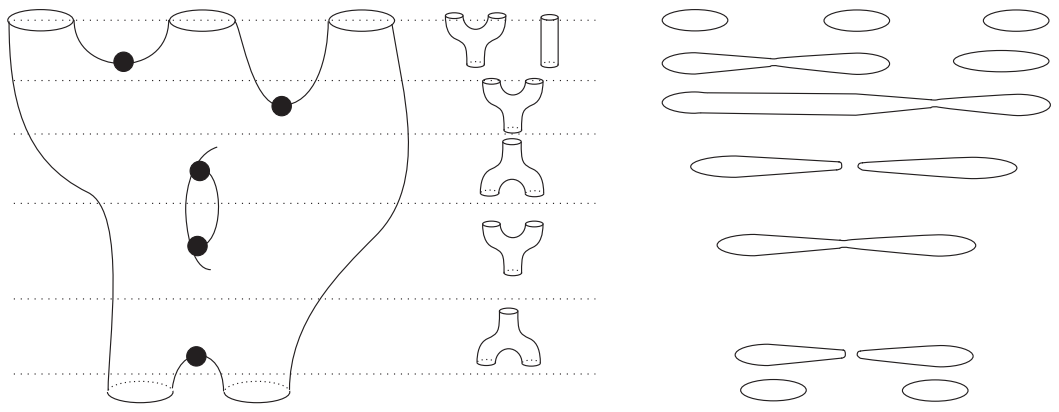

The diagram also shows how the critical points are exactly the points at which closed strings are created and annihilated. A similar algorithm can be applied in **OCob**. The net outcome is that $M$ is represented by a word $w(M)$ in the generators. We must show that $w(M)$ can be reduced to normal form. This can be done by using the relations to 'bubble' all the multiplications to the top, and all the comultiplications to the bottom. Also twist maps can be bubbled up and down in this way, using the twist relations (See Kock's book [50] for details). In this way one reduces $w(M)$ into normal form, thus establishing the sufficiency of the relations. □

This immediately establishes:
**Theorem 3.3.1 (⇐) (I)** : *A Frobenius algebra $A$ defines an open string theory $Z$ : **OCob** $\to$ **Vect**. A commutative Frobenius algebra $B$ defines a closed string theory $Z$ : **2Cob** $\to$ **Vect**.*

*Proof.* Define $Z$ : **OCob** $\to$ **Vect** by $Z(S^1) = A$, $Z(\,) = \mu$, $Z(\,) = \Delta$, $Z(\,) = \eta$ and $Z(\,) = \epsilon$. Lemma 3.3.6 tells us that (1) this is sufficient to define $Z$, since **OCob** is generated by these cobordisms, and (2) $Z$ is well-defined, because all the relations between , , and that one needs are guaranteed to be present amongst the $\mu$, $\Delta$, $\eta$ and $\epsilon$ which define $A$. A similar argument holds for **2Cob**. □

The importance of the formulation of a Frobenius algebra (Theorem 3.3.4) in terms of multiplication, comultiplication, unit and counit is now clear : these maps are in one-to-one correspondence with the local behaviour of a surface near a critical point.

[11] That is, a function $f : M \to I$ all of whose critical points are nondegenerate. If $M$ has a boundary, it is required that $f^{-1}(\partial I) = \partial M$.

### 3.3.6 Proof of Theorem 3.3.1 (II)

This proof is taken from Bakalov and Kirillov's book [18]. We consider only closed strings (open strings should be similar) and work in the Atiyah picture, so we are trying to assign vectors $Z(M) \in Z(\partial M)$ to manifolds with boundary $M$, in a way which is consistent with cutting up the manifold in all possible ways. Now Lemma 3.3.6 showed how a 2d manifolds could be decomposed into pieces of the form:

(a) disks : 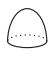 (b) cylinders : 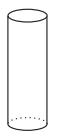 (c) trinions :

The advantage of working in the Atiyah picture is that it is not necessary to consider the time reversed (upside down) pictures, since in this picture $Z(M)$ depends only on $M$ as a manifold and not on $M$ as a cobordism (the word cobordism does not feature in Atiyah's definition)[12]. We will need the following classical lemma:

**3.3.7 Lemma. (Hatcher and Thurston ) [43]** *Any two ways to cut a 2-manifold $M$ into disks, cylinders and trinions can be related by isotopy of $M$ and a sequence of 'simple moves':*

(i) 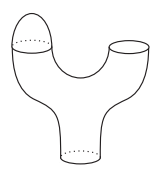 = 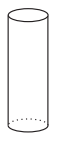 (ii) 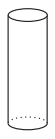 =  (i.e. removal) (iii) 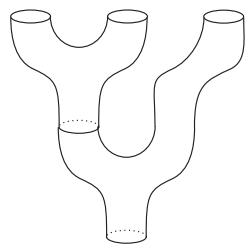 = 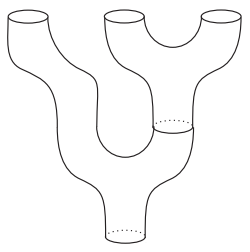 (iv) 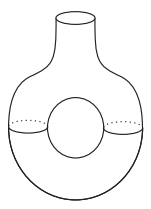 = 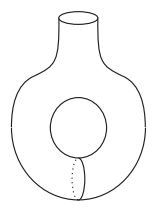

Now we have our second proof:
**Theorem 3.3.1 (⇐) (II)** : *A commutative Frobenius algebra $B$ defines a closed string theory $Z$ :* ***2Cob*** *→* ***Vect****.*

*Proof.* In the Atiyah picture, we must assign vector spaces $Z(\Sigma)$ to 1-manifolds $\Sigma$, isomorphisms of vector spaces $Z(f) : Z(\Sigma_1) \xrightarrow{\sim} Z(\Sigma_2)$ to diffeomorphisms $f : \Sigma_1 \xrightarrow{\sim} \Sigma_2$, and vectors $Z(\partial M)$ to manifolds with boundary $M$. We define $Z(S^1) = B$, and $Z(S^1 \sqcup \cdots \sqcup S^1) = B \otimes \cdots \otimes B$. There are really only two kinds of orientation-preserving diffeomorphisms between circles, $f : S^1 \xrightarrow{\sim} S^1$ the identity and $g : S^1 \xrightarrow{\sim} S^{1*}$ the diffeomorphism onto the circle with opposite orientation. Let $Z(f) = \mathrm{id}$, and $Z(g) : B \to B^*$ be given via the duality maps coming from the Frobenius form. Finally, let $Z(\text{disk}) = 1_B$, $Z(\text{cylinder}) = \mathrm{id}_B \in B \otimes B^*$, and $Z(\text{trinion}) = \mu \in B^* \otimes B^* \otimes B$, and let $Z(M)$ be given by chopping up $M$ into disks, cylinders and trinions in an arbitrary way, and then calculating the resultant $Z(M) \in Z(\partial M)$ using the gluing rule. Lemma 3.3.7 shows that this does not depend on the decomposition of $M$ used. □

### 3.3.7 Proof of Theorem 3.3.1 (III)

Consider an embedded version of **2Cob** called **TubeCob**. Similarly to **OCob**, objects are nonnegative integers $n$ representing the number of input or output circles, and morphisms $M : n_1 \to n_2$ are submanifolds of $\mathbb{R}^2 \times [0,1]$ with $\partial M = S_{n_1} \times \{0\} \cup S_{n_2} \times \{1\}$, where $S_k$ is a collection of $k$ radius $\frac{1}{3}$ circles centered at $(0,1), \cdots, (0,k)$ in the plane. The key point is that morphisms are considered equivalent only up to isotopy. Unlike **2Cob**, **TubeCob** is thus a nonsymmetric braided monoidal category, because one cannot deform the over braid into the under one:

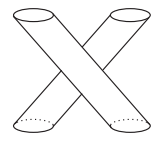 ≠ 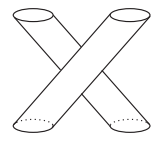 (3.48)

It is reasonably clear, however, that this is the *only* difference between **2Cob** and **TubeCob**. Suppose we define a TQFT for **TubeCob** as a braided monoidal functor $Z$ : **TubeCob** → **Vect**, where the braiding on **Vect** is the symmetric one, $v \otimes w \to w \otimes v$. Then $Z($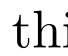$) = Z($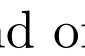$)$, so that this kind of $TQFT$ is completely equivalent to our previous notion.

#### Tube diagrams

Now we would like to find an algorithm which takes $M$ as input and outputs a word $w(M)$ in the generators of **TubeCob**. In the case of **2Cob** and **OCob**, knowledge of a Morse function (time) was sufficient[13]. This is not enough here because a Morse function will not detect the braidings. Here we propose another method, which amusingly is nothing but Feynman's famous time-slicing argument for deriving the path integral, applied to our present situation where the spatial sections can merge and split apart.

[12] It is hard to deny that this is certainly a technical advantage of the Atiyah picture.
[13] This is not entirely true. Even in **2Cob**, a Morse function alone could not detect the symmetry maps

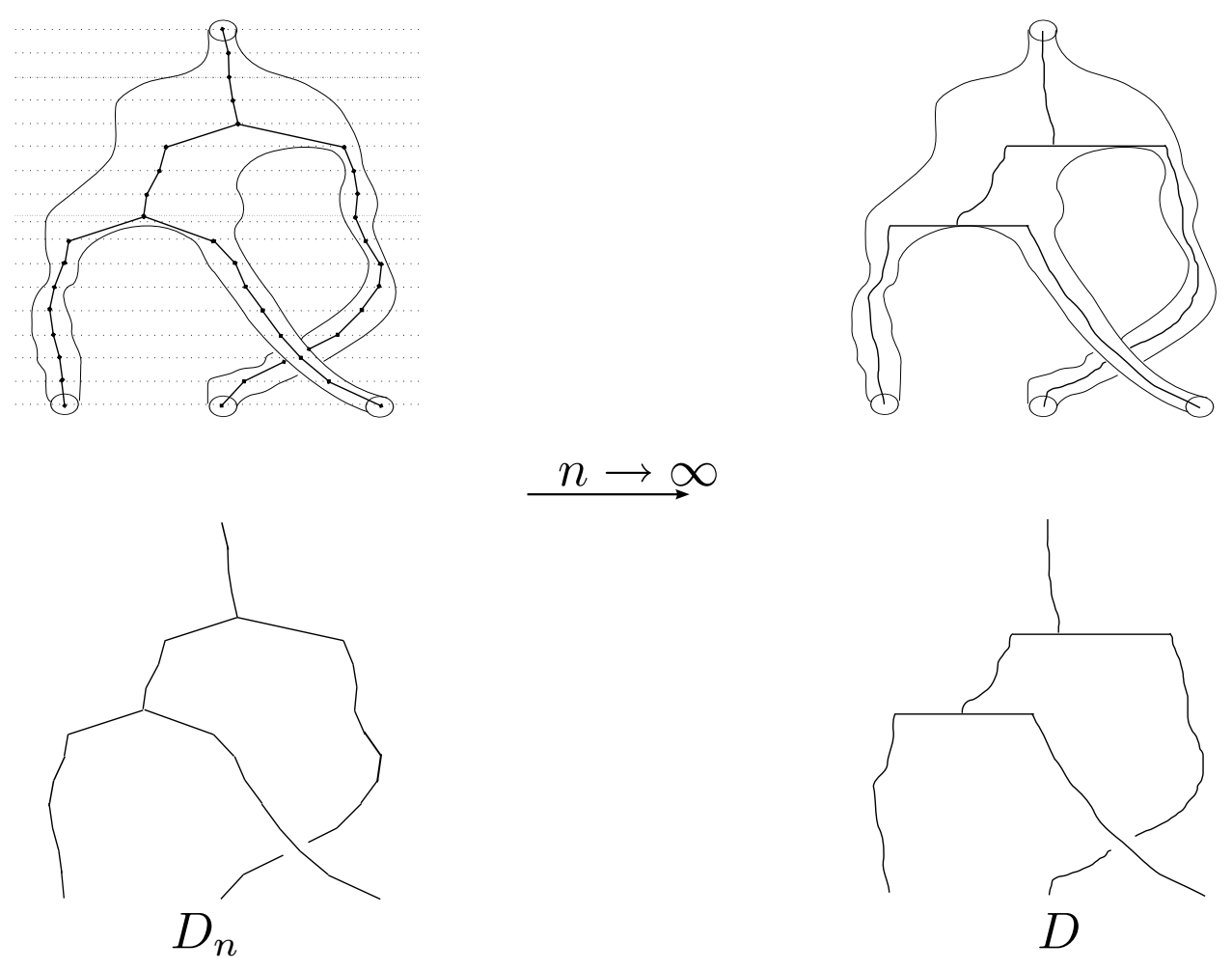

Figure 3.1: A tube picture being sliced up to form a tube diagram.

Let us agree to call a concrete morphism in **TubeCob** (that is, a concrete submanifold of $\mathbb{R}^3$) a *tube picture*. We shall now describe a method which associates to each tube picture $M$ a planar *tube diagram*[14].

Given the tube picture $M$, draw $n$ planes parallel to the $xy$ plane and having heights between 0 and 1 (see Fig. 3.1). For each $i$, the $i$'th plane generically intersects $M$ in a family of simple closed curves $\sigma$ ($\sigma$ might be a single point). For each plane $i$ and each $\sigma$ construct the point $p^i_\sigma$, the geometric center of the planar curve defined by $\sigma$, by averaging over $\sigma$. Now connect points $p^i_\sigma$ and $p^{i+1}_{\sigma'}$ lying on consecutive planes by a straight line in 3d if and only if $\sigma$ is connected to $\sigma'$ in the segment of $M$ lying between the planes $i$ and $i+1$. This procedure associates to the tube picture $M$ a discretized diagram composed of line segments living in 3d space. Define $D_n$ as the projection of this diagram onto the $xz$ plane. If two line segments cross over or under each other in space, indicate this on their projections in the planar diagram. The *tube diagram* is defined as $D = \lim_{n\to\infty} D_n$. Once we have $D = D(M)$, it is easy to split it into segments containing elementary generators, and hence we can form a word $w(M) = w(D(M))$. The strategy of the third proof of Theorem 3.3.1 is to classify the behaviour of $w(M)$ under smooth deformations of $M$. These correspond to a number of elementary 'graph moves', which turn out to be identical to those defining a Frobenius algebra, thus completing the proof.

## 3.4 The Closed/Open category and D-branes

### 3.4.1 Introduction.

We now consider a cobordism category **TopString** in which both open and closed strings make their appearance. A TQFT in this sense is a functor $Z$ : **TopString** $\to$ **Vect**, and is a baby model of string field theory. The idea is to understand the algebraic consequences of the open/closed transitions in **TopString**. These kind of ideas were originally studied in the context of open/closed conformal field theory by Lewellen [54]. Moore and Segal [62, 63] have shown recently how many of the open questions in string field theory can be explicitly solved in the simpler topological setting, and it is their approach we shall follow here (an alternative approach is presented in [46, 52]). Others have studied how **TopString** is related to string topology, that is, the study of algebraic and topological properties of loop spaces [21, 22, 77][15].

[14]This is a personal construction. Since we do not have embedding maps, we don't have access to the cores like Turaev did [89]. This means we have to find the cores ourselves. It might look artificial, but in a certain sense it is more natural than Turaev's cores since these diagrams work hand-in-hand with the critical values of the surface, whereas the ordinary cores don't.

We will define the notion of a $D$-brane in this setting, and then continue to classify them completely, under certain conditions. D-branes are notorious constructs in real string theory whose precise mathematical definition remains slightly elusive and mysterious. Nevertheless, one approach is to consider them as conformally invariant boundary conditions in 2d open and closed conformal field theory. The idea is that the closed strings define a 'string background' which places constraints on the boundary conditions. Since D-branes carry vector-bundle type geometry through the fields which live on them, a first step towards classifying them is through K-theory. Recall that (topological) K-theory is the rough-and-ready algebraic approach to classifying vector bundles on a topological space $X$ : simply define the abelian group $K_0(X)$ as the canonical group obtained from the semigroup which contains one generator $\langle V \rangle$ for each type of vector bundle on $X$, with relations $\langle V \rangle + \langle W \rangle = \langle V \oplus W \rangle$. Of course there is infinitely more to K-theory than this no-nonsense description, but this approach is enough to motivate why, for conformal field theories that do actually have a spacetime interpretation, it is believed that the K-theory of the spacetime is a rough indication of the types of D-branes which may live there. Of course it is only a *rough* description - and one may refine it by introducing twisted K-theory, elliptic cohomology, derived categories, and so on [83, 82].

Our main result will be to show that, in the case of a semisimple closed string algebra $B$, the possible D-branes are classified by $K^0(\text{spec}(B))$.

### 3.4.2 TopString

We consider a D-brane as a 'surface on which an open string ends'. Let us suppose that the finite index set $I$ labels the different kind of D-branes there are. So for instance, the identity in **OCob** would now be interpreted as the propagation of a string with one end living on $a \in I$ and the other end living on $b \in I$:

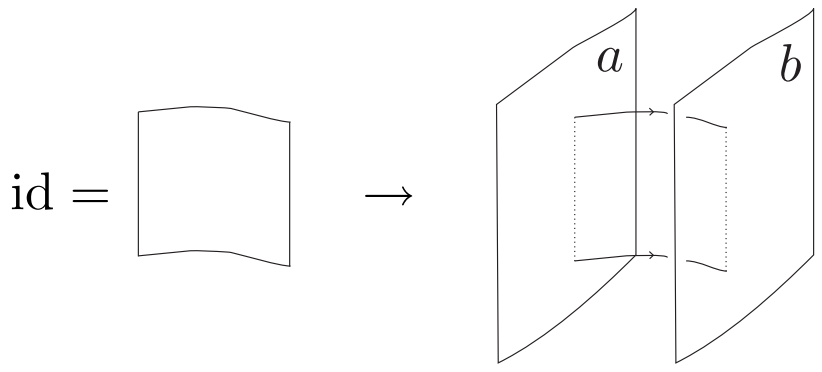

The boundary is thus partitioned into three parts : the incoming string, the outgoing string (each of which has its endpoints labeled from $I$), and the 'free' boundaries (which we have drawn as a dotted line) which live on the D-branes, and which are themselves (1 dimensional) cobordisms from the endpoints of the incoming string to the endpoints of the outgoing string.

Let us define a *sequence of strings* as finite sequence $(a_1, a_2, \ldots, a_n)$ where $a_i$ is either equal to a closed string $S^1$ or $a_i = I_{(a,b)}$, an open string labeled by its D-brane endpoints. Now suppose $\Sigma_1$ and $\Sigma_2$ are sequences of strings. We define an *open-closed cobordism* $M : \Sigma \to \Sigma'$ as a compact surface $M$ together with embeddings[16]

$$\Sigma \xrightarrow{i} \partial M \xleftarrow{i'} \Sigma' \tag{3.49}$$

where $i$ is orientation preserving and $i'$ is orientation reversing, and such that:

- $i(\Sigma)$ and $i'(\Sigma')$ are disjoint, and
- the free boundary $\partial_f M := \partial M - i(\Sigma) - i(\Sigma')$ must form a 1d cobordism $i(\partial(\Sigma)) \to i'(\partial(\Sigma'))$ which only connects labels of the same type.

Cobordisms are considered equivalent up to diffeomorphism preserving the input and output maps. Together they form a category **TopString**. The tensor product, on objects, is given by concatenation of strings, and on arrows by disjoint union. There are also symmetry cobordisms which swap block-wise the elements of the sequences of strings. These constructions present **TopString** as a symmetric monoidal category.

[15] The idea is that the geometric realization of **TopString** is a loop space.

[16] $\partial M$ consists of circles. By an embedding of an interval into a circle we simply mean an embedding of the interval into a subset of the circle.

**TopString** has open and closed string sectors, as well as maps connecting these sectors. The closed string sector is identical to **2Cob**, but the open string sector is not quite the same as **OCob** because we are considering the strings abstractly, and not embedded in the plane. Thus there is a twist cobordism which was not present before:

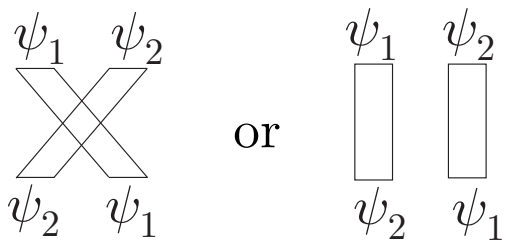

Here $\psi_1, \psi_2$ are just shorthand for the images of the intervals under the embedding map $i$ (we shall use $\psi$'s for labeling open strings, and $\phi$'s for labeling closed strings), but under $Z$ are of course vectors living inside vector spaces. Also we are assuming for now that there is only one type of D-brane, whose label we will omit in drawings. The Y-piece is still not commutative:

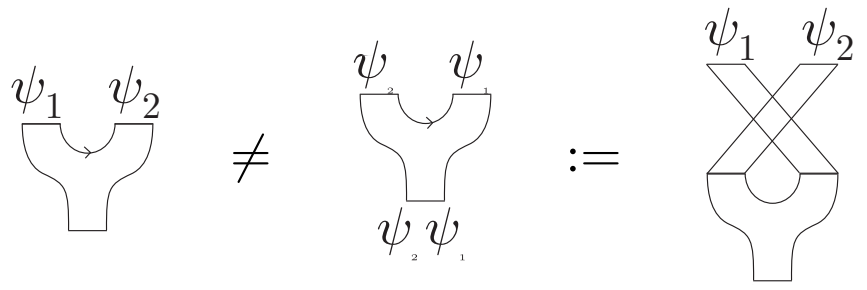

This is because an orientation preserving diffeomorphism must preserve the cyclic order of the labels. However, the horseshoe does become symmetric, since there is no output boundary which gets in the way:

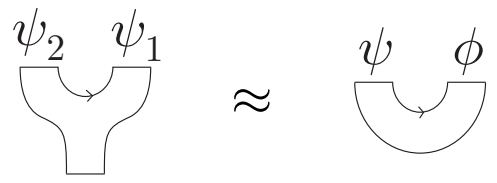

Thus an topological string theory $Z$ : **TopString** $\rightarrow$ **Vect** gives us a symmetric (in general noncommutative) Frobenius algebra $A$ from the open string sector, and a commutative Frobenius algebra $B$ from the closed string sector.

### 3.4.3 The open/closed frontier.

However, the algebras $A$ and $B$ have more structure, because there are open/closed transitions inside **2Cob**. We have the forward and reverse 'pennywhistle':

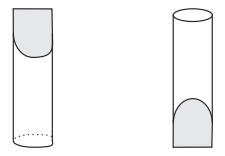

(3.50)

The shaded regions represent the inner surface of the cylinder - drawings like (3.50) can be confusing! These translate into maps $i^* := Z(\,) : A \rightarrow B$ and $i_* := Z(\,) : B \rightarrow A$ and we now examine their properties.

The first relation is that it makes no difference if one multiplies two closed strings, and then converts them to open strings, or if one first converts the closed strings to open strings, and then multiplies them as open strings:

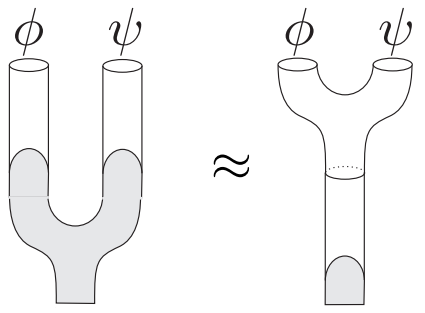

(3.51)

That is, $i_*(\phi \cdot \psi) = i_*(\psi) \cdot i_*(\phi)$, so that $i_* : B \rightarrow A$ is an algebra homomorphism.
The second relation is that $i_*$ is unit-preserving, $i_*(1_B) = 1_A$:

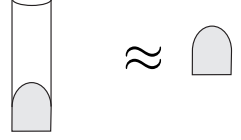

(3.52)

The third relation says that, although the multiplication in the open sector is not commutative, the hybrid multiplication is:

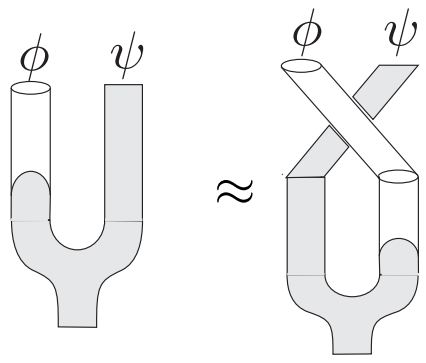

(3.53)

Here we have applied the 'hybrid' open-and-closed twist. Thus $i_*(\phi)\psi = \psi i_*(\phi)$, so that $i_* : B \to Z(A)$. The fourth relation is that the notation $i_*$ and $i^*$ was not accidental. The following diagram,

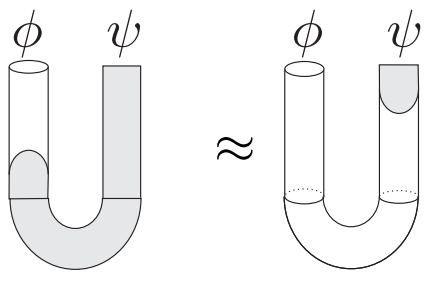

(3.54)

shows that $(i_*(\phi), \psi)_A = (\phi, i^*(\psi))_B$. This means that $i_*$ and $i^*$ are adjoint maps with respect to their Frobenius pairings.
The fifth relation has to do with the remarkable map $\pi : A \to A$:

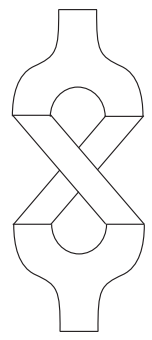

(3.55)

$\pi$ appears to nothing to do with the closed string sector. Yet the following remarkable topological manipulations make nonsense of this belief:

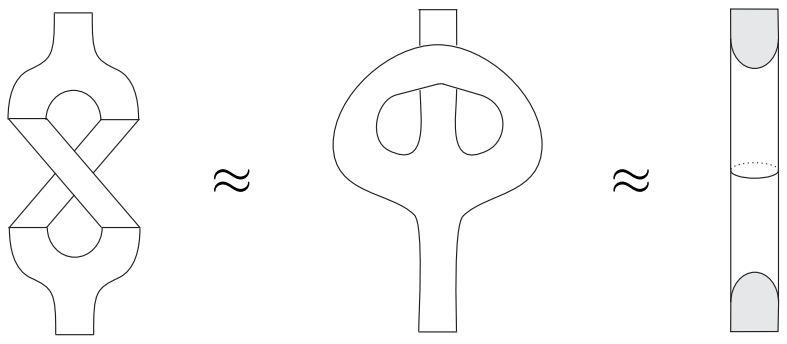

(3.56)

Thus the fifth relation is that $\pi$ *factorizes through the closed sector*,

$$\pi = i_* i^* \tag{3.57}$$

. This interesting relation has been denoted by Moore as the *Cardy condition* [19], after the physicist J.L. Cardy who analyzed the constraints on the possible boundary states in a class of rational conformal field theories. The Cardy condition places a very severe constraint on the interaction between the closed and open string algebras.

Lewellen obtained relations (3.51)-(3.56) during a study of open/closed conformal field theory [54][17], where he argued that they were sufficient. Moore and Segal specialized them to the topological setting [62, 63], as well as Lazariou [46]. Thus we have the following folk-theorem (the author has never seen an explicit proof):

**3.4.1 'Theorem'.** *To give an open and closed topological field theory $F$ : **TopString** $\to$ **Vect** is to give*

*(a) A commutative Frobenius algebra $B$.*

*(b) A (normally noncommutative) Frobenius algebra $A$.*

*(c) A unit-preserving homomorphism $i_* : B \to Z(A)$, such that*

$$i^* i_* = \pi := \mu_A \sigma_{A,A} \Delta_A \tag{3.58}$$

*where $i^*$ is the adjoint to $i_*$ under the Frobenius pairings on $A$ and $B$.*

[17] There are extra relations in conformal field theory. These turn out to be tautologies in the topological case [63].

### 3.4.4 Classification of open string algebras

We want to understand what kind of open strings can live in a given closed string background. Suppose we are given the closed string algebra $B$. Our task is to classify the possible compatible $A$'s. It turns out that if $B$ is semisimple there is a complete classification and a good spacetime interpretation. One may view the theory as a baby model of string field theory. On the other hand, if $B$ is not semisimple, then there are solutions which do not possess a spacetime interpretation. Nevertheless, they have been considered relevant to branes on Calabi-Yau manifolds [63].

As we saw in Section 3.3.4, the most general semisimple commutative Frobenius algebra takes the form

$$B = \bigoplus_{i=1}^{n} \mathbb{C} a_i \tag{3.59}$$

where $a_i a_j = \delta_{ij} a_j$, and the only degree of freedom is in the choice of the traces, $\epsilon_i := \epsilon(a_i) \neq 0$. If the TQFT is unitary then the $\epsilon_i$ must be real and positive.

In this case, one can show that the general solution for $A$ is simply a direct sum of matrix algebras, one for each dimension of $B$:

$$A = \bigoplus_{i=1}^{n} \mathrm{Mat}_{k_i}(\mathbb{C}) \tag{3.60}$$

where the dimensions of the matrix algebras $k_i$ are arbitrary integer choices. The second choice one needs to make is to choose a square root of each $\epsilon_i$, which we shall write as $\sqrt{\epsilon_i}$ with the choice of sign understood. If we decompose $\psi \in A$ as

$$\psi = \begin{pmatrix} \psi_1 & 0 & 0 & \cdots \\ 0 & \psi_2 & 0 & \cdots \\ \vdots & \vdots & \vdots & \ddots \end{pmatrix} \tag{3.61}$$

then the open string trace $\epsilon_A$ and the maps $i^*$ and $i_*$ are defined by:

$$\epsilon_A(\psi_i) = \sqrt{\epsilon_i}\,\mathrm{Tr}(\psi_i) \tag{3.62}$$

$$i^*(\psi_i) = \frac{\mathrm{Tr}(\psi_i)}{\sqrt{\epsilon_i}} a_i \tag{3.63}$$

$$i_*(a_i) = 0 \oplus \cdots \oplus \mathrm{id}_i \oplus \cdots \oplus 0. \tag{3.64}$$

This is the complete solution of our algebraic problem, in the case when $B$ is semisimple.

### 3.4.5 Spacetime interpretation

We want to give a geometric interpretation of this solution. This is done via the Gelfand-Naimark theorem, which associates topological spaces to commutative algebras such as $B$, and vice versa. Given a topological space $X$, the space of functions $C(X)$ over $X$ is a commutative algebra. On the other hand, given a commutative algebra $B$, one defines $\mathrm{Spec}(B)$ as the space of homomorphisms from $B$ into the complex numbers,

$$\mathrm{Spec}(B) := \{\chi : B \stackrel{\text{linear}}{\to} \mathbb{C}, \chi(a_1 a_2) = \chi(a_1)\chi(a_2)\} \tag{3.65}$$

$\mathrm{Spec}(B)$ is itself a topological space, and in fact the space of functions over $\mathrm{Spec}(B)$ is isomorphic to the original algebra $B$! Thus there is a one-to-one correspondence between commutative algebras and topological spaces[18]. What is the spacetime (topological space) associated with $B$? It is easy to see that

$$\mathrm{Spec}(B) = \{\chi_1\} \cup \{\chi_2\} \cup \cdots \cup \{\chi_n\} \tag{3.66}$$

with $\chi_i(a_j) = \delta_{ij}$. So spacetime consists of $n$ discrete points! We see that the open string algebra $A$ is a choice of matrix algebra for each point in spacetime, and can thus be interpreted as the algebra of endomorphisms of a *vector bundle* over spacetime (this is precisely the kind of relation we expected). More precisely, if we

[18]More precisely, there is a one-to-one correspondence between compact Hausdorff spaces and commutative unital $C^*$ algebras.

let $W \to \mathrm{Spec}(B)$ be the vector bundle over spacetime where a copy of $\mathbb{C}^{k_i}$ sits above each spacetime point $\chi_i$ (equivalently $a_i$):

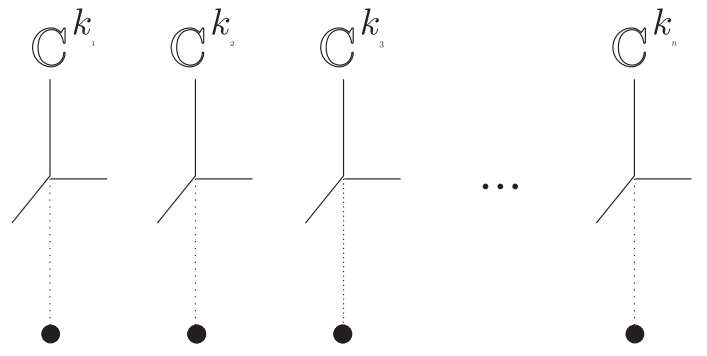

Then $A = \mathrm{Hom}(W, W)$ where we interpret Hom as meaning vector bundle homomorphisms, i.e. they split up into matrix algebras over each point.

### 3.4.6 Adding D-branes

This was only a first step. We now consider the case when there are multiple types of D-branes. The open strings $I_{ab}$ have their ends living on D-branes:

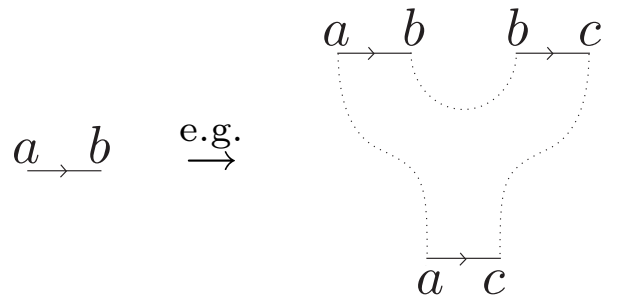

The TQFT transforms these into vector spaces $A_{ab}$ and maps $A_{ab} \otimes A_{bc} \to A_{ac}$. Thus instead of a single vector space $A$, we have a whole system of them, together with a whole system of multiplications, comultiplications, etc. They are subject to labeled versions of the same geometric relations we had before. For instance, the Frobenius relation becomes:

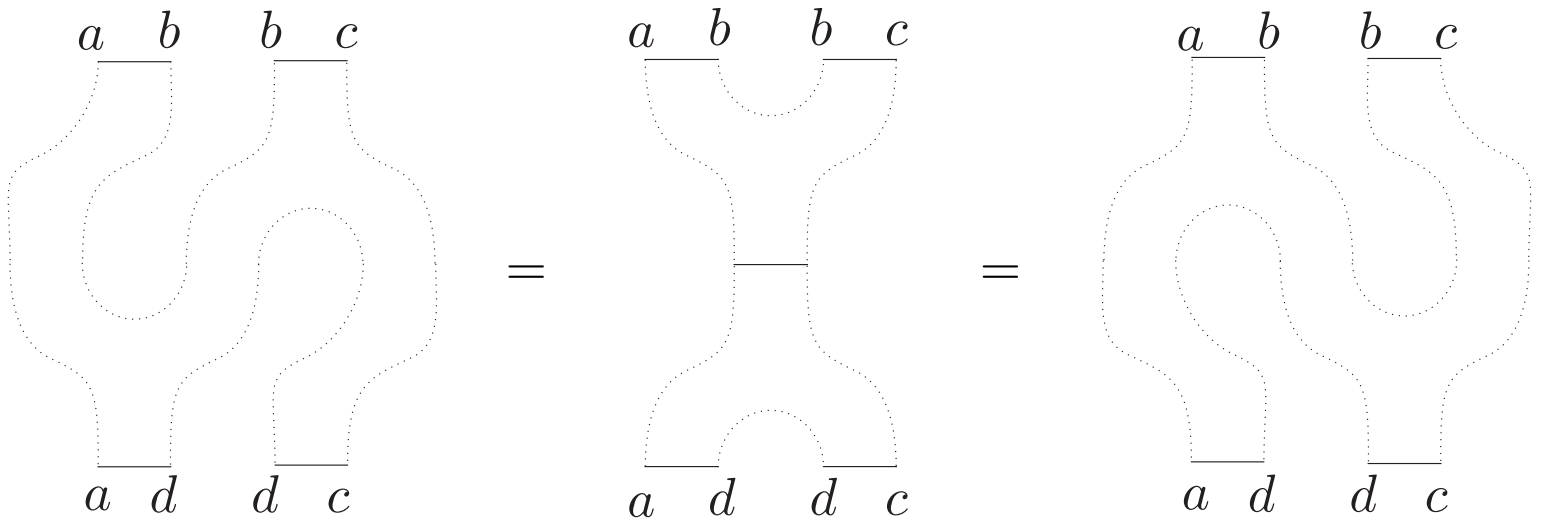

Similarly, all the open/closed relations (3.51)-(3.56) go over into labeled versions. Despite the apparent complexity, nothing essentially new has been added, and the solution is a direct generalization of what we had for only one type of D-brane. For each $a \in I$, we freely choose $n$ vector spaces $\mathbb{C}^{k_i^{(a)}}$ to be associated with the D-brane $a$ and which live over the $n$ spacetime points. We interpret this as a vector bundle over spacetime $W_a \to X$ for each $a \in I$. Finally we choose an overall sign for each spacetime point, as before. Then we have

$$A_{ab} = \mathrm{Hom}(W_a, W_b) \qquad \text{(as vector bundles).} \tag{3.67}$$

The trace $\epsilon_{A_{ab}}$ and maps $i_*$ and $i^*$ are given by the same formula as before, (3.62) - (3.64). Let us summarize : for each D-brane $a$, we choose $n$ natural numbers $k_1^{(a)}, \ldots, k^(a)_n$. We also choose a sign for each spacetime point. This is equivalent to choosing, for each D-brane, $n$ integers. We notice that this is precisely the algebraic $K$-theory of the closed string algebra $B$! That is,

$$K_0(B) = K^0(\mathrm{Spec}(B)) = \mathbb{Z} \oplus \cdots \oplus \mathbb{Z}. \tag{3.68}$$

One defines algebras $B$ and $B'$ to be *Morita equivalent* if $K_0(B) = K_0(B')$, and the Morita theorem tells us that this is same as saying that $B$ and $B'$ have equivalent categories of representations. We have established:

**3.4.2 Theorem (Moore and Segal).** *Let the closed string algebra $B$ be semisimple. Then the admissible boundary conditions for $B$ are precisely*

$$\Lambda = K^0(\mathit{Spec}(B)) = K_0(B) = \{ \text{ isomorphism classes of objects in } B\text{-mod}\}. \tag{3.69}$$

Of course, the appearance of $K^0$ should have been anticipated once we started generating 'vector bundles' i.e. matrix algebras. Moreover we are clearly using only the most rudimentary form of '$K$-theory' here. But we shall see in Sec. 5.1 how one can carry over this result to fully fledged conformal field theory!

# Chapter 4

# 2d gauge TQFT's, triangulations, and Yang-Mills theory

We will show in this chapter how one may obtain elegant constructions for TQFT's by considering principal $G$-bundles on manifolds. This model was first introduced by Dijkgraaf and Witten [25]. It can also be computed from triangulations, where it is interpreted as a lattice gauge theory. The only data that is necessary is a semisimple, finite-dimensional algebra $A$, but the usual case is that $G$ is a finite group and $A = \mathbb{C}[G]$. Generalizing these ideas to a continuous gauge group, we shall obtain topological Yang-Mills theory, from which ordinary Yang-Mills can be regarded as a perturbation. As Segal has pointed out [80], the surprising thing about the Dijkgraaf-Witten model is that a great part of quantum field theory can be seen as the study of various generalizations of this toy model in which finite groups are replaced by arbitrary Lie groups.

## 4.1 A functorial approach to gauge theory

### 4.1.1 Connections and gauge transformations.

We recall here some definitions in order to set the notation. The arena for gauge theory is a principal $G$-bundle $P \xrightarrow{\pi} M$, where $G$ is a compact Lie group. In our notation, $G$ acts from the right on $P$, $p \to p \cdot g$. In what follows, we shall often make the assumption that $P$ is a trivial bundle. A *gauge transformation* is a diffeomorphism $P \to P$ leaving the basepoints invariant, $\pi(f(p)) = f(p)$, and preserving the right action, $f(p \cdot g) = f(p) \cdot g$. These transformations form the (infinite dimensional) gauge group,

$$GA(P) = \{f : P \to P,\ \pi(f(p)) = f(p),\ f(p \cdot g) = f(p) \cdot g\}. \tag{4.1}$$

There is a natural correspondence between $GA(P)$ and the those maps $C(P, G)$ from $P$ to $G$ which convert the right action into the adjoint action,

$$C(P, G) := \{\tau : P \to G,\ \tau(p \cdot g) = g^{-1}\tau(p)g\}. \tag{4.2}$$

If $f \in GA(P)$, define $\tau : P \to G$ by the relation $f(p \cdot g) = p \cdot \tau(p)$ (this does not depend upon the choice of $p$). Similarly, if $\tau \in C(P, G)$, define $f : P \to P$ by $f(p) = p \cdot \tau(p)$. One easily establishes that these maps define an (anti)-isomorphism of groups between $GA(P)$ and $C(P, G)$.

In the case when $P$ is a trivial bundle we can regard gauge transformations as maps from the base manifold $M$ to the structure group $G$. A map $\tau \in C(P, G) : P \to G$ is uniquely defined by its restriction to the image of a local section $s : M \supset U \to P$. For, using the section one may assign group elements $g$ to each $p \in \pi^{-1}(U)$ by setting $p = s(x) \cdot g(x)$, $x \in U$. Hence, given an arbitrary map $\hat{\tau} : s(U) \to G$, we may recover $\tau$ by setting $\tau(p) = \tau(s(x) \cdot g(x)) = g^{-1}\hat{\tau}(x)g$. Thus if $P$ is trivial (so that there exists global section),

$$GA(P) \approx C(P, G) \approx C(M, G) := \{s : M \to G\}. \tag{4.3}$$

A connection on $P$ is a smooth decomposition of the tangent space at every point $p \in P$ into a vertical part $V_p = \mathrm{Ker}(\pi_*)_p$ and a horizontal part $H_p$, $T_pP = V_p \oplus H_p$, which is compatible with the right action $R_g$ of $G$ on $P$. Equivalently, a connection $A$ is a $G$-equivariant one form with values in the lie algebra $\mathfrak{g}$ of $G$. $G$-equivariance means that

$$A(\xi_p) = \xi_p \tag{4.4}$$
$$A(R_{g*}X) = \mathrm{ad}(g^{-1})A(X) \tag{4.5}$$

The horizontal subspaces $H_p$ are then the kernels of $A$. Normally one does not wish to distinguish connections $A$ and $A'$ if they differ by a gauge transformation, i.e. if $f_*(H_p) = H'_{f(p)}$ for some $f \in GA(P)$.

### 4.1.2 The connection as a functor

The connection $A$ assigns to each path $\sigma$ in $M$ from $x$ to $y$ and a point $p \in P_x$, a curve $\tilde{\sigma}_p$ in $P$ which begins at $p$ and ends in the fibre above $y$, and which proceeds always in the horizontal direction determined by $H_p$. This curve is known as the *lift* of $\sigma$, and sets up a bijection $\Gamma_A(\sigma) : P_x \xrightarrow{\sim} P_y$ called the *holonomy* of $A$ along $\sigma$. Define $\mathbf{Path}(M)$ as the category of smooth paths in $M$ : the objects are points in $M$ and the arrows are paths from $x$ to $y$, up to parameterization. A connection $A$ on $P$ then gives rise to a functor

$$\Gamma_A : \mathbf{Path}(M) \to \mathbf{Vect} \tag{4.6}$$

which sends $x$ to $\pi^{-1}(x)$ and $x \xrightarrow{\sigma} y$ to the holonomy of $A$ along $\sigma$. Conversely, given such a functor $\Gamma$, one can reconstruct the connection $A$ providing the functor satisfies important properties such as parameterization invariance [11, 17]. This is an important tool in the loop representation of quantum gravity.

The resemblance of $\Gamma_A$ to a TQFT is striking, and has been pointed out by Picken [72]. The holonomy functor $\Gamma_A$ can be viewed as a kind of *embedded* 1d TQFT, similar to an embedded 2d TQFT $F : \mathbf{TubeCob} \to \mathbf{Vect}$. An interesting question is whether one can define a quantum field theory action which gives rise to $\Gamma_A$; for instance the Kontsevich integral construction is a candidate [72].

Suppose $P$ has been trivialized with a global section $\chi : M \to P$, and that $\sigma : x \to y$ is a curve in $M$ with its lift chosen such that $\tilde{\sigma}(0) = \chi(x)$. Then this defines a unique $g \in G$ satisfying $\tilde{\sigma}(1) = \chi(y) \cdot g$. Moreover if $\sigma'$ is a curve from $y$ to $z$ with $\tilde{\sigma}'(0) = \chi(y)$, then $\widetilde{(\sigma' \circ \sigma)}(1) = \chi(z) \cdot gg'$. Thus when $P$ is trivial (4.6) is equivalent to specifying a functor

$$\Gamma : \mathbf{Path}(M) \to \mathbf{G}, \tag{4.7}$$

where $\mathbf{G}$ is viewed as the category with only one object and whose morphisms are a copy of $G$. In this picture, gauge transformations act, via (4.3), by conjugating the holonomies,

$$\Gamma(\sigma) \to (s \cdot \Gamma)(\sigma) := s(y)^{-1}\Gamma_A(\sigma)s(x), \tag{4.8}$$

This is nothing other than a natural transformation,

$$\mathbf{Path}(M) \underset{\Gamma'}{\overset{\Gamma}{\Downarrow s}} \mathbf{G}\,. \tag{4.9}$$

The connection is defined to be *flat* if it assigns the same group element to homotopic paths, i.e. $\Gamma(\sigma) = \Gamma(\sigma')$ if $\sigma'$ is homotopic to $\sigma$ (as paths from $x$ to $y$):

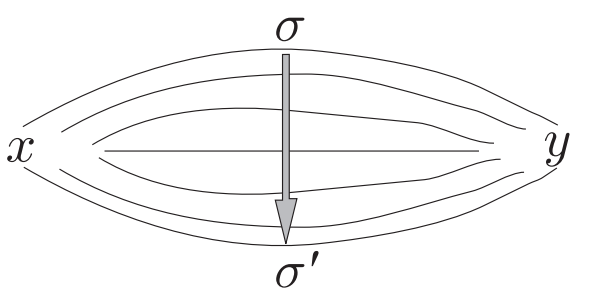

This notion allows us to think of $\mathbf{Path}(M)$ as a *2-groupoid* (see Appendix A). The objects are points in $M$ and the morphisms are paths in $M$ as before. The 2-morphisms between two paths $\sigma : x \to y$ and $\sigma' : x \to y$ are the homotopies taking $\sigma$ into $\sigma'$. This view of gauge theory is rather powerful as it provides a natural way to understand 'higher-dimensional' connections such as the $B$ field arising in string theory [73, 12].

For our purposes, it is enough to define $\widehat{\mathbf{Path}}(M)$ as the groupoid of paths in $M$ considered up to homotopy. Note that

$$\mathrm{Hom}_{\widehat{\mathbf{Path}}(M)}(x, x) = \pi_1(M, x) \tag{4.10}$$

for any $x \in M$. A connection $\Gamma : \mathbf{Path}(M) \to G$ is flat precisely when it factors through the projection functor $\mathbf{Path}(M) \to \widehat{\mathbf{Path}}(M)$.

Consider the moduli space $\mathcal{M}$ of all flat, gauge inequivalent, $G$ connections on $M$. The discussion above shows that

$$\mathcal{M} = \mathrm{Hom}(\pi_1(M), G)/G, \tag{4.11}$$

where the quotient action of $G$ generates the gauge equivalence classes of connections, (4.8). The spaces $\mathrm{Hom}(\pi_1(M), G)$ and $\mathcal{M}$ have a tremendously rich geometrical structure, which is highly nontrivial even if $M$ and $G$ are rather simple! This is an intriguing subject which can be studied from a variety of angles. The classic reference is Atiyah and Bott [8].

### 4.1.3 Triangulated spaces

We wish to apply the language of gauge theory to triangulated spaces, or in physics terminology, to engage in lattice gauge theory[1]. A 'triangulated space' is most elegantly defined via simplicial sets and singular complexes, language familiar from cohomology. Since we wish to specialize, ultimately, to two dimensions, we shall employ here a 'working definition' of a triangulated space, which we shall occasionally call a *simplicial 2-graph*, as a nod to the existence of the more coherent theory.

**4.1.1 Definition.** *A triangulated space, or a* simplicial 2-graph, *consists of*

- *a set $V$ of vertices*
- *a set $E$ of edges*
- *a set $T$ of triangles*

*with maps*

$$V \underset{d_0}{\overset{d_1}{\leftleftarrows}} E \overset{d_2,\, d_1,\, d_0}{\Lleftarrow} T \tag{4.12}$$

*where $d_i$ means 'leave out the $i$th vertex', satisfying*

$$\begin{aligned} d_1 d_1 &= d_1 d_2 && (4.13)\\ d_0 d_2 &= d_1 d_0 && (4.14)\\ d_0 d_0 &= d_0 d_1 && (4.15) \end{aligned}$$

The idea behind the definition is that the edges are oriented in such a way that $d_0$ and $d_1$ of an edge denote its source and target respectively:

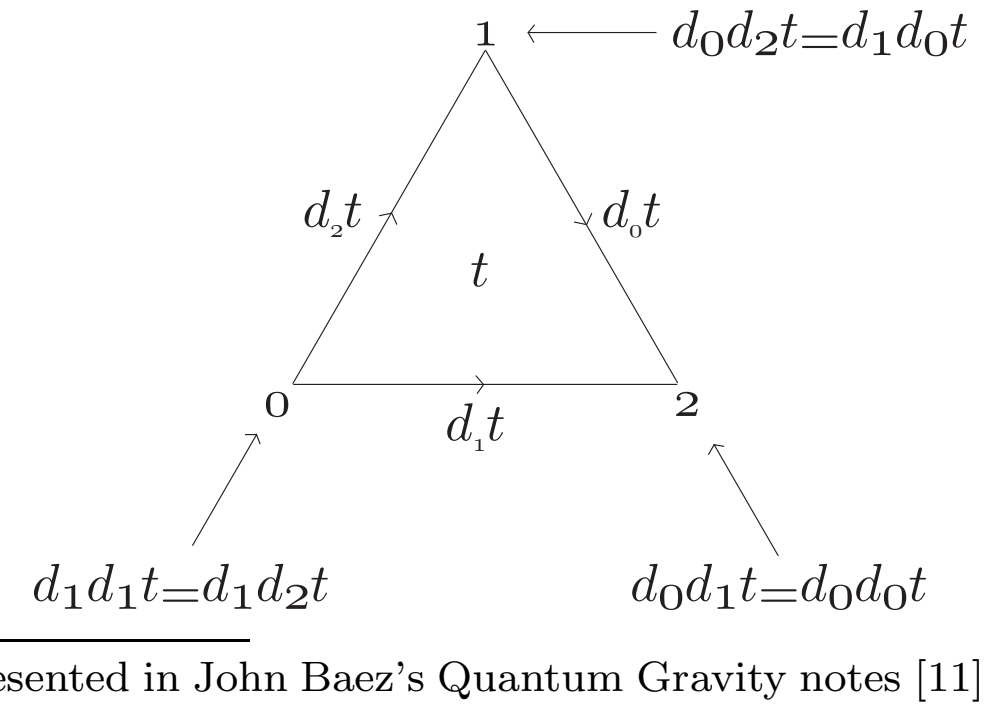

[1]This material is most elegantly presented in John Baez's Quantum Gravity notes [11], a series which is highly recommended.

A triangulated space is a discrete model of a smooth surface. To any surface $M$, one can assign a simplicial 2-graph, and moreover, to any simplicial 2-graph, one can create a smooth surface $M$, the geometric realization of the simplicial data. The holes in the surface come from edges which do not arise as the boundaries, via $d_i$, of triangles.

Let us now investigate what principal bundles and connections look like on a triangulated space $K$. The vertices play the role of spacetime. All principal bundles over $K$ are trivial since we can always find a global section (a map from the vertices to the fibers). A path in $K$ is a sequence $e_1 e_2 \cdots e_n$ of edges from $E$, such that $d_1(e_i) = d_0(e_{i+1})$ for all $i = 1 \ldots (n-1)$. A connection is just an assignment of a group element to each edge,

$$\Gamma : E \to G. \tag{4.16}$$

A gauge transformation is an assignment of a group element to each vertex,

$$s : V \to G. \tag{4.17}$$

These act on connections in the same way as in Eqn. (4.8).

In order to define what it means for a connection to be flat, we need the notion of homotopy of paths in a triangulated space. One does this by passing to the geometric realization of $K$, where homotopy has its natural topological meaning, and then pulling this definition back to the algebraic data. The result is that two paths $P_1$ and $P_2$ are homotopic if and only if $P_2$ can be obtained from $P_1$ by a process of 'sliding edges over triangles' from $T$:

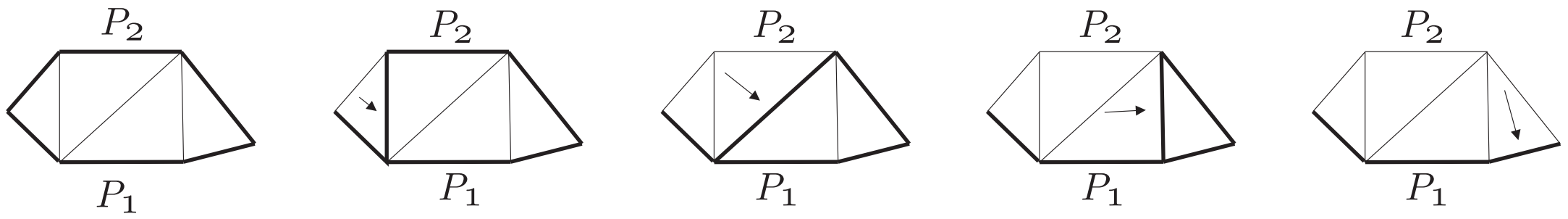

This will not be possible when the necessary triangles are not present in $T$. We define a flat connection on a triangulated space as one which is invariant under simplicial homotopy. If we define the 'simplicial homotopy' group $\pi_1(K)$ as the group of homotopy equivalence classes of paths (from some chosen fixed point), then a connection $\Gamma$ is flat if and only if its value for closed paths factors through the projection $P \to [P]$. In other words, the group multiplication must close around a contractible loop. The moduli space $\mathcal{M}$ of flat connections on $K$ is given by the same formula as in Eqn. (4.11).

For closed manifolds $K_g$ of genus $g$, one can often provide an explicit formula for $\mathcal{M}$, since $\pi_1(K_g)$ is known to be a group with $2g$ generators $a_1, b_1, \ldots, a_g, b_g$ satisfying the single relation

$$\prod_{i=1}^{g} [a_i, b_i] = 1 \tag{4.18}$$

where $[a, b]$ is the commutator $aba^{-1}b^{-1}$.

**Example : the torus**

What is the moduli space of flat $\mathbb{Z}/2$ connections[2] i? Here is a triangulation of the torus:

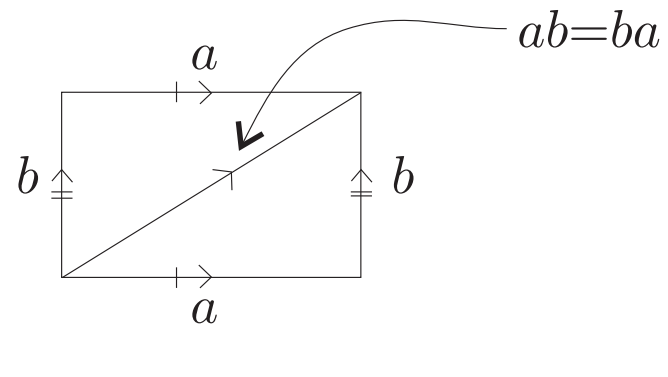

One may choose $a$ and $b$ freely, but flatness requires that the diagonal edge be equal to $ab$ $(= ba)$. Thus

$$|\mathcal{M}_g| = 4. \tag{4.19}$$

[2]Strictly speaking, there is only one connection on a principal $G$-bundle on the torus, so the term 'connection' really means 'bundle' here. See Section 4.2.

## 4.2 Dijkgraaf-Witten model

The Dijkgraaf-Witten [25] model is a TQFT (in any number of dimensions) associated to a finite group $G$, which functions as the gauge group[3]. In two dimensions, it is a discrete model for $BF$ theory (also known as topological Yang-Mills) while in three dimensions it is a discrete model for Chern-Simons theory. Dijkgraaf and Witten also demonstrated how to define deformed, twisted versions of the theory using a cohomology class $\alpha \in H^4(G, \mathbb{Z})$. The finer details of this construction were provided by Freed and Quinn [40], and an excellent introduction, explaining how quantum groups arise in the three-dimensional version of the model, can be found in Freed's lectures [33].

The fields in the finite gauge group model are principal $G$-bundles $P$ over a manifold $M$. At first glance it is somewhat strange to consider a nonlocal object like a principal bundle as a field. Observe though that a finite $G$-bundle is composed of 'sheets', locally of order $|G|$, identified together in various ways. This is very similar to the branches a holomorphic function $f$ on the complex plane. Also, principal bundles *glue together* like fields, in the sense that one needs only to match the restrictions of two bundles on their common boundary in order to glue them together.

For example, here are the two ($\mathbb{Z}/2$-principal bundles) on the circle $S^1$:

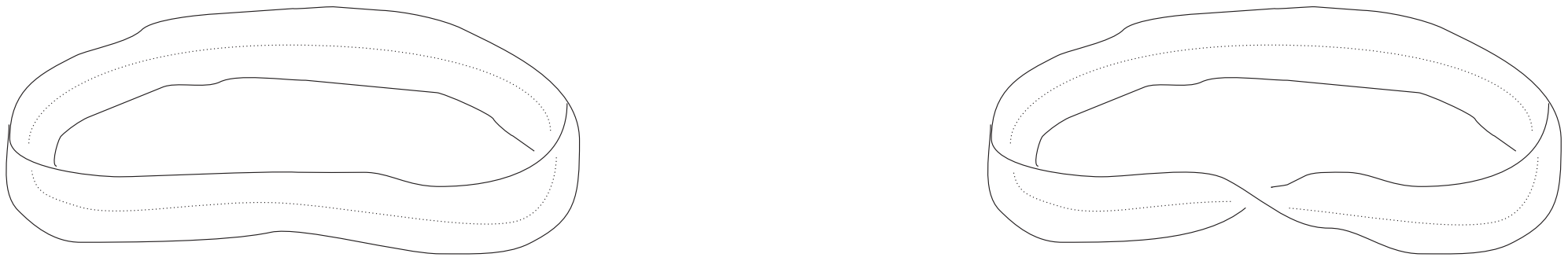

(4.20)

We want to talk about the Dijkgraaf-Witten model in gauge theory language, and hence would like to maintain the idea that the holonomies of connections are the fundamental variables. But if the gauge group is finite, then there is no freedom in the vertical (fiber) direction. For a given principal $G$-bundle, there is only one connection, and it is always flat in the topological sense. Thus the only degree of freedom is to vary the principal bundle $P$ itself. Allowing the topological type of the principal bundle to vary is a component of all gauge theories (eg. solitons). In the finite gauge group model, it is the *only* degree of freedom.

Recall once more the time evolution formula for gauge quantum field theory,

$$\langle \hat{A}_2 | U | \hat{A}_1 \rangle = \int_{A|_{\Sigma_1}=A_1}^{A|_{\Sigma_2}=A_2} \mathcal{DA} \exp iS[A], \tag{4.21}$$

which expresses the probability for a wave function $\hat{A}_1$ sharply peaked around the connection $A_1$ on $\Sigma_1$ to evolve to the wave function $\hat{A}_2$ sharply peaked around the connection $A_2$ on $\Sigma_2$. We are instructed to sum over all connections $A$ on $M$ such that $A|_{\Sigma_1} = A_1$ and $A|_{\Sigma_2} = A_2$, weighting each $A$ with the action $e^{iS[A]}$. More precisely, one should reformulate this procedure such that the sum only runs over gauge-equivalence classes of connections.

The Dijkgraaf-Witten topological field theory $Z$ is defined in the same spirit. For a closed $(d-1)$ manifold $\Sigma$, let $C_\Sigma$ be the space of fields on $\Sigma$,

$$C_\Sigma = \{G\text{-bundles } P \to \Sigma\}. \tag{4.22}$$

Regarding fields as principal bundles marks a subtle difference with ordinary quantum field theory, where the space of fields forms a set. In contrast, the space of 'all principal bundles over $\Sigma$' cannot technically be a set[4]. On the other hand, this allows us to view the underlying symmetries of the theory as extra structure. We should regard $C_\Sigma$ as a category (more precisely, a groupoid), whose objects are $G$-bundles and whose morphisms are maps which commute with the $G$-action (each such morphism must be an isomorphism since

[3]The word 'gauge group' is being used in the traditional physics way here, as the group $G$ which acts on the principal bundle. This is not the group $GA(P)$ from (4.1).

[4]Otherwise one runs into set theory conundrums such as Russell's paradox.

the $G$-action is free and transitive). This gives a picture of $C_\Sigma$ which we can draw as:

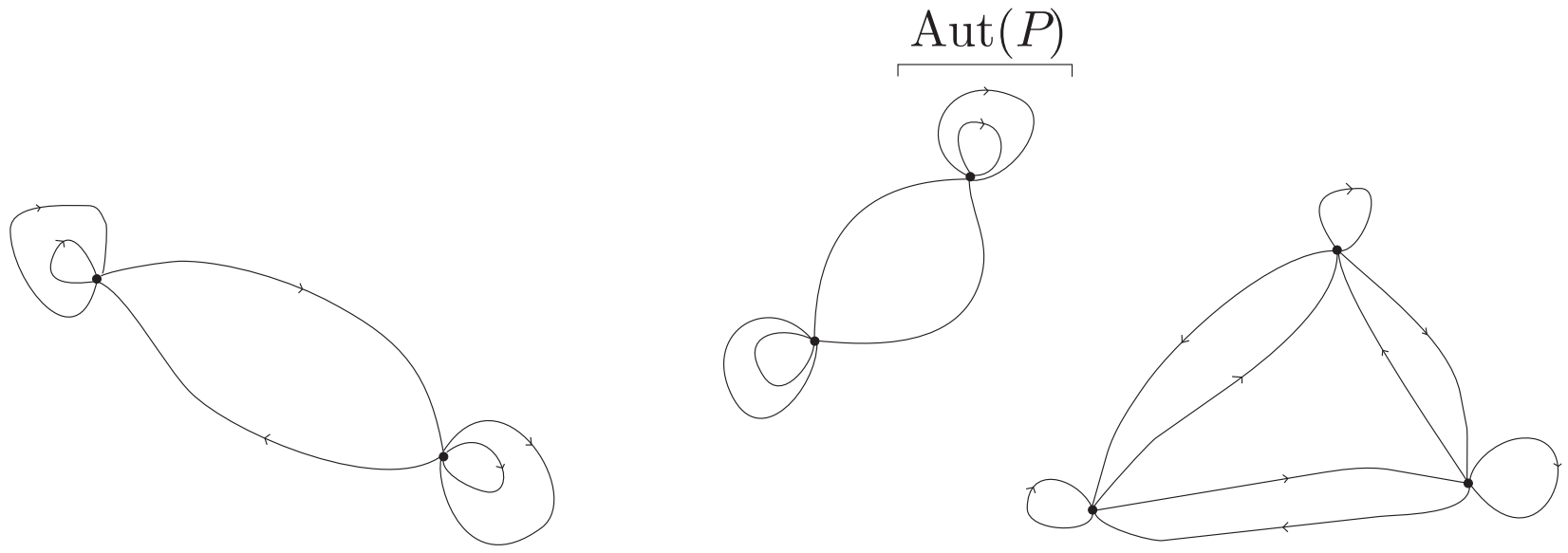

(4.23)

Now let $\overline{C_\Sigma}$ be the (finite) set of equivalence classes of fields in $C\Sigma$. Then we set

$$Z(\Sigma) := \{\text{functions } f : \overline{C_\Sigma} \to \mathbb{C}), \tag{4.24}$$

the space of wavefunctions of the fields. $Z(\overline{C_\Sigma})$ has a canonical basis of functions $\{\hat{P}\}$ whose value on $\overline{P}$ is 1 and zero otherwise.

For $M$ a $d$-manifold with boundary, we extend the notation $\overline{C_M}$ to refer to the finite set of isomorphism classes of $G$-bundles over $M$. In the Atiyah picture, to complete the description of $Z$ we need to define $Z(M) \in Z(\partial M)$, and this is done as follows:

$$Z(M)(Q) = \sum_{\overline{P} \in C_M : P|_{\partial M} = Q} \frac{1}{\text{Aut}(P)}. \tag{4.25}$$

The sum is understood to be over all equivalence classes of bundles $\overline{P}$ which have a representative $P$ which restricts to $Q$ on $\partial M$, weighting each $P$ with the reciprocal of the size of its automorphism group. If $\partial M = \Sigma_1 \sqcup \cdots \sqcup \Sigma_k$ then $Q|\partial M$ is specified by its restriction $Q_1, \ldots, Q_k$ to each $\Sigma_i$, so that we can view $Z(M)$ as a tensor $Z(M)_{Q_1 Q_2 \cdots Q_k}$, or as a $k$-point function $\langle Q_1 Q_2 \cdots Q_k \rangle_M$. The 2-point functions on the cylinder $\Sigma \times I$ define a metric on $Z(\Sigma)$:

$$g_{QR} = \langle Q\, R \rangle_{\Sigma \times I} =$$
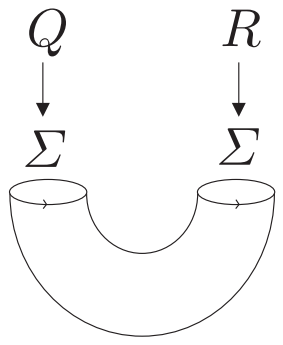

. (4.26)

From this we define the inverse metric $g^{QR} = g^{-1}_{QR}$. Suppose we wish to view a $d$-manifold $M$ with boundary $\partial M = \Sigma_1 \sqcup \Sigma_2 \sqcup \cdots \sqcup \Sigma_k$ as a cobordism $M : \Sigma_{\text{in}} \to \Sigma_{\text{out}}$, by choosing $l$ of the components of $\partial M$ as input boundaries and the remaining $k - l$ components as output boundaries, i.e.

$$\Sigma_{\text{in}} = \Sigma_{j_1} \sqcup \cdots \sqcup \Sigma_{j_l}, \tag{4.27}$$
$$\Sigma_{\text{out}} = \Sigma_{j_{l+1}} \sqcup \cdots \sqcup \Sigma_{j_k}. \tag{4.28}$$

Then we simply use the metric to raise indices appropriately,

$$Z_{Q_1, \cdots, Q_k} \to Z_{Q_1, \cdots, Q_k}^{Q_{j_{l+1}} \cdots Q_{j_k}} = g^{Q_{j_{l+1}} R_1} \cdots g^{Q_{j_{l+k}} R_k} Z_{Q_1, \cdots, Q_k}, \tag{4.29}$$

which allows us to define

$$Z(M) : Z(\Sigma_{j_1}) \otimes \cdots \otimes Z(\Sigma_{j_l}) \to Z(\Sigma_{j_{l+1}}) \otimes \cdots \otimes Z(\Sigma_{j_k}). \tag{4.30}$$

In other words we define $Z(M) : Z(\Sigma_{\text{in}}) \to \Sigma_{\text{out}}$ by

$$\langle \bar{P}_2 | Z(M) | \bar{P}_1 \rangle = \sum_{P|_{\Sigma_1} = P_1}^{P|_{\Sigma_2} = P_2} \frac{g(\partial M)}{|\text{Aut}(P)|}, \tag{4.31}$$

where $g(\partial M)$ is a sequence of metric-raising operations containing one factor for each component of $\Sigma_{\text{out}}$. Equation (4.31) takes precisely the same form as (4.21), with the ill-defined measure $\mathcal{DA}$ basically being replaced by the well-defined measure $1/|\text{Aut}(P)|$; this is what makes finite group gauge theory tractable.

This completes the description of the Dijkgraaf-Witten model. Let us now compute further the space of fields, as well as the group $\text{Aut}(P)$. Fix a point $x \in M$ and a set $X$ equipped with a free and transitive $G$-action. There is no harm if we restrict our attention to $G$-bundles $P \to M$ such that $P_x = X$. Also, fix an element $p_0 \in X$. Parallel transport around a loop $\sigma$ based at $x$ will send $p_0 \to p_0 \cdot g$. Taken together, these holonomies determine a homomorphism $\phi_P : \pi_1(M, x) \to G$. This homomorphism in turn completely determines the bundle, which we can thus write as $P_\phi$. Note that if we had chosen a different reference point $p_0' = p_0 \cdot h$ the holonomy would have been conjugated $g \to h^{-1}gh)$. Thus there is a right action of $G$ on $\text{Hom}(\pi_1(M, X), G)$, defined by conjugating the image of $\phi$,

$$(\phi \cdot g)(\sigma) = g^{-1}\phi(\sigma)g, \tag{4.32}$$

and $P_{\phi_1}$ is isomorphic to $P_{\phi_2}$ when $\phi_2 = \phi_1 \cdot g$ for some $g$ (i.e. they lie in the same orbit). Thus

$$\mathcal{P}_\Sigma = \text{Hom}(\pi_1(M, x), G)/G. \tag{4.33}$$

Taken together with (4.11), we have shown that there is a $1-1$ correspondence between:

$$\text{Flat } G\text{-connections on a triangulation of } M \overset{\text{1-1}}{\longleftrightarrow} G\text{-bundles over } M. \tag{4.34}$$

In our setup, a automorphism of the bundle $f : P \to P$ is uniquely determined by what $f$ does to $p_0$, i.e. by the $g \in G$ such that $f(p_0) = p_o \cdot g$. For then one may define $f$ on $X$ by setting $f(p_0 \cdot h) := f(p_0) \cdot h = p_0 \cdot gh$, and extend $f$ to the whole of $P$ by continuity. However, not all choices of $g$ are allowed - continuity requires that $f$ must be consistent if one travels around a loop $\sigma$:

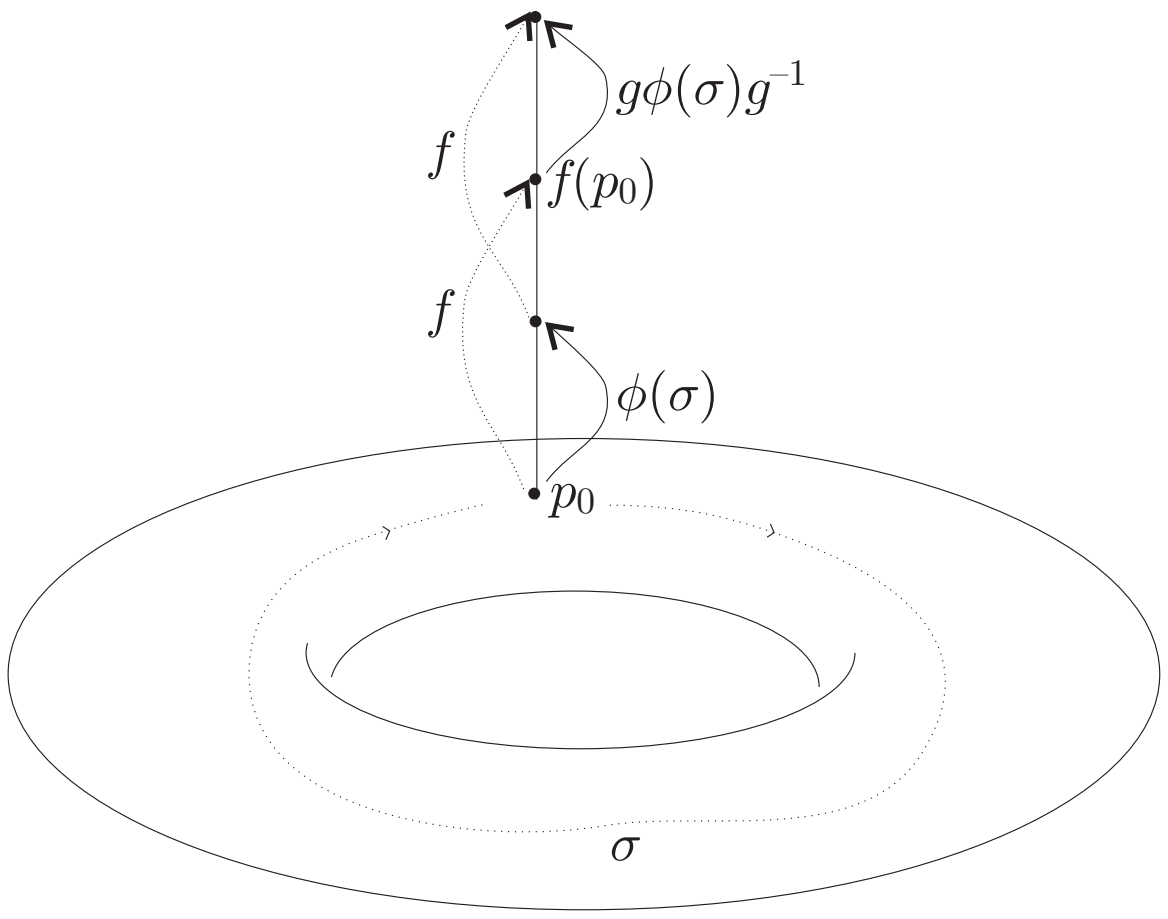

That is, $f(p_0 \cdot \phi(\sigma)) = f(p_0) \cdot g^{-1}\phi(\sigma)g$. Thus we require that $g$ commutes with the image of $\phi_P$, i.e.

$$\text{Aut}(P_\phi) = \{g \in G : g\phi(\cdot) = \phi(\cdot)g\}. \tag{4.35}$$

Notice the simplification achieved; when $G$ is continuous, the collection of gauge transformations $\text{Aut}(P)$ is an infinite dimensional group, while in our case $\text{Aut}(P)$ is actually a subgroup of $G$.

The topological invariant $Z(M)$ that the model assigns to a closed $d$-manifolds $M$ is, by (4.25), simply the number of $G$-bundles on $M$, each weighted by the reciprocal of the size of its automorphism group. We have seen how isomorphism classes of $G$-bundles on $M$ are in 1-1 correspondence with the number of orbits in $\text{Hom}(\pi_1(M, x), G)$ under the $G$-action (4.32). Eqn. (4.35) shows that under this correspondence, $\text{Aut}(P_\phi)$ is realized as the stability subgroup $\text{Stab}(\phi)$. On the other hand, by a general property of $G$-sets we must have

$$|G| = |\text{Orbit}(\phi)||\text{Stab}(\phi)|. \tag{4.36}$$

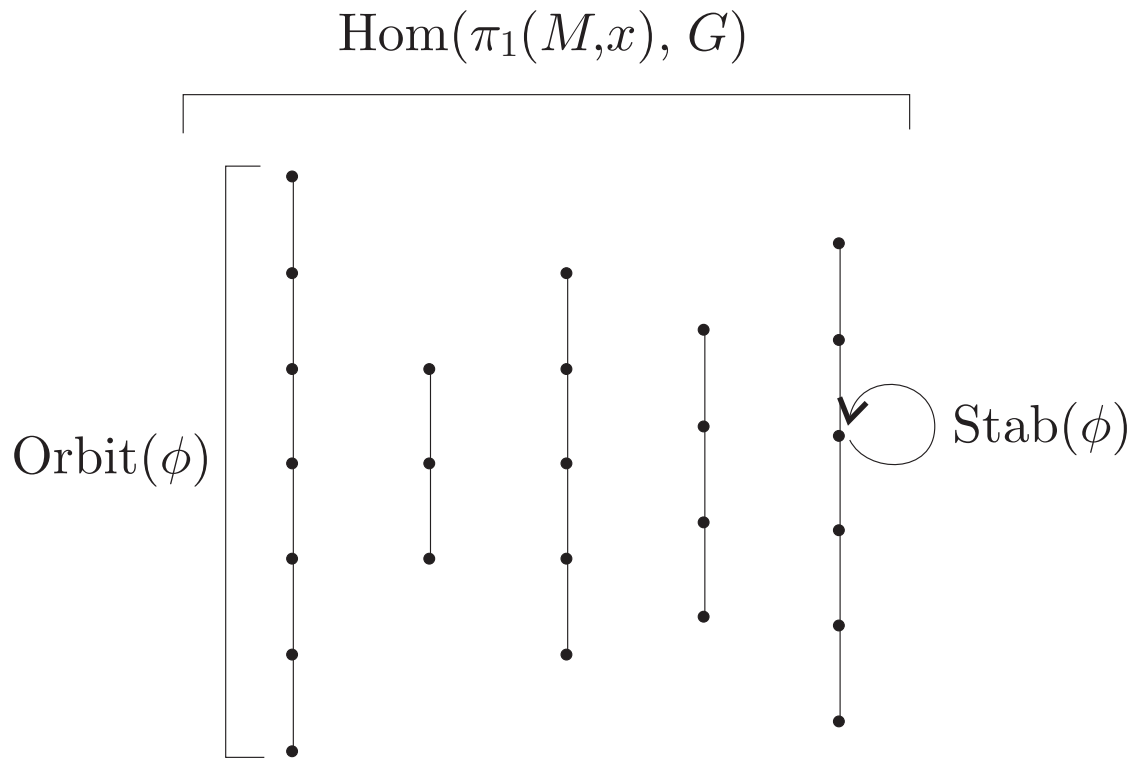

Figure 4.1: The space $\mathrm{Hom}(\pi_1(x, M), G)$ decomposes under the adjoint action into orbits which are the equivalence classes of bundles $P \to M$.

This allows us to compute (see Fig. 4.1)

$$\begin{aligned} Z(M) &= \sum_{\text{Orbits}} \frac{1}{\mathrm{Stab}(\phi)} \\ &= \sum_{\text{Orbits}} \frac{|\mathrm{Orbit}(\phi)|}{|G|} \\ &= \frac{1}{|G|} \sum_{\text{Orbits}} |\mathrm{Orbit}| \\ &= \frac{|\mathrm{Hom}(\pi_1(M,x),G)|}{|G|}. \end{aligned} \tag{4.37}$$

This is obviously a topological invariant of $M$, as expected. Also observe that had we not chosen the weighting 1Aut($P$) in (4.25), we would not have been able to write $Z(M)$ explicitly in terms of $|\mathrm{Hom}(\pi_1(M,x),G)|$. For instance, if we had weighted each bundle with unit weight, then we would have had

$$Z(M) = |\mathrm{Hom}(\pi_1(M,x),G)/G|, \tag{4.38}$$

and this is combinatorially a far more intricate number to calculate[5].

The line of reasoning in (4.37) can be repeated to allow us to calculate the vectors $Z(M) \in Z(\partial M)$ when $M$ has a boundary. Thus we may compute (4.25) via the homotopy groups as:

$$Z(M)(Q) = \frac{1}{|G|}|\{\phi \in \mathrm{Hom}(\pi_1(M,x),G)\,:\,\phi \text{ `restricts' to a bundle isomorphic to } Q \text{ on } \partial M.\}|. \tag{4.39}$$

To say that $\phi$ 'restricts' to a bundle on $\partial M$ means that we choose special points $y_i$ for each component $\Sigma_1$ of $\partial M$, $i = 1 \dots n$, as well as paths $\sigma_i : x \to y_i$ in $M$ for each $i$, so that we obtain a restriction map

$$\mathrm{Hom}(\pi_1(M,x),G) \xrightarrow{\partial} \mathrm{Hom}(\pi_1(\Sigma_1,y_1),G) \times \cdots \times \mathrm{Hom}(\pi_1(\Sigma_n,y_n),G). \tag{4.40}$$

Choosing different $x$, $y_i$ or $\sigma_i$ changes this restriction map only by conjugation, and since this does not change the associated bundles on $\Sigma_i$, the restriction operation is well-defined on the level of bundles.

Let us once more summarize the main idea of the Dijkgraaf-Witten model. It describes the quantum evolution of *wave functions of principal bundles*. The probability of one bundle to evolve into another is determined by the weighted number of bundles which interpolate between the two. In particular, the probability is zero if there is no way to interpolate from the input to the output bundle.

[5]On the other hand, these numbers are precisely the dimensions of the Hilbert spaces in the 3d version of the theory - see Sec. 5.6 for a comparison.

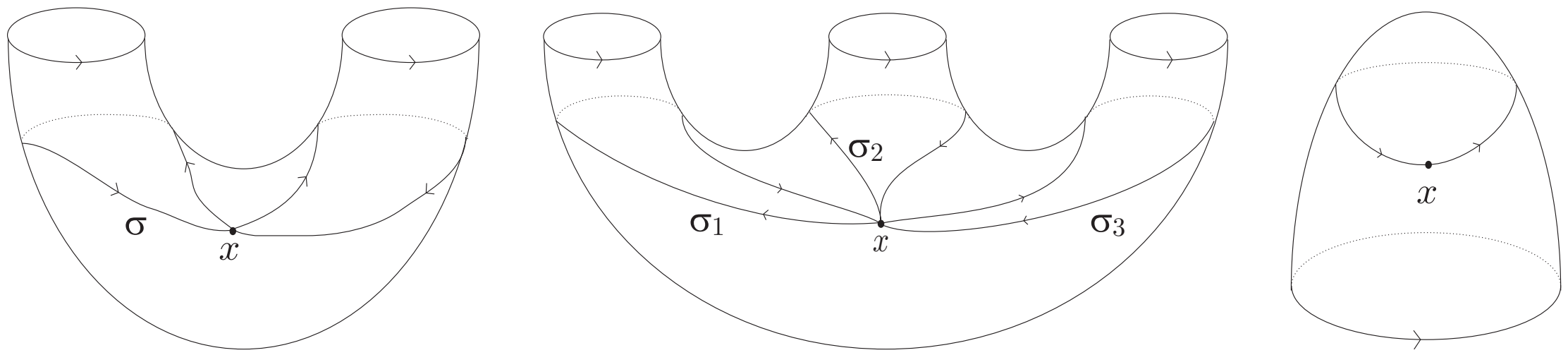

(a) (b) (c)

Figure 4.2: (a) For the bent cylinder there is only a single nontrivial loop $\sigma$, whose orientation matches one of the circles but is reversed with respect to the other. (b) The trinion has three fundamental loops $\sigma_1$, $\sigma_2$, $\sigma_3$ satisfying the relation $\sigma_1\sigma_2\sigma_3 = 1$. (c) The unit disc has trivial homotopy.

## 4.3 The Dijkgraaf-Witten model in 2d and $C_{\text{class}}(G)$

We shall now compute the model explicitly in two dimensions, and show that the underlying Frobenius algebra is the space $C_{\text{class}}(G)$ of class functions on $G$ (see Section 3.3.3). If nothing else, the Dijkgraaf-Witten model thus provides a beautiful geometric interpretation of $C_{\text{class}}(G)$ in terms of the quantum evolution of principal bundles over 2d surfaces.

We know from Chapter 2 that to compute a 2d TQFT it is enough to compute the Hilbert space $Z(S^1)$, the pairing $Z(\;)$, the multiplication map $Z(\;)$ and the unit map $Z(\;)$. Fix a point $x \in S^1$. Then $\pi_1(S^1, x) = \mathbb{Z}$ so that to give a map $\phi : \pi_1(S^1, x) \to G$ is to give an arbitrary $g \in G$. Conjugation of $\phi$ corresponds to conjugation of $g$ so that by (4.33),

$$Z(S^1) = C_{\text{class}}(G). \tag{4.41}$$

From now on we shall represent a bundle $Q$ over $S^1$ by its corresponding conjugacy class $\alpha_Q$. To compute $Z(\;)$ we should first compute the metric $Z(\;)$. Fig. 4.2a shows that if $Q$ and $R$ (represented by $\alpha$ and $\beta$ respectively) are bundles living on the two boundary circles of , then by (4.39) we have

$$g_{\alpha\beta} \equiv Z_{\alpha\beta} = \delta_{\alpha\beta^{-1}} \frac{|\alpha|}{|G|}. \tag{4.42}$$

The inverse metric is $g^{\alpha\beta} = \delta^{\alpha\beta^{-1}} \frac{|G|}{|\alpha|}$. Fig. 4.2b shows that the three-point function is given by

$$Z_{\alpha\beta\gamma} = \frac{1}{|G|} |\{g_1 \in \alpha, g_2 \in \beta, g_3 \in \gamma : g_1 g_2 g_3 = 1\}. \tag{4.43}$$

Thus by (4.30) we can compute the multiplication $e_\alpha \otimes e_\beta \to Z^\gamma_{\alpha\beta} e_\gamma$ as

$$Z^\gamma_{\alpha\beta} = g^{\kappa\gamma} Z_{\alpha\beta\kappa} = |\{h \in \beta : gh^{-1} \in \alpha\}|. \tag{4.44}$$

Finally Fig. 4.2c shows that the unit computes as $1 = \{e\}$, the conjugacy class of the identity. The same figure tells us that the Frobenius form computes as

$$\epsilon(\alpha) = \frac{1}{|G|} \delta_{\alpha,1}, \tag{4.45}$$

since the conjugacy class of the identity is trivial. Comparing (4.42), (4.44) and (4.45) with (3.30) and (3.31) shows that $Z(S^1)$ carries the same Frobenius algebra structure as $C_{\text{class}}(G)$.

## 4.4 The triangulation construction in 2d

### 4.4.1 Remarks on Lattice TQFT's.

There is an alternative formulation of the Dijkgraaf-Witten model in terms of lattice gauge theory. This construction works in any number of dimensions, but we shall specialize to two dimensions. These ideas originated with the work of Fukama, Hosono and Kawai (abbreviated FHK) [34], although we shall, for the most part, follow John Baez's elegant reformulation [13].

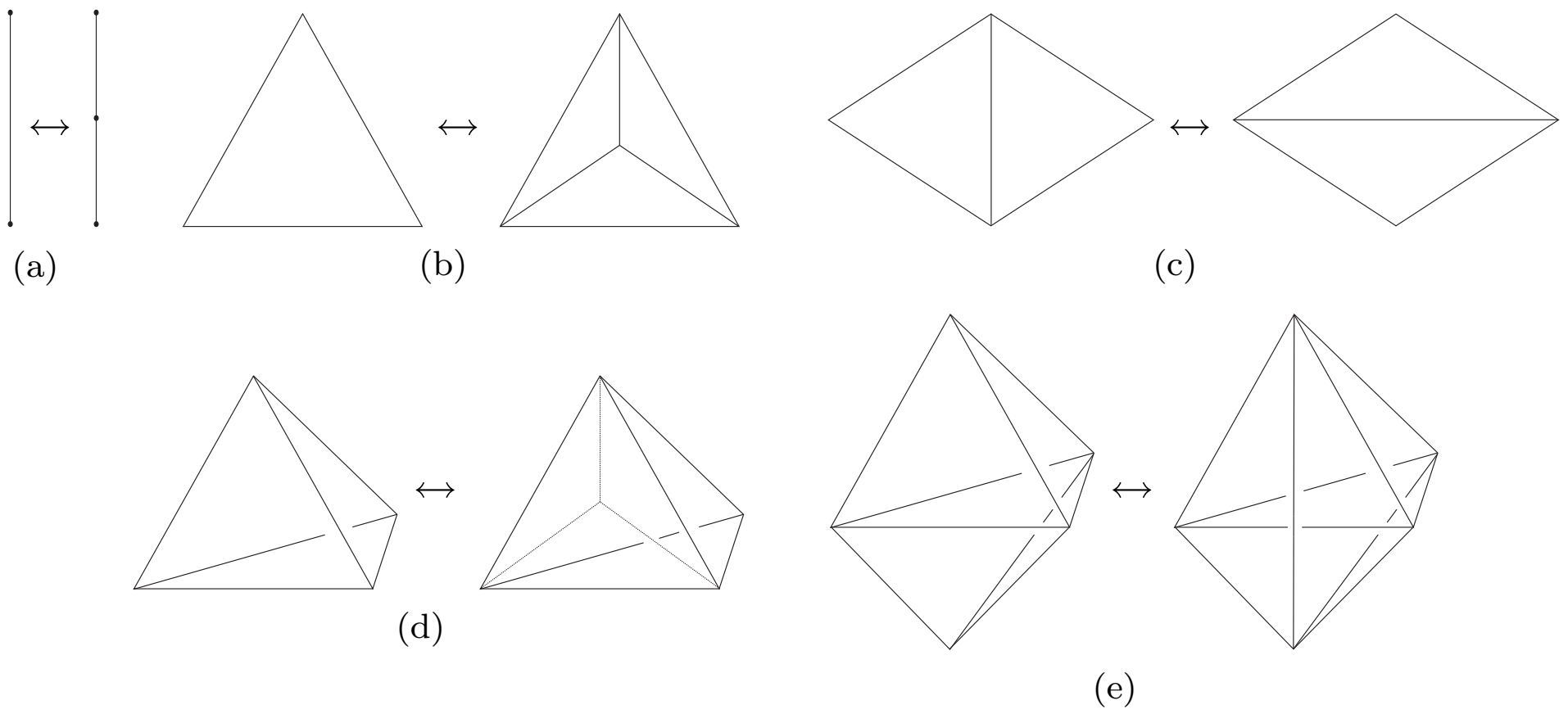

Figure 4.3: The Pachner moves in dimensions 1, 2 and 3. (a) : Subdivision of a line in one dimension. (b) and (c) : The 1-3 and 2-2 moves in two dimensions. (d) and (e) : The 2-3 move and 1-4 move in three dimensions.

The notion of producing TQFT's (in any dimensions) from combinatorial data associated to triangulations is a very important one, and has been studied by many authors. Characterizing manifolds and the cobordism category algebraically is an extremely difficult problem and apparently tractable only for $n = 2$ (although see [49] for progress on $n = 3$). On the other hand, the local 'Pachner moves' [71] necessary to pass between any two triangulations of a manifold are known for all dimensions, and correspond to viewing an $n$-simplex from the back and from the front (see Fig. 4.3 ). This observation allows one to construct a full-blown TQFT $Z$ in two stages. Firstly construct a *lattice* TQFT $\tilde{Z}$ using algebraic data which are invariant under the triangulation moves. Although $\tilde{Z}$ is independent of the triangulation in the bulk, it will in general still depend on the triangulation on the boundary. The second stage is to remove this boundary triangulation dependence - a process which can be elegantly formulated in category theory terms as the *colimit* under coarse-graining of the triangulation [102]:

$$\begin{array}{ccc} \text{Algebraic object } A \text{ invariant under Pachner moves} & \longrightarrow & \text{Lattice TQFT } \tilde{Z} \\ & & \big\downarrow \text{colimit} \\ \text{Algebraic object } B \text{ which resembles } \mathbf{nCob} & \longrightarrow & \text{TQFT } Z \end{array} \tag{4.46}$$

The advantage gained is that $A$ is a far simpler object than $B$, and the difference in complexity grows rapidly as $n$ increases. The price to pay is that not all TQFT's can be obtained in this way. In addition, the category of piecewise linear manifolds (constructed from triangulations) is not equivalent to the category of ordinary manifolds for $n \geq 6$ [10], although this is hardly a physical problem if we are modest about the dimension of spacetime! Finally, the output TQFT $Z$ often remains difficult to compute explicitly, even though its existence is assured.

For instance, in two dimensions the Pachner moves turn out to require a semisimple algebra $A$ (which doesn't have to be commutative). The process $\tilde{Z} \to Z$ corresponds to taking the center $Z(A)$ of the algebra (do not confuse the two meanings of $Z$!). Here $Z(A)$ is viewed as a commutative Frobenius algebra by inheriting the canonical form $(a, b) = \mathrm{Tr}(L_a L_b)$ from $A$:

$$\begin{array}{ccc} \text{Semisimple algebra } A & \longrightarrow & \text{Lattice TQFT } \tilde{Z} \\ & & \big\downarrow \text{take center } A \to Z(A) \\ \text{Commutative Frobenius algebra B} & \longrightarrow & \text{TQFT } Z \end{array} \tag{4.47}$$

By Wedderburn's theorem [5], a semisimple algebra $A$ is isomorphic to a direct sum of matrix algebras,

$$A \simeq M_{n_1} \oplus M_{n_2} \oplus \cdots \oplus M_{n_k}. \tag{4.48}$$

Thus to specify $A$ is simply to specify a list of $k$ integers $(n_1, n_2, \ldots, n_k)$. On the other hand, specifying $B$ is more difficult, although in the semisimple case it amounts to a choice of $k$ nonzero complex numbers $(\lambda_1, \lambda_2, \ldots, \lambda_k)$ (see Sec. 3.3.4). Evidently for $n = 2$ there is not a large difference in complexity between the two approaches.

A sneak preview of Chapter 4 is in order. In three dimensions, the Pachner moves turn out to require a 'semisimple 2-algebra' $A$ (a 2-algebra is the categorification of an algebra). On the other hand, to specify a full-blown TQFT requires a modular category $B$, which can be regarded as a 2-algebra equipped with nontrivial data such as a braiding (the analogue of commutativity in 2d), a twist, and duality morphisms. Thus we already see the complexity gap between $A$ and $B$ growing. In fact, $B$ can be regarded, as a braided monoidal category, as the center of the monoidal category $A$! In other words, the entire diagram (4.47) has been categorified in a straightforward fashion - a remarkable tribute to the power of categorical language:

$$\begin{array}{ccc} \text{Semisimple 2-algebra } A & \longrightarrow & \text{Lattice TQFT } \tilde{Z} \\ & & \Big\downarrow \scriptstyle{\text{take center of monoidal category } A \to Z(A)} \\ \text{Modular category B} & \longrightarrow & \text{TQFT } Z \end{array} \tag{4.49}$$

Returning to two dimensions, we make contact with gauge theory by choosing $A = C(G)$ (or equivalently $\mathbb{C}[G]$). The Peter-Weyl theorem tells us that $C(G)$ is semisimple since it decomposes as a direct sum of matrix algebras over the irreducible representations $\rho_1, \rho_2, \ldots, \rho_n$ of $G$,

$$\mathbb{C}[G] \simeq \mathrm{End}(V_{\rho_1}) \oplus \mathrm{End}(V_{\rho_2}) \oplus \cdots \oplus \mathrm{End}(V_{\rho_n}). \tag{4.50}$$

Since the output TQFT is $Z(C(G))$, we see that it is a lattice theory realization of the Dijkgraaf-Witten model.[6].

### 4.4.2 FHK construction

#### The lattice TQFT $\tilde{Z}$

Suppose we are given a semisimple algebra $A$ with basis $\{e_i\}$. For $A = \mathbb{C}[G]$ we can take the basis to be $\{e_h\}$ for $h \in G$. Define a metric on $A$ by

$$g_{ij} := \mathrm{Tr}(L_{e_i} L_{e_j}) = g_{ji}. \tag{4.51}$$

Here $L_a : A \to A$ is the multiplication map $x \to ax$. Since $A$ is semisimple, it is a direct sum of matrix algebras (4.48), so that $g_{ij}$ is the normal trace pairing of matrices. Thus $g_{ij}$ has an inverse which we denote by $g^{ij}$. In fact, it is not hard to show the converse (see [34]), i.e. that

$$\text{A is semisimple (a finite sum of matrix algebras)} \Leftrightarrow g_{ij} \text{ in (4.51) is nondegenerate.} \tag{4.52}$$

Now define the algebra coefficients $m_{ij}^k$ by

$$e_i e_j = m_{ij}^k e_k. \tag{4.53}$$

Note that we can express the metric in terms of these coefficients as

$$g_{ij} := m_{ik}^l m_{jl}^k = g_{ji} \tag{4.54}$$

We define the three-point functions $m_{ijk}$ which are symmetric under cyclic permutations of the indices by

$$m_{ijk} := g_{lk} m_{ij}^l = \mathrm{Tr}(L_{e_i} L_{e_j} L_{e_k}). \tag{4.55}$$

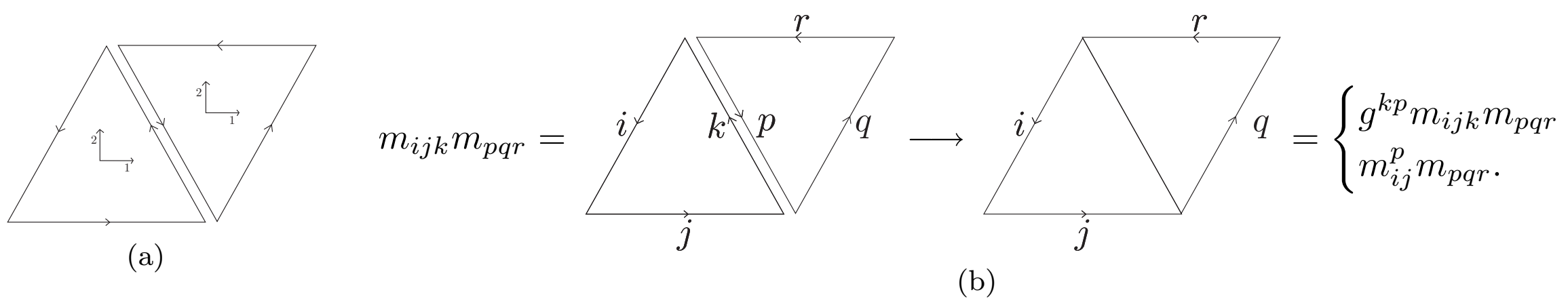

Figure 4.4: (a) Use the orientation $(1, 2)$ on $M$ to define an orientation on each edge by laying vector 1 parallel to the edge with vector 2 pointing inwards. Note common edges are oriented in opposite ways. (b) Gluing triangles together.

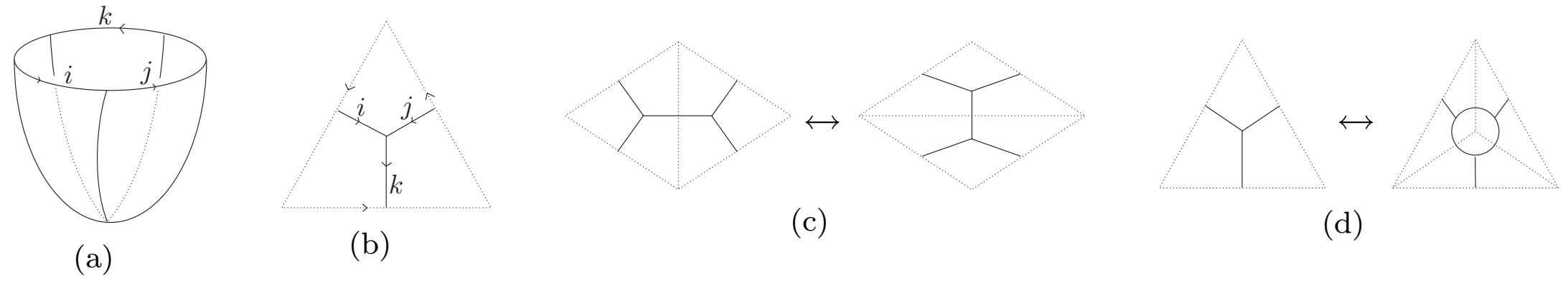

Figure 4.5: (a) Triangulation of unit disc. (b) Passing to the duals reinterprets the algorithm in terms of Feynman diagrams. (c) The 2-2 move is the Frobenius law ('crossing symmetry'). (d) The 1-3 move is the bubble move.

Now let $M$ be an oriented surface with boundary $\partial M$ (having the induced orientation). Triangulate $M$, assigning a cyclic orientation to the edges of each triangle by the 'inward-pointing' procedure, as in Fig. 4.4a. Assign indices (group elements) to each edge of the triangle, so that each triangle is assigned a symbol $m_{ijk}$. Now glue triangles together by contracting their indices with the metric, as in Fig. 4.4b. The output is a string of tensor symbols $m^{\cdots}_{\cdots}$, whose indices can be divided into unpaired lower indices which come from boundary edges, and pairs of upper and lower indices which are contracted with each other. We define the $n$-point function $\tilde{m}_{a_1 a_2 \cdots a_k}$ to be the resulting tensor. For example, to the triangulation in Fig. 4.5a we assign

$$\tilde{m}_{ijk} = m_{ip}^{c} m_{j}^{pd} m_{kcd} \tag{4.56}$$

For $A = \mathbb{C}[G]$ we calculate the metric as

$$g_{ij} = |G| \delta_i^{j^{-1}} \quad \leftrightarrow \quad g^{ij} = \frac{1}{|G|} \delta_i^{j^{-1}}, \tag{4.57}$$

Here we are employing the $\delta$ notation adapted for groups, i.e.

$$\delta_i^{jk} = \begin{cases} 1 & \text{if } i = jk \\ 0 & \text{otherwise.} \end{cases} \tag{4.58}$$

Similarly the algebra coefficients compute as

$$m_{ijk} = |G| \delta_{ijk}^{1}, \quad m_{ij}^{k} = \delta_{ij}^{k}, \quad m_i^{jk} = \frac{1}{|G|} \delta_i^{kj}, \quad m^{ijk} = \frac{1}{|G|^2} \delta_1^{ikj}. \tag{4.59}$$

That is, the three point function $m_{ijk}$ equals $|G|$ if the product of the group elements around the oriented edges is the unit, and is zero otherwise. Raising an index corresponds to reversing the orientation of an edge, and introducing a factor of $\frac{1}{|G|}$ (This can also be interpreted via the dual of the triangulation, see Fig. 4.5b ). Thus we are implementing precisely the philosophy of lattice gauge theory outlined in Sec. 4.1.3. The point is that we are *summing over flat connections on* $M$.

[6]The choice $A = C(G)$ is not very restrictive since 'most' algebras of the form (4.48) can be generated as the matrix algebras on the irreducible representations of some $G$. When $A$ is not of the form $C(G)$ then we are apparently doing gauge theory with a 'forbidden' gauge group!

Now let $\Sigma$ be a triangulated 1-manifold (a union of circles). To each edge of $\Sigma$ associate one copy of $A^*$, and tensor these together to form

$$\tilde{Z}(\partial\Sigma) = \bigotimes_{\text{no. of edges}} A^*. \tag{4.60}$$

We see that the tensor $\tilde{m}_{a_1 a_2 \cdots a_k}$ assigned to $M$ can be interpreted as a vector $Z(M) \in Z(\partial M)$,

$$\tilde{Z}(M)(e_1 \otimes e_2 \otimes \cdots \otimes e_k) = \tilde{m}_{a_1 a_2 \cdots a_k}. \tag{4.61}$$

We now show that $\tilde{Z}(M)$ is independent of the triangulation in the bulk of $M$, by establishing its invariance under the two Pachner moves in Fig. 4.3 (b) and (c). Note from Fig. 4.5 (c) and (d) that the dual of the 2-2 move becomes the familiar 'crossing symmetry' while the dual of the 1-3 move refers to the vanishing of '1-loop diagrams'. The 2-2 move requires that

$$m_{ij}^p m_{pk}^l = m_{jk}^p m_{ip}^l, \tag{4.62}$$

which is just the associative law for multiplication $(e_i e_j)e_k = e_i(e_j e_k)$. Similarly the 1-3 move requires that

$$g_{ij} = C_{ik}^l C_{jl}^k, \tag{4.63}$$

which, by (4.54), requires that the metric on $A$ be given by (4.51). Hence by (4.52) $A$ must be semisimple. Finally note that the existence of the unit $1_A \in A$ can also be inferred from the invariance under the Pachner moves. For the vanishing of the 1-loop functions allows us to define

$$1_A := $$ 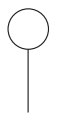 $$= g^{ij} m_{ij}^k e_k \tag{4.64}$$

which, taken together with (4.62) and (4.63), gives $1_A e_i = e_i = e_i 1_A$ for all $i$.

Suppose we wish to view $M$ as a cobordism $M : \Sigma_1 \to \Sigma_2$ by partitioning its boundary into input and output circles $\partial M = \Sigma_1 \sqcup \Sigma_2$. Then we employ the same trick as in (4.29), and use the metric to raise the indices of the edges on the output circles. To summarize, we have shown that there is a 1-1 correspondence between

$$\text{Lattice TQFT's } \tilde{Z} \text{ which sum over flat connections via (4.52)} \longleftrightarrow \text{Semisimple algebras } A. \tag{4.65}$$

**The procedure $\tilde{Z} \to Z$**

The assignment of $\tilde{Z}(M) : \tilde{Z}(\Sigma_1) \to \tilde{Z}(\Sigma_2)$ to a triangulated cobordism $M : \Sigma_1 \to \Sigma_2$ behaves in many ways like a functor. Firstly notice that $\tilde{Z}$ respects composition in the sense that

$$\tilde{Z}(M'M) = \tilde{Z}(M')\tilde{Z}(M). \tag{4.66}$$

where the input triangulation of $M'$ matches the output triangulation of $M$:

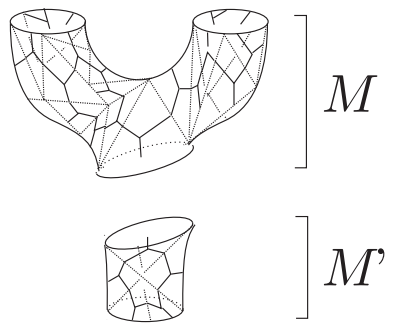

(4.67)

There is no natural identity map, since one can triangulate the input and output circles of a cylinder $\Sigma \times I$ in many ways. But consider those cylinders $\Sigma \times I$ whose input and output circles have identical triangulations. Then the corresponding linear operators are idempotents:

$$\begin{aligned} \tilde{Z}(\Sigma \times I)\tilde{Z}(\Sigma \times I) &= \tilde{Z}(\Sigma \times I)(\Sigma \times I) \qquad (\tilde{Z} \text{ is functorial}) \\ &= \tilde{Z}(\Sigma \times I) \qquad \text{(Triangulation invariance in the bulk)} \end{aligned} \tag{4.68}$$

We now define the functor $Z$ by restricting to the images of these idempotents. Firstly, set

$$Z(\Sigma) = \mathrm{Im}\tilde{Z}(\Sigma \times I) \subseteq \tilde{Z}(\Sigma) \tag{4.69}$$

We shall presently show that, up to isomorphism, this does not depend on the triangulation of $\Sigma$. First notice that one can consistently restrict the entire theory $\tilde{Z}$ to the subspaces $Z(\Sigma)$. This is because if $M : \Sigma_1 \to \Sigma_2$ is an arbitrary triangulated cobordism, then one can perform the following remarkable computation,

$$\begin{aligned} \tilde{Z}(M)Z(\Sigma) &= \tilde{Z}(M)\mathrm{Im}\tilde{Z}(\Sigma_1 \times I) && \text{(Defn.)} \\ &\subseteq \mathrm{Im}\tilde{Z}(M)\tilde{Z}(\Sigma \times I) \\ &= \mathrm{Im}\tilde{Z}(M(\Sigma \times I)) && \text{(By (4.66))} \\ &= \mathrm{Im}\tilde{Z}(M) && \text{(By (4.68))} \\ &= \mathrm{Im}\tilde{Z}((\Sigma_2 \times I)M) && \text{(By (4.66))} \\ &\subseteq \mathrm{Im}\tilde{Z}(\Sigma_2 \times I) \\ &= Z(\Sigma_2) && \text{(Defn.)} \end{aligned}$$

This computation is far better expressed graphically:

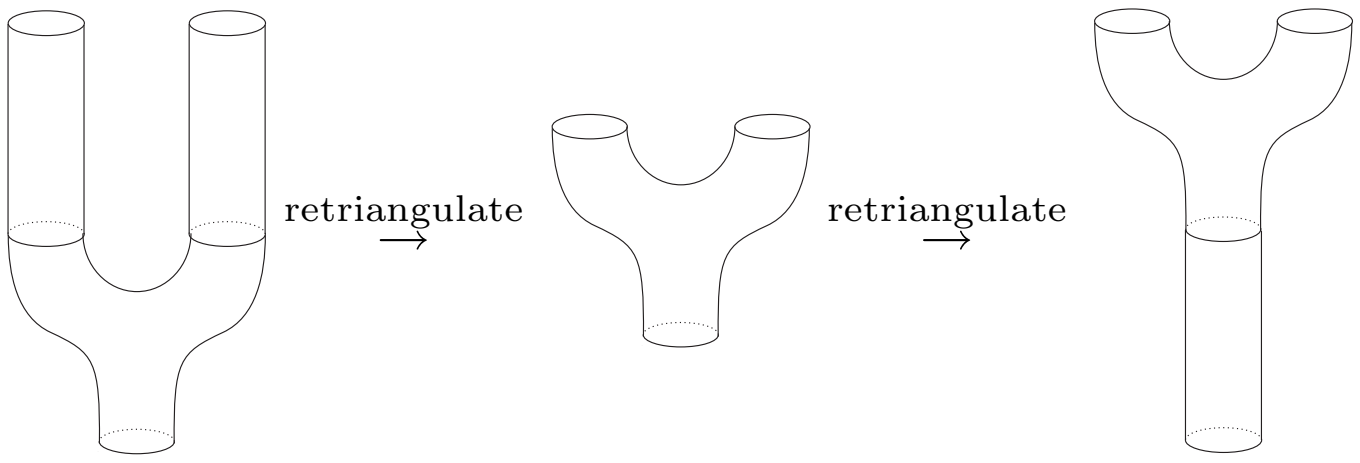

It allows us to define $Z(M) : Z(\Sigma_1) \to Z(\Sigma_2)$ as

$$Z(M) = \tilde{Z}(M)|_{Z(\Sigma_1)}. \tag{4.70}$$

Observe that $Z$ is a true functor from (triangulated) **2Cob** to **Vect**; it inherits $Z(M'M) = Z(M')Z(M)$ from $\tilde{Z}$, while $Z(\Sigma \times I) = \mathrm{id}$ since $Z(\Sigma \times I) = \tilde{Z}(\Sigma \times I)|_{Z(\Sigma)} = 1_{Z(\Sigma)}$ because $\tilde{Z}(\Sigma \times I)$ is a projection onto $Z(\Sigma)$.

Finally, suppose $\Sigma'$ and $\Sigma''$ are two different triangulations of the same underlying 1d manifold $\Sigma$. Then there is a canonical isomorphism

$$\alpha : Z(\Sigma') \xrightarrow{\sim} Z(\Sigma'') \tag{4.71}$$

given by $\alpha = Z(\Sigma \times I)$ where the triangulation on $\Sigma \times I$ is arbitrary, subject to the constraint that it matches $\Sigma'$ on its input boundary and $\Sigma''$ on its output boundary ($\alpha$ is an isomorphism since its inverse is given by the cylinder with triangulation reversed). Thus all traces of triangulation dependence have been removed and $Z$ is a full-blown 2d TQFT.

### Computing $Z$

We now calculate $Z(\Sigma) = \mathrm{Im}\tilde{Z}(\Sigma \times I)$ for $\Sigma = S^1$ and show that it is equal to the center of $A$. We triangulate $S^1 \times I$ as in Fig. 4.6a. We also show the Poincaré dual in Fig. 4.6b. Fattening out the Feynman diagram into a ribbon, as in Fig. 4.6c, shows that the map $Z(\Sigma \times I)$ is precisely the enigmatic map $\pi : A \to A$, so important in the discussion of open/closed string theory from Sec. 3.4![7] In other words, *we are trying to generate closed string algebra $B$ starting with the open string algebra $A$.*

We shall presently prove that $\mathrm{Im}(\pi) = Z(A)$ from first principles, using the elegant graphical argument presented in [13]. Observe though, that we can use the open/closed correspondence to establish that $\mathrm{Im}(\pi) \subseteq Z(A)$. Namely, we know that, if $B$ is the closed string theory corresponding to $A$, then the Cardy condition (3.57) holds:

$$\pi = i_* i^* \tag{4.72}$$

. where $i^* : A \to B$ and $i_* : B \to A$, as in Chapter 2. There we also proved that $i_*$ maps into the center of $A$, completing the argument (note that this does not use the semisimplicity of $A$).

[7] This remarkable fact has not yet been appreciated in the literature.

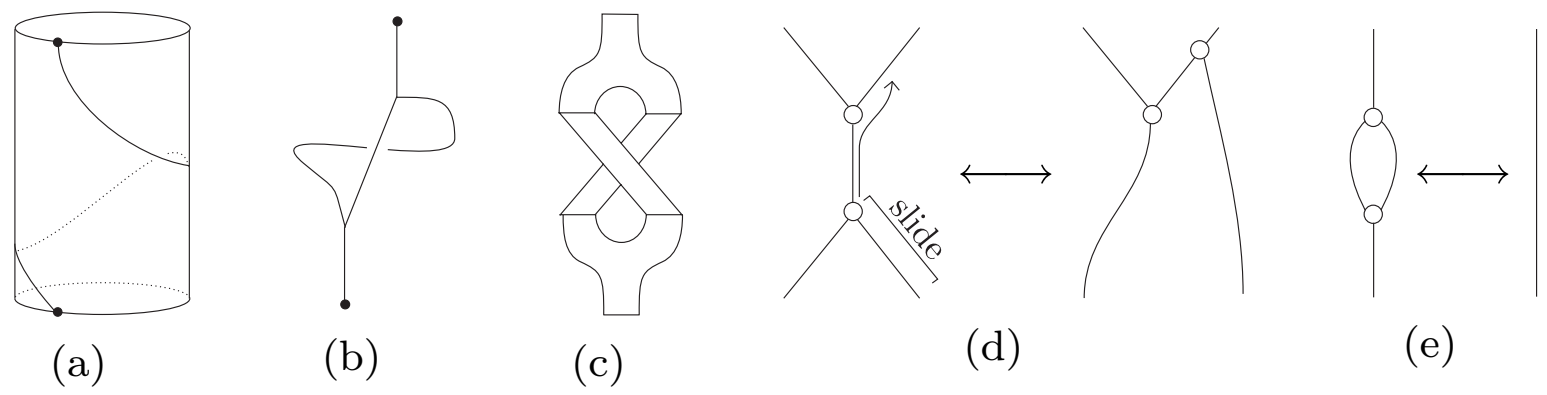

Figure 4.6: (a) A triangulation of the cylinder. (b) The Poincaré dual. (c) Thickening to obtain the map $\pi$. (d) The 2-2 move as a 'sliding' move. (e) The 1-3 move as a bubble move.

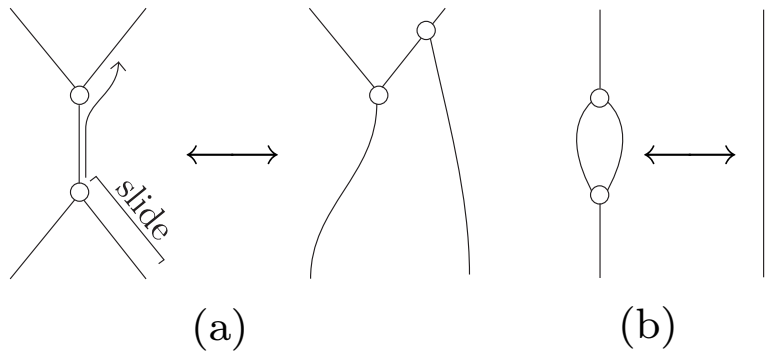

Figure 4.7: The triangulation moves from Fig. 4.5 interpreted as graph moves. (a) The frobenius move allows one to slide 'vertices along the wires'. (b) The bubble move allows one to cancel bubbles.

Now for the graphical argument. We shall operate at the level of the dual diagrams, and interpret these using graph moves in the graphical calculus manner familiar from Chapter 1. The 2-2 Frobenius move in Fig. 4.5c tells us that we can 'slide' interaction vertices (which we will draw as small circles) freely along the lines, as in Fig. 4.7a. Similarly the bubble move in Fig. 4.5d allows us to cancel any bubbles which may appear, as in Fig. 4.7b.
Here is a graphical proof that $\pi(A) \subseteq Z(A)$:

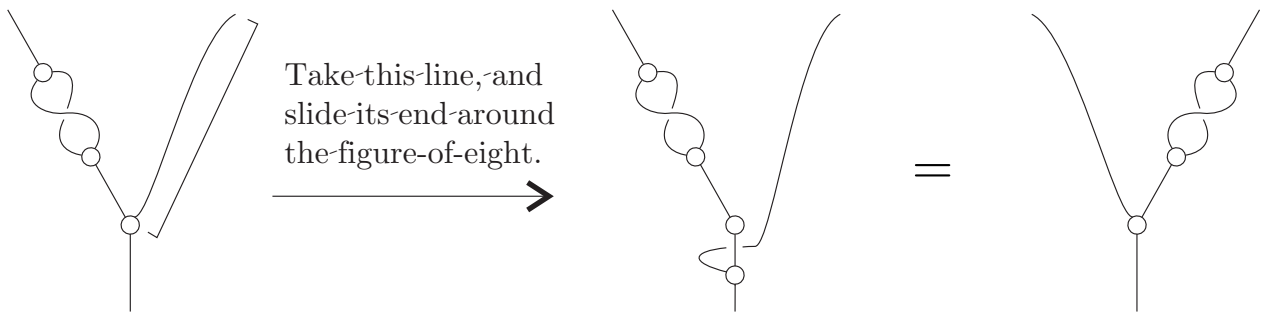

Note that we only used the sliding law (so that it holds for any Frobenius algebra). Here is a proof that if $a \in Z(A)$, then $\pi(a) = a$, i.e. $Z(A) \subseteq \mathrm{Im}(\pi)$:

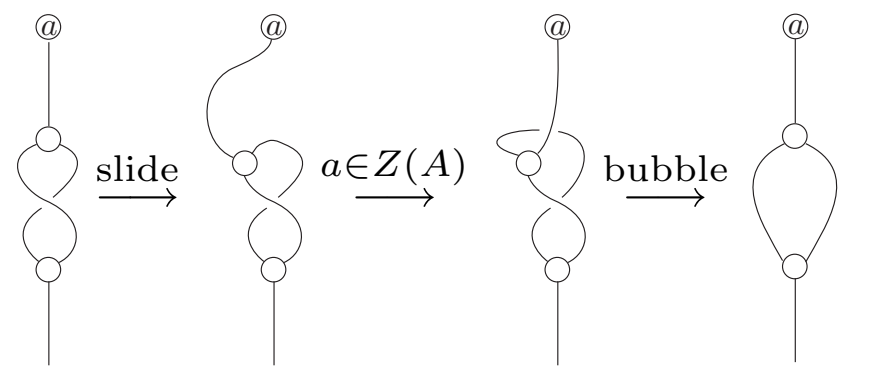

This establishes that $Z(S^1) = Z(A)$ .[8]

### 4.4.3 Computing the topological invariants $Z(M)$

When $A = \mathbb{C}[G]$, then it breaks up as a sum of matrix algebras over the irreducible representations $\rho$ of $G$,

$$\mathbb{C}[G] \simeq \mathrm{End}(V_{\rho_1}) \oplus \mathrm{End}(V_{\rho_2}) \oplus \cdots \oplus \mathrm{End}(V_{\rho_n}). \tag{4.73}$$

The center $Z(A)$ is the span of the identity matrices $1_\rho$ on each representation $\rho$. The trace on $A$ then restricts to $Z(A)$ as

$$\epsilon(1_\rho) = \mathrm{Tr}(L_{1_\rho}) = \dim(\rho)^2. \tag{4.74}$$

[8]Again, do not confuse the two meanings of $Z$.

Note that it is the square of the dimension, since we are taking the trace in the regular representation. In particular, Frobenius trace $\epsilon$ evaluates on the identity $1 = \sum_p 1_\rho$ as

$$\epsilon(1) = \sum_\rho \dim(\rho)^2 = |G|. \tag{4.75}$$

In the original Dijkgraaf-Witten model expressed in terms of principal bundles, we computed the trace on $Z(\mathbb{C}[G])$ (see (4.42)) as

$$\epsilon(g) = \frac{1}{|G|}\delta_{g,1} \quad \Rightarrow \quad \epsilon(1) = \frac{1}{|G|}. \tag{4.76}$$

Comparing (4.75) and (4.76) shows that the trace in the lattice version of the Dijkgraaf-Witten model has been rescaled:

$$\epsilon_{\text{lattice}} = |G|^2 \epsilon_{\text{principal bundles}}. \tag{4.77}$$

Bearing this in mind, let us now proceed using the lattice trace $\epsilon_{\text{lattice}}$. The basis dual to $\{1_\rho\}$ is $\{1^\rho = \frac{1}{\dim(\rho)^2} 1_\rho\}$. The handle element $\omega$ computes as

$$\omega = \sum_\rho 1_\rho 1^\rho = \sum_\rho \frac{1}{\dim(\rho)^2} 1_\rho. \tag{4.78}$$

Hence we conclude from the general results in Chapter 2 that the topological invariants $Z_g$ assigned to closed genus $g$ manifolds $M$ by the lattice version of the Dijkgraaf-Witten model are given by

$$Z_g^{\text{lattice}} = \epsilon_{\text{lattice}}(\omega^g) = \sum_\rho \frac{1}{(\dim \rho)^{2g-1}}. \tag{4.79}$$

If we rescale the trace $\epsilon(\cdot) \to \lambda\epsilon(\cdot)$ then we know from Chapter 2 that the invariants $Z_g$ scale as $Z_g \to \lambda^{1-g} Z_g$. Thus,

$$Z_g^{\text{principal bundles}} = |G|^{2g-2} \sum_\rho \frac{1}{(\dim \rho)^{2g-1}}. \tag{4.80}$$

Comparing this with the topological result (**??**) gives us an interesting formula for counting homomorphisms from $\pi_1(M, x)$ into $G$, expressed in terms of the representations of $G$:

$$\frac{|\text{Hom}(\pi_1(M,x), G)|}{|G|} = |G|^{2g-2} \sum_\rho \frac{1}{(\dim \rho)^{2g-1}}. \tag{4.81}$$

### 4.4.4 Rational 2d topological field theories

I would like to point out here how the concept of *rationality* familiar from conformal field theory has its counterpart even in our simple topological setting. Namely, we know that a semisimple algebra is a direct sum of matrix algebras,

$$A \simeq M_{n_1} \oplus M_{n_2} \oplus \cdots \oplus M_{n_k}, \tag{4.82}$$

and we saw above that the commutative Frobenius algebra $Z(A)$ obtained from the FHK construction is the algebra spanned by the identity matrices $\{1_i\}$, with trace $\epsilon(1_i) = n_i^2$. One sees immediately that *the FHK construction can only produce commutative Frobenius algebras (2d TQFT's) whose ratio of traces in the idempotent basis are rational numbers.* This is true even if we allow ourselves to rescale the trace by a global constant.

## 4.5 Area-dependent theories

As stated in the introduction, we would like to generalize our topological field theories to the case when the cobordisms have extra structure, namely a volume form. Thankfully, any two volume forms on $M$ with the same total volume are related by a diffeomorphism of $M$, so that we are not introducing extraordinary new difficulties into our theory : only an extra number $t$, the total volume, for each cobordism. That is, the effect of this change is simply that the operator $Z_M^{(t)}$ associated to $M$ depends on a number $t > 0$ which is the volume of $M$, and gluing of cobordisms naturally adds volumes, hence we require

$$Z_{M_1 \circ M_2}^{(t_1+t_2)} = Z_{M_1}^{t_1} \circ Z_{M_2}^{t_2}. \tag{4.83}$$

A nasty consequence is that the cobordisms-with-volume no longer form a category, since there is no identity. Nevertheless, they do at least form a 'category without identities', $\mathbf{nCob}^\omega$, where the $\omega$ reminds us that we are dealing with volume dependent theories. Specializing to two dimensions, a topological quantum field theory of this extended sort is a functor $Z : \mathbf{2Cob}^\omega \to \mathbf{Vect}$.

An interesting change is that the functoriality of $Z$ no longer implies that the vector spaces $Z(\Sigma)$ are finite dimensional. Hence, we should be a little careful, and we shall take them to be locally convex and complete topological vector spaces. On each space $Z(\Sigma)$ we have a semigroup of operators (the time translation operators) $\{U_t\}_{t>0}$ coming from the cylinder cobordisms. The argument of Theorem 3.3.5 now proves that $U_t$ is of trace class. One can show that the semigroup $\{U_t\}$ defines a 'rigging' of $Z(\Sigma)$, that is, two complete topological vector spaces $\breve{Z}(\Sigma)$ and $\hat{Z}(\Sigma)$ with maps

$$\breve{Z}(\Sigma) \to Z(\Sigma) \to \hat{Z}(\Sigma) \tag{4.84}$$

which are injective with dense images. We shall not elaborate on this type of analysis here, but rather mention the analogue of the theorem for finite dimensional Frobenius algebras, in the sense that,

**4.5.1 Theorem.** A two dimensional area-dependent theory is the same thing as a commutative topological algebra $A$ with a non-degenerate trace $\beta : A \to C$ and a trace-class approximate unit, i.e. a family $\{\epsilon_t\}_{t>0}$ in $A$ such that (see [79])

(a) $\epsilon_t \to 1$ as $t \to 0$

(b) $\epsilon_s \epsilon_t = \epsilon_{s+t}$, and

(c) multiplication by $\epsilon_t$ is a trace-class operator $A \to A$

## 4.6 Yang-Mills in 2d

### 4.6.1 Formalities

If $G$ is a compact Lie group then it is natural to generalize $C[G]$ to the ring $\mathcal{F}_G$ of smooth $L^2$ functions on $G$, under convolution, i.e.

$$f \star g(x) = \int_G dy f(y) g(xy^{-1}) \tag{4.85}$$

where $dy$ is the Haar measure on $G$. For the trace we take

$$\beta(f) = f(1) \tag{4.86}$$

The ring $\mathcal{F}_G$ does not have a unit, for its natural unit would be the Dirac delta function at the identity element of $G$. The most obvious choice of approximate unit is to take $\epsilon_t$ to be the heat kernel, that is, the fundamental solution of the equation,

$$\frac{\partial \epsilon_t}{\partial t} = \triangle \epsilon_t \tag{4.87}$$

where $\triangle$ is the Laplacian constructed from the canonical metric associated to the Lie group $G$. Thus $\epsilon_t$ is the smooth function to which $\delta$ diffuses in time $t$. We have assembled all the appropriate data for Theorem 4.5.1.

We recall that the $L^2$ functions on $G$ are spanned by the characters on $G$, i.e. the class functions $\chi(gxg^{-1}) = \chi(x)$. Moreover, to each character one can associate an irreducible representation $V$. In this way, exactly as in the finite gauge group setting, one shows that the invariant for a closed surface of genus $g$ and area $t$ is

$$Z(M_g) = \sum_V \frac{e^{-t\lambda_V}}{(\mathrm{dim}V)^{2g-2}} \tag{4.88}$$

where $\lambda_V$ is the eigenvalue of the Casimir operator on $V$. Equation 4.88 is a remarkable formula which deserves closer attention. It tells us that to calculate the topological invariants associated to a gauge group $G$ is an exercise in representation theory and in summation of an infinite series. For $G = SU(2)$, one can complete the calculation to obtain

$$Z(M_g) = 2(2\pi^2)^{1-g}\zeta(2g-2) \tag{4.89}$$

The appearance of the zeta function may seem remarkable to some mathematicians, but physicists are used to this kind of thing when evaluating partition functions! As for the finite gauge group model, what we would *really* like is an interpretation of this quantum field theory in terms of gauge theory and path integrals. It is to this subject that we now turn.

### 4.6.2 The action

It was first shown by Migdal [60] (for a good exposition, see [97]) that in two dimensions, the pure Yang-Mills gauge theory with any gauge group is exactly soluble, and can be described by a triangulation invariant theory. As mentioned in the introduction, one key feature of two dimensions is that there are no propagating degrees of freedom - there are no gluons. This does not make the theory trivial, but does mean that we must investigate the theory on spacetimes of nontrivial topology or with Wilson loops to see degrees of freedom. Since there are so few degrees of freedom one might suspect that there is a very large group of local symmetries. Indeed, this theory has a much larger invariance group than just the group of gauge transformations $GA(P)$. It is invariant under area-preserving diffeomorphisms of the base space, SDiff$(M)$.

One should remember that Yang-Mills theory is not a generally covariant theory since it relies for its formulation on a fixed Riemannian or Lorentzian metric on the spacetime manifold $M$. The gauge fields are connections $A$ on $P$, and the Yang-Mills action is given by

$$S[A] = -\frac{1}{4e^2}\int_M \mathrm{Tr}(F \wedge \star F) \tag{4.90}$$

where $F$ is the curvature of $A$ and Tr is the trace in the fundamental representation of $G$. The action $S[A]$ is gauge-invariant so it can be regarded as a function on the space of connections on $M$ modulo gauge transformations. The group Diff$(M)$ acts on this space, but the action is not diffeomorphism invariant. However, in two dimensions $\star F$ is a lie-algebra valued scalar function $f$ on $M$. Since the metric on $M$ provides us with an area form $\omega$, we may in fact write

$$F = f\omega. \tag{4.91}$$

This simplifies the action considerably to the form

$$S[A] = -\frac{1}{2e^2}\int_M \omega \mathrm{Tr} f^2. \tag{4.92}$$

Since $f$ is a scalar , one sees immediately that the group of symmetries is the group which preserves the volume element $\omega$, that is, the group of area preserving diffeomorphisms of $M$, SDiff$(M)$. Remember that this means we are preserving the total area,

$$t = \int_M \omega. \tag{4.93}$$

Notice that the action is invariant under $e^2 \to \lambda e^2, \omega \to \frac{\omega}{\lambda}$ for any real number $\lambda$. This shows that the partition function

$$Z(M, e^2, t) = \int_M \mathcal{D}A \exp(-S[A]) \tag{4.94}$$

is really only a function of the combination $e^2 t$.

One might wonder if there is not a way to reformulate the theory so as to regain *complete* topological invariance. Looking at the action it seems that this would involve the limit $e \to 0$, although it is not clear how to take this limit. One can perform a standard but clever trick in which this limit is made tractable. Namely, introduce a new theory which has dynamical variables $A$ *and* $\phi$, where $\phi$ is a lie-algebra valued scalar function on $M$, by the action (note how $e^2$ enters in the numerator)

$$S'[A, \phi] = -\frac{e^2}{2} \int_M d\omega \mathrm{Tr}\phi^2 - i \int_M \mathrm{Tr}\phi F \tag{4.95}$$

Note that terms such as $\mathrm{Tr}\phi F$ provide their own area form (since it is a 2-form) and hence do not need an area form to integrate against. The partition function is

$$Z'[M, e^2 t] = \int \mathcal{D}\phi \mathcal{D}A \exp(\frac{e^2}{2} \int \mathrm{Tr}\phi^2) \exp(i \int d\omega \mathrm{Tr}(\phi F) \tag{4.96}$$

Then integrate out $\phi$, using the path integral analogue of the formula,

$$\int_{\infty}^{+\infty} \frac{dx}{\sqrt{2\pi}} \exp(-\frac{e^2}{2} x^2 - ixy) = \exp(-\frac{y^2}{2e^2}). \tag{4.97}$$

Once this has been done, the partition function takes the same form as (4.94). Hence $Z'$ is actually our old theory in disguise, although in a form in which it is much easier to take the limit $e \to 0$. In this limit, we obtain a new theory (still dependent on two fields $A$ and $\phi$),

$$S''[A, \phi] = -i \int_M \mathrm{Tr}(\phi F) \tag{4.98}$$

Just like Chern-Simons theory, this theory does not need an area form for its formulation and is hence (apparently) manifestly topologically invariant. One would assume that, if one has a solution for the full (area-dependent) theory, then one can take the limit $e \to 0$ to obtain the solution for this 'topological theory'. We shall indeed investigate this later in this section.

### 4.6.3 Quantization

We recall briefly the basic notions of quantum field theory in the Schrödinger picture[9]. This view is a direct generalization of quantum mechanics of point particles to the case of fields. Just as quantum mechanics deals with wavefunctions $\psi(x)$ which give the probability for finding the particle at a certain point $x$ in space, quantum field theory deals with wavefunction*als* $\Psi[\phi]$ which give the probability for finding the *field* in a certain configuration $\phi$. A basis for the quantum mechanical Hilbert space are the position eigenkets $|x\rangle$, which represent systems in which the chance of finding the particle is sharply peaked around $x$. Similarly a basis for the quantum field theoretical Hilbert space are the field eigenkets $|\phi(\cdot)\rangle$, which represent systems where the chance of finding the field is sharply peaked around a certain field configuration $\phi(\cdot)$. The position operator in quantum mechanics is represented by multiplication by $x$, while the momentum operator is represented by differentiation with respect to $x$. Similarly, in quantum field theory, the field measurement operator is represented by multiplication by an entire field $\phi(\cdot)$, while the field momentum operator is represented by functional differentiation. This type of language might be the most fundamental way of looking at quantum field theory, but few physicists are familiar with it since in this viewpoint, the particles have apparently disappeared! In fact there is a complete equivalence between this picture and the ordinary picture in terms of eigenkets of the momentum operator, which are called 'particles'. The Schrödinger picture is also called the coherent state formalism in quantum optics, where it is a useful representation. To this list one may add Yang Mills theory in two dimensions, which is most naturally solved in the Schrödinger picture.

In Yang-Mills theory in $d$ dimensions, the physical Hilbert space associated with a $d-1$ manifold $M$ (endowed with a $G$-bundle P) is always $L^2(\mathcal{A}/GA(P))$, where $\mathcal{A}$ is the space of connections on $P$. When spacetime is two dimensional, space is one dimensional and hence all connections are flat. Recall from section 4.24 that this means that the physical Hilbert space is given by

$$L^2\left(\mathrm{Hom}(\pi_1(M), G)/G\right) \tag{4.99}$$

[9]For an excellent reference, see [42]. In some sense, this picture is, philosophically speaking, more fundamental since the *fields* are central, and not the*particles*.

Now if space is one dimensional, it must be a union of circles. The fundamental Hilbert space to determine is the space assigned to one circle, since the rest are obtained by tensor product. Now for a circle $\Sigma$, $\pi_1(\Sigma)$ is freely generated by the path which wraps once around $\Sigma$. Hence a homomorphism in $\mathrm{Hom}(\pi_1(\Sigma), G)$ is simply an arbitrary choice of an element in $G$, so that (4.99) reduces to the space of conjugation invariant (or class invariant) functions on the group $G$. This is a remarkably elegant characterization of the physical Hilbert space, and is worthy of further reflection. Specifically, one should review carefully how the space of wave *functionals* of connections ended up becoming a space of ordinary *functions* on the group $G$. One must not make the mistake however of thinking that connections have been replaced by points in $G$!

We have thus reduced the problem of Yang-Mills in two dimensions to analysis on the compact Lie group $G$. The Hilbert space is the space of $L^2$-class functions on $G$, and the inner product is

$$\langle f|g\rangle = \int_G dx f^*(x) g(x) \tag{4.100}$$

where $dx$ is the standard Haar measure on $G$. The Peter-Weyl theorem tells us that the Hilbert space is spanned by the characters $\chi_R$ of the irreducible representations $R : G \to \mathrm{End}(V)$ of $G$,

$$\chi_R(x) = \mathrm{Tr} R(x) \tag{4.101}$$

One may show that the Hamiltonian is diagonalized in this basis, and that it is the quadratic Casimir $C_2(R)$ of the representation $R$,

$$H = \sum_R \frac{e^2}{2} L C_2(R), \tag{4.102}$$

where $L$ is the length of the circle.

Armed with the Hamiltonian, we can now explicitly solve the theory by writing down the linear maps $Z(M, e^2 t)$ associated to cobordisms $M$ of area $t$. This is because the Hamiltonian gives us the time evolution operator $U(t) = \exp(-Ht)$, and we may compute the amplitudes for surfaces by decomposing them into simpler units. The basic idea [61] is that the cylinder of area $t$ gives rise to an operator

$$\mathrm{Cylinder} = \sum_R |R\rangle\langle R| \exp(-\frac{e^2 t}{2} C_2(R)) \tag{4.103}$$

Similarly, the pair of pants gives rise to an operator

$$\mathrm{Pants} = \sum_R |R\rangle \otimes |R\rangle \otimes |R\rangle \frac{\exp(-\frac{e^2 t}{2} C_2(R))}{\mathrm{dim} R} \tag{4.104}$$

and the cap gives rise to the operator

$$\mathrm{Cap} = \sum_R \mathrm{dim} R \langle R| \exp(-\frac{e^2 t}{2} C_2(R)) \tag{4.105}$$

Using these basic cobordisms, we can calculate the linear operator $Z$ for any cobordism. In particular, we recover the topological invariants (4.88),

$$Z(M_g) = \sum_R \frac{e^{-tC_2(R)}}{(\mathrm{dim} V)^{2g-2}} \tag{4.106}$$

This is a direct generalization of our finite gauge group results from Section 4.4.3.

# Chapter 5

# Three dimensional TQFT's and higher dimensional categories

## 5.1 Modular categories

Recall from Chapter 1 the notion of a semisimple ribbon category $(\mathcal{C}, \oplus, \otimes)$ as an $\mathbb{C}$-linear rigid braided monoidal category with compatible twist, which also has direct sums (i.e. it is an abelian category), such that every object is a finite sum of simple objects. As usual we index the simple objects of $\mathcal{C}$ by the set $I$. Note that there is an involution on $I$ given by taking the dual of a simple object, $i \to i^* \to i^{**} \simeq i$. Recall also the definition of the Grothendieck ring $K(\mathcal{C})$ of $\mathcal{C}$ as the algebra[1] one obtains by decategorifying $\mathcal{C}$, i.e. by taking the algebra spanned by the equivalence classes $\langle i \rangle$ of simple objects $i \in I$ with addition $\langle i \rangle + \langle j \rangle = \langle i \oplus j \rangle$ and multiplication $\langle i \rangle \langle j \rangle = \langle i \otimes j \rangle$. Note that this ring is commutative since $\mathcal{C}$ is braided.

When thinking about semisimple ribbon categories, the best examples to keep in mind are those arising from conformal field theory. Here $\mathcal{C}$ is the category of representations of the chiral algebra, and the simple objects are the primary chiral vertex operators (the space of conformal blocks). The braiding on $\mathcal{C}$ accounts for the presence of braid group statistics in two dimensions. The twist on a simple object $x$ is the fractional part of the conformal weight,

$$\theta_x = e^{-2\pi i \Delta_x} \mathrm{id}_x. \tag{5.1}$$

The dual $i^*$ of an object $i$ corresponds to charge conjugation. The Grothendieck ring of $\mathcal{C}$ is the fusion ring of the conformal field theory.

For $i, j \in I$, define a matrix $s$ by defining its matrix elements $s_{ij} \in \mathbb{C}$ by the following elementary link diagram[2]:

$$s_{ij} =$$

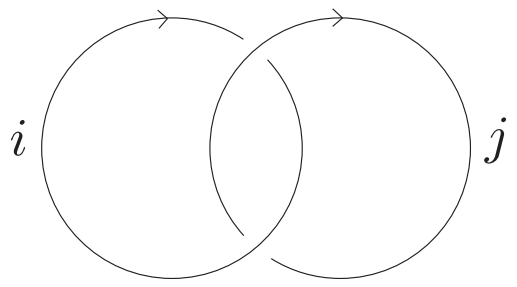

(5.2)

In conformal field theory, $s$ is (up to normalization) the modular $S$-matrix which implements the modular transformation $\tau \to -\frac{1}{\tau}$. We define a *modular category* as a semisimple ribbon category $\mathcal{C}$ such that:

(a) $\mathcal{C}$ has only a finite number of isomorphism classes of simple objects : $|I| < \infty$.

(b) The matrix $s$ is invertible.

It is called a 'modular' category because in this case one can define a projective action of the group $SL_2(\mathbb{Z})$ (which is known as the modular group) on certain objects in $\mathcal{C}$ - a concept well known in conformal field theory.

[1] We shall often abuse language and refer to the 'Grothendieck ring' as an algebra when we really mean the algebra $K(\mathcal{C}) \otimes_{\mathbb{Z}} \mathbb{C}$. This is common practice.

[2] It is hoped that the reader remembers how to evaluate such diagrams. If not, their memory will be refreshed by returning to Chapter 2.

The statement that $s$ is invertible means that the pairing $(\langle i\rangle, \langle j\rangle)$ on the commutative Grothendieck algebra $K(\mathcal{C})$ is nondegenerate. Thus $K(\mathcal{C})$ looks suspiciously like a Frobenius algebra, but it is not quite[3], because the pairing is not associative, $(\langle i\rangle\langle j\rangle, \langle k\rangle) \neq (\langle i\rangle, \langle j\rangle\langle k\rangle)$, as the following topological argument shows:

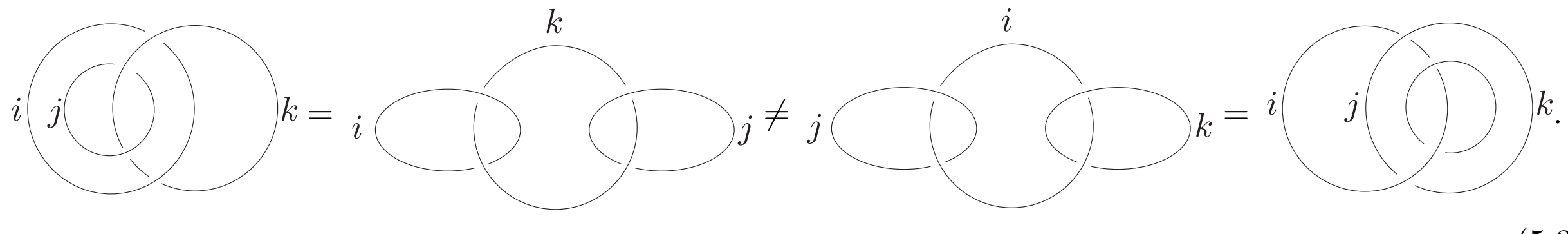

(5.3)

Nevertheless, one of the aims of this chapter is to argue that, just as commutative Frobenius algebras provide the algebraic data for 2d TQFT's, similarly modular categories provide the algebraic data for 3d TQFT's. This last statement is highly non-trivial; for instance there is no theorem stating that '*all* 3d TQFT's come from modular categories'; this is an active area of research. Still, it is helpful to think of a modular category as a categorification of a Frobenius algebra.

The projective action of the modular group $SL_2(\mathbb{Z})$ works as follows. Recall from Chapter 1 the definition of the twist and the quantum dimension of a simple object $i \in I$:

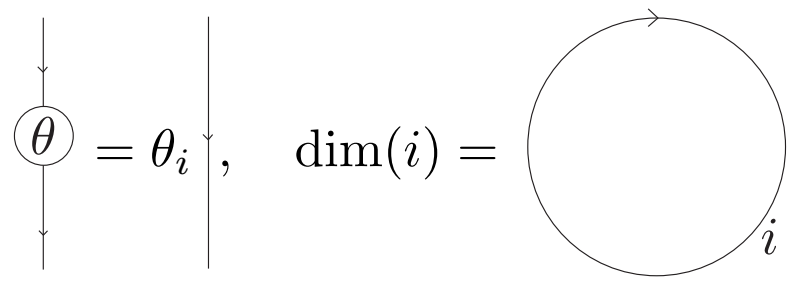

(5.4)

The twist $\theta : i \to i$ must be a scalar multiple of the identity since $i$ is a simple object. Now define the numbers [4]

$$p^{\pm} = \sum_{i \in I} \theta_i^{\pm 1} \dim(i)^2, \quad D = \sqrt{p^+ p^-}, \quad \zeta = (\frac{p^+}{p^-})^{\frac{1}{6}}. \tag{5.5}$$

Also define the matrices

$$S_{ij} = \frac{1}{D} s_{ij} \quad \text{(S matrix)} \tag{5.6}$$
$$T_{ij} = \delta_{ij}\theta_i \quad \text{(matrix of conformal weights),} \tag{5.7}$$
$$C_{ij} = \delta_{ij^*} \quad \text{(charge conjugation matrix)} \tag{5.8}$$

Then one can show [18] from topologically manipulating diagrams that the following relations hold :

$$(ST)^3 = \zeta^3 s^2, \quad S^2 = C, \quad CT = TC, \quad C^2 = 1. \tag{5.9}$$

On the other hand, the modular group $SL_2(\mathbb{Z})$ is generated by the $2 \times 2$ matrices

$$S = \begin{pmatrix} 0 & -1 \\ 1 & 0 \end{pmatrix}, \quad T = \begin{pmatrix} 1 & 1 \\ 0 & 1 \end{pmatrix}, \tag{5.10}$$

satisfying the relations

$$(ST)^3 = s^2, \quad S^4 = 1. \tag{5.11}$$

Comparing (5.9) and (5.11) shows that the matrices $S$ and $T$ constructed from the modular category data give a projective representation of $SL_2(\mathbb{Z})$.

This rather curious and unexpected property can be explained by the fact that a modular category gives rise to a 3d TQFT. We know from Sec. 3.1.3 that any TQFT $Z$ : **nCob** $\to$ **Vect** defines an action of the mapping class groups $\Gamma(\Sigma)$ of closed $(n-1)$ manifolds $\Sigma$ on the vector spaces $Z(\Sigma)$ of the theory. The mapping class group of the torus is equal to $SL_2(\mathbb{Z})$ (see Fig. 5.1), so that the projective[5] representation (5.9) is just one of a whole series of mapping class group representations generated from $\mathcal{C}$.

[3]At least, it does not seem so to me, despite a statement in the literature to the contrary [74]. As I understand it [53], the

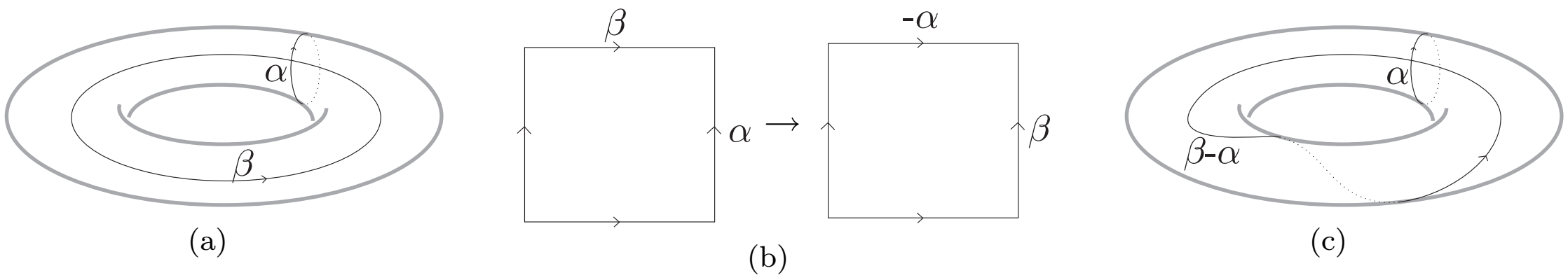

Figure 5.1: (a) The two nontrivial loops on the cylinder. (b) The modular transformation S : $\tau \to -\frac{1}{\tau}$ (c) The Dehn twist $T : (\alpha, \beta) \to (\alpha, \beta - \alpha)$.

### 5.1.1 Conformal field theories living inside modular categories

The main examples of modular categories come from the category of representations of quantum groups at roots of unity and the chiral algebras associated to rational conformal field theories. The representation categories $\mathcal{M}$ of these latter chiral algebras form the *Moore-Seiberg data* [64] of the rational conformal field theory. This data is not enough to specify a *full* conformal field theory with local correlation functions - for this one needs to choose consistent conformal boundary conditions for Riemann surfaces with boundary ('open strings'). An important recent development, which employs in a central manner all the concepts outlined in this thesis, is the work of Jürgen Fuchs *et. al.* [35, 36, 37, 38, 39]. They have shown, remarkably, how the collection of different types of boundary conditions for a given conformal field theory is in one-to-one correspondence with *symmetric, special Frobenius algebras $B$ living inside $\mathcal{M}$.* A 'Frobenius algebra inside $\mathcal{M}$' is nothing but an object $\mathcal{A}$ inside $\mathcal{M}$ armed with a product $A \otimes A \to 1$, coproduct $1 \to A \otimes A$, unit $1 \to A$ and counit $A \to 1$ satisfying the Frobenius conditions (3.21). To be symmetric means that the multiplication maps should behave naturally with respect to duality, and to be special is nothing but our familiar notion of *semisimplicity* (see Fig. 4.6e) inside $\mathcal{M}$:

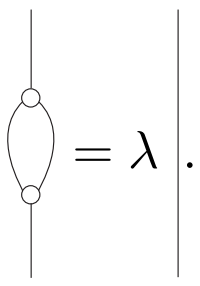

Moreover, for a given type of boundary condition $B \in \mathcal{M}$, the collection of boundary conditions of this type is in one-to-one correspondence with $B$-modules, that is, the category of representations of $B$. This is an eery and astonishing extension of Moore and Segal's classification of boundary conditions for open strings in the topological setting (Theorem 3.4.2), and is certain to generate much interest in the future.

### 5.1.2 Verlinde formula inside modular categories

Let $N_{ij}^k$ be the fusion rules of $\mathcal{C}$, that is the algebra $K(\mathcal{C})$ has multiplication given by

$$x_i x_i = \sum_k N_{ij}^k x_k. \tag{5.12}$$

We can view this as defining $|I|$ multiplication matrices $N_i$ representing left multiplication by $x_i$. Also let $C(I)$ be the algebra of functions on $I$,

$$C(I) = \{f : I \to \mathbb{C}\} \tag{5.13}$$

---

quest is still on to achieve the 'right' definition of a modular category which unifies the approaches of Turaev [89] and Tillman [88].

[4] Interestingly, Vafa [93] has used ideas from conformal field theory to prove that in any modular category, $\zeta$ and $\theta_i$ are roots of unity. That is, the twist always untwists itself after finitely many steps.

[5] One can of course replace $T$ by $\frac{T}{\zeta}$ to get a true representation, but this cannot be done for the actions of the mapping class groups of higher genus surfaces. This is because a modular category gives rise, initially, to an *anomalous* TQFT, although the anomalies are later removed.

Note that $C(I)$ has a basis of delta functions $\{\epsilon_i\}$ which we will normalize as $\epsilon_i(j) = \frac{\delta_{ij}}{\dim(i)}$. These multiply in a diagonal fashion,

$$\epsilon_i \epsilon_j = \delta_{ij} \frac{1}{\dim(i)} \epsilon_i. \tag{5.14}$$

We will now show that the $S$ matrix 'diagonalizes the fusion rules' of $K(\mathcal{C})$ by constructing an isomorphism $\mu : K(\mathcal{C}) \xrightarrow{\sim} F(I)$.

**5.1.1 Theorem.** *Let $\mu : K(\mathcal{C}) \xrightarrow{\sim} F(I)$ be defined by*

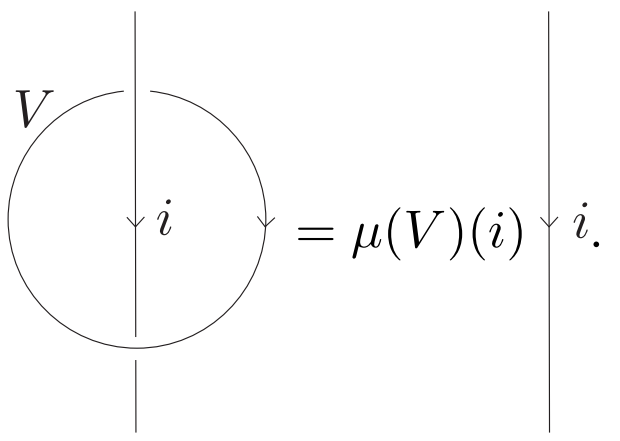

(5.15)

*Then $\mu$ is an isomorphism. Moreover, the change of basis in $K(\mathcal{C})$ from $\{x_i\} \to \{e_i := \mu^{-1}(\epsilon_i)\}$ is given precisely by the $S$ matrix,*

$$x_i = \sum_j S_{ij} e_j. \tag{5.16}$$

*Proof.* Firstly let us calculate the morphism $\alpha : V_i \to V_i$ in (5.15) when $V$ is a simple object $V_i$. Since $\mathrm{End}(V_i) = \mathbb{C}$, $\alpha$ must be equal to $a_{ij}\mathrm{id}_{V_i}$ for some $a_{ij} \in \mathbb{C}$. Then we proceed by closing the loops on both sides:

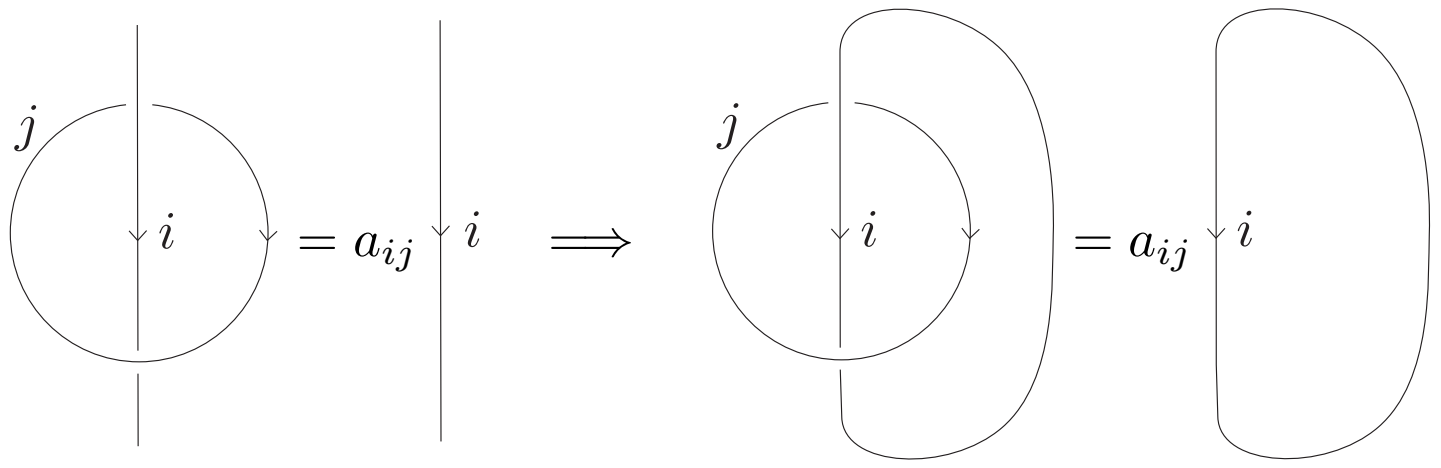

The left-hand side of the latter equality is equal to $s_{ij}$ while the right-hand side is equal to $a_{ij}\dim(V_i)$. Thus $a_{ij} = \mu(V_j)(i) = \frac{1}{\dim(V_i)}\tilde{s}_{ij}$. Comparing with the normalization of the $\epsilon_i$ (below (5.13)) and (5.6) this establishes that it is the $S$ matrix which performs the change of basis,

$$\mu(x_i) = \sum_i S_{ij} \epsilon_j. \tag{5.17}$$

The $S$ matrix is invertible so this is indeed a change of basis. It only remains to show that $\mu : K(\mathcal{C}) \to C(I)$ is an algebra homomorphism, i.e. that $\mu(\langle V_j \otimes V_k \rangle) = \mu(\langle V_j \rangle)\mu(\langle V_k \rangle)$. But this follows from the fact that the loops can be slid under each other:

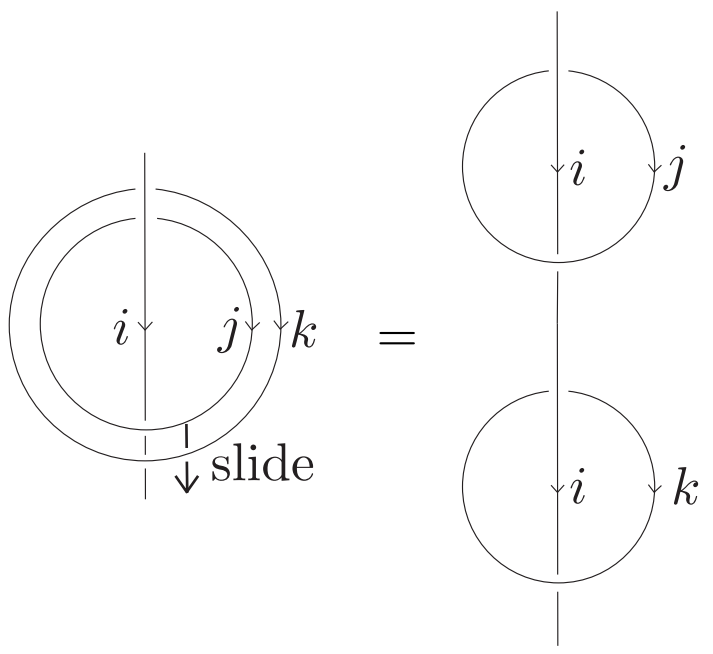

(5.18)

□

From (5.17) we see that the $S$ matrix diagonalizes the multiplication operators $N_i$,

$$SN_iS^{-1} = D_i \tag{5.19}$$

where $D_i$ is the diagonal matrix $(D_i)_j^k = \delta_{jk}\frac{1}{\dim(j)}S_{ij}$. Whenever we have an algebra whose multiplication diagonalizes according to (5.19), we can recover the original algebra coefficients $N_{ij}^k$ from the matrices $D_i$ since $N_i = S^{-1}D_iS$, or

$$N_{ij}^k = \sum_r \frac{S_{ir}S_{jr}S_{k^*r}}{S_{0r}}. \tag{5.20}$$

This is known as the *Verlinde formula*, and can be viewed as a kind of Fourier transformation on the algebra. It is quite surprising because it implies that the right hand side is a nonnegative integer, something which was not obvious from the definition of the $S$ matrix. An interesting example taken from Ferguson [28] is the modular category associated to the quantum group $U_q(\mathfrak{sl}(2))$ where $q$ is an $n$th root of unity. This category has the same fusion rules $N_{ij}^k$ as the ordinary representations of $SU(2)$,

$$V_i \otimes V_j = V_{|i-j|} \oplus V_{|i-j|+1} \oplus \cdots \oplus V_{i+j}. \tag{5.21}$$

These rules are independent of the level $n$. However, the $S$ matrix in this category can be computed to be:,

$$S_{ij} = \sqrt{\frac{2}{n+2}} \sin\left(\frac{\pi(i+1)(j+1)}{l+2}\right). \tag{5.22}$$

It is not at all obvious that we can recover the $N_{ij}^k$ in (5.21) from the Verlinde formula (5.20).

## 5.2 3d TQFT's from modular categories

### 5.2.1 3d manifolds from links

We have spent considerable time developing the theory of semisimple ribbon categories, of which modular categories are an example. The data which lurks in such categories is both *topological* (through the braiding and the twist) and *algebraic* (through the interaction of $\otimes$ with $\oplus$). The topological information is inherently three dimensional, since braids and twists can always be straightened in higher dimensions, and are trivial in lower dimensions.

Thus we have seen in Chapter 1 how ribbon categories give rise to invariants of oriented knots and links, by simply labeling the components with objects from the category. Recall that these links should really be viewed as ribbons, which is the same as to say that all knots must come equipped with a *framing*, a nonzero transverse vector field along the knot (this traces out a ribbon, see Fig. 5.2). The correspondence between framed knots and three dimensional topology is very strong. The fundamental fact underlying the entire program of obtaining 3d TQFT's from modular categories is the following classical theorem of Lickorish and Wallace [55, 56]:

**5.2.1 Theorem of Lickorish and Wallace** . *Every closed, orientable connected 3d manifold $M$ can be obtained by surgery of $S^3$ along a framed link $L \subset S^3$.*

The surgery operation is defined as follows. Let $T_i$ be a small tubular neighborhood of the $i$th component of $L$. Each $T_i$ is a solid torus having boundary $\partial T_i$ equal to an ordinary torus. $\partial T_i$ has 1-cycles $\alpha_i$ and $\beta_i$ as shown in Fig. 5.2c, where the framing of the link determines a canonical way for $\beta_i$ to wrap around $\partial T_i$. For each $i$, we remove the solid torus $T_i$ and glue back a standard torus, twisted in such a way that $\alpha_i \to -\beta$ and $\beta_i \to \alpha$. This gives us a new 3d manifold $M_L$ obtained from surgery along $L$, and the theorem states that all 3d manifolds can be obtained in this way. This is quite remarkable, but is best understood as a statement that *links have enough complexity to generate all 3d manifolds.* Indeed, many non-isomorphic links $L$ generate the same 3d manifold $M_L$. Thankfully, it is known precisely when this occurs, due to a theorem of Kirby [48].

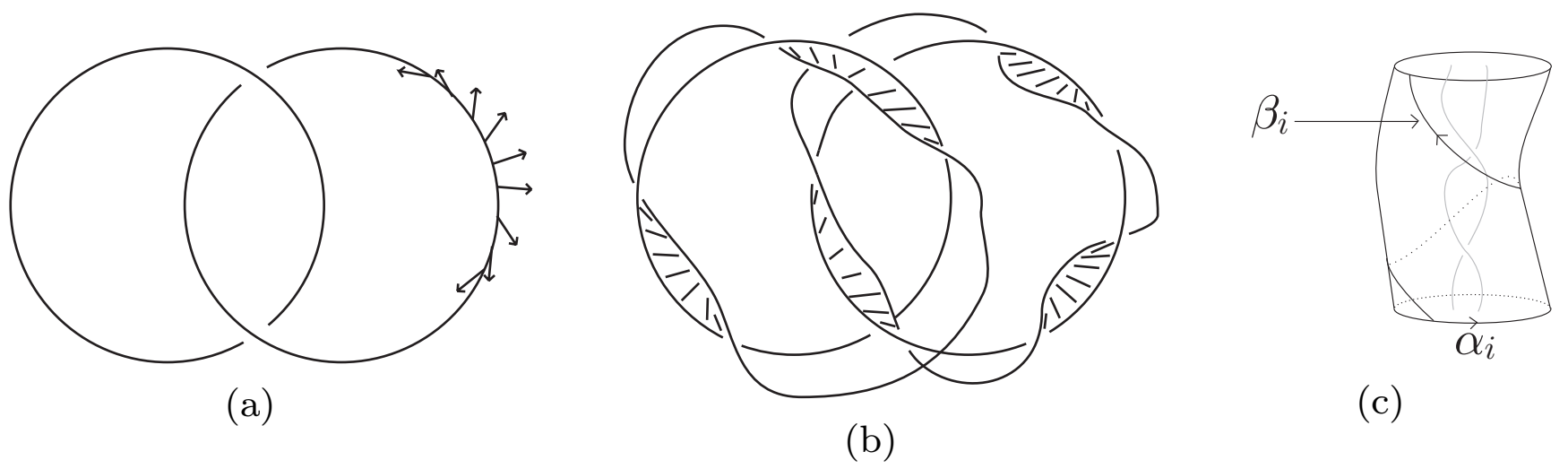

Figure 5.2: (a) A framed link. (b) The equivalent ribbon. (c) The ribbon link $L_i$ and its tubular neighborhood, which defines the one cycles $\alpha_i$ and $\beta_i$.

**5.2.2 Kirby Calculus** . *$M_L \simeq M_{L'}$ if and only if $L$ and $L'$ differ by a finite sequence of the following two elementary moves:*

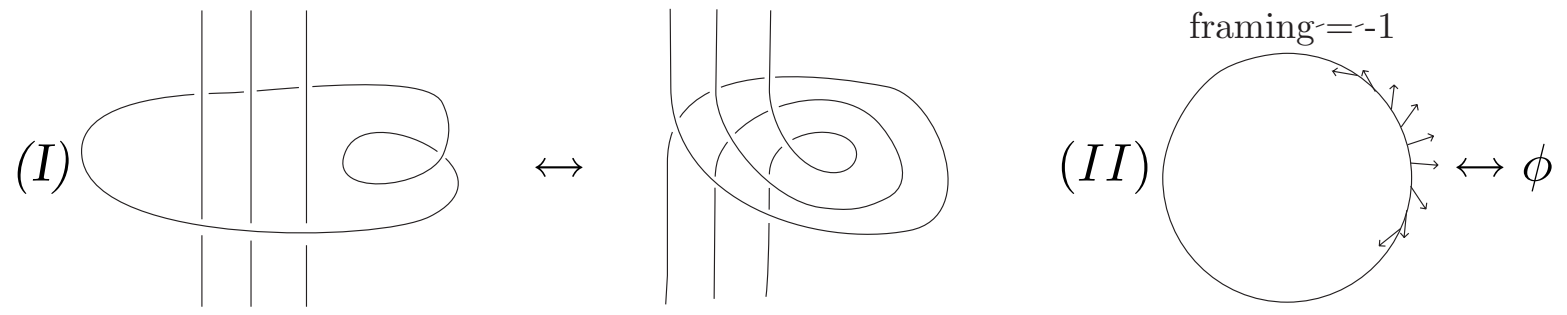

Thus we have a 1-1 correspondence between

$$\left\{ \begin{array}{c} \text{Invariants of framed links which do not change} \\ \text{under the Kirby moves} \end{array} \right\} \overset{\text{1-1}}{\longleftrightarrow} \{\text{Invariants of 3-manifolds}\}.$$

In particular, Reshetikhin and Turaev [89, 90] showed how to define precisely such an invariant $Z(M)$ from the data in a modular category. Moreover these invariants $Z(M)$ for closed manifolds $M$ can be extended into a full blown topological field theory $Z$ : **3Cob** $\to$ **Vect**, which we shall presently describe. The Reshetikhin-Turaev invariant is

$$Z(M_L) = D^{-|L|-1} F(L) (\frac{p^+}{p^-})^{\sigma(L)/2}. \tag{5.23}$$

That is, present $M$ as surgery on a framed link $L$. Then calculate the above number from the link, where $D$ and $p^{\pm}$ are defined in (5.5), $\sigma(L)$ is the *writhe* of the link (the no. of full twists), and $F(L)$ is a weighted summation over the link invariants associated to all possible labelings $i \in I$ of the strands:

$F(L) = F\left(\right.$ 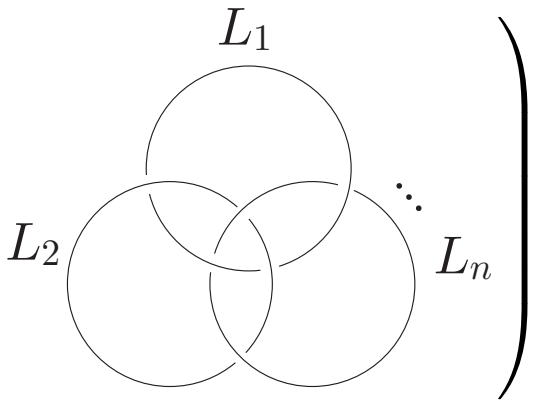

$\left.\right) = \sum_{i_1, i_2, \cdots, i_n \in |I|} \dim(i_1)\dim(i_1)\cdots\dim(i_n)$ 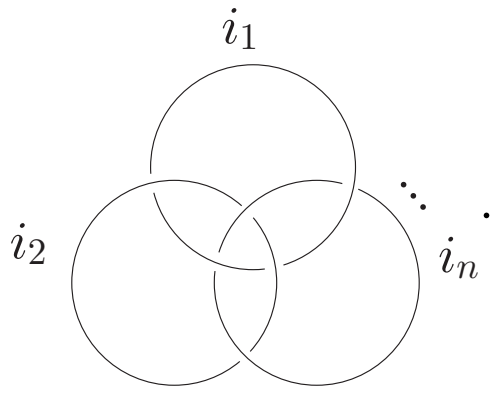

. (5.24)

Turaev and Reshetikhin showed that (5.23) is an invariant of 3d manifolds by showing that it is invariant under the Kirby moves. To the present author's knowledge, there is still no known way to obtain this invariant from homological or homotopical methods[6]. Formula (5.23) can be viewed as a *(finite) path integral over the simple objects in $\mathcal{C}$*. It behooves us to pause and reflect on its comparison with Witten's partition function from Chern-Simons theory (1.4),

$$Z(M) = \int_{\mathcal{A}} \mathcal{D}A \exp \frac{ik}{4\pi} \int_M \mathrm{Tr}(A \wedge dA + \frac{2}{3} A \wedge A \wedge A). \tag{5.25}$$

[6] At least, not in a way which works for all modular categories. We shall presently construct the modular category associated to the finite group model. We know from Chapter 3 that in this case, $Z(M) = \frac{1}{|G|}|\mathrm{Hom}(\pi_1(M,x),G)|$ which is clearly expressed in terms of homotopy information.

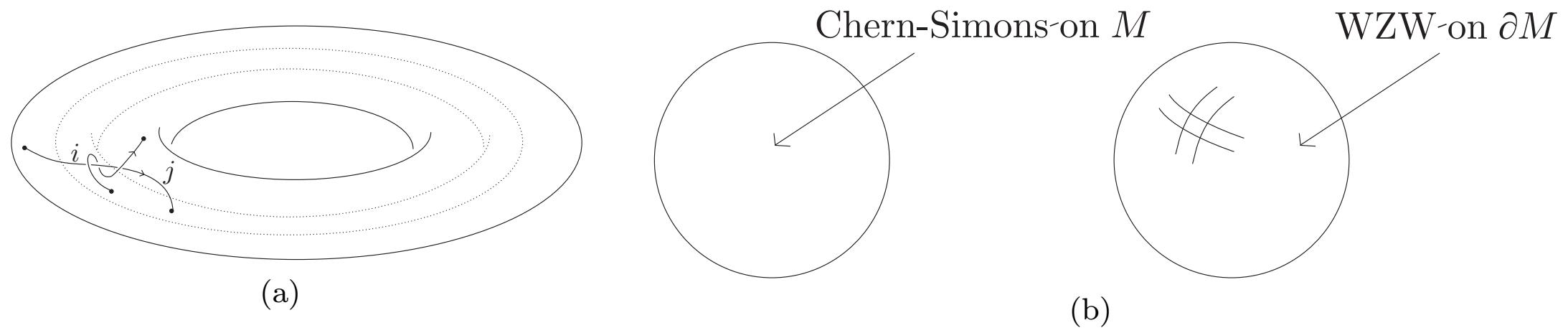

Figure 5.3: (a) A $\mathcal{C}$-marked cobordism inside the solid torus having boundary two torii. The objects $V_i$ and vectors on the points on the boundary reflect the links 'poking' through the surface. (b) A Chern-Simons theory on a 3-manifold $M$ (eg. the unit ball) restricts to a WZW theory on $\partial M$ (eg. $S^2$).

We summarize the salient features of these two approaches in the following table:

| Witten | Reshetikhin-Turaev |
|---|---|
| Manifold $M$ | Link $L$ representing $M$ |
| Compact Lie group and an integer level eg. $SU(2), k$ | Modular category eg. $U_q(\mathfrak{su}_2)$ where $q^k = 1$ |
| Infinite sum over connections $A$ | Finite sum over simple objects $V_i$ |
| Weighting each by the action $S[A]$ | Weighting each by the factors $\dim(V_i)$ F(L) |

A rigorous understanding of the relationship between these two approaches remains the central problem in the field. We shall resume this discussion in Sec. 5.4.

### 5.2.2 Constructing $Z$ : 3Cob $\to$ Vect

We shall construct a 3d TQFT $Z$ by presenting everything in sight in terms of links, and then use the modular category $\mathcal{C}$ to get numbers and vectors from these links. We shall only construct an ordinary 3d TQFT, that is, 'simply' a functor $Z$ : **3Cob** $\to$ **Vect**. This simplifies much of the analysis found in [89, 18]. One can (and should) go much further and construct an $\mathcal{C}$*-extended* 3d TQFT, which in this context means one which assigns vector spaces to $\mathcal{C}$-marked surfaces $\Sigma$ and linear operators to $\mathcal{C}$-marked cobordisms $M : \Sigma_1 \to \Sigma_2$. A $\mathcal{C}$-marked surface is one which has a finite number of marked points $p_1, \ldots, p_k$, each of which is given the data of a nonzero tangent vector $v_i$, and object $W_i \in \mathcal{C}$ and a sign $\epsilon_i = \pm$. These marked points are referring to the intersection with the boundary $\partial M$ of some link inside $M$, each of whose strands is labeled with an object of $\mathcal{C}$ (see Fig. 5.3a).

The physical relationship between ordinary TQFT's and $\mathcal{C}$-extended TQFT's is that the former describes only the dynamics of the fields, while the latter describes observables as well. To see this, consider Chern-Simons theory. Witten showed that Chern-Simons on a three-manifold $M$ with boundary $\partial M$ is essentially characterized by a corresponding two dimensional Wess-Zumino conformal field theory living on $\partial M$ (see Fig. 5.3b). In particular, Chern-Simons expectation values for Wilson lines ending at $k$ points in the boundary is described by the associated Wess-Zumino theory on the boundary with $k$ punctures carrying the representations of the free Wilson lines . We shall study extended TQFT's from a different viewpoint in Sec. ***.

We now define the ordinary TQFT $Z$ associated with $\mathcal{C}$, in the Atiyah picture. Firstly define the fundamental object (familiar from representation theory for compact groups)

$$H = \sum_i V_i \otimes V_i^*. \tag{5.26}$$

The only vector spaces around are the Hom spaces in $\mathcal{C}$, so the entire approach is based on them. We use the notation

$$\langle V \rangle = \mathrm{Hom}(\mathbb{1}, V). \tag{5.27}$$

For a genus $g$ closed 2d surface $\Sigma$, set

$$Z(\Sigma_g) = \langle H^g \rangle. \tag{5.28}$$

This means that $Z(-\Sigma_g) = Z(\Sigma_g)$, but this has always been the case in our constructions (see (4.51)). As always, what we do is define a nondegenerate pairing $(Z(\Sigma_g) \otimes Z(-\Sigma_g) \to \mathbb{C}$, which gives us an identification $Z(-\Sigma) \simeq Z(\Sigma)^*$. For $\phi : 1 \to H^g, \psi : 1 \to H^g$, define $(\phi, \psi) \in \mathbb{C}$ by

$$(\phi, \psi) = \frac{1}{D^g}\left(1 \xrightarrow{\phi\otimes\psi} H^g \otimes H^g \xrightarrow{(\eta\otimes\mathrm{id})\otimes\cdots\otimes(\eta\otimes\mathrm{id})} H^g \otimes H^g \xrightarrow{e_H\otimes\cdots\otimes e_H}\right). \tag{5.29}$$

Here $\eta|_{V_i\otimes V_i^*} = \dim(V_i)^{-1}\mathrm{id}$, and $e_H$ is the evaluation map $H \otimes H \to 1$ which comes from the ordinary evaluation map $e : H^* \otimes H \to 1$ and the canonical isomorphism $H \simeq H^*$. We will not comment on this pairing, other than to say that a proof that it is nondegenerate and symmetric can be found in Turaev's book [89].

Now we need a way to present 3d manifolds with boundary, using surgery on links. Let $L_g$ denote an uncolored coupon with $g$ uncolored 'handle' links attached:

$L_g =$ 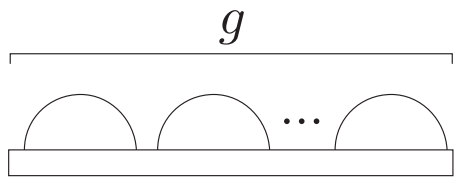

Observe that if $M$ is a closed manifold with $L_g$ living inside it, then removing a tubular neighborhood $T$ of $L_g$ from $M$ creates a manifold $M \setminus T$ having boundary a genus $g$ torus.

This gives us a way to create any arbitrary manifold with boundary by first performing surgery along a link in $S^3$ to get a closed manifold $M$, and then removing an $L_g$ from $M$. We define a *special link* $X = (L, L')$ as a diagram in $S^3$ composed of coupon type links $L_g^a, a = 1 \dots p$ and ordinary links $L'$, which are not allowed to be attached to the coupons. For example:

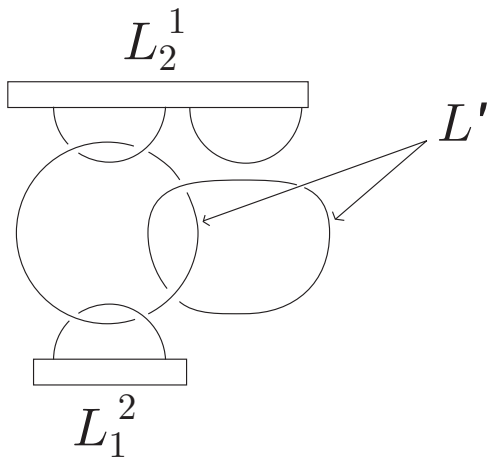

From a special link $X = (L, L')$ we can form a manifold with boundary $M_X$ as explained above by setting

$$M_X = M_{(L,L')} = M_{L'} \setminus \cup_a M_{L_g^a} \tag{5.30}$$

$M_X$ is a 3d manifold with $p$ boundary components of genus $g_a$, $a = 1 \dots p$. We state without proof (see [89]) that all manifolds with boundary can be obtained in this way.

We now define

$$Z(M_X) \in Z(\partial M_X)^* = \otimes_{a=1}^p \langle H^{g_a}\rangle^* \tag{5.31}$$

and then use the pairing (5.29) to 'lower the indices' as usual so that $Z(M_X) \in Z(\partial M_X)$ (with an abuse of notation). To define $Z(M_X)$, view a collection $(\phi_1, \dots, \phi_p) \in \langle H^{g_1}\rangle \times \cdots \times \langle H^{g_p}\rangle$ as defining a colouring $X_{\phi_1\cdots\phi_p}$ of the coupons in $X$, and then sum over all labelings of the remaining strands (as in (5.24), except that now we have labeled coupons):

$$Z(M_X)(\phi_1 \otimes \cdots \phi_p) = F(X_{\phi_1\cdots\phi_p}). \tag{5.32}$$

For example, the manifold $M_X$ defined by the special link above has a boundary consisting of a torus and a genus 2 torus, and $Z(M_X) \in Z(\partial M_X)$ is defined by

$Z(M_X)(\phi_1 \otimes \phi_2) = F$ 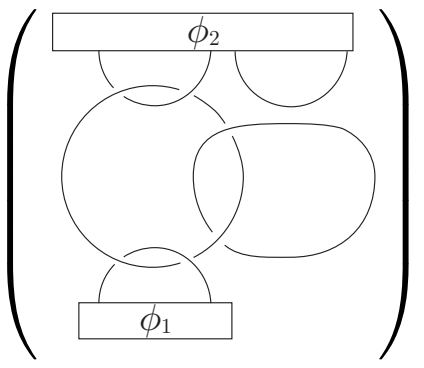 (5.33)

This completes the description of the TQFT $Z : \mathbf{3Cob} \to \mathbf{Vect}$. A proof that it behaves correctly under cutting and pasting, as well as all the other properties, can be found in [89].

We end this section by considering the dimensions of the Hilbert spaces $Z(\Sigma_g)$ that the theory assigns to genus $g$ surfaces $\Sigma_g$. From (5.28) we have that

$$Z(\Sigma_g) = \left\langle (\sum_i V_i \otimes V_i^*)^g \right\rangle = \left\langle (\sum_i \oplus_j N_{ii^*}^j V_j)^g \right\rangle \tag{5.34}$$

By repeated use of the Verlinde formula (5.20), this can be expressed directly in terms of the *quantum* dimensions of the simple objects,

$$\dim Z(\Sigma_g) = D^{2g-2} \sum_i \frac{1}{\dim(V_i)^{2g-2}}. \tag{5.35}$$

This result is also commonly known as the Verlinde formula. Alternatively, it can be computed from first principles (as in all TQFT's) by evaluating the topological invariant $Z(\Sigma_g \times S^1)$.

## 5.3 Quantum doubles of monoidal categories

In this section we outline an elegant method of creating ribbon categories, starting from a general monoidal category.

### 5.3.1 Finding braidings

Suppose we are given a monoid $X$ (a set with a multiplication operation and an identity element). Then we can construct another monoid $Z(X)$, the center of $X$:

$$Z(X) = \{x \in X : xa = ax \text{ for all } ain X\}. \tag{5.36}$$

$Z(X)$ is a commutative monoid. Since a monoidal category is the categorification of a monoid, we can ask if there is an analogous construction on a monoidal category $\mathcal{C}$. The output should be a braided monoidal category $Z(\mathcal{C})$, since a braiding is the categorification of commutativity.

The answer is affirmative, and we follow here the exposition given by Müger [66] as well as Street [85]. Let $\mathcal{C}$ be a strict monoidal category. For a fixed object $x$, we are going to 'seek out' those collections of morphisms $\{x \otimes y \to y \otimes x : y \in \mathcal{C}\}$ with other objects in $\mathcal{C}$ which might qualify as a braiding from $x$'s perspective. So define a *half braiding* for an object $x$ as a function $f_x$ which for every $y \in \mathcal{C}$ associates a map $f_x(y) \in \mathrm{Hom}_{\mathcal{C}}(x \otimes y, y \otimes x)$,

$f_x(y) =$ 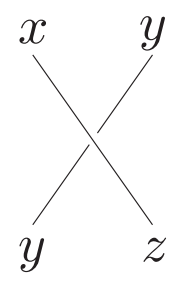

satisfying:

(naturality)

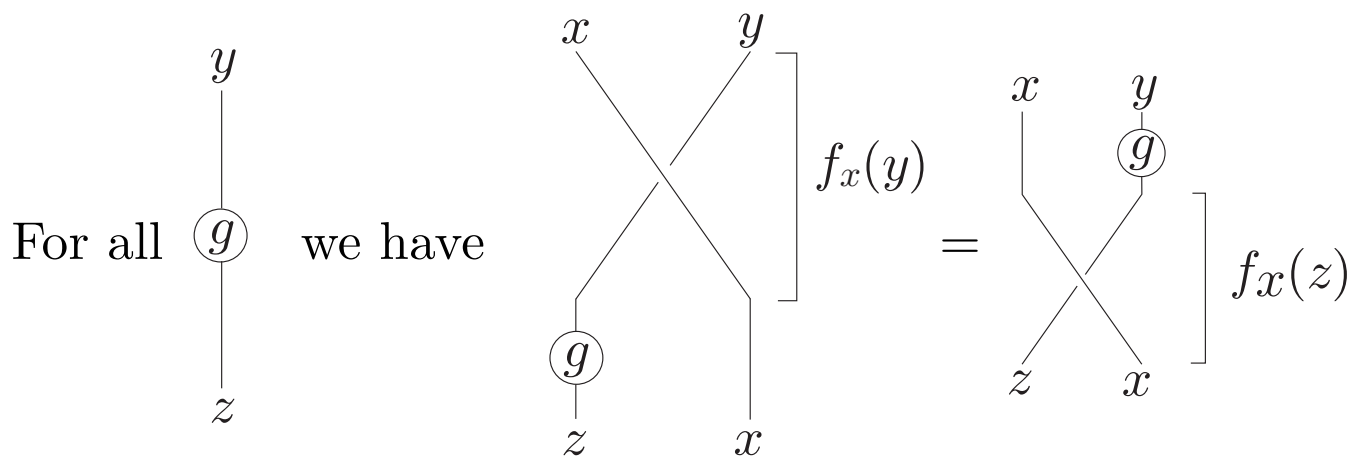

(bilinear)

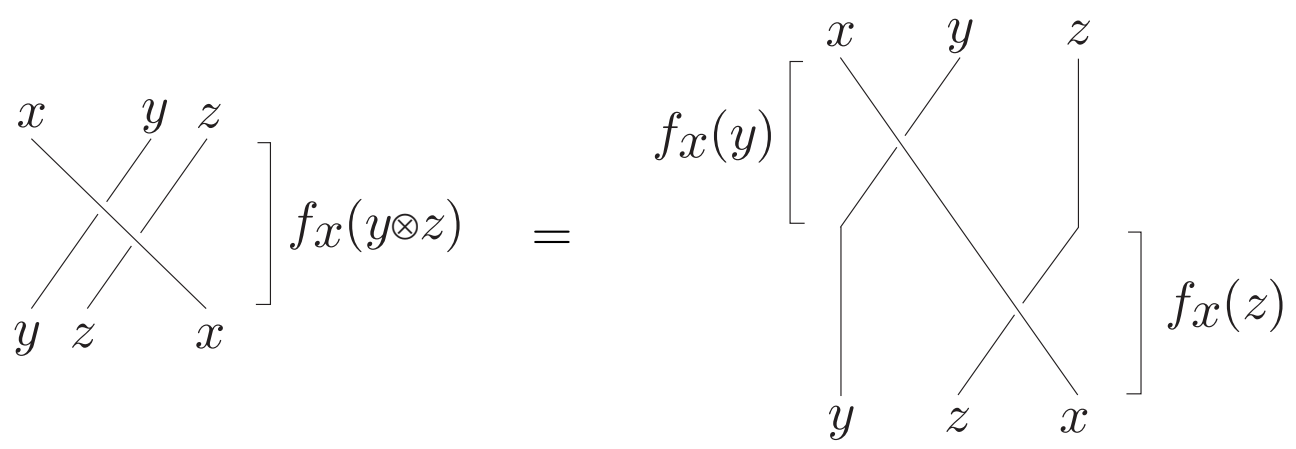

We define the objects of $Z(\mathcal{C})$ by

$$\mathrm{Obj}(Z(\mathcal{C})) = \{(x, f_x) : x \in \mathrm{Obj}(\mathcal{C}), f_x \text{ is a half braiding for } x\}. \tag{5.37}$$

A morphism $s : (x, f_x) \to (y, f_y)$ is defined to be a morphism $s \in \mathrm{Hom}(x, y)$,

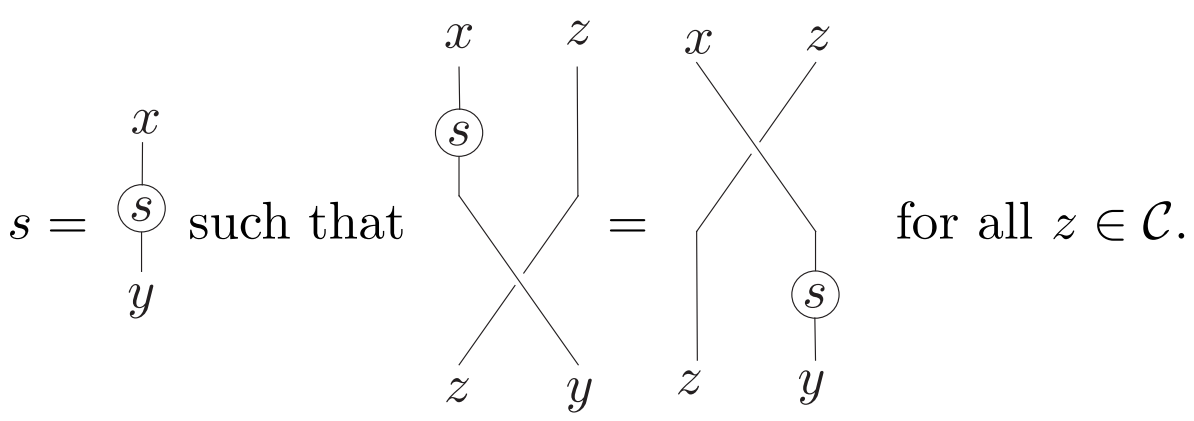

It is then easy to see that $Z(\mathcal{C})$ forms a category with composition induced from $\mathcal{C}$. Moreover it is a braided monoidal category, with:

unit $1_{Z(\mathcal{C})} = (1_{\mathcal{C}}, f_1)$ (where $f_1(x) = id_x$ for all $x$).

tensor product $(x, f_x) \otimes (y, f_y) = (x \otimes_{\mathcal{C}} y, f_{x\otimes_{\mathcal{C}} y})$ where braiding $\sigma_{((x,f_x),(y,f_y))} = f_x(y) : x \otimes y \to y \otimes x$.

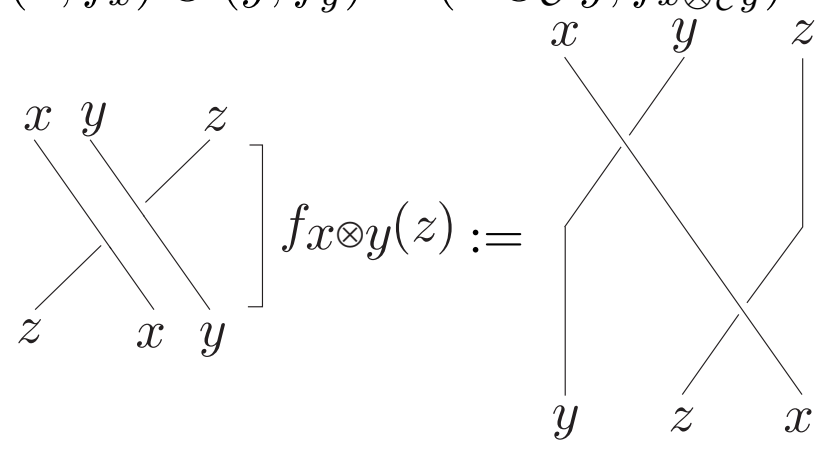

This completes the construction of $Z(\mathcal{C})$. Note that if $\mathcal{C}$ had a braiding to begin with, then $Z(\mathcal{C})$ does not in general reproduce $\mathcal{C}$, since there may be 'other ways' to braid $\mathcal{C}$, which $Z(\mathcal{C})$ takes into account. Roughly speaking, there are more objects in $Z(\mathcal{C})$, but there are less morphisms, since the 'non-covariant' morphisms have been left out.

### 5.3.2 Finding twists

We can use the same idea to find twist maps $\theta_x : x \to x$. Namely, if $\mathcal{C}$ is a braided monoidal category, then define $\mathcal{C}^Z$ as the category of automorphisms of $\mathcal{C}$. The objects are pairs $(x, \theta_x)$ where $\theta_x : x \to x$ is an automorphism of $x$ in $\mathcal{C}$. A morphism $f : (x, \theta_x) \to (y, \theta_y)$ is a morphism $f : x \to y$ in $\mathcal{C}$ such that

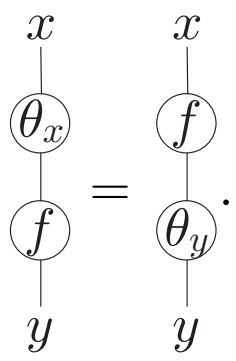

We could define the tensor product on $\mathcal{C}^Z$ by $(x, \theta_x) \otimes (y, \theta_y) = (x, \theta_x \otimes \theta_y)$, but if $\mathcal{C}$ is braided it is better to define it as

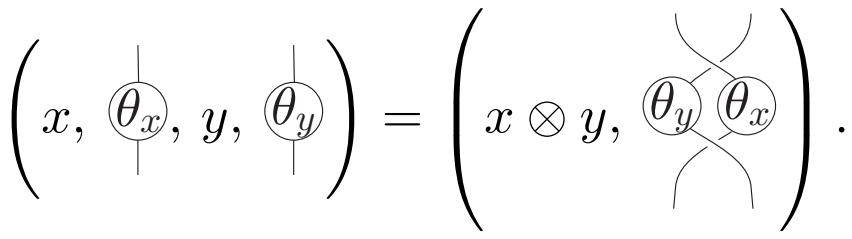

Compare with (2.36). Moreover, the braiding on $\mathcal{C}$ ascends to a braiding on $\mathcal{C}^Z$ in the obvious way. Finally, we can set the twist map $\theta_{(x,\theta_x)}$ of an object $(x, \theta_x)$ to be equal to $\theta_x$ itself. It is not hard to see that $\mathcal{C}^Z$ is now a balanced (braided with a compatible twist) monoidal category.

### 5.3.3 Restricting to a ribbon category

The last ingredient we need to construct a ribbon category is *duality*. Dual objects cannot be created using the above technique, since duality requires for each object $x$, the existence of a dual object $x^*$ *and* the pairing and copairing maps, $e_x : x^* \otimes x \to 1$ and $i_x : 1 \to x \otimes x^*$. Indeed, suppose we started with a general monoidal category $\mathcal{C}$ and then formed $Z(\mathcal{C})^Z$. One can then show that an object $(x, f_x, \theta_x)$ of $Z(\mathcal{C})^Z$ has a dual (inside $Z(\mathcal{C})^Z$) if and only if $x$ has a dual inside $\mathcal{C}$, and such that an additional constraint involving the compatibility of the braiding and the duality is satisfied. Thus we are forced to use *restriction* as the only way to generate duality, and while we're about it, we might as well restrict ourselves to ribbon subcategories (where the duality is compatible with the twist).

Thus suppose $\mathcal{C}$ is a balanced monoidal category, and let $\mathcal{N}(\mathcal{C})$ be the full subcategory of $\mathcal{C}$ consisting of those objects $x$ which have a dual $x^*$ satisfying

$$\theta_{x^*} = (\theta_x)^* \tag{5.38}$$

One can show that $\mathcal{N}(\mathcal{C})$ is in fact closed under the balanced monoidal structure, so that it is a ribbon category.

We have now pieced together all the steps we need. Starting with a general monoidal category $\mathcal{C}$, we can generate a ribbon category $\mathcal{D}(\mathcal{C})$ called the *quantum double* of $\mathcal{C}$, defined as

$$\mathcal{D}(\mathcal{C}) = \mathcal{N}(Z(\mathcal{C})^Z). \tag{5.39}$$

If $\mathcal{C}$ was $\mathbb{C}$-linear and semisimple, then so is $\mathcal{D}(\mathcal{C})$.

## 5.4 Modular categories from path integrals

Daniel Freed has offered a deep and profound explanation for the appearance of modular categories in 3d TQFT's (see [30, 31, 32] and especially [33] for a good overview). Freed's approach is nothing short of a revolutionary new view of a $d$ dimensional quantum field theory as *a machine which produces a tower of algebraic structures for each dimension $d-k$, each obtained by a path integral over fields taking values in its predecessor.*

This philosophy unifies, in principle, the following three incarnations of a 3d TQFT:

(a) A topologically invariant action, $A \to \exp iS[A] \in U(1)$, where $A$ is a field on a 3d manifold $X$ (eg. Chern-Simons action).

(b) An explicit functor $Z : \mathbf{3Cob} \to \mathbf{Vect}$ (eg. the principal bundle defn. of the finite group model)

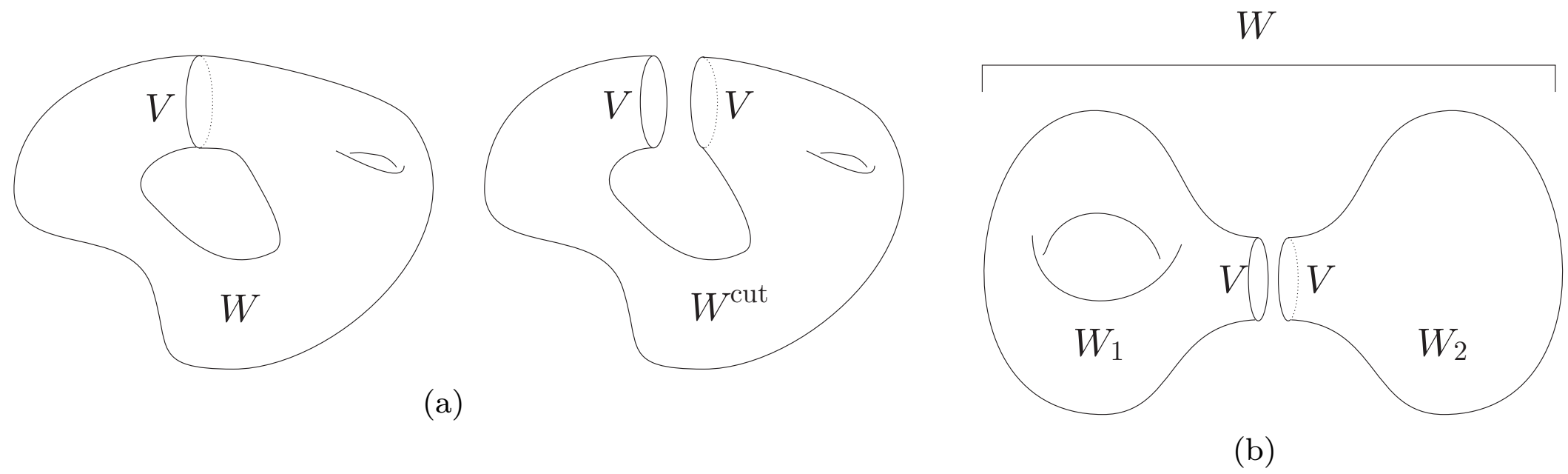

Figure 5.4: (a) Cutting a spacetime $W$ along $V$ to obtain $W^{\text{cut}}$ (compare with Sec. 3.1.1). (b) Gluing together two disconnected manifolds $W_1$ and $W_2$ along $V$ to form $W$.

(c) A modular category (eg. the Reshetikhin-Turaev viewpoint).

Freed's method begins at the uppermost level (an action), and proceeds downwards. This is certainly the most natural from a physics viewpoint, and it is also the most elegant (compare the definition of a modular category with the Chern-Simons action). Lie groups and three dimensional manifolds enter in an intrinsic fashion right from the very beginning, and are not artificially 'manufactured' from surgery on links, etc.

One starts with a *classical field theory*, which satisfies various properties, the most important being *locality* (computable from cutting and pasting) and *symmetry* (behaves naturally under diffeomorphisms). Let us adopt the notation that $X, Y, T$ (resp. $X', Y', T'$) refer, respectively, to closed 3, 2 and 1 manifolds (resp. with boundary), and let $W$ be variable.

A classical field theory consists of a space of 'fields' $C_W$ for each $W$. These are usually local functions on $W$ if we interpret 'function' liberally enough. There should be a way to restrict fields to the boundary, $C_{W'} \to C_{\partial W'}$. For each field $\psi \in C_W$, there is defined an 'action' $\exp iS[\psi]$, which takes values in various structures according to the following table:

Classical Field Theory

$$\phi \in C_X \longrightarrow \exp iS[\phi] \text{ (a number)} \tag{5.40}$$

$$\gamma \in C_Y \longrightarrow \exp iS[\gamma] \text{ (a 1d Hilbert space) } \ni \phi|_{\partial X'} \tag{5.41}$$

$$\alpha \in C_T \longrightarrow \exp iS[\alpha] \text{ (a 1d 2-Hilbert space) } \ni \gamma|_{\partial Y'}$$

A '2-Hilbert space' is the categorification of a Hilbert space i.e. it is a $\mathbb{C}$-linear abelian, monoidal category $\mathcal{C}$ with an inner product $\mathcal{C} \times \mathcal{C} \to \mathcal{L}$, where $\mathcal{L}$ is the category of Hilbert spaces. It can be thought of as a category whose objects are formal sums of vector spaces (eg. a category of representations is a 1d 2-Hilbert space). Normally the assigning of $\gamma$ and $\alpha$ to one dimensional Hilbert and 2-Hilbert spaces respectively is determined in a natural geometric fashion from the behaviour of the original (number-valued) action $\exp iS[\phi]$. [7] This assignment is subject to:

(a) Symmetry. If $f : W_1 \to W_2$ is a diffeomorphism, then it should induce a map on the fields $f^* : C_{W_2} \to C_{W_1}$ such that the action is preserved:

$$\exp iS_{W_1}[f^*\psi] = \exp iS_{W_2}[\psi]. \tag{5.42}$$

(b) Locality. If we cut a spacetime $W$ along $V$ to obtain $W^{\text{cut}}$ as in Fig. 5.4a, and if $\psi^{\text{cut}}$ denotes the pullback field on $W^{\text{cut}}$, then we require that:

$$\exp iS_{W^{\text{cut}}}[\psi^{\text{cut}}] = \exp iS_W[\psi]. \tag{5.43}$$

This completes the geometrical input of a classical field theory. The data (5.41) can be thought of as a collection of bundles over the spaces of fields,

$$\begin{array}{ccc} L_X & L_Y & L_T \\ \big\downarrow{\scriptstyle \text{number}} & \big\downarrow{\scriptstyle \text{1d Hilbert space}} & \big\downarrow{\scriptstyle \text{1d 2-Hilbert space}} \\ C_X & C_Y & C_T \end{array} \tag{5.44}$$

whose fibers are given by numbers, 1d Hilbert spaces, and 1d 2-Hilbert spaces respectively.

[7] Freed showed as motivation that the Wess-Zumino-Witten action provides a 2d example of the above data.

The next step is to quantize the theory by a tower of path integrals. The important point is that this is an automated process - all the information of the theory resides in the initial data, viz. the spaces of fields $C_W$ and the actions $\exp iS[W]$.

The quantum theory integrates over the space of fields $C_W$ to assign an algebraic object to $W$. We are going to get:

$$\begin{array}{ccl} & \textsc{Quantum Field Theory} & \\ X & \longrightarrow & Z(X) \text{ (a number)} \qquad (5.45) \\ X' & \longrightarrow & Z(X') \in Z(\partial X') \text{ (a vector in a Hilbert space)} \\ Y & \longrightarrow & Z(Y) \text{ (a Hilbert space)} \qquad (5.46) \\ Y' & \longrightarrow & Z(\partial Y') \text{ (a Hilbert space in a 2-Hilbert space)} \\ T & \longrightarrow & Z(T) \text{ (a 2-Hilbert space)} \\ & \vdots & \qquad (5.47) \end{array}$$

This works as follows. For 3d manifolds:

$$\begin{aligned} X \quad &\longrightarrow \quad \int_{\phi \in C_X} \exp iS[\phi] = \text{a number.} \\ X' \quad &\longrightarrow \quad \int_{\phi \in C'_X} \exp iS[\phi] \\ &= \int_{\gamma \in C_{\partial X'}} \underbrace{\int_{\phi \in C_{X'} : \phi|_{\partial X'} = \gamma} \exp iS[\phi]}_{\text{vector in } \exp iS[\gamma]} = \text{a vector in } Z(\partial X) \end{aligned}$$

For 2d manifolds:

$$\begin{aligned} Y \quad &\longrightarrow \quad \int_{\gamma \in C_Y} \underbrace{\exp iS[\phi]}_{\text{1d Hilbert space}} = \text{a Hilbert space.} \\ Y' \quad &\longrightarrow \quad \int_{\gamma \in C'_Y} \exp iS[\gamma] \\ &= \int_{\alpha \in C_{\partial Y'}} \underbrace{\int_{\gamma \in C_{Y'} : \gamma|_{\partial Y'} = \alpha} \exp iS[\phi]}_{\text{1d Hilbert space in } \exp iS[\alpha]} = \text{a Hilbert space in } Z(\partial Y) \end{aligned}$$

For 1d manifolds (circles):

$$T \quad \longrightarrow \quad \int_{\alpha \in C_T} \underbrace{\exp iS[\phi]}_{\text{1d 2-Hilbert space}} = \text{a 2-Hilbert space.}$$

Alternatively, one may view the spaces $Z(W)$ as the space of sections of the bundles in (5.44),

$$Z(W) = \text{ space of sections of the bundle } L_W \to C_W. \tag{5.48}$$

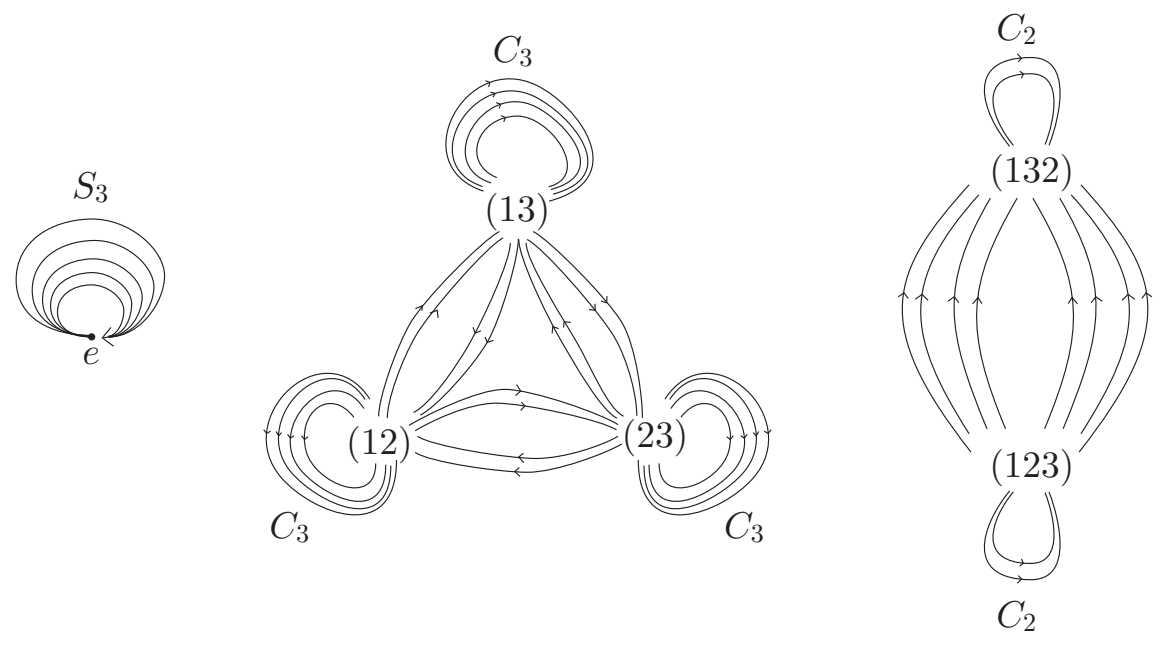

Figure 5.5: (a) $\mathfrak{C}_{S^1}$ for $G = S_3$, the symmetric group on three letters. Here $C_3 \simeq \mathbb{Z}_3$ is the cyclic group of order 3.

The main point about presenting the $n$-Hilbert spaces as path integrals over fields taking values in $(n-1)$-Hilbert spaces is that this allows one to carry over the gluing and symmetry laws (5.42) and (5.43) to these spaces. Thus not only do we have the usual topological field theory result (see Sec. 3.1.1) that (if we interpret $W_i$ and $V$ in Fig. 5.4b as 3- and 2-manifolds) the vectors $Z(W_1)$ and $Z(W_2)$ living inside the Hilbert space $Z(V)$ glue together to form the number $Z(W)$,

$$Z(W) = (Z(W_1), Z(W_2))_{Z(V)}, \tag{5.49}$$

but the same expression above can be interpreted as a gluing law for the quantum Hilbert spaces themselves; i.e. if we interpret $W_i$ and $V$ as 2- and 1-manifolds, then the monsters $Z(W_1)$ and $Z(W_2)$ living inside the 2-Hilbert space $Z(V)$ glue together to form the Hilbert space $Z(W)$ (recall that the inner product on a 2-Hilbert space outputs a Hilbert space). Expression (5.49) is then nothing but the Verlinde formula (5.20).

Using the symmetry and gluing laws leads to algebraic structure on the higher Hilbert spaces $Z(W)$. If we start with an $n$-dimensional field theory, then a codimension $p$ manifold $W_p$ will output a $p$-Hilbert space $Z(W_p)$. These higher categories hold more and more structure, a kind of compensation mechanism to make up for the loss in dimension. Amazingly, we shall see that the 2-Hilbert space of the circle, $Z(S^1)$, has the structure of a modular category.

## 5.5 The modular category of the 3d finite group model

Let us now illustrate how Freed's approach works in the case of the finite group Dijkgraaf-Witten model from Sec. 4.2. Our goal is to see that $Z(S^1)$ is a modular category.

We know what the spaces of fields are; for each closed manifold $W$, the space of fields $C_W$ is actually the *groupoid* of $G$-bundles over $W$. These spaces are large and a little bit awkward to work with. Since a $G$-bundle over $S^1$ is determined up to conjugation by its holonomy around the circle, a more convenient model is the full subcategory $\mathfrak{C}_{S^1}$ of $C_{S^1}$, whose objects are elements $x \in G$ and with a morphism labeled by $x \xrightarrow{g} gxg^{-1}$ for every pair of elements $(x, g)$ (see Fig. 5.5a). $\mathfrak{C}_{S^1}$ is a picture of the action of $G$ on itself by conjugation. Note that the set of objects isomorphic to $x$ is the conjugacy class $[x]$ of $x$, and $\mathrm{Aut}(x) = Z_x$, the centralizer of $x$.

Having defined the space of fields, we now define the action. It is *trivial*! That is,

$$\begin{aligned} \phi \in C_X \quad &\longrightarrow \quad \exp iS[\phi] = 1 \\ \gamma \in C_Y \quad &\longrightarrow \quad \exp iS[\gamma] = \mathbb{C} \\ \alpha \in C_{S^1} \quad &\longrightarrow \quad \exp iS[\alpha] = \mathcal{L} \end{aligned}$$

One sees that this satisfies the necessary gluing and symmetry axioms. In the quantized theory, we must

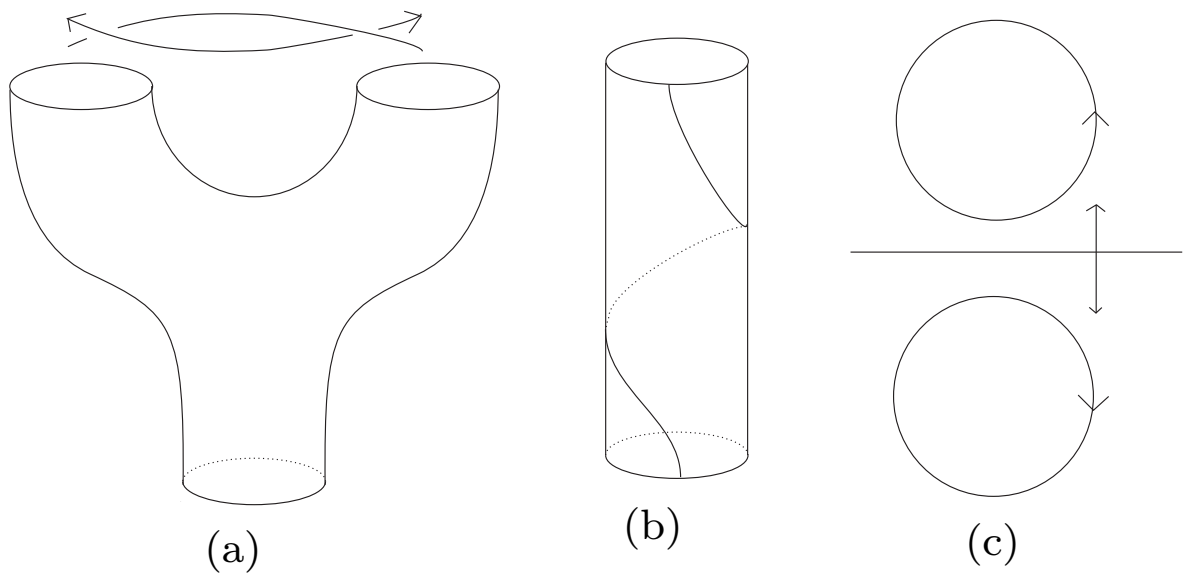

Figure 5.6: (a) The braiding diffeomorphism of the pair of pants. (b) The Dehn twist of the cylinder. (c) The reflection diffeomorphism of the circle.

integrate over the fields, or equivalently, determine the space of sections (5.48) of the 'trivial bundle'

$$\begin{array}{c} L \\ \big\downarrow {\scriptstyle \text{fiber} = \mathcal{L}} \\ \mathfrak{C}_{S^1} \end{array} \quad , \tag{5.50}$$

which has constant fibers equal to $\mathcal{L}$. The appropriate notion of a 'section' of $L$ is really a *functor* $\mathcal{W} : \mathfrak{C}_{S^1} \to \mathcal{L}$, that is, a collection of Hilbert spaces $W_x$ for each $x \in G$, together with maps $\rho_{x \xrightarrow{g}} : W_x \xrightarrow{\sim} W_{gxg^{-1}}$ for each arrow $x \xrightarrow{g}$ in $\mathfrak{C}_{S^1}$, which compose together like $\mathfrak{C}_{S^1}$. Note that the dimension of $W_x$ must be fixed over the conjugacy class $[x]$ of $x$, but it can jump across conjugacy classes. For each section $\mathcal{W}$ we can form the direct sum

$$W_{\mathcal{W}} = \oplus_x W_x \tag{5.51}$$

which is now simultaneously a

$$\begin{array}{rcl} \text{representation of } G & : & g \cdot w_x = \rho(x \xrightarrow{g}) w_x \in W_{gxg^{-1}} \text{ for } w_x in W_x \\ G\text{-graded vector space} & : & g \cdot W_x = W_{gxg^{-1}} \end{array}$$

We can call $W_{\mathcal{W}}$ an *equivariant vector bundle* over $G$, and the collection of such vector bundles defines

$$Z(S^1) = \mathrm{Vect}_G(G) = \{\text{equivariant vector bundles } W \to G\}. \tag{5.52}$$

The standard inner product in Ł is $(V_1, V_w) = \bar{V}_1 \otimes V_2$, and we extend this to $Z(S^1)$ by summing over the fibers. The game is now to see what extra structure exists on $Z(S^1)$, by using the gluing and symmetry laws. As usual, we can view a 2d surface $Y'$ with boundary as a cobordism from the 'input' circles to the output circles, and then use the inner product on $Z(S^1)$ to turn the Hilbert spaces $Z(Y')$ into a functor (the categorification of a linear map) $Z(Y') : \partial Y'_{\text{in}} \to \partial Y'_{\text{out}}$. In this way we obtain a tensor product $\otimes : Z(S^1) \boxtimes Z(S^1) \to Z(S^1)$ from the pair of pants, making $Z(S^1)$ a monoidal category. Here we have used the gluing property.

On the other hand, we can also use the symmetry law (5.42). Consider the braiding diffeomorphism on the pair of pants in Fig. 5.6a. This will correspond to a natural transformation $\square_1 \otimes \square_2 \to \square_2 \otimes \square_1 l$; that is, we get a *braiding* on $Z(S^1)$. Next consider the Dehn twist diffeomorphism of the cylinder in Fig. 5.6b. This will correspond to a natural transformation id $\to$ id of the identity functor on $Z(S^1)$; in other words we get a *twist*. We can also use diffeomorphisms at the level of the circle itself. Namely, let $f : S^1 \to S^1$ be an orientation reversing diffeomorphism of the circle as in Fig. 5.6c. This will correspond to a **-involution* on $Z(S^1)$ - a map $Z(S^1) \to Z(S^1)^{\mathrm{op}}$ which squares to the identity - which gives us duality[8] One uses geometric reasoning like this to show that $Z(S^1)$ satisfies all the axioms of a modular category.

[8]This is not strictly true. We get dual objects $x^*$, but not the maps $e : x^* \otimes x \to 1$ and $1 \to x \otimes x^*$. Clearing up this issue is an interesting topic [88]. Ross Street defines this as $*$-autonomy as opposed to autonomy [86].

## 5.6 The quantum double of a finite group

We have shown how to take the quantum double of a monoidal category $\mathcal{C}$ in order to get a modular category $\mathcal{D}(\mathcal{C})$. Since monoidal categories are (generally) representations of Hopf algebras (see Appendix B), and since modular categories are (generally) representations of quantum groups, one might ask if there is a direct construction on a Hopf algebra $H$ to give a quantum group $D(H)$. This is in fact the original approach of Drinfeld, and it is called the quantum double construction. This commutes with our categorical construction in the sense that

$$\mathrm{Rep}(D(H)) \simeq \mathcal{D}(\mathrm{Rep}(H)). \tag{5.53}$$

In particular, we shall apply this construction to compute the quantum double of the algebra $C(G)$ of functions on $G$. $C(G)$ has a basis $\{\delta_g\}_{g\in G}$ consisting of delta functions

$$\delta_g(x) = \delta_{g,x} = \begin{cases} 1 & \text{if } g = x \\ 0 & \text{otherwise.} \end{cases}$$

It is in fact a Hopf algebra since it has:

| | |
|---|---|
| multiplication | $\delta_g \delta_h = \delta_{g,h}\delta_g$ |
| unit | $1 = \sum_{g\in G} \delta_g$ |
| comultiplication | $\Delta(\delta_g) = \sum_{g_1 g_2 = g} \delta_{g_1} \otimes \delta_{g_2}$ |
| counit | $\epsilon(\delta_g) = \delta_{g,e}$ |
| antipode | $S(\delta_g) = \delta_{g^{-1}}$. |

Here we are drawing results from [18]. The quantum double $D(C(G))$, which we shall write as $D(G)$, is, as a vector space, equal to $C(G)\otimes \mathbb{C}[G]$. It is a Hopf algebra with

| | |
|---|---|
| multiplication | $(\delta_g \otimes x)(\delta_h \otimes y) = \delta_{gx,xh}(\delta_g \otimes xy)$ |
| unit | $1 = 1_{C(G)} \otimes e$ |
| comultiplication | $\Delta(\delta_g \otimes x) = \sum_{g_1 g_2 = g}(\delta_{g_1} \otimes x) \otimes (\delta_{g_2} \otimes x)$ |
| counit | $\epsilon(\delta_g \otimes x) = \delta_{g,e}$ |
| antipode | $S(\delta_g \otimes x) = \delta_{x^{-1}g^{-1}x} \otimes x^{-1}$. |

The Hopf algebra $D(G)$ is quasitriangular with R-matrix

$$R = \sum_{g\in G} (\delta_g \otimes e) \otimes (1 \otimes g)$$

Observe that $C(G)$ and $\mathbb{C}[G]$ embed in $D(G)$ as algebras and that $D(G)$ is their semidirect product,

$$D(G) = C(G) \rtimes \mathbb{C}[G]. \tag{5.54}$$

Consider the category $\mathrm{Rep}(D(G))$ of representations of $D(G)$ (this is the same as $\mathcal{D}(\mathrm{Rep}(G))$). Now, a representation of $C(G)$ as an algebra is nothing but a $G$-graded vector space, since the $\delta_g$ are projectors. On the other hand, a representation of $\mathbb{C}[G]$ is nothing but a representation of $G$. Taken together, we see that a representation of $D(C(G))$ is a $G$-graded representation of $G$. *Thus Rep(D(G)) is precisely equal to $Vect_G(G)$, the 2-Hilbert space $Z(S^1)$ of the circle computed in Sec. 5.5!* This is a remarkable result. It implies that the passage from the 2d finite group theory to the 3d finite group theory corresponded to *taking the quantum double of the category of representation of G.*

Looking at Fig. 5.5, the following explicit characterization of $\mathrm{Rep}(D(G))$ (taken from [18]) makes intuitive sense. The irreducible representations $V_{\bar{g},\pi}$ are labeled by pairs $(\bar{g},\pi)$, where $\bar{g}$ is a conjugacy class in $G$ and $\pi$ is an irreducible representation of the centralizer $Z(g)$. The duality, braiding and twist structures from (5.6)-(5.8) are calculated to be:

$$V^*_{(\bar{g},\pi)} \simeq V_{(g^{-1},\pi^*)} \text{ (so that } C_{(\bar{g},\pi),(\bar{g'},\pi'} = \delta_{(\bar{g},\pi),(g^{-1},\pi^*)}) \tag{5.55}$$

$$S_{(\bar{g},\pi),(\bar{g'},\pi')} = \frac{1}{|Z(g)||Z(g')|} \sum_{h\in G: hg'h^{-1}\in Z(g)} \mathrm{tr}_\pi(hg'^{-1}h^{-1})\mathrm{tr}_{\pi'}(h^{-1}g^{-1}h) \tag{5.56}$$

$$\theta_{(\bar{g},\pi)} = \frac{\mathrm{tr}_\pi(g)}{\mathrm{tr}_\pi(e)} \tag{5.57}$$

This explicit data offers us an interesting opportunity to 'close the loop' of all the constructions we have developed in the last two chapters. That is, using the above data we can get a TQFT $Z : \mathbf{3Cob} \to \mathbf{Vect}$ via the Reshetikhin-Turaev construction in Sec. 5.2.2. Let us, for example, compute the dimension of the Hilbert space of a genus $g$ surface via the general Verlinde formula for 3d TQFT's obtained from modular categories (5.35):

$$\dim Z(\Sigma_g) = D^{2g-2} \sum_i \frac{1}{\dim(V_i)^{2g-2}}. \tag{5.35}$$

We know from the geometrical construction that this should be the number of principal $G$-bundles $P \to \Sigma_g$, which we calculated in (4.33):

$$|C_{\Sigma_g}| = |\mathrm{Hom}(\pi_1(M,x),G)/G|. \tag{4.33}$$

Now $\dim(V_i) = \frac{S_{0i}}{S_{00}}$, and plugging in the values of $S$ in (5.56) gives the quantum dimension of the simple objects as

$$\dim(V_{(\bar{h},\pi)}) = \frac{|G|\dim(\pi)}{|Z(h)|}, \tag{5.58}$$

where on the right hand side we mean the ordinary (i.e. not quantum) dimension of $\pi$. Plugging this into the Verlinde formula (5.35) above, and recalling that $D^2 = \sum_i \dim(V_i)^2$, apparently gives an interesting formula for the number of homomorphisms *up to conjugation* of the fundamental group of $\Sigma_g$ into $G$:

$$|\mathrm{Hom}(\pi_1(\Sigma_g,x),G)/G| = \left( \sum_{(\bar{h}',\pi')} \frac{|G|\dim(\pi')}{|Z(h')|} \right)^{g-1} \sum_{(\bar{h},\pi)} \left( \frac{|Z(h)|}{|G|\dim(\pi)} \right)^{2g-2}. \tag{5.59}$$

This is to be compared with the formula we obtained from topological invariants for closed genus $g$ surfaces *from the 2d version of the theory* in (4.81):

$$\frac{|\mathrm{Hom}(\pi_1(\Sigma_g,x),G)|}{|G|} = |G|^{2g-2} \sum_\rho \frac{1}{(\dim\rho)^{2g-1}}. \tag{4.81}$$

Recall that $\rho$ is here a representation of $G$, not of a centralizer subgroup. Staring at these formulas makes us pause for reflection. We see that there is a relationship between the topological invariants of a theory in $n$ dimensions, and the dimensions of the Hilbert spaces of the same theory in $n+1$ dimensions. Presumably (5.59) can be analysed using the explicit presentation given in (4.18) of $\pi_1(\Sigma_g,x)$ as being generated by the $2g$ generators $a_1, b_1, \ldots, a_g, b_g$ subject to the relations

$$\prod_{i=1}^{g} [a_i, b_i] = 1. \tag{4.18}$$

However, to the present author's knowledge, no analysis explcitly relating the abovementioned formulas has appeared in the literature[9].

We end this chapter by briefly mentioning that one can in fact *twist* the finite group model. For example, in two dimensions this corresponds, on the algebraic level, to defining a new multiplication on $\mathbb{C}[G]$ by

$$g \star h = c(g,h)gh \tag{5.60}$$

where $c(g,h) : G \times G \to \mathbb{C}^\times$. Associativity then tells us that $c$ is a cohomology class in $H^2(G,\mathbb{C}^\times)$. On the geometric (principal bundles) level, this corresponds to choosing an integral cohomology class in $H^3(BG,\mathbb{Z})$. In three dimensions one must choose (algebraically) a cohomology class $c \in H^3(G,\mathbb{C}^\times)$, which corresponds to twisting the associator map in the monoidal category. On the geometric side, one must choose a cohomology class $c \in H^4(BG,\mathbb{Z})$. See [13, 25, 30].

[9] Although, see [23], [51], [65], Dijkgraaf and Witten's original paper [25], as well as the recent [99].

# Appendix A

# Appendix A - Categories and Functors

**Defn 1** *A* category $\mathcal{C}$ *consists of the following data:*

*(a) a class* $|C$ *whose elements are the objects of* $\mathcal{C}$*;*

*(b) a set* $Hom_{\mathcal{C}}(A, B)$ *of morphisms from* $A$ *to* $B$ *for every ordered pair* $(A, B)$ *of objects in* $\mathcal{C}$*;*

*(c) a composition function*

$$Hom_{\mathcal{C}}(B, C) \times Hom_{\mathcal{C}}(A, B) \to Hom_{\mathcal{C}}(A, C)$$

*for any objects* $A, B$ *and* $C$ *in* $\mathcal{C}$*;*

*(d) an identity morphism* $\mathrm{id}_A$ *in* $Hom_{\mathcal{C}}(A, A)$ *for every object* $A$ *in* $\mathcal{C}$*.*

*This data must satisfy the following conditions:*

*(a) identity morphisms are left and right units for composition :*

$$\mathrm{id}_B \circ f = f = f \circ \mathrm{id}_A.$$

*(b) composition of morphisms is associative :*

$$h \circ (g \circ f) = (h \circ g) \circ f.$$

We shall often refer to objects $A$ in categories as $A \in \mathcal{C}$ when we should really talk about objects $A$ in $|C|$ (which is not necessarily a set). In a category $\mathcal{C}$, we say that two objects $A$ and $B$ are *isomorphic* when there are morphisms $f : A \to B$ and $g : B \to A$ such that $gf = \mathrm{id}_A$ and $fg = \mathrm{id}_B$. A *groupoid* is a category where all the morphisms are invertible. Also, given a category $\mathcal{C}$ we can form a new category $\mathcal{C}^{\mathrm{op}}$ which has the same objects as $\mathcal{C}$ but with all morphisms reversed.

**Defn 2** *Let* $\mathcal{C}$ *and* $\mathcal{D}$ *be two categories. By a* functor $F : \mathcal{C} \to \mathcal{D}$ *we mean a rule that assigns to each object* $A$ *in* $\mathcal{C}$ *an object* $F(A)$ *in* $\mathcal{D}$*, and to each morphism* $f : A \to B$ *in* $\mathcal{C}$ *a morphism* $F(f) : F(A) \to F(B)$ *in* $\mathcal{D}$*, such that:*

*(a)* $F(\mathrm{id}_A) = \mathrm{id}_{F(A)}$ *for each object* $A$ *in* $\mathcal{C}$*,*

*(b)* $F(g \circ f) = F(g) \circ F(f)$ *for any two composable morphisms* $f$ *and* $g$ *in* $\mathcal{C}$*.*

**Defn 3** *Suppose* $F$ *and* $G$ *are two functors* $\mathcal{C} \underset{G}{\overset{F}{\Downarrow}} \mathcal{D}$ *from* $\mathcal{C}$ *to* $\mathcal{D}$*. A* natural transformation $\mathcal{C} \underset{G}{\overset{F}{\Downarrow \alpha}} \mathcal{D}$ *is a collection of morphisms*

$$\{\alpha_A : F(A) \to G(A)\}$$

*indexed by the objects of $\mathcal{C}$, such that for any morphism* $A \xrightarrow{f} B$ *in $\mathcal{C}$, the following square commutes in $\mathcal{D}$:*

$$\begin{array}{ccc} F(A) & \xrightarrow{\alpha_A} & G(A) \\ {\scriptstyle F(f)}\big\downarrow & & \big\downarrow{\scriptstyle G(f)} \\ F(B) & \xrightarrow[\alpha_B]{} & G(B) \end{array}$$

*If all the components $\alpha_A$ are isomoprhisms, then we say that $\alpha$ is a* natural isomorphism *from $F$ to $G$.*

Thinking a bit about natural transformations, and especially the picture 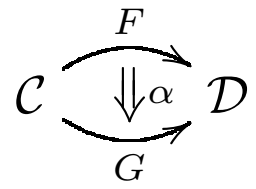, suggests that they are a hint of *higher dimensional structure*. This leads us to define 2-categories, very roughly, as follows (See [10] for further details).

**Defn 4** *A* 2-category $\mathfrak{C}$ *consists of objects $A$, 1-morphisms* $A \xrightarrow{f} B$ *and 2-morphisms* 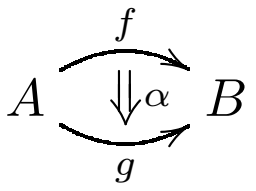, *where the objects and 1-morphisms behave as in an ordinary category, and we have the following two composition rules for 2-morphisms. (a)* Verical composition*:*

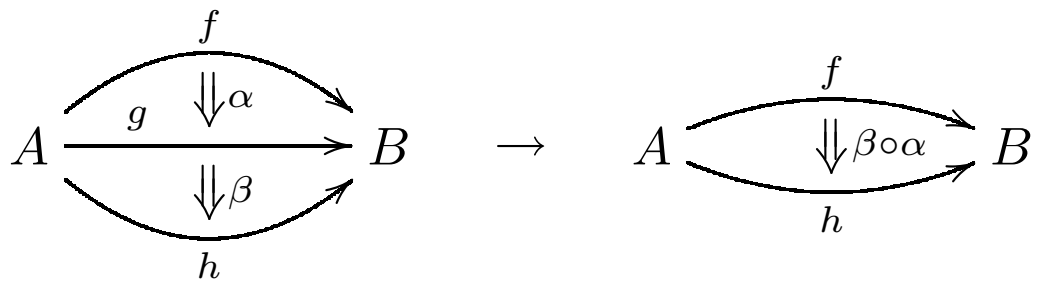

*and* horizontal composition*:*

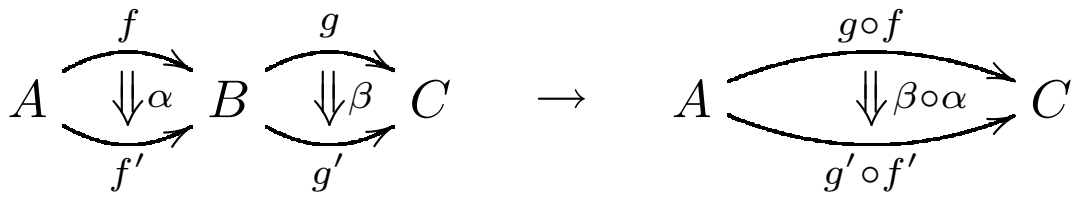

*which are associative, and such that the* interchange law *holds, i.e. that the following diagram is independent of which way it is interpreted:*

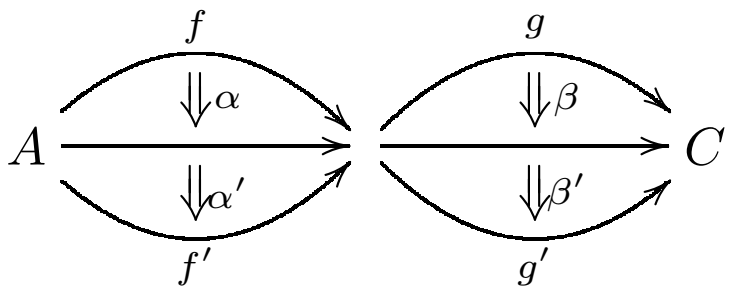

The most familiar example of a 2-category would be `Cat`, the 2-category whose objects are categories, whose 1-morphisms are functors, and whose 2-morphisms are natural transformations. A 2-groupoid is a 2-category in which all 1-morphisms and all 2-morphisms are invertible.

# Appendix B

# Appendix B - Hopf Algebras and Quantum Groups

We saw in Sec. 5.3 that one way to form modular categories is to start with a rigid (i.e. having duality) monoidal category and apply the quantum double construction. The question is, what data do we need to form rigid monoidal categories? The answer is that one needs a *Hopf algebra*[1].

If $G$ is a group, then we know how to form the representation category $\mathrm{Rep}(G)$. This category is monoidal since if $\rho_1$ and $\rho_2$ are representations then their tensor product is also a representation,

$$(\rho_1 \otimes \rho_2)(g) := \rho_1(g) \otimes \rho_2(g) \in \mathrm{End}(V_1 \otimes V_2). \tag{B.1}$$

This seemingly innocuous fact is really quite remarkable since the same trick will not work for representations of algebras (it will not be linear under scalar multiplication). On the other hand we know that representations of $G$ are the same as representations of $\mathbb{C}[G]$. Closer reflection reveals that the tensor product on the representations of $\mathbb{C}[G]$ are really using a kind of comultiplication map $\Delta : \mathbb{C}[G] \otimes \mathbb{C}[G] \to \mathbb{C}[G]$ defined by

$$\Delta(g) = g \otimes g \quad \text{(and extended to the rest by linearity).} \tag{B.2}$$

The tensor product of $\rho_1$ and $\rho_2$, considered as reps of $\mathbb{C}[G]$, is then given by

$$\rho_1 \otimes \rho_2(a) = (\rho_1 \otimes \rho_2)\Delta(a). \tag{B.3}$$

Let us also consider the duality structure in $\mathrm{Rep}(G)$. We know that for every representation $\rho$ on $V$, there is a dual representation $\rho^*$ on $V^*$ given by

$$\rho^*(g)(f)(x) = f(\rho(g^{-1})x). \tag{B.4}$$

We see that the inverse map $g \to g^{-1}$ enters in an essential way. Consideration of $\mathrm{Rep}(G)$ leads us to define a *Hopf algebra* as an algebra which possesses all the necessary structure so that its category of representations is a rigid monoidal category. Namely, $(A, \mu, \eta, \Delta, \epsilon, S) =$ 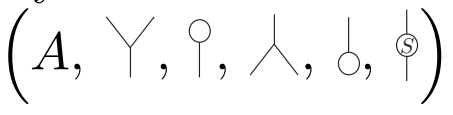 is a Hopf algebra if:

(a) 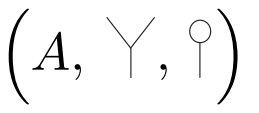 is an algebra : 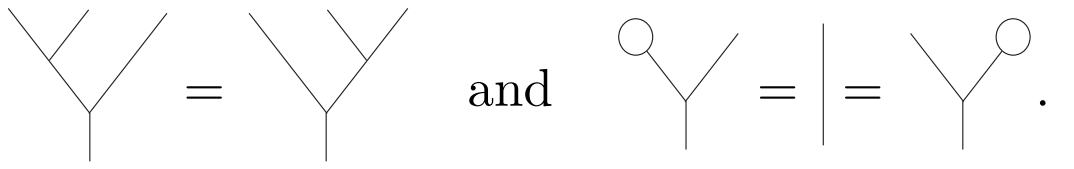

(b) 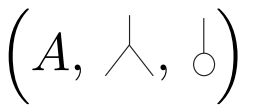 is a coalgebra : 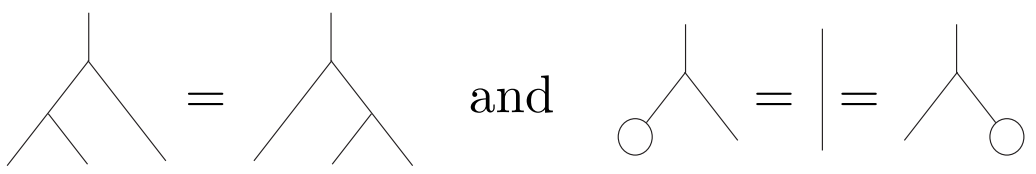

with the compatibility conditions[2] 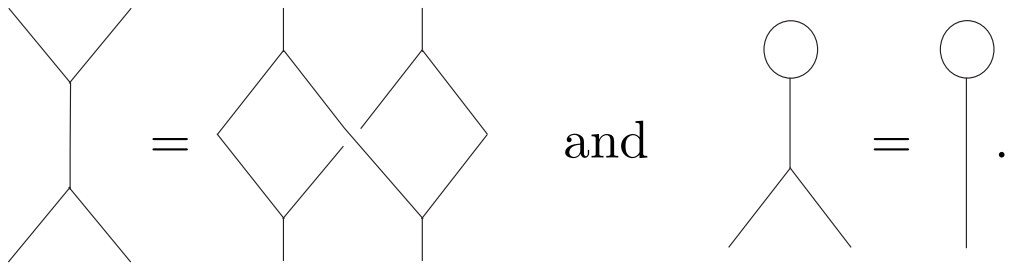

[1] We are restricting ourselves here to the $\mathbb{C}$-linear setting.

(c)

is an invertible map called the antipode : 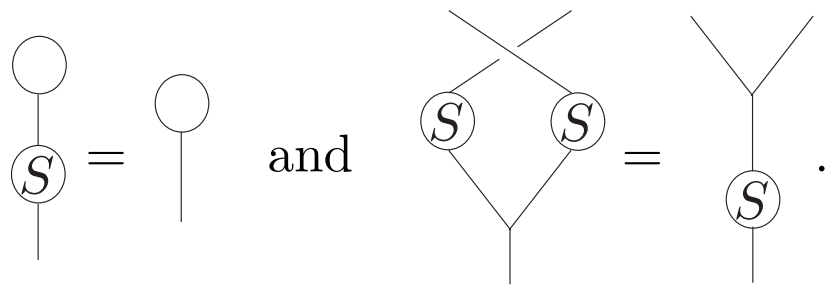

These definitions are to be compared with that of a Frobenius algebra (3.21). The first two conditions say that $A$ is a *bialgebra*, and the presence of the antipode makes it a Hopf algebra[3].

Similarly, in order for Rep($H$) to have a braiding it is necessary that $H$ be a *quasitriangular* Hopf algebra. That is, there exists an invertible element $R =$ 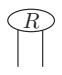 $\in H \otimes H$ such that:

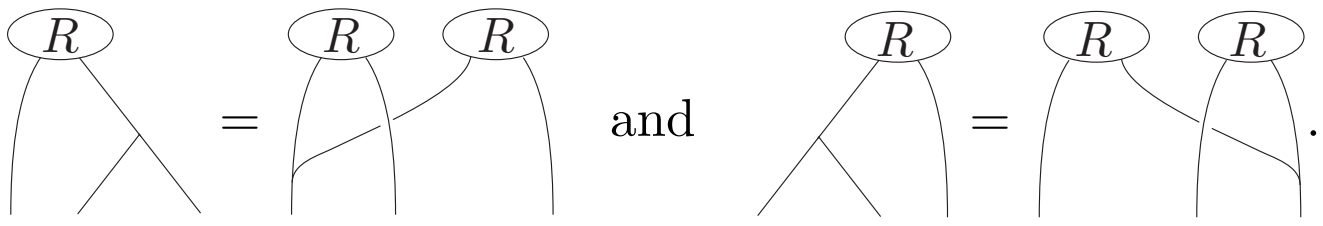

(B.5)

Finally, in order for Rep($H$) to have a twist, it should have a *grouplike* element $\theta \in H$ (i.e. $\Delta(\theta) = \theta \otimes \theta$) such that

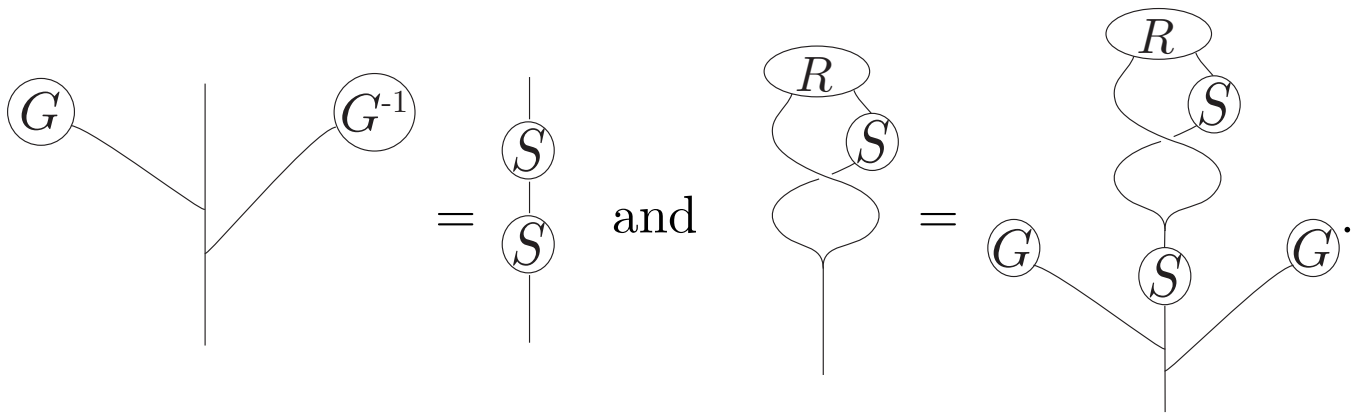

(B.6)

Such a monster is called a ribbon Hopf algebra, for obvious reasons. The best known constructions of them are quantum groups, which are deformations of the universal enveloping algebra $U(\mathfrak{g})$ of a Lie algebra $\mathfrak{g}$. For the group $SU(2)$, this is the algebra $U(\mathfrak{sl}_2)$ generated by $\{x, y, h\}$ with relations $[h, x]$

$$[h, x] = 2x \tag{B.7}$$
$$[h, y] = -2y \tag{B.8}$$
$$[x, y] = h \tag{B.9}$$

where $[a, b] = ab - ba$. Its Hopf algebra structure is determined by

$$\Delta(a) = a \otimes 1 + 1 \otimes a, \quad \epsilon(a) = 0, \quad S(a) = -a, \tag{B.10}$$

where $a$ is in $\{x, y, h\}$ and extended to formal products in the obvious way since these are all algebra homomorphisms. Its ribbon and quasitriangular structure is trivial,

$$R = 1 \otimes 1, \quad G = 1. \tag{B.11}$$

This algebra admits a one parameter family of deformations, $U_q(\mathfrak{sl}_2)$ depending on a complex pararmeter $q$, and as $q \to 1$ we recover $U(\mathfrak{sl}_2)$. These are called quantum groups. Namely, set $U(\mathfrak{sl}_2)$ as the algebra generated by $\{x, y, h\}$ subject to the relations

$$[h, x] = 2x \tag{B.12}$$
$$[h, y] = -2y \tag{B.13}$$
$$[x, y] = \frac{q^h - q^{-h}}{q - q^{-1}}. \tag{B.14}$$

Of course we have to find a way to interpret (B.14). This is done as follows. Define $U'(\mathfrak{sl}_2)$ as the completion of $U(\mathfrak{sl}_2)$ in the topology of convergence in every finite dimensional representation, and let the closure of the subalgebra generated by $h$ be $\mathcal{H}$. Note that the right hand side of (B.14) is in $\mathcal{H}$ for all $q \neq 0$, because every finite dimensional representation is spanned by eigenvectors of $h$ with integer eigenvalues (so that the exponents and division can be carried out safely). Consider the algebra spanned by free products of $x$ and $y$ and elements of $\mathcal{H}$, with the topology inherited from $\mathcal{H}$. The above equations generate an ideal in this algebra, and the quotient by the closure of that ideal is the algebra $U_q(\mathfrak{sl}_2)$.

[2]In all the diagrams on this page, it is understood that the braidings refer to the ordinary switch map of vector spaces $v \otimes w \to w \otimes v$.

[3]To see that these conditions encode the requirements on $H$ such that Rep($H$) is a Hopf algebra, see [14].

Choose a square root $\sqrt{q}$. The Hopf algebra structure of $U_q(\mathfrak{sl}_2)$ is given by

$$\begin{aligned} \Delta(x) &= x \otimes \sqrt{q}^h + \sqrt{q}^{-h} \otimes x \qquad & S(x) &= -qx \\ \Delta(y) &= y \otimes \sqrt{q}^h + \sqrt{q}^{-h} \otimes y \qquad & S(y) &= -q^{-1}y \end{aligned} \tag{B.15}$$

and the rest as in the $U(\mathfrak{sl}_2)$ case. The quasitriangular structure is quite a bit tricker. Define quantum integers and the quantum factorials as

$$[n]_q = \frac{q^n - q^{-n}}{q - q^{-1}}, \quad [n]_q! = [n]_q[n-1]_q \cdots [1]_q \tag{B.16}$$

Then

$$R = \sqrt{q}^{h \otimes h} \sum_{n=0}^{\infty} \frac{\sqrt{q}^{n(n+1)}(1-q^{-2})^n}{[n]_q!} x^n \otimes y^n. \tag{B.17}$$

$$G = q^h. \tag{B.18}$$

The infinite sum for $R$ converges because only finitely many terms are nonzero, as one can show by considering the ordinary reps of $SU(2)$. Taken together, this defines the quantum group $U_q(\mathfrak{sl}_2)$. Away from a $q$ a root of unity, its representations are the same as for ordinary $SU(2)$, so that there are infinitely many simple objects. At a root of unity, the quantum factorials $[n]_q$ equal zero at some point, and the whole system of representations is truncated, so that there are only finitely many simple objects. This is when the category becomes modular.

# Bibliography

[1] A. Lauda, *Frobenius algebras and open string topological field theories in the plane*, `math.QA/0508349`.

[2] A. Lauda, *Frobenius algebras and ambidextrous adjunctions*, `math.QA/0502550`.

[3] L. Abrams, *Two-dimensional topological quantum field theories and Frobenius algebras*, J. Knot Th. Ramif. **5**, 569-587 (1996).

[4] L. Abrams, *Modules, comodules and cotensor products over Frobenius algebras*, J. Alg. **219**, 201-213 (1999).

[5] T. W. Hungerford, *Algebra*, 8th Ed., Springer-Verlag 1997.

[6] A. Arageorgis, J. Earman and L. Ruetsche, *Weyling the time away: the non-unitary implementability of quantum field dynamics on curved spacetime*, from *Studies in the History and Philosophy of Modern Physics*, **33** (2002), 151-184.

[7] M. Atiyah, *Topological quantum field theory*, Publications Mathmatiques de l'IHS, 68 p. 176-186 (1988).

[8] M. Atiyah and R. Bott, *The Yang-Mills Equations over Riemann Surfaces*, Phil. Trans. R. Soc. Lond. A **308**, 523-615 (1982).

[9] J. Baez, *Quantum Quandaries : A Category-Theoretic Perspective*, `http://math.ucr.edu/home/baez/quantum/quantum.ps`, also to appear in *Structural Foundations of Quantum Gravity*, eds. Steven French, Dean Rickles and Juha Saatsi, Oxford University Press.

[10] J. Baez and J. Dolan, *Higher-dimensional algebra and topological quantum field theory*, Jour. Math. Phy. 36(1995), 6073-6105.

[11] J. Baez, Quantum Gravity notes, `http://math.ucr.edu/home/baez/qg-fall2004/connection.pdf`

[12] J. Baez, *Higher Yang-Mills Theory*, `arXiv:hep-th/0206130`

[13] J. Baez, Quantum Gravity seminar, Fall 2004, May 15 notes, `http://math.ucr.edu/home/baez/qg-fall2004/f04week09.pdf`.

[14] J. Baez, Quantum Gravity seminar, Fall 2004, April 17 notes, `http://math.ucr.edu/home/baez/qg-spring2001/4.17.01.pdf`.

[15] J. Baez, Quantum Gravity seminar, Winter 2005, `http://math.ucr.edu/home/baez/qg-winter2005/`

[16] J. Baez, Higher Gauge Theory, Homotopy Theory and n-Categories, `http://math.ucr.edu/home/baez/calgary/`.

[17] J. Barrett, *Holonomy and path structures in general relativity and Yang-Mills theory*, Int. J. Theor. Phys., **30** (1991), 1171-1215.

[18] B. Bakalov and A. Kirillov, Jr. *Lectures on Tensor Categories and Modular Functors*, AMS University Lecture Series Vol 21 (2001). Also available at `http://www.math.sunysb.edu/ kirillov/tensor/tensor.html`.

[19] J.L. Cardy, *Boundary conditions, fusion rules and the Verlinde formula*, Nucl. Phys. **B324** 581 (1989).

[20] *Clay Mathematics Institute Millenium Problems*, `www.claymath.org/millenium/`

[21] R. L. Cohen, A.A. Voronov, *Notes on String Topology*, http://de.arxiv.org/PS_Cache/math/pdf/0503/0503625.pdf (2005)

[22] R. L. Cohen, *Homotopy and geometric perspectives on string topology*, http://www.abdn.ac.uk/ wpe006/conference/skye/present/Cohen_13p.pdf. (2005)

[23] A. Coste, Terry Gannon, Phillippe Ruelle, *Finite group modular data*, Nucl. Phys. B 581 [PM] (2000) pg. 679-717.

[24] R. Dijkgraaf, *A geometric approach to two dimensional conformal field theory.* Ph.D. Thesis, University of Utrecht, 1989.

[25] R. Dijkgraaf and E. Witten, *Topological gauge theories and Group Cohomology*, Comm. Math. Phys., **129**, 393-429 (1990).

[26] B. Durhuus and T. Jonsson, *Classification and construction of unitary topological quantum field theories in two dimensions*, Jour. Math. Phys. 35 (1994), 5306-5313.

[27] B. Eckmann and P. Hilton, *Group like structures in categories*, Math. Ann. **145** (1962), 227-255.

[28] K. Ferguson, *Link invariants associated to TQFT's with finite gauge groups.* Jour. Knot. Th. and its Ram., Vol. 2 No. 1 (1993) pp. 11-36.

[29] M. Finkelberg, *An equivalence of fusion categories*, Geom. Funct. Anal. **6**, 249-267.

[30] D.S. Freed and F. Quinn, *Chern-Simons Theory with Finite Gauge Group*, Commun. Math. Phys. 156, 435-472 (1993).

[31] D.S. Freed, *Higher Algebraic Structures and Quantization*, Commun. Math. Phys. 159, 343-398(1994).

[32] D.S. Freed, *Locality and Integration in Topological Field Theory*, in Group Theoretical methods in Physics, Volume 2, eds. M.A. del Olmo, M. Santander and J. M. Guilarte, Ciemat, 1993, 35-54, also at http://arxiv.org/abs/hep-th/9209048.

[33] D.S. Freed, *Quantum Groups from Path Integrals*, http://arxiv.org/abs/q-alg/9501025.

[34] M. Fukama, S. Hosono, H. Kawai, *Lattice topological field theory in two dimensions*, Commun. Math. Phys. 161 (1994) 157-176 .

[35] J. Fuchs, I. Runkel and C. Schweigert, *TFT construction of RCFT correlators I : Partition Functions*, Nucl.Phys. B646 (2002) 353-497.

[36] J. Fuchs, I. Runkel and C. Schweigert, *TFT construction of RCFT correlators II: Unoriented world sheets* , Nucl.Phys. B678 (2004) 511-637.

[37] J. Fuchs, I. Runkel and C. Schweigert, *TFT construction of RCFT correlators III: Simple currents* , Nucl.Phys. B694 (2004) 277-353 .

[38] J. Fuchs, I. Runkel and C. Schweigert, *TTFT construction of RCFT correlators IV: Structure constants and correlation functions* , Nucl.Phys. B715 (2005) 539-638.

[39] J. Fuchs, I. Runkel and C. Schweigert, *TFT construction of RCFT correlators V: Proof of modular invariance and factorisation* , hep-th/0503194.

[40] D.S. Freed and F. Quinn, *Chern-Simons Theory with Finite Gauge Group*, Comm. Math. Phys. 156, no. 3 (1993), 435472.

[41] M. Gromov, *Pseudo-holomorphic curves in symplectic manifolds*, Inv. Math. **82** (1985) 307-347.

[42] B. Hatfield, *Quantum Field Theory of Point Particles and Strings*, Addison- Wesley, New York, 1992.

[43] A. Hatcher and W. Thurston, *A presentation for the mapping class group of a closed orientable surface*, Topology **19** (1980), 221-237.

[44] V. Jones, *A new knot polynomial and von Neumann algebras*, Notices A.M.S. 33, p.219-225, (1986).

[45] A. Joyal, R. Street, *An Introduction to Tannaka Duality and Quantum Groups*, Part II of 'Category Theory , Proceedings, Como 1990' Lecture Notes in Math. **1488** (Springer-Verlag Berlin, Heidelberg 1991) 411-492. Or : http://www.maths.mq.edu.au/ street/CT90Como.pdf.

[46] C.I. Lazaroiu, *On the structure of open-closed topological field theory in two dimensions*, hep-th/0009042 (2000).

[47] L.H. Kauffman, *Knots and Physics*, Teaneck, NJ, World Scientific Press, 1991 (K & E Series on Knots and Everything, vol. 1)

[48] See eg. Prasolov, V.V. and Sossinsky, A.B., *Knots, links, braids and 3-manifolds. An introduction to the new invariants in low-dimensional topology.* Translations of Mathematical Monographs **154**, Amer. Math. Soc., Providence, RI, (1997).

[49] T. Kerler, *Towards an Algebraic Characterization of 3-dimensional Cobordisms*, Contemp. Math., 318 (2003) 141-173. Also available at arXiv:math.GT/0106253.

[50] J. Kock, *Frobenius Algebras and 2D Topological Quantum Field Theories*, Cambridge University Press (2004), Cambridge.

[51] T. H. Koornwinder *et al*, J. Phys. A: Math. Gen. 32 8539-8549 (1999).

[52] C.I. Lazaroiu, *Instanton amplitudes in open-closed topological field theory*, hep-th/0011257 (2000).

[53] A. Lauda, private communication.

[54] D. Lewellen, *Sewing constraints for conformal field theories on surfaces with boundaries*, Nucl. Phys. **B372** 654 (1992).

[55] A.D. Wallace, *Modifications and cobounding manifolds*, Can. J. Math. **12** (1960) 503.

[56] B.R. Lickorish, *A representation of orientable combinatorial 3-manifolds*, Ann. of Math. (2) 76 (1962) 531-540.

[57] S. Majid, *Reconstruction Theorems and Rational Conformal Field Theories*, Int. J. Mod. Phys. A. 6 (1991) 4359-4374.

[58] S. MacLane, *Categories for the Working Mathematician*, Springer-Verlag, Berlin, 1972.

[59] Quotation from Peter May, in honour of one of the co-founders of category theory, Saunders Maclane (1910-2005). Taken from http://math.ucr.edu/home/baez/maclane.html.

[60] A. Migdal, Zh. Eksp. Teor. Fiz. **69** (1975) 810 (Sov. Phys. Jetp. **42** 413).

[61] S. Cordes, G. Moore and S. Ramgoolam, *Lectures on 2d Yang-Mills Theory, Equivariant Cohomology and Topological Field Theories*, Nucl. Phys. Proc. Suppl. 41 (1995) 184-244.

[62] G. Moore, *Some comments on Branes, G-flux, and K-theory,* hep-th/0012007, (2000).

[63] G. Moore, *D-Branes, RR-Fields and K-Theory I*, audio and slides available at http://online.itp.ucsb.edu/online/mp01/moore1/, (2001).

[64] G. Moore, N. Seiberg, *Classical and Quantum Conformal Field Theory*, Commun. Math. Phys. 123, 177-254 (1989).

[65] M. Müger, *On the structure of modular categories*, Proc. Lond. Math. Soc. 87, 291-308 (2003).

[66] M. Müger, special lecture given at University of California Riverside, see http://math.ucr.edu/home/baez/qg-spring2001/.

[67] T. Nakayama, *On Frobeniusean algebras. I*, Ann. of Math. (2), 40:611-633, (1939).

[68] T. Nakayama, *On Frobeniusean algebras. II*, Ann. of Math. (2), 42:1-21, (1941).

[69] C. Nesbitt, *On the regular representations of algebras.*, Ann. of Math. (2), 39, (1938).

[70] R. Oeckl, *Quantum Geometry and Quantum Field Theory*, Ph.D. Thesis, University of Cambridge, (2000). Available at http://www.matmor.unam.mx/ robert/papers/phd_thesis.pdf.

[71] U. Pachner, *Konstruktionsmethoden und das Kombinatorische Homöomorphiproblem für Triangulationen Kompakter Semilinearer Mannigfaltigkeiten*, Abh. Math. Sem. Hamb. **57** (1987) 69-86.

[72] R. Picken, *Reflections on Topological Quantum Field Theory*, Rep. Math. Phys., **40**, 295-303, (1997). Also available at `http://arxiv.org/PS_cache/q-alg/pdf/9707/9707002.pdf`.

[73] R. Picken, *TQFT and gerbes*, Alg. and Geom. Topology, **4**, (243-272) (2004).

[74] H.B. Posthuma, *Quantization of Hamiltonian loop group actions*, Phd Thesis, University of Amsterdam (2003). (See page 36). Available at `http://remote.science.uva.nl/ npl/proefschrift.pdf`.

[75] S. Abramsky and B. Coecke, *A categorical semantics of quantum protocols*, `quant-ph/0402130`.

[76] F. Quinn, *Lectures on axiomatic topological quantum field theory.* In : Geometry and quantum field theory (Park City, UT, 1991), 323-453, IAS/Park City Math. Ser., 1, Amer. Math. Soc., Providence, RI, (1995).

[77] N.A. Baas, R.L. Cohen, A. Ramírez, *The topology of the category of open and closed strings* `http://de.arxiv.org/PS_cache/math/pdf/0411/0411080.pdf` (2004)

[78] A.S. Schwarz, *Lett. Math. Phys.* **2** 247 (1978).

[79] G. Segal, Lectures at Stanford University on topological field theory, (partly) available at `http://www.cgtp.duke.edu/ITP99/segal/`.

[80] G. Segal, *The definition of conformal field theory*, from 'Topology, Geometry and Quantum Field Theory : Proceedings of the 2002 Oxford Symposium in Honour of the 60th birthday of Graeme Segal', London Mathematical Society Lecture Note Series (No. 308).

[81] G. Segal, *Elliptic Cohomology*, Séminaire Bourbaki 695 (1988) 187-201.

[82] G. Segal, it Topological structures in string theory, Phil. Trans. R. Soc. Lond. A (2001) **359**, 1389-1398.

[83] S. El-Showk, *Classifying D-brane RR-charge in low energy string theory using K-theory*, M.Sc. thesis, University of Amsterdam, available at `http://www.netherrealm.net/ sheer/uva/master_thesis/thesis_main.pdf`.

[84] R. Street and A. Joyal, *An introduction to Tannaka duality and quantum groups*, Part II of 'Category Theory, Proceedings, Como 1990' (Editors A. Carboni, M.C. Pedicchio and G. Rosolini), Lecture Notes in Math. **1488** (Springer-Verlag Berlin, Heidelberg 1991) 411-492. Also available at `http://www.maths.mq.edu.au/ street/CT90Como.pdf`.

[85] R. Street, *The quantum double and related constructions*, Journal of Pure and Applied Algebra, 132 (1998) 195-206.

[86] R. Street, *Frobenius monads and pseudomonoids*, Jour. Math. Phys., Vol 45. No. 10 (Oct 2004).

[87] A. Joyal and R. Street, *The geometry of tensor calculus I*, Adv. Math. 88 (1991), no. 1, 55-112.

[88] U. Tillman, *L-linear structures for k-linear categories and the definition of a modular functor*, J. London. Math. Soc. (2) 58 (1998) 208-228.

[89] V.G. Turaev, *Quantum Invariants of Knots and 3-Manifolds*, de Gruyter Studies in Mathematics 18, (1994).

[90] N.Y. Reshetikhin and V.G. Turaev, *Invariants of 3-manifolds via link polynomials and quantum groups*, Invent. Math. **103** (1991), 547-597.

[91] V.G. Turaev, *Homotopy field theory in dimension 2 and group-algebras*, `math.QA/9910010`

[92] V.G. Turaev, *Homotopy field theory in dimension 3 and crossed group-categories*, `math.GT/0005291`

[93] C. Vafa, *Toward classification of conformal field theories*, Phys Lett. **B 206** (1988), 421-426.

[94] C. Weibel, *An introduction to homological algebra*, Cambridge University Press (1995).

[95] E. Witten, *Quantum field theory and the Jones polynomial*, Comm. Math. Phys. **121** (1989), 351-399.

[96] E. Witten, *Two-dimensional gravity and intersection theory on moduli space*, Surveys of Diff. Geom. 1 (1991) 243-310.

[97] E. Witten, *On quantum gauge theories in two dimensions*, Commun. Math. Phys. 141, 153 (1991).

[98] E. Witten, *Chern-Simons theory as a string theory*, Prog. Math. 133 (1995), 637-678.

[99] S. Willerton, *The twisted Drinfeld double of a finite group via gerbes and finite groupoids*, `arXiv::math.QA/0503266`.

[100] U. Bunke, P. Turner and S. Willerton, *Gerbes and homotopy quantum field theories*, Alg. and Geom. Top. Volume 4 (2004) 407-437

[101] F. Wilczek, in *Fractional Statistics and Anyon Superconductivity*, edited by F. Wilczek, World Scientific Publishing Co., Singapore 1990.

[102] D.N. Yetter, *Triangulations and TQFT's*, taken from Quantum Topology, edited by Louis H. Kauffman and Rany A. Baadhio, Series on Knots and Everything Vol. 3, World Scientific, (1993).